\documentclass[preprint,3p,times]{elsarticle}

\usepackage[section]{placeins}
\usepackage{multirow}
\usepackage{makecell}
\usepackage{longtable}
\usepackage{booktabs}
\usepackage{supertabular}
\usepackage{supertabular}
\usepackage{color}
\usepackage{ulem}
\usepackage{colortbl}

\usepackage[table]{xcolor}
\usepackage{setspace}
\usepackage{bm}

\usepackage{soul}
\soulregister\cite7
\soulregister\citep7
\soulregister\citet7
\soulregister\ref7
\soulregister\pageref7

\usepackage{cancel}

\usepackage[bf]{caption}

\usepackage[ruled,linesnumbered,lined,commentsnumbered]{algorithm2e}

\newcommand{\tabincell}[2]{\begin{tabular}{@{}#1@{}}#2\end{tabular}}

\usepackage{amssymb}
\usepackage{amsmath}
\newtheorem{theorem}{Theorem}
\newtheorem{lemma}{Lemma}

\newtheorem{definition}{Definition}

\journal{Elsevier}

\begin{document}

\begin{frontmatter}

%\title{Improving the Z-type WENO schemes by introducing the 
%posteriori adaptive order-preserving mapping}
%\title{On improving the Z-type WENO schemes by introducing LOP 
%mapping}
%\title{Improved Z-type WENO schemes with locally order-preserving 
%mapping}
%\title{Application of Locally Order-Preserving Mapping in Improving 
%the WENO-Z-type Schemes}
\title{A general improvement in the WENO-Z-type schemes}
% running title: [A GML-GMWENO-X scheme]

\author[a]{Ruo Li}
\ead{rli@math.pku.edu.cn}

\author[b,c]{Wei Zhong\corref{cor1}}
\ead{zhongwei2016@pku.edu.cn}
%\ead{zhongwei@nint.ac.cn}

\cortext[cor1]{Corresponding author}

\address[a]{CAPT, LMAM and School of Mathematical Sciences, Peking
University, Beijing 100871, China}

\address[b]{School of Mathematical Sciences, Peking University,
Beijing 100871, China}

\address[c]{Northwest Institute of Nuclear Technology, Xi'an 
710024, China}

\begin{abstract}
  
  A new type of finite volume WENO schemes for hyperbolic problems 
  was devised in \cite{MOP-WENO-ACMk} by introducing the 
  order-preserving (OP) criterion. In this continuing work, we 
  extend the OP criterion to the WENO-Z-type schemes. We firstly 
  rewrite the formulas of the Z-type weights in a uniform form from 
  a mapping perspective inspired by extensive numerical 
  observations. Accrodingly, we build the concept of the locally 
  order-preserving (LOP) mapping which is an extension of the 
  order-preserving (OP) mapping and the resultant improved 
  WENO-Z-type schemes are denoted as LOP-GMWENO-X. There are four 
  major advantages of the LOP-GMWENO-X schemes superior to the 
  existing WENO-Z-type schemes. Firstly, the new schemes can amend 
  the serious drawback of the existing WENO-Z-type schemes that most 
  of them suffer from either producing severe spurious oscillations 
  or failing to obtain high resolutions in long calculations of  
  hyperbolic problems with discontinuities. Secondly, they can 
  maintain considerably high resolutions on solving problems with 
  high-order critical points at long output times. Thirdly, they can 
  obtain evidently higher resolution in the region with 
  high-frequency but smooth waves. Finally, they can significantly 
  decrease the post-shock oscillations for simulations of some 2D 
  problems with strong shock waves. Extensive benchmark examples are 
  conducted to illustrate these advantages.

\end{abstract}

%%% Local Variables:
%%% mode: latex
%%% TeX-master: "article"
%%% End:


\end{frontmatter}

\section{Introduction}
\label{secIntroduction}
During the past thirty years, the weighted essentially 
non-oscillation (WENO) schemes \cite{ENO-Shu1988,ENO-Shu1989,WENO-LiuXD,WENO-JS,WENOoverview,AMR-WENO-Liu-CiCP,HGKS-WENO-Ji-CiCP,Ada-WENO-Chen-CiCP2021,DG-WENO-Zhu-CiCP2021}, as well as the 
recently-published low-dissipation shock-capturing ENO-family 
schemes, dubbed TENO \cite{TENO-Fu-JCP2016,TENO-Fu-JCP2018,TENO-Fu-CPC2019-01,TENO-Fu-CPC2019-02,TENO-Fu-JSC2019,TENO-Fu-CiCP2019,TENO-Fu-CMAME2021,TENO-Haimovich-CF2017}, is a major area of interest within the field of high-resolution 
numerical simulation for the following hyperbolic conservation laws
\begin{equation}
    \begin{array}{ll}
      \dfrac{\partial \mathbf{u}}{\partial t} +
      \displaystyle\sum\limits_{\alpha = 1}^{d} \dfrac{\partial
      \mathbf{f}_{\alpha}(\mathbf{u})}{\partial x_{\alpha}} = 0,
      & x_{\alpha} \in \mathbb{R}, t > 0,
\end{array}
\label{governingEquation}
\end{equation}
where $\mathbf{u} = (u_{1}, u_{2}, \cdots, u_{m}) \in \mathbb{R}^{m}$
are the conserved variables and $\mathbf{f}_{\alpha}: \mathbb{R}^{m} 
\rightarrow \mathbb{R}^{m}$, $\alpha = 1,2,\cdots,d$ are the 
Cartesian components of flux.

By using the information of all candidate substencils of the 
original essentially non-oscillation (ENO) scheme \cite{ENO1987JCP71,
ENO1987V24, ENO1986, ENO1987JCP83} through a convex combination, Liu 
et al. \cite{WENO-LiuXD} proposed the $(r + 1)$th-order WENO scheme. 
Later, Jiang and Shu \cite{WENO-JS} improved it by introducing a 
new measurement of the smoothness of a solution over a particular 
substencil, say, local smoothness indicator (LSI), into the 
$(2r + 1)$th-order one, dubbed WENO-JS. WENO-JS is the most widely 
used one among the family of the WENO schemes since it was proposed 
as it can maintain the ENO property near discontinuities and obtain 
the designed convergence rate of accuracy in most smooth regions. 
However, it was commonly known that \cite{WENO-M,WENO-PM,WENO-IM,
WENO-RM260,WENO-MAIMi,WENO-ACM,WENO-Z,article_WENO-MS_JSC2019} 
WENO-JS fails to obtain the optimal accuracy near critical points of 
order $n_{\mathrm{cp}}$, and this was originally discovered by 
Henrick et al. \cite{WENO-M}. Here, $n_{\mathrm{cp}}$ stands for the 
order of the critical point for the function $f$ that satisfies 
$f' = 0, \cdots, f^{(n_{\mathrm{cp}})} = 0, f^{(n_{\mathrm{cp}} + 1)}
\neq 0$. In the same article, Henrick et al. performed the 
truncation error analysis and derived the necessary and sufficient 
conditions on the nonlinear weights of the WENO schemes to achieve 
the designed fifth-order convergence in smooth regions of the 
solution. Then, they designed a mapping function to the original 
weights of the WENO-JS scheme resulting in the mapped weights 
satisfying these conditions. The resultant scheme, denoted as 
WENO-M, successfully recovered the designed convergence orders even 
at or near the critical points. It is since the introduction of 
Henrick et al. \cite{WENO-M} that various versions of mapped WENO 
schemes, such as WENO-PM$k$ \cite{WENO-PM}, WENO-IM($k,A$) 
\cite{WENO-IM}, WENO-PPM$n$\cite{WENO-PPM5}, WENO-RM($mn0$) 
\cite{WENO-RM260}, WENO-MAIM$i$ \cite{WENO-MAIMi}, WENO-ACM 
\cite{WENO-ACM}, and etc., have been devoloped by devising different 
mapping functions under the similar principles of WENO-M. Over the 
past decade, there has been an increasing amount of literature on 
the long output time simulations of mapped WENO schemes 
\cite{WENO-PM,WENO-IM,WENO-RM260,WENO-AIM,WENO-MAIMi,MOP-WENO-ACMk,
MOP-WENO-X}. A key issue of WENO-M is that its resolution decreases 
dramatically when solving problems with discontinuities for long 
output times, and this drawback was first noticed and successfully 
fixed by Feng et al. \cite{WENO-PM}. After that, a series of mapped 
WENO schemes, as reported in \cite{WENO-PM,WENO-IM,WENO-RM260,
WENO-MAIMi}, have been developed to address this potential loss of 
accuracy properly. However, in these same articles, it was 
illustrated that these schemes have caused another serious problem 
that they produced severe spurious oscillations for long output time 
calculations because of the lack of robustness. Taken together, it 
is rather difficult for the previously published mapped WENO schemes 
to avoid spurious oscillations while obtain high resolutions at the 
same time for long output times simulations. The essential reason of 
such phenomena has been revealed in a recently published article 
\cite{MOP-WENO-ACMk} in which the core concept of 
\textit{oder-preserving (OP) mapping} was innovatively proposed 
resulting in the OP-Mapped WENO schemes \cite{MOP-WENO-ACMk,
MOP-WENO-X} that can preserve high resolutions and meanwhile prevent 
spurious oscillations no matter in short or long output time 
simulations. We refer to the corresponding literature for more 
details.

The necessary and sufficient conditions on the nonlinear weights for 
optimal order of convergence proposed by Henrick et al. in 
\cite{WENO-M} was simplified to a sufficient condition by Borges et 
al. in \cite{WENO-Z} where they proposed another version of 
nonlinear weights by introducing the global smoothness indicator 
(GSI) of higher-order. In \cite{WENO-Z}, the GSI was computed via a 
linear combination of the LSIs on substencils of the WENO-JS scheme 
and it was used to define the new nonlinear weights leading to the 
WENO-Z scheme. It was analyzed theoretically and examined 
numerically that \cite{WENO-Z} the nonlinear weights of WENO-Z 
satisfy the sufficient condition for optimal order of convergence 
and hence it can recover the designed convergence order of accuracy 
properly by choosing a suitable tunable parameter $p = 2$. As it 
only used available and previously unused information of WENO-JS, 
the extra computational cost of WENO-Z compared to WENO-JS is very 
small and even negligible in practice. Naturally, WENO-Z is much 
cheaper than WENO-M. By obeying the similar principles proposed by 
Borges et al. \cite{WENO-Z}, extensive nonlinear weights with 
various GSIs \cite{WENO-eta,WENO-eta-02,WENO-P,WENO-Zplus,
WENO-ZplusI,WENO-ZA,P-WENO,MWENO-P,WENO-D_WENO-A,WENO-NIP}, denoted 
as Z-type weights in this paper, have been developed, and we 
collectively call these resultant schemes as the ``WENO-Z-type'' 
schemes.

Despite the success mainly for short output time calculations, the 
WENO-Z-type schemes also suffer from either producing spurious 
oscillations or failing to achieve high resolutions on solving 
hyperbolic problems with discontinuities at long output times. 
The same issue with respect to the family of the mapped WENO schemes 
has received considerable attention over the past ten years 
\cite{WENO-PM,WENO-IM,WENO-RM260,WENO-RM-Vevek2018,WENO-AIM,
WENO-MAIMi,MOP-WENO-ACMk}. However, there has been little discussion 
about the WENO-Z-type schemes on this topic so far. Of course, it is 
worthy of scholarly attention to examine the performance of the 
WENO-Z-type schemes for long output time simulations, and we mainly 
focus our attention on this theme in this paper.

We firstly present the profiles of the implicit relationship between 
the nonlinear weights of the WENO-JS scheme and the Z-type weights, 
denoted as IMR (standing for \textit{implicit mapping relation}). 
Then, a uniform form of the formulas of the Z-type weights is 
provided from a mapping perspective. We develop the concept of 
\textit{locally order-preserving (LOP) mapping} and introduce it to 
get the improved Z-type weights. We conduct extensive numerical 
experiments of 1D linear advection equation with various intial 
conditions for long output times to demonstrate the advantages of 
the improved WENO-Z-type schemes. Moreover, some benchmark tests of 
1D and 2D Euler systems have been calculated to show the good 
performance of these schemes.

The rest of this paper proceeds as follows. In Section 
\ref{secWENO-Z-type}, by way of preliminaries, we provide a 
brief description of the WENO-JS \cite{WENO-JS} and WENO-Z \cite{WENO-Z} schemes. In Section \ref{LOP-GMWENO-X}, the improved 
WENO-Z-type schemes are constructed by introducing the locally 
order-preserving mapping. Several typical numerical tests are also 
performed to demonstrate the major advantages of these improved 
WENO-Z-type schemes in this section. In Section \ref{NumericalExperiments}, some more benchmarck examples are provided to 
show the remarkable performance and some additional enhancements of 
the new schemes. Concluding remarks are given in Section \ref{secConclusions}.

%%% Local Variables:
%%% mode: latex
%%% TeX-master: "article"
%%% End:

\section{Preliminaries}
\label{secWENO-Z-type}
\subsection{The finite volume methodology}
For simplicity but without loss of generality, we restrict our 
attention to the following one-dimensional scalar hyperbolic 
conservation law
\begin{equation}
\dfrac{\partial u}{\partial t} + \dfrac{\partial f(u)}{\partial x} 
= 0, \quad x_{l} < x < x_{r}, t > 0.
\label{eq:1D-scalar}
\end{equation}
Within the framework of the finite volume method, the computational 
domain is discretized into non-overlapping cells. In this section, 
we only focus on the uniform mesh cells and hence the domain is 
discretized into smaller uniform cells $I_{j} = [x_{j-1/2}, 
x_{j+1/2}]$ with the width $\Delta x = (x_{r} - x_{l})/N$, the 
interfaces $x_{j\pm1/2} = x_{j} \pm \Delta x/2$, and the cell 
centers $x_{j} = (x_{j-1/2} + x_{j+1/2})/2$. Then, we can transform 
Eq. \eqref{eq:1D-scalar} into the following semi-discretized form 
after some simple mathematical manipulations
\begin{equation}
\dfrac{\mathrm{d}\bar{u}_{j}(t)}{\mathrm{d}t} \approx -\dfrac{1}
{\Delta x}\bigg( \hat{f}\big(u_{j+1/2}^{\mathrm{L}},u_{j+1/2}^{\mathrm{R}}\big) - \hat{f}\big(u_{j-1/2}^{\mathrm{L}},u_{j-1/2}^{\mathrm{R}}\big) \bigg).
\label{eq:discretizedFunction}
\end{equation}
Let $\bar{u}(x_{j}, t)=\dfrac{1}{\Delta x}\int_{x_{j-1/2}}^{x_{j+1/2}
}u(\xi,t)\mathrm{d}\xi$ be the cell average of $I_{j}$. In Eq.
\eqref{eq:discretizedFunction}, $\bar{u}_{j}(t)$ is the numerical 
approximation to $\bar{u}(x_{j}, t)$ and $\hat{f}\big(u_{j \pm 1/2}^{
\mathrm{L}}, u_{j \pm 1/2}^{\mathrm{R}}\big)$ is the numerical flux 
used to replace the physical flux function $f(u)$ at the cell 
boundaries $x_{j \pm 1/2}$. In this paper, the global Lax-Friedrichs 
flux $\hat{f}(a,b) = \frac{1}{2}\big[f(a) + f(b) -\alpha(b - a)\big]$
with $\alpha = \max_{u} \lvert f'(u) \rvert$ is employed. The values 
of $u_{j \pm 1/2}^{\mathrm{L}/\mathrm{R}}$ will be reconstructed by 
the WENO reconstructions in this paper. It is trivial to know that 
$u_{j+1/2}^{\mathrm{R}}$ is symmetric to $u_{j+1/2}^{\mathrm{L}}$ 
with respect to $x_{j+1/2}$. Thus, just for the sake of brevity, we 
only describe the approximate procedure for $u_{j+1/2}^{\mathrm{L}}$ 
in the following discussion and we also drop the superscript ``L'' 
without causing any confusion.

In our calculations below, the ODE system Eq.\eqref{eq:discretizedFunction} is marched in time by using the third-order, 
TVD, Strong Stability Preserving (SSP) Runge-Kutta method \cite{ENO-Shu1988,SSPRK1998,SSPRK2001}.

\subsection{WENO-JS}\label{subsec:WENO-JS}
In the classical fifth-order WENO-JS scheme, a 5-point global 
stencil, hereafter named $S^{5}$, is used. $S^{5}$ is divided into 
three 3-point substencils, say, $\{S_{0}, S_{1}, S_{2}\}$ with 
$S_{k} = \{I_{j+k-2}, I_{j+k-1}, I_{j+k}\}, k = 0, 1, 2$. On each 
substencil $S_{k}$, the corresponding third-order approximation of 
$u_{j+1/2}$ is computed by
\begin{equation}
\begin{array}{l}
\begin{aligned}
&u_{j+1/2}^{0} = \dfrac{1}{6}(2\bar{u}_{j-2} - 7\bar{u}_{j-1}
+ 11\bar{u}_{j}), \\
&u_{j+1/2}^{1} = \dfrac{1}{6}(-\bar{u}_{j-1} + 5\bar{u}_{j}
+ 2\bar{u}_{j+1}), \\
&u_{j+1/2}^{2} = \dfrac{1}{6}(2\bar{u}_{j} + 5\bar{u}_{j+1}
- \bar{u}_{j+2}).
\end{aligned}
\end{array}
\label{eq:approx_ENO}
\end{equation}
Then, the fifth-order approximation of $u_{j+1/2}$ on the global 
stencil $S^{5}$ is defined by
\begin{equation}
u_{j + 1/2} = \sum\limits_{k = 0}^{2}\omega_{k}u_{j + 1/2}^{k},
\label{eq:approx_WENO}
\end{equation}
where the nonlinear weights $\omega_{k}$ for the classical WENO-JS 
scheme proposed by Jiang and Shu \cite{WENO-JS} are given as
\begin{equation} 
\omega_{k}^{\mathrm{JS}} = \dfrac{\alpha_{k}^{\mathrm{JS}}}{\sum_{l =
 0}^{2} \alpha_{l}^{\mathrm{JS}}}, \quad \alpha_{k}^{\mathrm{JS}} = \dfrac{d_{k}}{(\epsilon + IS_{k})^{2}}, \quad k = 0,1,2,
\label{eq:weights:WENO-JS}
\end{equation} 
and the coefficients $d_{0} = 0.1, d_{1} = 0.6, d_{2} = 0.3$ are the 
ideal weights of $\omega_{k}$ since they generate the central 
upstream fifth-order scheme for the global stencil $S^{5}$. The 
small positive number $\epsilon$ is used to prevent the denominator 
from becoming zero. Following \cite{WENO-JS}, the smoothness 
indicator $IS_{k}$ used to measure the regularity of the $k$th 
polynomial approximation $\hat{u}^{k}(x_{j})$ on the substencil 
$S_{k}$ is defined by
\begin{equation*}
IS_{k} = \sum_{l=1}^{2}\Delta x^{2l - 1}\displaystyle\int_{x_{j-1/2}}^{x_{j+1/2}}\bigg(\dfrac{\mathrm{d}^{l}\hat{u}^{k}(x)}{\mathrm{d}x^{l}}\bigg)^{2}\mathrm{d}x.
\label{eq:IS_k:def}
\end{equation*}
Accordingly, the smoothness indicators take on the following 
intuitive form
\begin{equation*}
\begin{array}{l}
\begin{aligned}
IS_{0} &= \dfrac{13}{12}\big(\bar{u}_{j - 2} - 2\bar{u}_{j - 1} + 
\bar{u}_{j} \big)^{2} + \dfrac{1}{4}\big( \bar{u}_{j - 2} - 4\bar{u}_
{j - 1} + 3\bar{u}_{j} \big)^{2}, \\
IS_{1} &= \dfrac{13}{12}\big(\bar{u}_{j - 1} - 2\bar{u}_{j} + \bar
{u}_{j + 1} \big)^{2} + \dfrac{1}{4}\big( \bar{u}_{j - 1} - \bar{u}_{
j + 1} \big)^{2}, \\
IS_{2} &= \dfrac{13}{12}\big(\bar{u}_{j} - 2\bar{u}_{j + 1} + \bar
{u}_{j + 2} \big)^{2} + \dfrac{1}{4}\big( 3\bar{u}_{j} - 4\bar{u}_{j 
+ 1} + \bar{u}_{j + 2} \big)^{2}.
\end{aligned}
\end{array}
\end{equation*}

\subsection{WENO-Z}\label{subsec:WENO-Z}
Borges et al. \cite{WENO-Z} proposed the WENO-Z scheme to recover 
the optimal convergence rate of accuracy at critical points. By 
introducing a GSI $\tau_{5}$, they suggested a new method to 
calculate the nonlinear weights
\begin{equation}
\omega_{k}^{\mathrm{Z}}=\dfrac{\alpha_{k}^{\mathrm{Z}}}{\sum_{l=0}^{2
}\alpha_{l}^{\mathrm{Z}}}, \quad \alpha_{k}^{\mathrm{Z}} = d_{k}
\bigg(1 + \Big(\dfrac{\tau_{5}}{IS_{k} + \epsilon} \Big)^{p} 
\bigg), \quad k = 0,1,2.
\label{eq:weights:Z}
\end{equation}
The original GSI here is
\begin{equation*}
\tau_{5} = \lvert IS_{0} - IS_{2} \rvert.
\label{eq:GSI:Z}
\end{equation*}
In Eq. \eqref{eq:weights:Z}, the parameters $IS_{k}$ and $\epsilon$ 
are the same as in the classical WENO-JS scheme, and $p$ is a 
tunable parameter. It was indicated by Borges et al. that 
\cite{WENO-Z}, if $p = 1$, the WENO-Z scheme only has fourth-order 
convergence rate of accuracy at critical points; if $p = 2$, it can 
achieve fifth-order convergence rate of accuracy. Unless indicated 
otherwise, we choose $p = 2$ in the present study.

The work of Borges et al. \cite{WENO-Z} inspired a significant 
increase in the study of WENO-Z-type schemes, like WENO-NS \cite{WENO-NS}, WENO-Z$\eta(\tau_{l})$ \cite{WENO-eta,WENO-eta-02}, 
WENO-Z+ \cite{WENO-Zplus}, WENO-Z+I \cite{WENO-ZplusI}, WENO-ZA \cite{WENO-ZA}, WENO-D and WENO-A \cite{WENO-D_WENO-A}, WENO-NIP \cite{WENO-NIP}, etc. For systematically reviews in detail, we refer the 
reader to the literature. Just for brevity in the presentation but 
without loss of generality, we mainly devote our attention to the 
WENO-Z, WENO-Z$\eta(\tau_{81})$ and WENO-A schemes in the rest of 
this paper. It should be pointed out that the new method proposed in 
this study below can easily be extended to other WENO-Z-type schemes.

%%% Local Variables:
%%% mode: latex
%%% TeX-master: "article"
%%% End:

\section{Design and properties of improved WENO-Z-type schemes with 
LOP mappings}
\label{LOP-GMWENO-X}
\subsection{A mapping perspective for Z-type weights}
\label{subsec:IMRs_Z-type}
From a mapping perspective, we can plot the profiles of 
$\omega_{k}^{\mathrm{JS}} \sim \omega_{k}^{\mathrm{X}}$ in practical 
calculations, where ``X'' stands for some kind of WENO-Z-type 
scheme, for example, X = Z, Z$\eta(\tau_{81})$, A, etc. As no 
explicit mapping functions are necessary for the plotting of 
$\omega_{k}^{\mathrm{JS}} \sim \omega_{k}^{X}$, we call it 
\textit{Implicit Mapping Relation (IMR)} for the sake of simplicity.

\begin{figure}[ht]
\flushleft
  \includegraphics[height=0.28\textwidth]
  {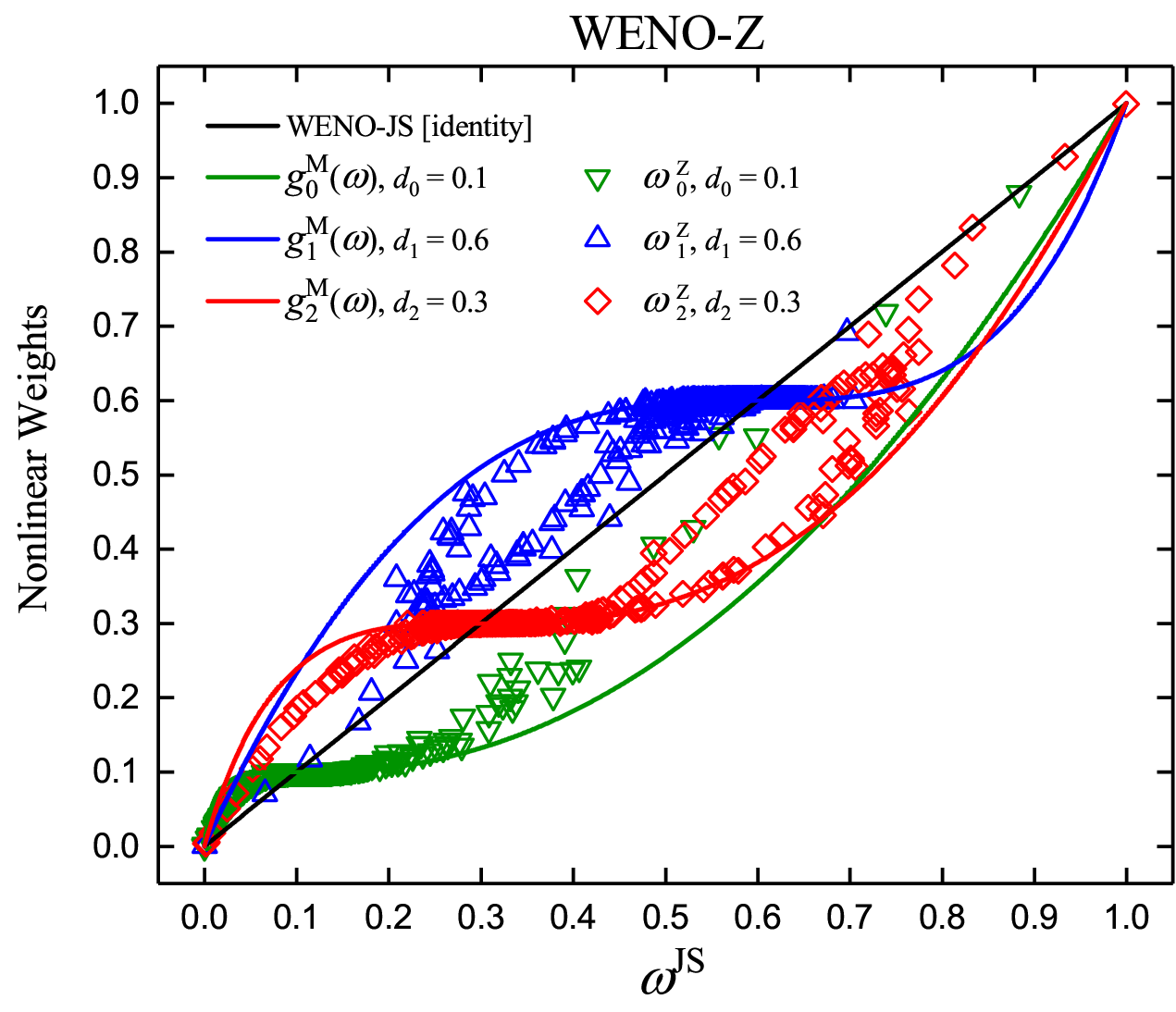}
  \includegraphics[height=0.28\textwidth]
  {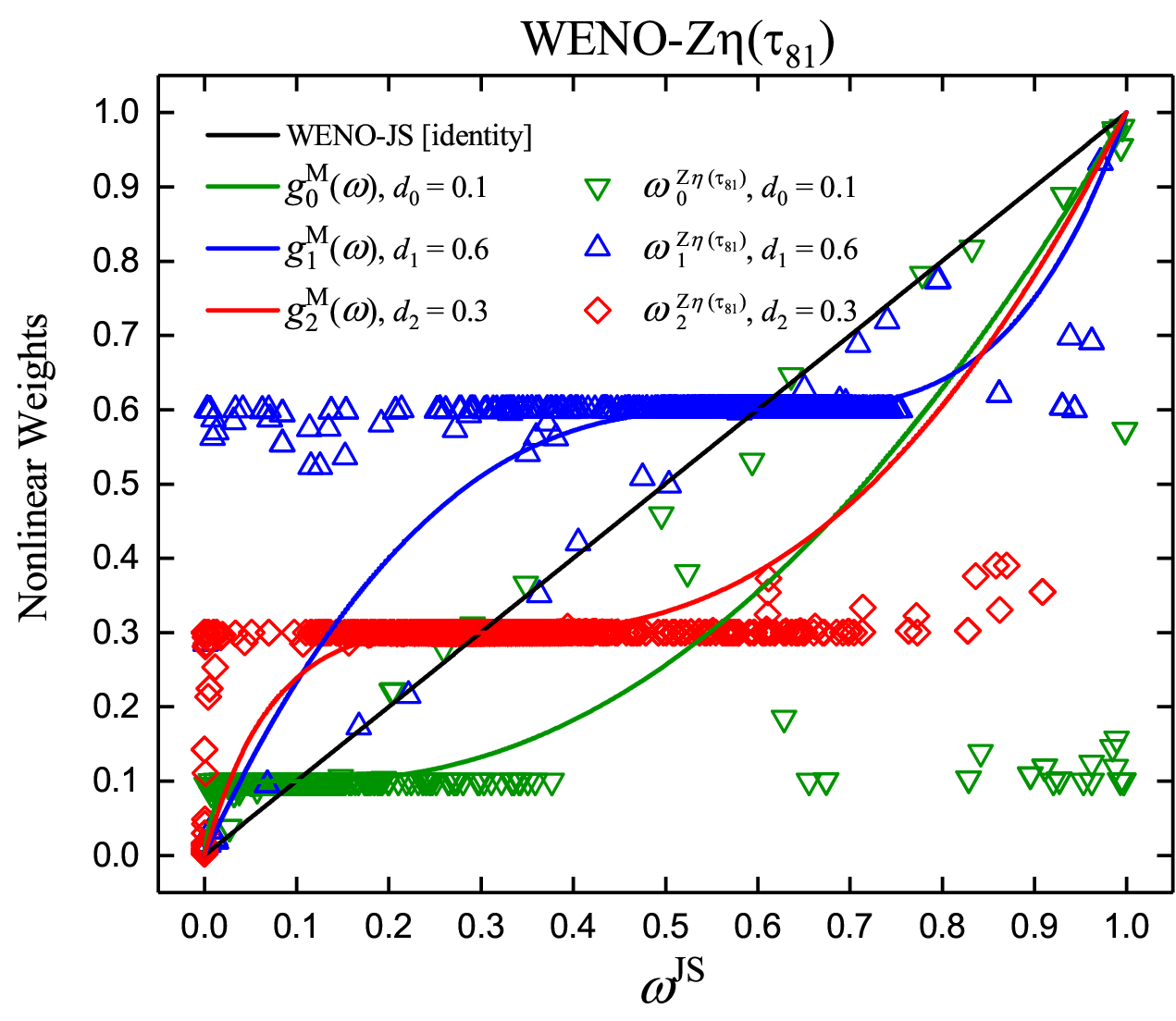}%\\
  \includegraphics[height=0.28\textwidth]
  {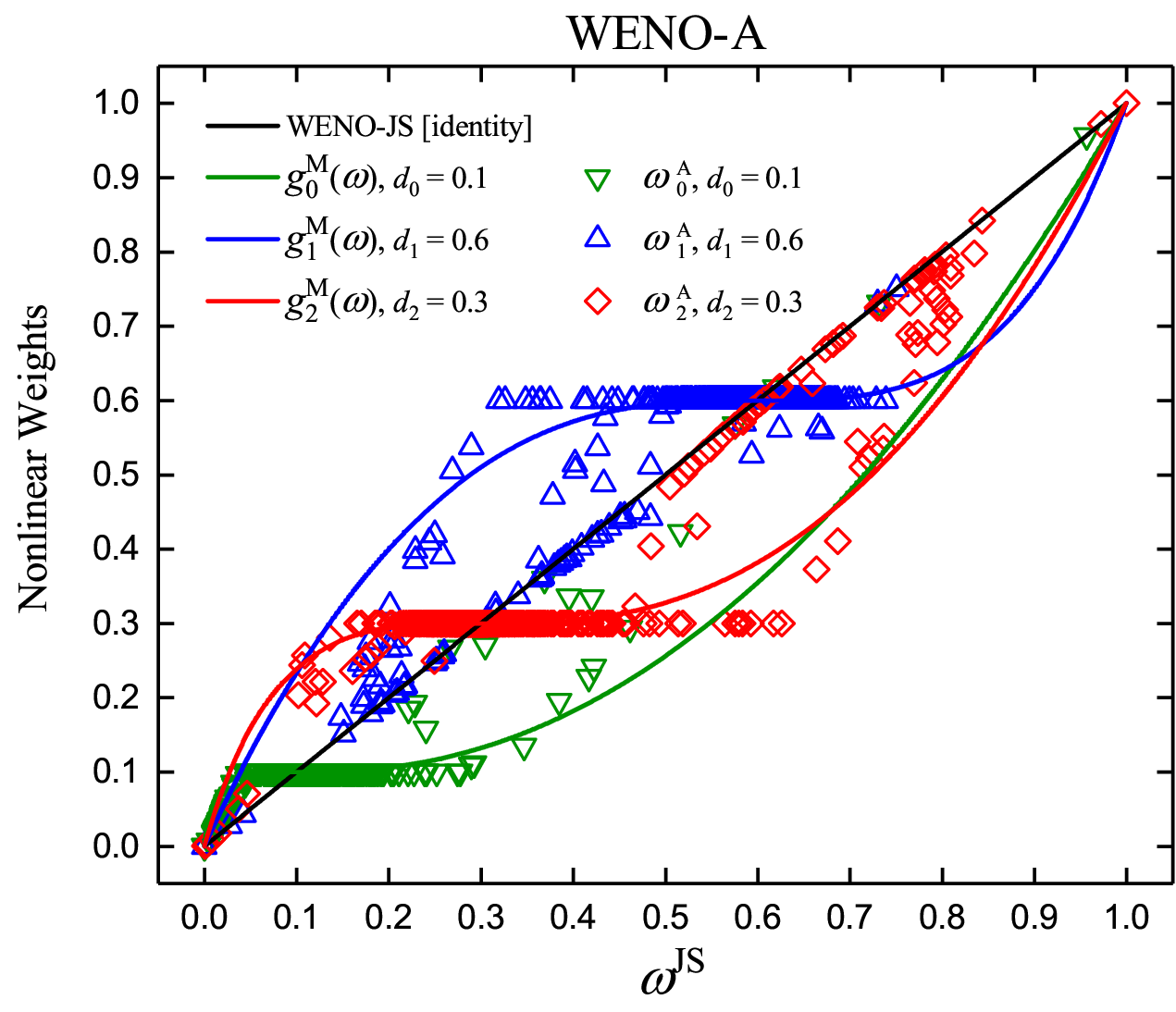}
\caption{The IMRs for the WENO-Z, WENO-Z$\eta(\tau_{81})$ and WENO-A 
schemes.}
\label{fig:IMRs:Z-type}
\end{figure}

Because the present study is primarily concerned with the 
performance of long-run calculations of the WENO-Z-type schemes, as 
examples, we plot the IMRs of the WENO-Z, WENO-Z$\eta(\tau_{81})$ 
and WENO-A schemes in Fig. \ref{fig:IMRs:Z-type} on solving the 
following initial-value problem for the linear advection equation at 
$t = 200$
\begin{equation*}
\begin{array}{ll}
\partial_{t} u+ \partial_{x} u = 0, & (x,t)\in[-1,1]\times (0,200],\\
u(x,0) = u_{0}(x), & x \in [-1, 1],
\end{array}
\label{eq:LAE1:Z-type}
\end{equation*}
where
\begin{equation*}
u_{0}(x) = \left\{
\begin{array}{ll}
1, & x \in [-1, 0], \\
0, & x \in (0,1].
\end{array}
\right.
\end{equation*}
To compare these IMRs with the traditional mappings of classical 
mapped WENO schemes, the designed mapping profiles of WENO-M, as 
well as the identity mapping of WENO-JS, are also plotted in Fig. \ref{fig:IMRs:Z-type}. Surprisingly, it is observed that there is an 
apparent similarity between the IMRs and the traditional mappings, 
that is, they both embrace the ``optimal weight interval'' that 
stands for the interval about $\omega_{k} = d_{k}$ over which the 
Z-type weight formulas or the explicit mapping functions attempt to 
use the corresponding optimal weight. It is widely reported that \cite{WENO-IM,WENO-AIM,WENO-MAIMi} the optimal weight interval is 
important for achieving the designed orders of convergence near the 
critical points and getting high resolutions near discontinuities. 
In spite of these advantages, the latest studies \cite{MOP-WENO-ACMk,MOP-WENO-X} have indicated that the optimal weight interval is 
harmful to preserving high resolutions and meanwhile avoiding 
spurious oscillations for long-run simulations with discontinuities, 
as it result in the increasing of the nonlinear weights for 
non-smooth stencils as well as the decreasing of the nonlinear 
weights for smooth stencils. So far, however, there has been little 
discussion about this topic for WENO-Z-type schemes. And thus, we 
set out to investigate it carefully in this study.

\subsection{The new WENO-Z-type schemes with locally 
order-preserving mappings}\label{subsec:Z-type_LOP}
Inspired by the observation above, we rewrite the Z-type weights in 
a general and meaningful form as follows
\begin{equation}
\omega_{k}^{\mathrm{X}} = \dfrac{\alpha_{k}^{\mathrm{X}}}{\sum_{j = 0}^{2}\alpha_{j}^{\mathrm{X}}}, \quad \alpha_{k}^{\mathrm{X}} = \Big( g^{\mathrm{GMX}} \Big)_{k}(\omega_{k}^{\mathrm{JS}}) = \lambda_{1}^{\mathrm{X}}(d_{k}) + \lambda_{2}^{\mathrm{X}}(d_{k})\cdot\omega_{k}^{\mathrm{JS}}, \quad k = 0,1,2.
\label{eq:general:Z-typeWeights}
\end{equation}
According to the specified Z-type weights, we can determine 
$\lambda_{1}^{\mathrm{X}}(d_{k})$ and 
$\lambda_{2}^{\mathrm{X}}(d_{k})$ easily. Considering the WENO-Z, 
WENO-Z$\eta(\tau_{81})$ and WENO-A schemes, we present their $\lambda_{1}^{\mathrm{X}}(d_{k})$ and $\lambda_{2}^{\mathrm{X}}(d_{k})$ in 
Table \ref{tab:lambda:Z-type}. We call $\Big( g^{\mathrm{GMX}} \Big)_{k}(\omega_{k}^{\mathrm{JS}})$ the \textit{generalized mapping function} of WENO-X.

\begin{table}[ht]
\renewcommand\arraystretch{1.5}
\footnotesize
\centering
\caption{$\lambda_{1}^{\mathrm{X}}(d_{k})$, 
$\lambda_{2}^{\mathrm{X}}(d_{k})$ for different Z-type weights.}
\label{tab:lambda:Z-type}
%\begin{tabular*}{lll}
%\begin{spacing}{1.25}
\begin{tabular*}{\hsize}
{@{}@{\extracolsep{\fill}}lll@{}}
\hline
WENO-X         & $\lambda_{1}^{\mathrm{X}}(d_{k})$ 
               & $\lambda_{2}^{\mathrm{X}}(d_{k})$\\
\hline
WENO-Z    & $d_{k}$ 
          & $\sum_{l=0}^{2}\alpha_{l}^{\mathrm{JS}}\cdot\tau_{5}^{2}$\\
WENO-Z$\eta(\tau_{81})$   
          & $d_{k}$ 
          & $\sum_{l=0}^{2}\alpha_{l}^{\mathrm{JS}}\cdot\dfrac{(
            IS_{k}^{\mathrm{JS}} + \epsilon)^{2}}{(\eta_{k} + 
            \epsilon)^{2}}\cdot\tau_{81}^{2}$ \\
WENO-A    & $d_{k}\cdot B$ 
          & $\sum_{l=0}^{2}\alpha_{l}^{\mathrm{JS}}\cdot 
          \Phi\tau_{5}^{p}\cdot(IS_{k}^{\mathrm{JS}} + 
          \epsilon)^{2-p}(1 - B).$\\
\hline
\multicolumn{3}{l}{\tabincell{l}{Here, $B = \mathrm{BOOL}\Bigg(\max
\bigg(1, \Phi\bigg(\dfrac{\tau_{5}}{IS_{k} + \epsilon}\bigg)^{p}
\bigg) = 1 \Bigg)$ and $\Phi = \min\Big\{ 1,\sqrt{\lvert IS_{0} - 
2IS_{1} + IS_{2} \rvert} \Big\}$.}} \\ 
\hline
\end{tabular*}
%\end{spacing}
\end{table}

In order to clarify our concerns and simplify the description, we 
simply state the following definitions.

\begin{definition}
Let $S^{5}$ denote the $5$-point global stencil centered around any 
location $x_{j}$, and $\omega_{0}, \omega_{1}, \omega_{2}$ denote 
the nonlinear weights with respect to the three $3$-point 
substencils of $S^{5}$, that is, $S_{0}, S_{1}, S_{2}$. Assume that
$\big(g^{\mathrm{GMX}}\big)_{k}(\omega), s=0, 1, 2$ is the 
generalized mapping function of the Z-type WENO-X scheme. If for any
$a, b \in \{0, 1, 2\}$, $\big( g^{\mathrm{GMX}}\big)_{a}(\omega_{a}) \geq \big(g^{\mathrm{GMX}} \big)_{a}(\omega_{a})$ holds true when $\omega_{a} > \omega_{b}$, and $\big( g^{\mathrm{GMX}}\big)_{a}(\omega_{a}) = \big(g^{\mathrm{GMX}} \big)_{b}(\omega_{b})$ holds true when $\omega_{a}=\omega_{b}$, then, the set of the 
generalized mappings \Big\{$\big( g^{\mathrm{GMX}}\big)_{k}(\omega), s=0, 1, 2$\Big\} is called \textbf{locally order-preserving (LOP)}.
\label{def:LOPM}
\end{definition}

\begin{definition} $x_{j}$ is a \textbf{non-OP point}, if one has
$m, n \in \{0, 1, 2\}$ and
\begin{equation}\left\{
\begin{array}{ll}
\begin{aligned}
&\big(\omega_{m}-\omega_{n}\big)\bigg(\big(g^{\mathrm{GMX}}\big)_{m}
(\omega_{m}) - \big(g^{\mathrm{GMX}}\big)_{n}(\omega_{n})\bigg) < 0, 
&\mathrm{if} \quad \omega_{m} \neq \omega_{n},\\
&\big(g^{\mathrm{GMX}}\big)_{m}(\omega_{m})\neq\big(g^{\mathrm{GMX}}
\big)_{n}(\omega_{n}), &\mathrm{if} \quad \omega_{m}=\omega_{n}.
\end{aligned}
\end{array}\right.
\end{equation}
Otherwise, $x_{j}$ is an \textbf{OP point}.
\label{def:OP-point}
\end{definition}

\begin{definition} Let any $a, b \in \{0, 1, 2\}$, then, the set of 
function $\mathbb{S}^{\mathrm{X}}$ is defined by
\begin{equation}
\begin{aligned}
\mathbb{S}^{\mathrm{X}}= \quad &\bigg\{\mathbf{LOP\_idx}(a,b,
\mathrm{X}):\mathbf{LOP\_idx}(a, b, \mathrm{X}) > 0\bigg\} \\
\bigcup &\bigg\{ 
\mathbf{LOP\_idx}(a, b, \mathrm{X}): \omega_{a}^{\mathrm{JS}} - 
\omega_{b}^{\mathrm{JS}} = g_{a}^{\mathrm{GMX}}\Big(\omega_{a}^{
\mathrm{JS}}\Big) - g_{b}^{\mathrm{GMX}}\Big(\omega_{b}^{\mathrm{JS}}
\Big) = 0 \bigg\},
\end{aligned}
\label{eq:postINDEX:SET}
\end{equation}
where
\begin{equation}
\mathbf{LOP\_idx}(a, b, \mathrm{X}) = 
\Big( \omega_{a}^{\mathrm{JS}} - \omega_{b}^{\mathrm{JS}} 
\Big)\Bigg( g_{a}^{\mathrm{GMX}}\Big(\omega_{a}^{\mathrm{JS}}\Big) - 
g_{b}^{\mathrm{GMX}}\Big(\omega_{b}^{\mathrm{JS}}\Big) \Bigg).
\label{eq:postINDEX}
\end{equation}
\label{def:postINDEX:SET}
\end{definition}

Trivially, according to Definition \ref{def:OP-point} and Definition
\ref{def:postINDEX:SET}, we have the following property.

\begin{lemma}
At $x_{j}$, for any $a, b = 0, 1, 2$ with $a \neq b$, if {\rm{\textbf{LOP\_idx}}}$(a, b, \mathrm{X}) \in \mathbb{S}^{\mathrm{X}}$, then $x_{j}$ is an \textit{OP} point to the WENO-X scheme. In other words, 
the set of mapping functions $\Big\{\big(g^{\mathrm{LOP-GMX}}\big)_{k}(\omega^{\mathrm{JS}}_{k}), s=0, 1, 2\Big\}$ is \textit{LOP} at $x_{j}$. Otherwise, if one has $a, b = 0, 1, 2$ with $a \neq b$, and
{\rm{\textbf{LOP\_idx}}}$(a, b, \mathrm{X}) \notin \mathbb{S}^{
\mathrm{X}}$, then $x_{j}$ is a \textit{non-OP} point to the WENO-X 
scheme.
\label{lem:postINDEX}
\end{lemma}

We devise a general method to build the improved WENO-Z-type schemes 
satisfying the \textit{LOP} generalized mappings by utilizing the 
\textbf{LOP\_idx} function given in Eq. \eqref{eq:postINDEX}, as 
shown in Algorithm \ref{alg:posteriori}.

\begin{algorithm}[htb]
\caption{A general method to construct \textit{LOP} mappings.}
\label{alg:posteriori}
%\SetKwData{Left}{left}\SetKwData{This}{this}\SetKwData{Up}{up}
\SetKwInOut{Input}{input}\SetKwInOut{Output}{output}
\Input{$s$, $d_{k}$, $\alpha^{\mathrm{JS}}_{k}$, 
$\omega^{\mathrm{JS}}_{k}$, $\big(g^{\mathrm{GMX}}\big)_{k}(\omega)$
}
\Output{$\Big\{\big(g^{\mathrm{LOP-GMX}}\big)_{k}(\omega^{
\mathrm{JS}}_{k}), s=0,1,\cdots,r-1\Big\}$, the new set of mapping 
functions that is \textit{LOP}}
\BlankLine
\emph{$\big( g^{\mathrm{GMX}} \big)_{k}(\omega)$ is a generalized 
mapping over $[0, 1]$, and the set of mapping functions 
$\Big\{\big(g^{\mathrm{GMX}} \big)_{k}(\omega), s = 0, 1,\cdots,r-1 
\Big\}$ is \textit{non-OP}}\;
set $s_{1}=0$\;
\While{$s_{1} \leq r - 2$}{
  set $s_{2} = s_{1} + 1$\;
  \While{$s_{2} \leq r - 1$}{
     $\kappa=\mathrm{\textbf{LOP\_idx}}(s_{1},s_{2},\mathrm{X})$\;
     \eIf{$\kappa \in \mathbb{S}^{\mathrm{X}}$}{
       $\lambda = 1$\;
     }{
       $\lambda = 0$\;
       \textbf{Break}\;
     } 
     $s_{2} ++$\; 
  }
  \If{$\lambda = 0$}{
    \textbf{Break}\;
  }
  $s_{1} ++$\;
}
set $s=0$\;
\While{$s \leq r - 1$}{
  \eIf{$\lambda = 1$}{
    $\big(g^{\mathrm{LOP-GMX}}\big)_{k}(\omega^{\mathrm{JS}}_{k})
    = \big(g^{\mathrm{GMX}}\big)_{k}(\omega^{\mathrm{JS}}_{k})$\;
  }{
    $\big(g^{\mathrm{LOP-GMX}}\big)_{k}(\omega^{\mathrm{JS}}_{k})
    = \alpha_{k}^{\mathrm{JS}}$. 
    \tcp{$\alpha_{k}^{\mathrm{JS}}$ is computed by 
    Eq.\eqref{eq:weights:WENO-JS}}
  }
  $s ++$\;
}
\end{algorithm}

\begin{theorem}
The set of mapping functions 
$\Big\{\big(g^{\mathrm{LOP-GMX}}\big)_{k}(\omega^{\mathrm{JS}}_{k})
, s=0, 1, 2\Big\}$ obtained through Algorithm 
\ref{alg:posteriori} is \textit{LOP}.
\label{theorem:g_MOP}
\end{theorem}
\textbf{Proof.} The proof can be divided into two cases:
$\lambda = 0$ and $\lambda = 1$. (1) We first prove the case of 
$\lambda = 0$. From a mapping perspective, the nonlinear weights of 
WENO-JS can be seen to be computed by an identity mapping, say, 
$\big( g^{\mathrm{JS}} \big)_{k}(\omega) = \omega,s= 0, 1, 2$. 
Following this treatment and according to Definition \ref{def:LOPM}, 
we can trivially prove that the set of mapping functions 
$\Big\{\big( g^{\mathrm{JS}} \big)_{k}(\omega), s=0, 1, 2\Big\}$ is 
\textit{LOP} with the widths of the optimal weight intervals to be 
zero. Then, according to Line 25 of Algorithm \ref{alg:posteriori}, we complete the proof of the case of $\lambda = 0$. (2) For the case 
of $\lambda = 1$, as $\kappa \in \mathbb{S}^{\mathrm{X}}$ (see Line 
6 of Algorithm \ref{alg:posteriori}), then, according to Lemma 
\ref{lem:postINDEX} and Line 23 of Algorithm \ref{alg:posteriori}, 
it is easy to obtain that $\Big\{\big(g^{\mathrm{LOP-GMX}}\big)_{k}(\omega^{\mathrm{JS}}_{k}), s=0, 1, 2\Big\}$ is \textit{LOP}.
$\hfill\square$ \\

Now, the improved Z-type weights satisfying the \textit{LOP} 
generalized mappings are computed by
\begin{equation}
\omega_{k}^{\mathrm{LOP-GMX}}=\dfrac{\alpha_{k}^{\mathrm{LOP-GMX}}}{\sum_{j=0}^{2}\alpha_{j}^{\mathrm{LOP-GMX}}}, \quad \alpha_{k}^{
\mathrm{LOP-GMX}} = \big(g^{\mathrm{LOP-GMX}}\big)_{k}(\omega^{
\mathrm{JS}}_{k}), \quad s = 0, 1, 2,
\label{eq:PoAOP-mapping}
\end{equation}
where $\big(g^{\mathrm{LOP-GMX}}\big)_{k}(\omega^{\mathrm{JS}}_{k})$ 
is given by Algorithm \ref{alg:posteriori}. The resultant scheme 
will be referred to as LOP-GMWENO-X.

\begin{figure}[ht]
\flushleft
  \includegraphics[height=0.28\textwidth]
  {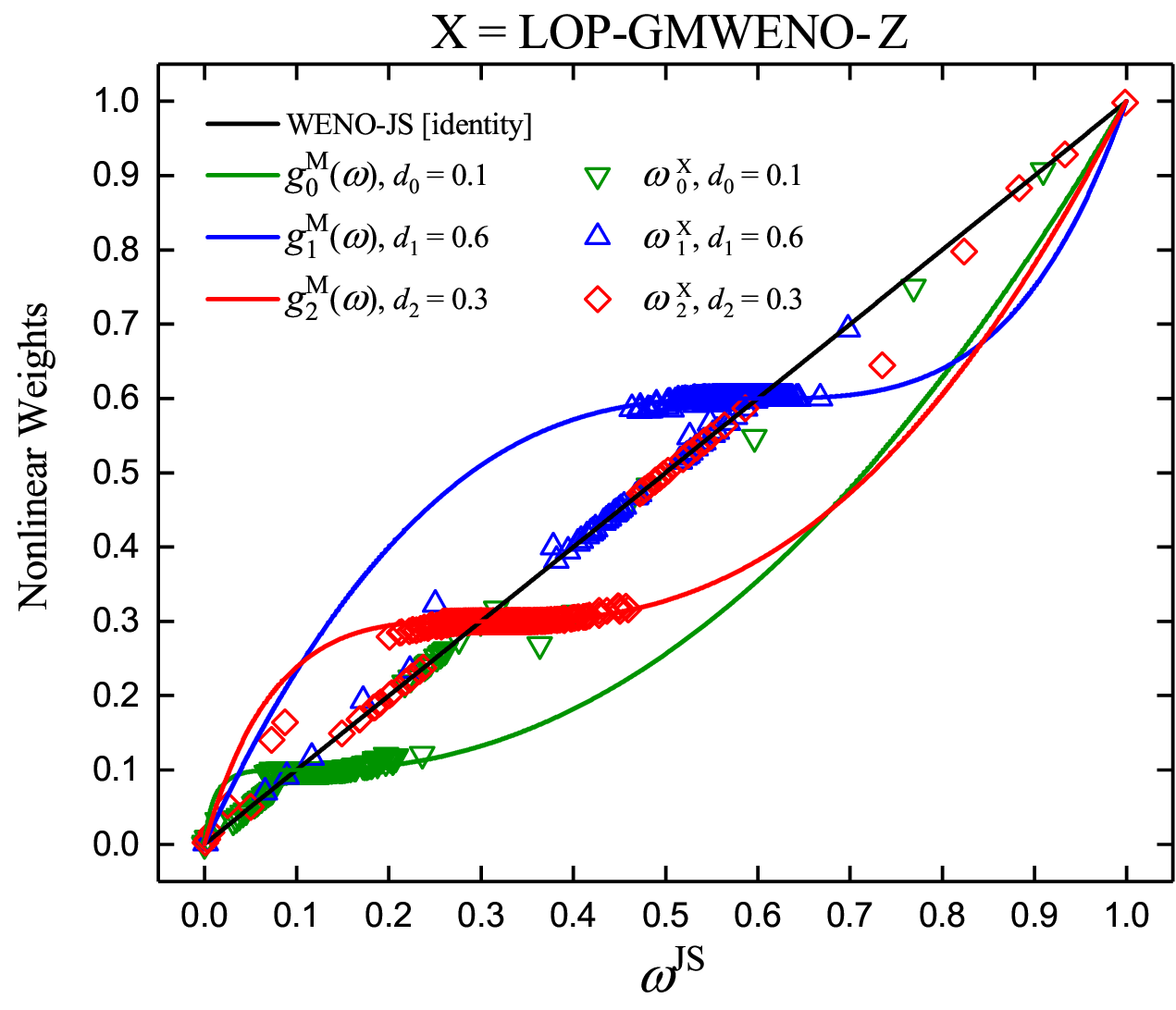}
  \includegraphics[height=0.28\textwidth]
  {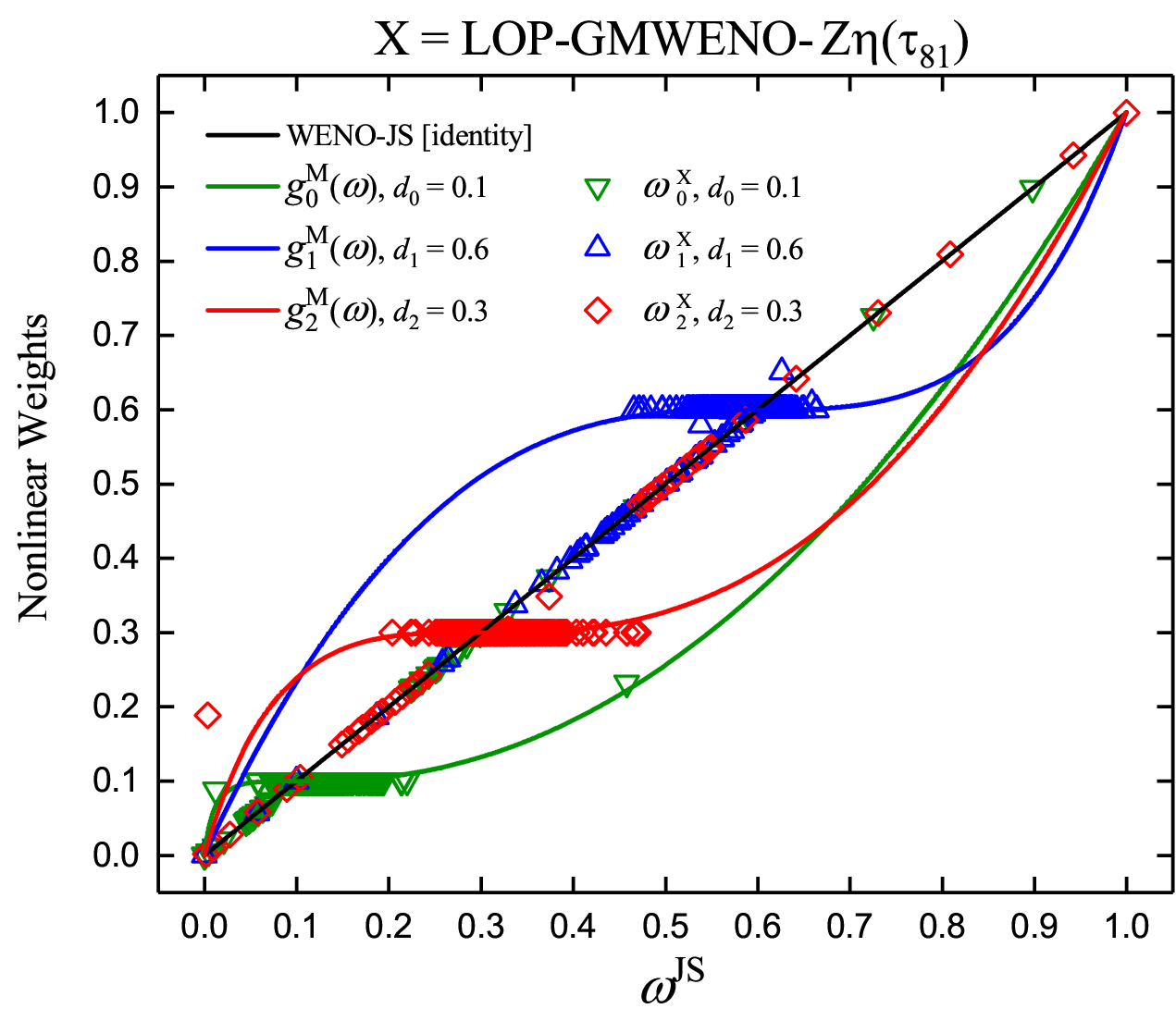}%\\
  \includegraphics[height=0.28\textwidth]
  {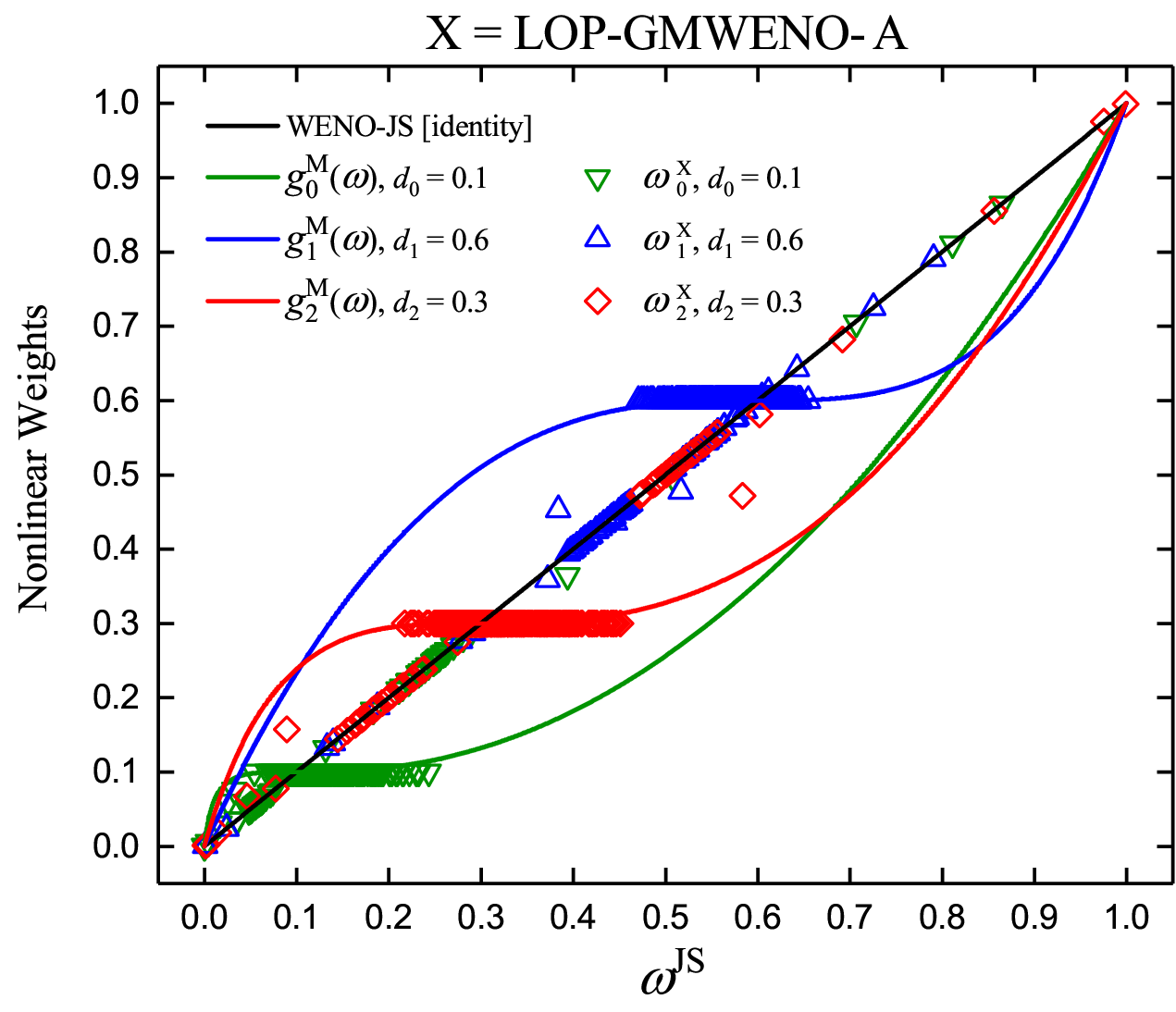}
\caption{The IMRs for the LOP-GMWENO-Z, LOP-GMWENO-Z$\eta(\tau_{81})$ and LOP-GMWENO-A schemes.}
\label{fig:IMR:GMZ}
\end{figure}

Similarly, we plot the IMRs of the LOP-GMWENO-Z, LOP-GMWENO-Z$\eta(\tau_{81})$ and LOP-GMWENO-A schemes in Fig. \ref{fig:IMR:GMZ} using 
exactly the same computing conditions as in Fig. \ref{fig:IMRs:Z-type}. It is clear that the LOP-GMWENO-X schemes also embrace apparent 
optimal weight intervals as their associated WENO-X schemes do. 
Moreover, we can intuitively observe a noticeable difference that 
many nonlinear weights of the LOP-GMWENO-X schemes drop onto the 
identity mappings. However, with the exception of WENO-A, this can 
not be observed for all other corresponding WENO-X schemes from 
Fig. \ref{fig:IMRs:Z-type}. According to Algorithm \ref{alg:posteriori}, the nonlinear weights dropping onto the identity 
mappings represent the non-OP points. For WENO-A, although some of 
its nonlinear weights also drop onto the identity mappings, no 
evidence was found for the fact that its IMR is LOP. Indeed, from 
Fig. \ref{fig:IMRs:Z-type} and Fig. \ref{fig:IMR:GMZ}, we can see 
that there are much less nonlinear weights of WENO-A dropping 
onto the identity mappings than those of LOP-GMWENO-A, and this 
may indirectly indicate that WENO-A can not identify the full non-OP 
points. In addition, more evidence will be presented numerically in 
the rest of this paper to support that the IMR of WENO-A is not LOP.

In order to examine the convergence properties of the LOP-GMWENO-X 
schemes, we perform a typical benchmark numerical test below. We 
will see that the LOP-GMWENO-X schemes can achieve almost exactly 
the same convergence properties as those of their associated WENO-X 
schemes.

\subsection{Convergence at critical points}\label{subsec:Critical}
There has been a growing number of publications 
\cite{WENO-IM,WENO-PM,WENO-RM260,WENO-PPM5,WENO-MAIMi,WENO-ACM,
MOP-WENO-ACMk,MOP-WENO-X,WENO-Z} focusing on the convergence at 
critical points since in \cite{WENO-M} it was pointed out that the 
fifth-order WENO-JS scheme suffers from the loss of accuracy and 
achieves only third-order convergence rate of accuracy at critical 
points of smooth solutions. 

As used in \cite{WENO-eta,WENO-eta-02,WENO-D_WENO-A}, the following 
test function is also considered in the present study to measure the 
accuracy of the LOP-GMWENO-X schemes
\begin{equation}
f(x) = x^{l}\exp(x), \quad x \in (-1, 1).
\label{eq:Critical:OnePoint}
\end{equation}
It is easy to check that Eq. \eqref{eq:Critical:OnePoint} has a 
critical point of order $n_{\mathrm{cp}} = l - 1$ at $x = 0$ where 
$f^{(j)}(0) = 0, j = 0, 1, \cdots, l-1$.

Numerical experiments on solving the test function above by setting 
$l = 2$ have been performed to compare the behaviors of the 
LOP-GMWENO-X schemes and their associated WENO-X schemes. For the 
purpose of comparison, we also presented the results of the the 
WENO-JS scheme and the WENO5 scheme using ideal linear weights 
(denoted as WENO5-ILW in this paper for brevity). 

Table \ref{tab:Convergence:Critical} shows the $L_{\infty}$ 
convergence behaviors for the considered schemes at the critical 
point $x = 0$. In terms of convergence rate of accuracy, we can see 
that: (1) the WENO5-ILW scheme has sixth-order accuracy, while as 
expected, the WENO-JS scheme only gets about third-order accuracy;
(2) the WENO-Z scheme and the associated LOP-GMWENO-Z scheme are 
both able to achieve about fifth-order accuracy; (3) the other two 
considered WENO-Z-type schemes, say, WENO-Z$\eta(\tau_{81})$ and 
WENO-A, and their associated LOP-GMWENO-Z$\eta(\tau_{81})$ and 
LOP-GMWENO-A schemes, have about sixth-order accuracy.

Furthermore, in terms of accuracy, we can observe that: (1) the 
WENO-JS scheme obtain the errors of 6 to 10 orders of magnitude 
larger than those of the the WENO-ILW scheme; (2) the WENO-Z scheme 
and the associated LOP-GMWENO-Z scheme get the errors of 3 to 4 
orders of magnitude larger than those of the the WENO-ILW scheme; 
(4) the other two considered WENO-Z-type schemes and their 
associated LOP-GMWENO-X schemes can obtain the errors with the same 
order of magnitude as those of the WENO-ILW scheme. 

In summary, it has been shown from this test that all the 
LOP-GMWENO-X schemes can achieve similar numerical errors and 
convergence orders as those of their associated WENO-X schemes.

\begin{table}[ht]
%\begin{footnotesize}
%\begin{scriptsize}
\begin{myFontSize}
\centering
\caption{Convergence rate of accuracy at the critical point.}
\label{tab:Convergence:Critical}
\begin{tabular*}{\hsize}
{@{}@{\extracolsep{\fill}}lllllllll@{}}
\hline
\space    &\multicolumn{2}{l}{\cellcolor{gray!35}{WENO5-ILW}}  
          &\multicolumn{2}{l}{\cellcolor{gray!35}{WENO-JS}}  
          &\multicolumn{2}{l}{\cellcolor{gray!35}{WENO-Z}}  
          &\multicolumn{2}{l}{\cellcolor{gray!35}{LOP-GMWENO-Z}}\\
\cline{2-3}  \cline{4-5}   \cline{6-7}  \cline{8-9}
$\Delta x$            & $L_{\infty}$ error & $L_{\infty}$ Order 
                      & $L_{\infty}$ error & $L_{\infty}$ Order 
              & $L_{\infty}$ error & $L_{\infty}$ Order
                & $L_{\infty}$ error & $L_{\infty}$ Order\\
\Xhline{0.65pt}
0.01                  & 7.01182E-13        & -
                      & 8.25085E-07        & -
                      & 1.58511E-09        & -
                      & 1.58511E-09        & -\\
0.005                 & 1.19831E-14        & 5.8707 
                      & 8.74036E-08        & 3.2388 
                      & 3.75866E-11        & 5.3982 
                      & 3.75866E-11        & 5.3982 \\
0.0025                & 1.84002E-16        & 6.0251 
                      & 9.20947E-09        & 3.2465 
                      & 8.72344E-13        & 5.4292 
                      & 8.72344E-13        & 5.4292 \\
0.00125               & 2.84843E-18        & 6.0134 
                      & 9.65225E-10        & 3.2542 
                      & 1.94104E-14        & 5.4496 
                      & 1.94104E-14        & 5.4496\\
0.000625              & 4.42179E-20        & 6.0094 
                      & 1.00658E-10        & 3.2614 
                      & 4.54841E-16        & 5.4583 
                      & 4.54841E-16        & 5.4583 \\
\hline
\space    &\multicolumn{2}{l}{\cellcolor{gray!35}{WENO-Z$\eta(\tau_{81})$}}  
          &\multicolumn{2}{l}{\cellcolor{gray!35}{LOP-GMWENO-Z$\eta(\tau_{81})$}} 
          &\multicolumn{2}{l}{\cellcolor{gray!35}{WENO-A}}  
          &\multicolumn{2}{l}{\cellcolor{gray!35}{LOP-GMWENO-A}}\\
\cline{2-3}  \cline{4-5}   \cline{6-7}  \cline{8-9}
$\Delta x$            & $L_{\infty}$ error & $L_{\infty}$ Order 
                      & $L_{\infty}$ error & $L_{\infty}$ Order 
              & $L_{\infty}$ error & $L_{\infty}$ Order
                & $L_{\infty}$ error & $L_{\infty}$ Order\\
\Xhline{0.65pt}
0.01                  & 7.77097E-13        & -
                      & 7.87300E-13        & -
                      & 7.87184E-13        & -
                      & 7.88374E-13        & -\\
0.005                 & 1.19831E-14        & 6.0190 
                      & 1.19831E-14        & 6.0378 
                      & 1.19831E-14        & 6.0376 
                      & 1.19831E-14        & 6.0398\\
0.0025                & 1.84002E-16        & 6.0251 
                      & 1.84002E-16        & 6.0251 
                      & 1.84002E-16        & 6.0251 
                      & 1.84002E-16        & 6.0251\\
0.00125               & 2.84843E-18        & 6.0134 
                      & 2.84843E-18        & 6.0134 
                      & 2.84843E-18        & 6.0134 
                      & 2.84843E-18        & 6.0134\\
0.000625              & 4.42179E-20        & 6.0094 
                      & 4.42179E-20        & 6.0094 
                      & 4.42179E-20        & 6.0094 
                      & 4.42179E-20        & 6.0094\\
\hline
\end{tabular*}
\end{myFontSize}
%\end{scriptsize}
%\end{footnotesize}
\end{table}

\subsection{Comparison on solving 1D linear advection equation with long output times}\label{subsec:Long-run}
\subsubsection{With discontinuities}
The most excellent performance of the proposed LOP-GMWENO-X schemes 
is that they can preserve high resolutions and in the meantime 
prevent spurious oscillations for large output times. To manifest 
this, we solve the 1D linear advection equation $u_{t} + u_{x} = 0$, 
$x \in (-1, 1)$ with the following two different initial conditions 
here.

Case 1. The initial condition is given by
\begin{equation}
\begin{array}{l}
u(x, 0) = \left\{
\begin{array}{ll}
\dfrac{1}{6}\big[ G(x, \beta, z - \hat{\delta}) + 4G(x, \beta, z) + G
(x, \beta, z + \hat{\delta}) \big], & x \in [-0.8, -0.6], \\
1, & x \in [-0.4, -0.2], \\
1 - \big\lvert 10(x - 0.1) \big\rvert, & x \in [0.0, 0.2], \\
\dfrac{1}{6}\big[ F(x, \alpha, a - \hat{\delta}) + 4F(x, \alpha, a) +
 F(x, \alpha, a + \hat{\delta}) \big], & x \in [0.4, 0.6], \\
0, & \mathrm{otherwise},
\end{array}\right. 
\end{array}
\label{eq:LAE:SLP}
\end{equation}
where $G(x, \beta, z) = \mathrm{e}^{-\beta (x - z)^{2}}, F(x, \alpha
, a) = \sqrt{\max \big(1 - \alpha ^{2}(x - a)^{2}, 0 \big)}$, and 
the constants are $z = -0.7, \hat{\delta} = 0.005, \beta = \dfrac{
\log 2}{36\hat{\delta} ^{2}}, a = 0.5$ and $\alpha = 10$. It is 
known that this problem consists of a Gaussian, a square wave, a 
sharp triangle and a semi-ellipse. 

Case 2. The initial condition is given by
\begin{equation}
\begin{array}{l}
u(x, 0) = \left\{
\begin{array}{ll}
0,   & x \in [-1.0, -0.8] \bigcup (-0.2, 0.2] \bigcup (0.8, 1.0], \\
0.5, & x \in (-0.6, -0.4] \bigcup (0.2, 0.4]  \bigcup (0.6, 0.8], \\
1,   & x \in (-0.8, -0.6] \bigcup (-0.4, -0.2] \bigcup (0.4, 0.6].
\end{array}\right. 
\end{array}
\label{eq:LAE:BiCWP}
\end{equation}

Case 1 and Case 2, dubbed SLP and BiCWP respectively, were carefully 
studied in \cite{MOP-WENO-ACMk}. Indeed, these two cases are widely 
used to test the performance of the WENO schemes for long output 
times, and we refer to \cite{WENO-IM,WENO-RM260,WENO-RM-Vevek2018,
WENO-AIM,WENO-MAIMi,MOP-WENO-ACMk,MOP-WENO-X} for more details.

To examine the convergence properties, both Case 1 and Case 2 are 
computed to the final time of $t = 2000$. The LOP-GMWENO-X schemes 
and their associated WENO-X schemes, as well as the WENO-JS and 
WENO-ILW schemes are considered. The CFL number is set to be $0.1$ 
here just for the purpose of keeping consistent with the articles 
\cite{WENO-PM,WENO-IM,WENO-MAIMi,MOP-WENO-ACMk,MOP-WENO-X} in which 
the long-run simulation of Case 1 was widely studied. Actually, we 
have also verified numerically that our conclusion is still 
well-supported even for a larger CFL number, such as 0.5.

The $L_{1}$-, $L_{2}$- and $L_{\infty}$- norms of numerical errors 
are computed by
\begin{equation}
\displaystyle
\begin{aligned}
%\begin{array}{l}
L_{1} = h \cdot \displaystyle\sum\nolimits_{i=1}^{N}\Big\lvert 
u_{i}^{\mathrm{exact}} - (u_{h})_{i} \Big\rvert, \quad 
L_{2} = \sqrt{\displaystyle\sum\nolimits_{i=1}^{N} h \cdot 
\Big(u_{i}^{\mathrm{exact}} - (u_{h})_{i} \Big)^{2}}, \quad
L_{\infty} = \max_{1\leq i \leq N}\lvert u_{i}^{\mathrm{exact}} - (u_{h})_{i} \rvert,
%\end{array}
\end{aligned}
\label{normsDefinition}
\end{equation}
where $N$ is the number of the cells and $h$ is the associated 
uniform spatial step size. $(u_{h})_{i}$ is the numerical solution 
and $u_{i}^{\mathrm{exact}}$ is the exact solution. 

Here, we present the $L_{1}$- and $L_{\infty}$-norm errors and 
orders of convergence with $N = 100, 200, 400, 800$, as shown in 
Table \ref{table:AccuracyTest:SLP:t2000} and Table \ref{table:AccuracyTest:BiCWP:t2000} for Case 1 and Case 2 respectively. 
Clearly, in terms of accuracy, for both Case 1 and Case 2, we can 
see that: (1) WENO-ILW outperforms all other schemes, while WENO-JS 
has the largest numerical errors; (2) the numerical errors of the 
LOP-GMWENO-X schemes are comparable to those of the associated 
WENO-X schemes in general; (3) closer inspection of the table shows 
that LOP-GMWENO-Z and LOP-GMWENO-A schemes perform better than the 
associated WENO-Z and WENO-A schemes for $N \geq 400$. In addition, 
the LOP-GMWENO-X schemes have significantly larger $L_{1}$-norm 
convergence orders of accuracy compared to the associated WENO-X 
schemes.

To better demonstrate the improvement of the proposed schemes, 
without loss of generality, we re-calculate both Case 1 and Case 2 
with a final time of $t = 200$ by using $N = 1600$ and $N = 3200$ 
respectively. For comparison purpose, we also compute Case 1 and 
Case 2 using WENO-JS and WENO-M.

\begin{figure}[!ht]
\centering
\includegraphics[height=0.33\textwidth]
{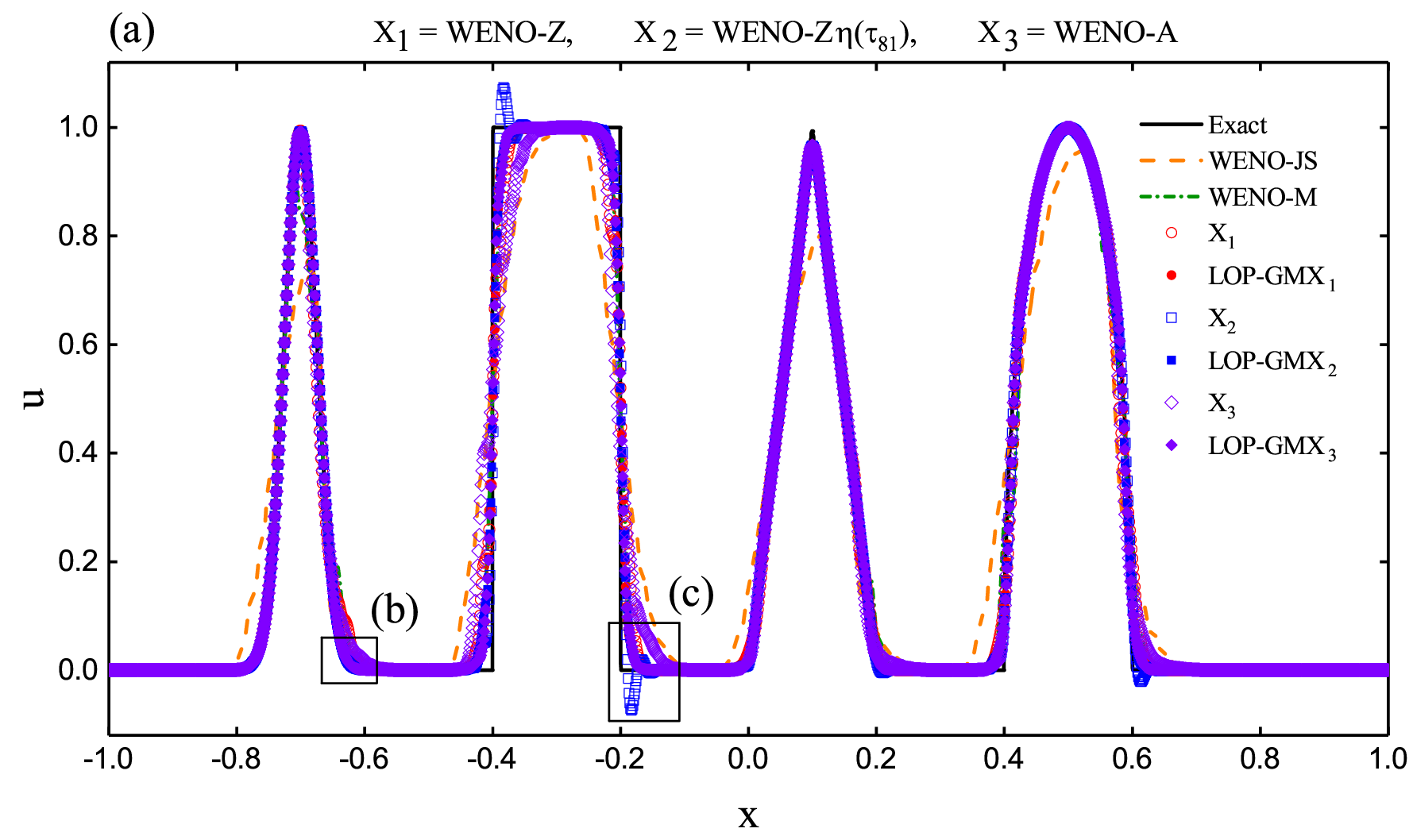}
\includegraphics[height=0.33\textwidth]
{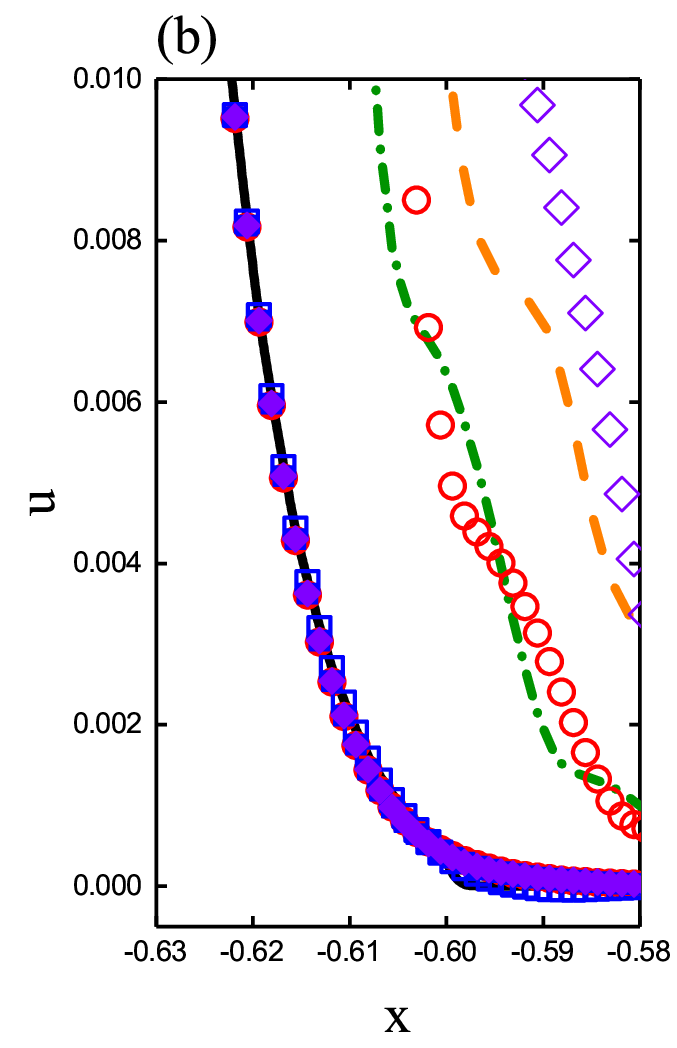}
\includegraphics[height=0.33\textwidth]
{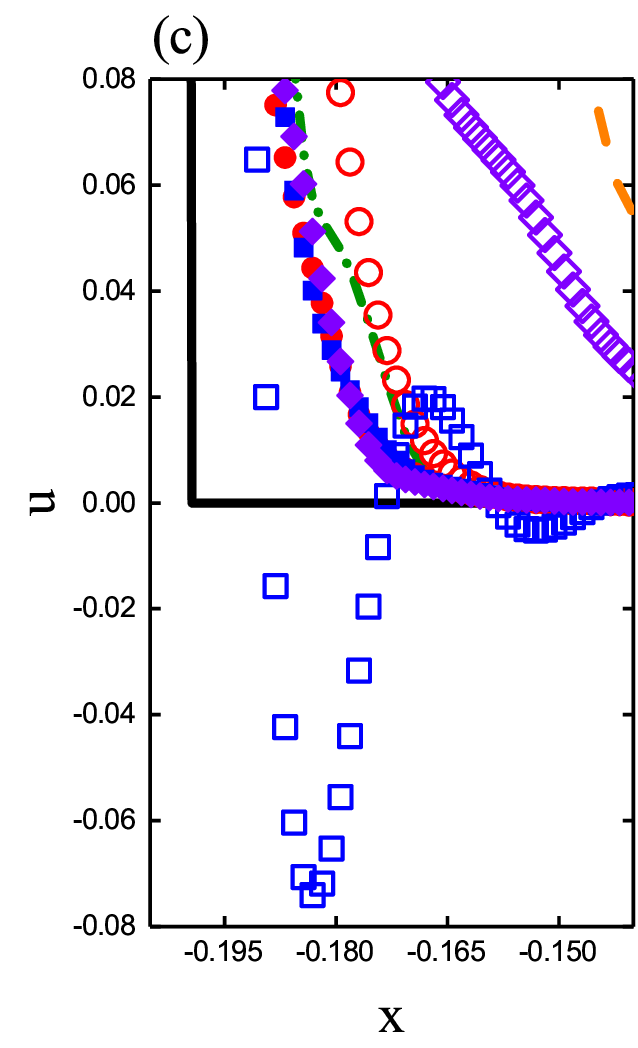}
\caption{Performance of the considered WENO schemes for SLP (Case 1) 
at output time $t = 200$ with a uniform mesh size of $N = 1600$.}
\label{fig:SLP:N1600:01}
\end{figure}

\begin{figure}[!ht]
\centering
\includegraphics[height=0.33\textwidth]
{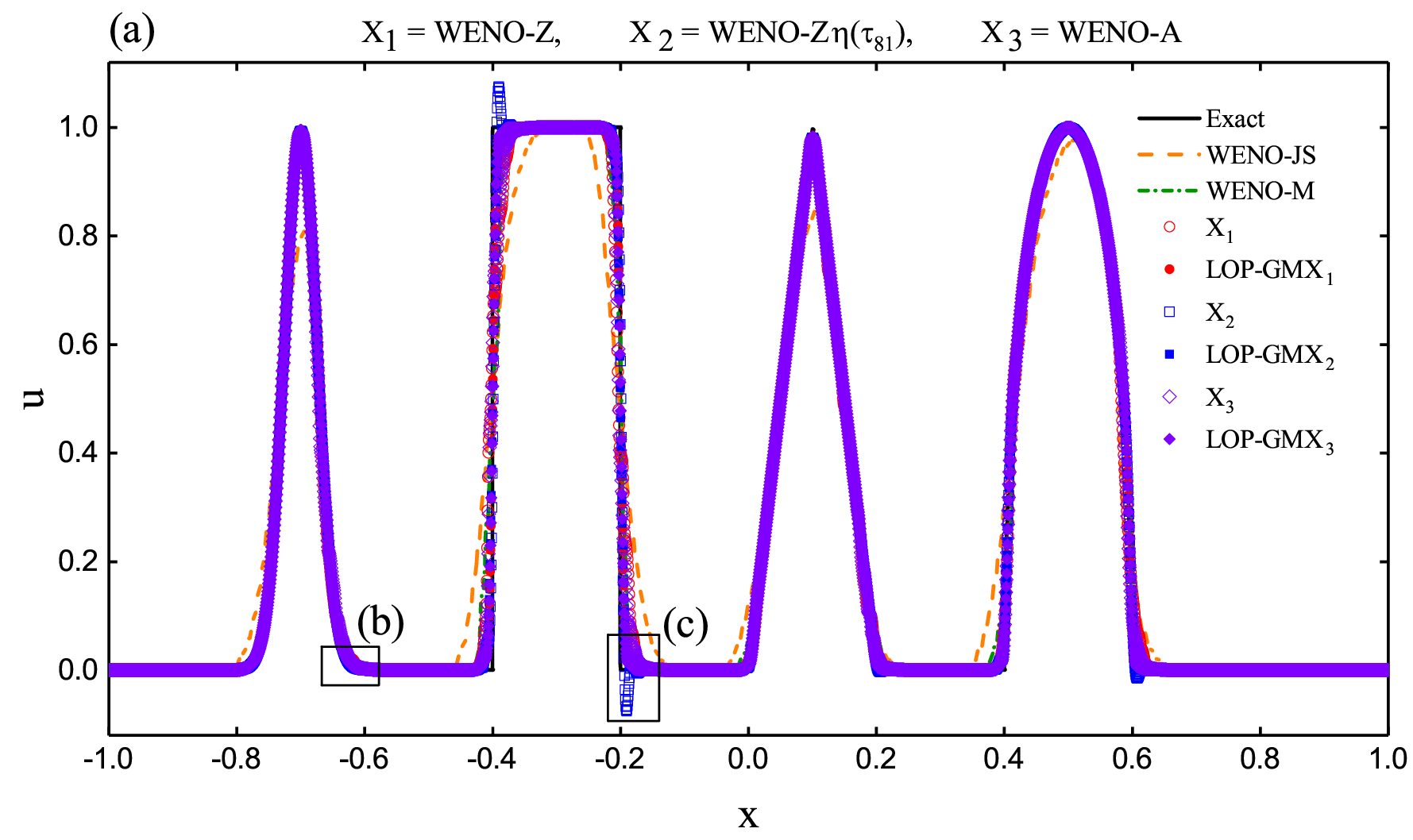}
\includegraphics[height=0.33\textwidth]
{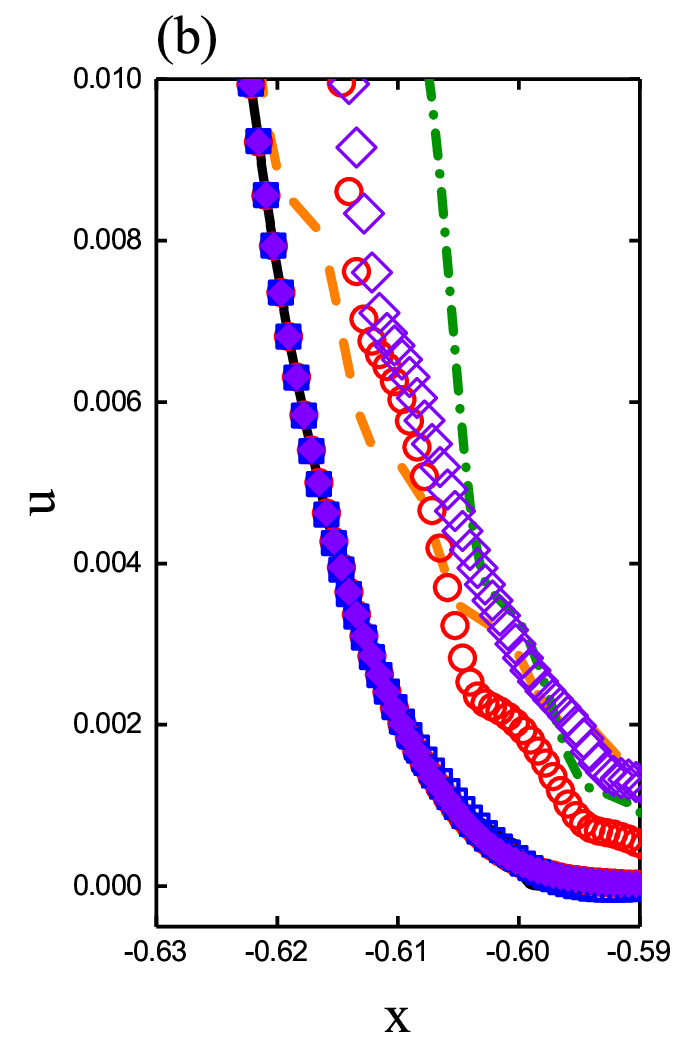}
\includegraphics[height=0.33\textwidth]
{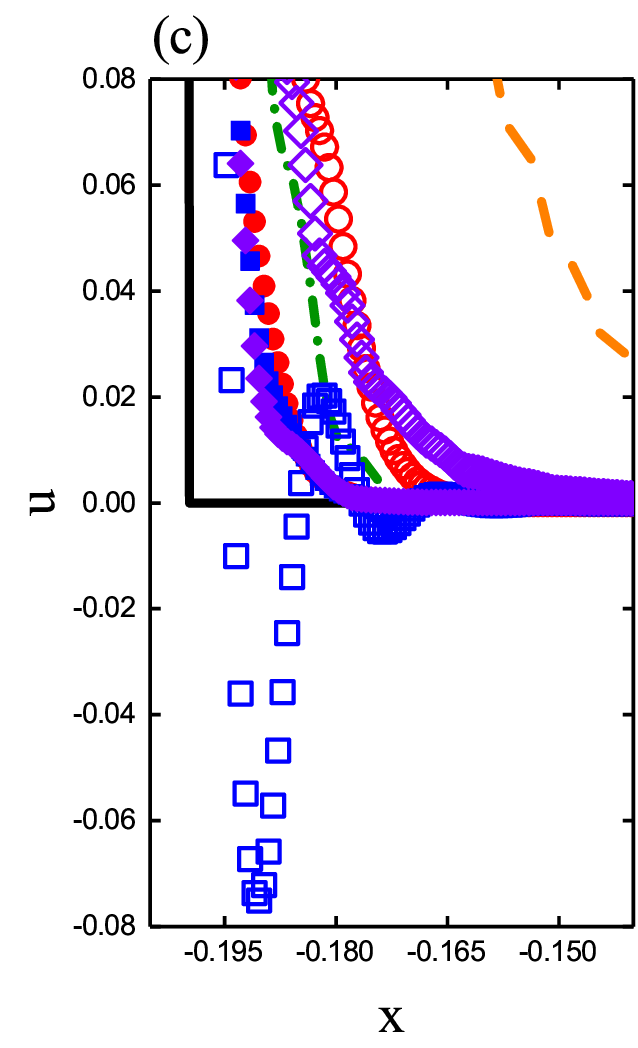}
\caption{Performance of the considered WENO schemes for SLP (Case 1) 
at output time $t = 200$ with a uniform mesh size of $N = 3200$.}
\label{fig:SLP:N3200:01}
\end{figure}

\begin{figure}[!ht]
\centering
\includegraphics[height=0.33\textwidth]
{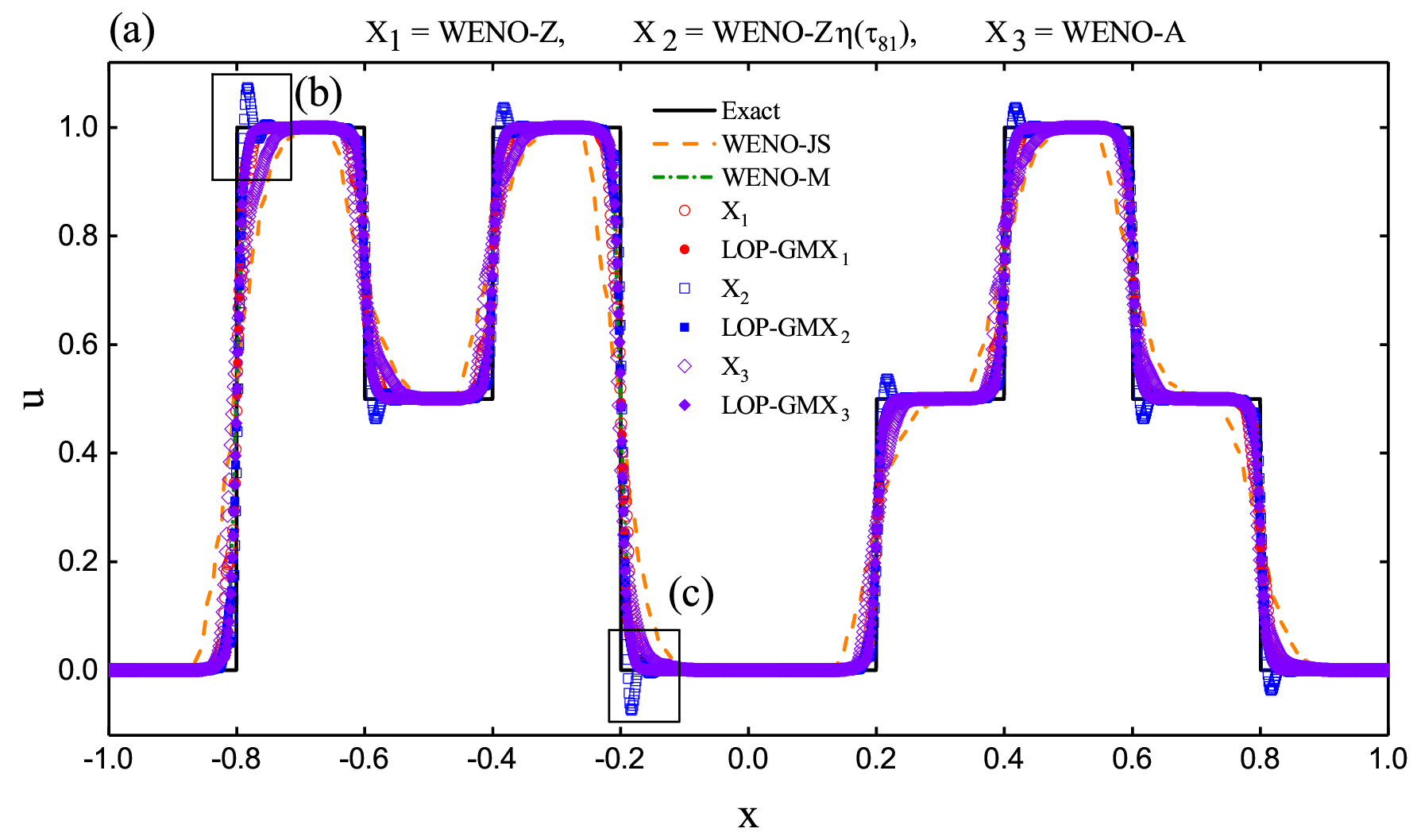}
\includegraphics[height=0.33\textwidth]
{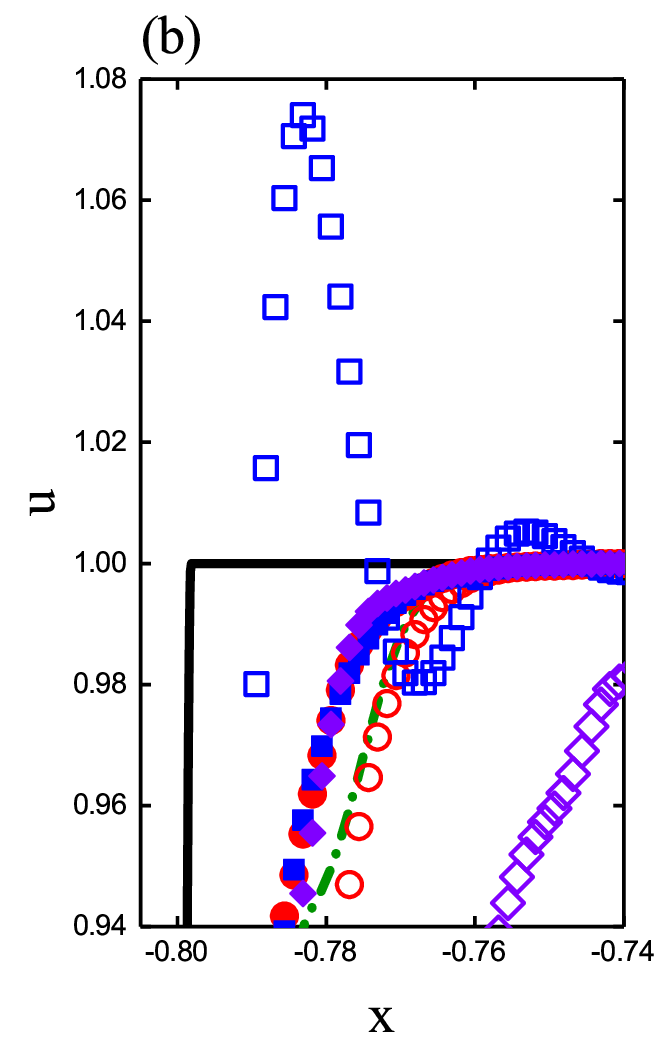}
\includegraphics[height=0.33\textwidth]
{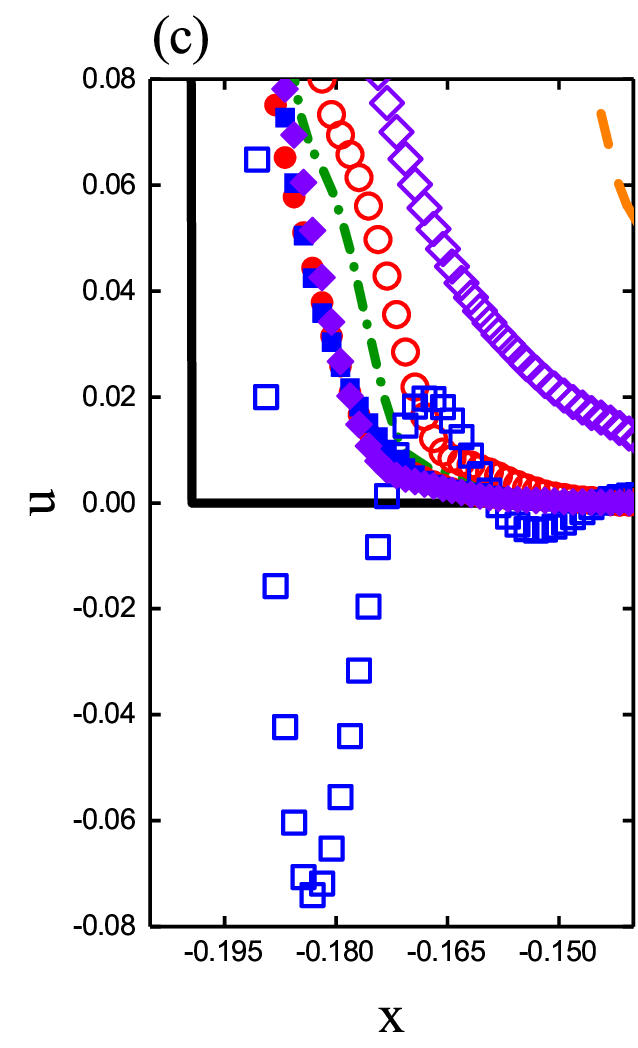}
\caption{Performance of the considered WENO schemes for BiCWP (Case 2) at output time $t = 200$ with a uniform mesh size of $N = 1600$.}
\label{fig:BiCWP:N1600:01}
\end{figure}

\begin{figure}[!ht]
\centering
\includegraphics[height=0.33\textwidth]
{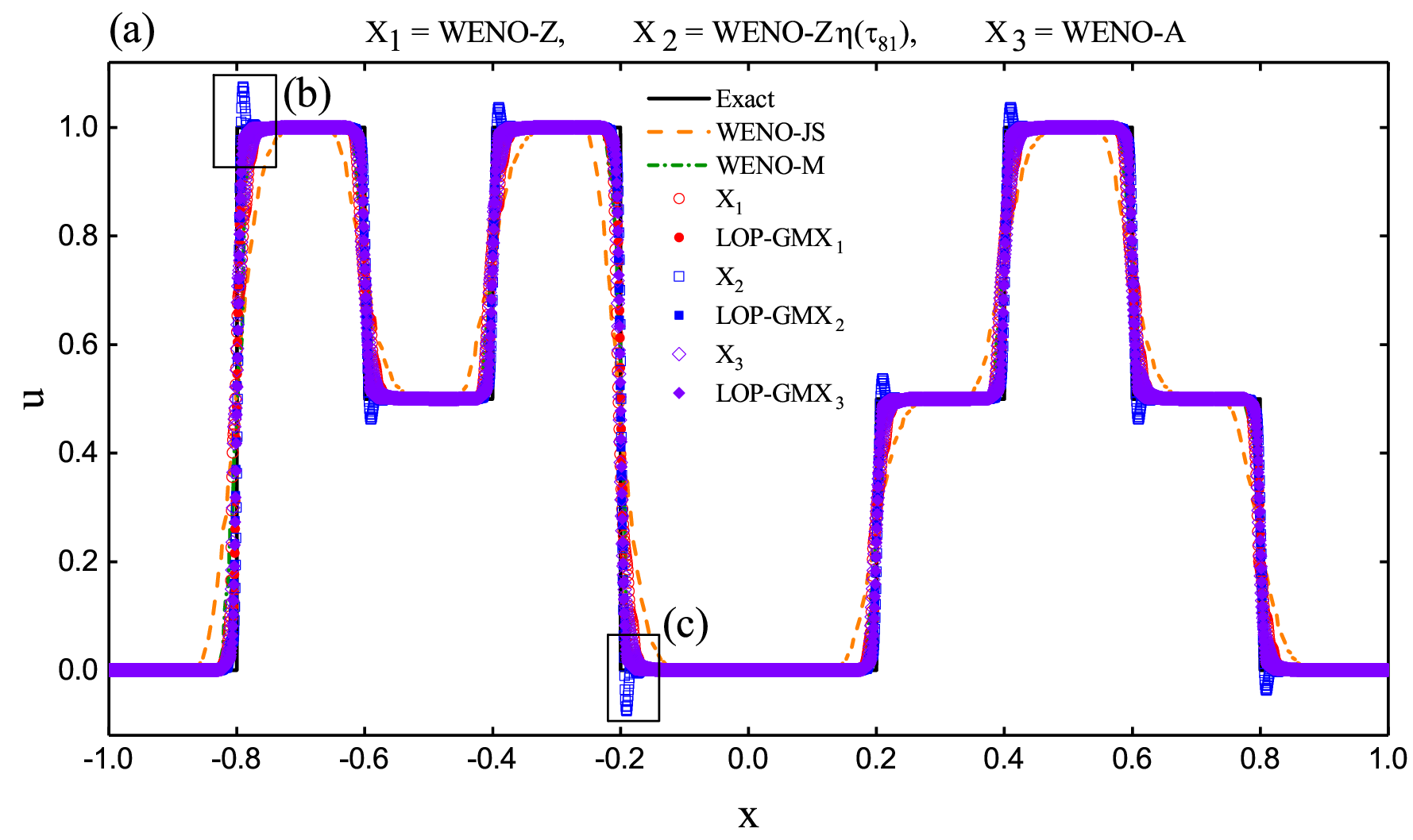}
\includegraphics[height=0.33\textwidth]
{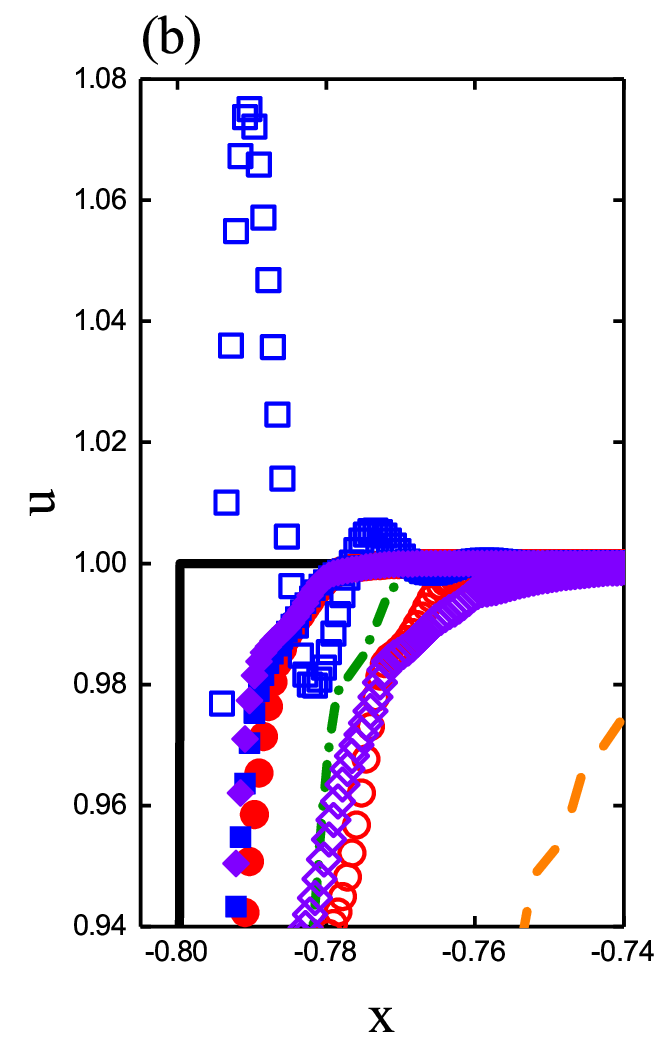}
\includegraphics[height=0.33\textwidth]
{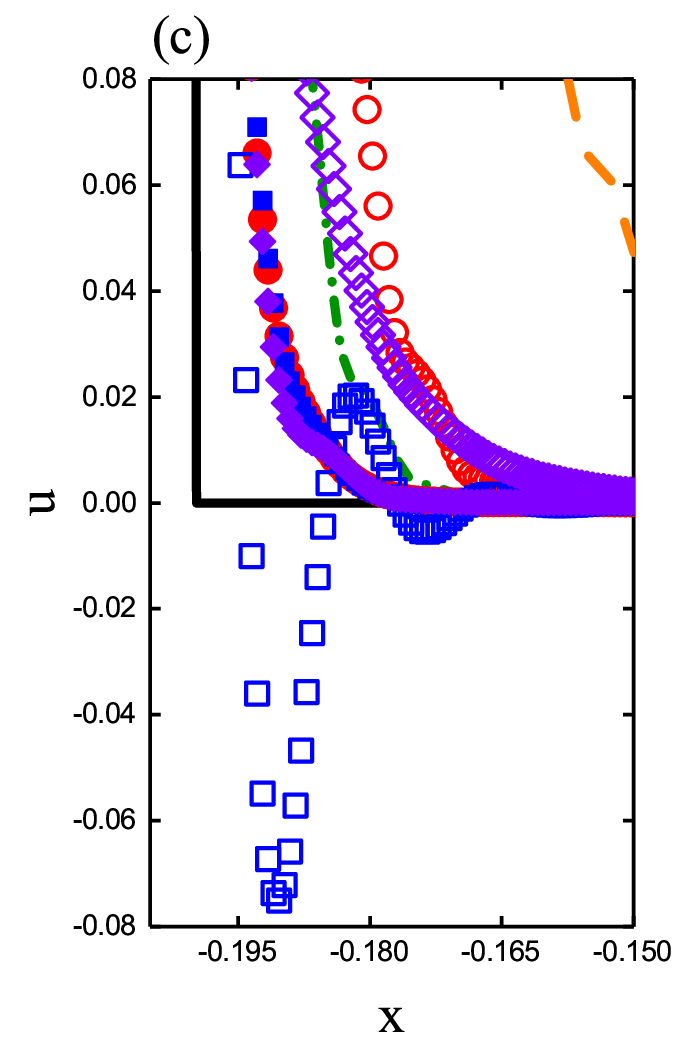}
\caption{Performance of the considered WENO schemes for BiCWP (Case 2) at output time $t = 200$ with a uniform mesh size of $N = 3200$.}
\label{fig:BiCWP:N3200:01}
\end{figure}

The solutions for $N = 1600, 3200$ with respect to Case 1 and Case 2 
are plotted in Fig. \ref{fig:SLP:N1600:01} and Fig. \ref{fig:SLP:N3200:01}, Fig. \ref{fig:BiCWP:N1600:01} and Fig. \ref{fig:BiCWP:N3200:01}, respectively. It is clear that, for both 
Case 1 and Case 2 regardless of the grid numbers: (1) the 
WENO-Z and WENO-A schemes have significantly lower resolution than 
the associated LOP-GMWENO-Z and LOP-GMWENO-A schemes that produce 
sharper solutions than the WENO-JS and WENO-M schemes; (2) the 
WENO-Z$\eta(\tau_{81})$ scheme results in severe spurious 
oscillations near discontinuities while the LOP-GMWENO-Z$\eta(\tau_{81})$ scheme can remove these oscillations successfully.

In summary, the following conclusion can be drawn from the present 
numerical tests that the use of the LOP property is helpful for the 
WENO-Z-type schemes to remove the spurious oscillations on the 
premise of achieving high resolutions for long output time 
simulations of hyperbolic problems with discontinuities. And this is 
a major contribution to the field of WENO-Z-type schemes in the 
present study.

\begin{table}[!ht]
%\begin{footnotesize}
%\begin{scriptsize}
\begin{myFontSize}
\centering
\caption{Numerical errors and convergence orders of accuracy on 
solving SLP, $t = 2000$.}
\label{table:AccuracyTest:SLP:t2000}
\begin{tabular*}{\hsize}
{@{}@{\extracolsep{\fill}}cllllllll@{}}
\hline
\space    &\multicolumn{4}{l}{\cellcolor{gray!35}{WENO5-ILW}}  
          &\multicolumn{4}{l}{\cellcolor{gray!35}{WENO-JS}}  \\
\cline{2-5}  \cline{6-9}
$N$                   & $L_{1}$ error      & $L_{1}$ order 
                      & $L_{\infty}$ error & $L_{\infty}$ order
                      & $L_{1}$ error      & $L_{1}$ order 
                      & $L_{\infty}$ error & $L_{\infty}$ order \\
\Xhline{0.65pt}
100                   & 4.70125E-01        & -
                      & 7.41389E-01        & -
                      & 6.33519E-01        & -
                      & 7.17815E-01        & - \\
200                   & 2.27171E-01        & 1.0493
                      & 5.14236E-01        & 0.5278
                      & 6.12899E-01        & 0.0477
                      & 7.99265E-01        & -0.1551\\
400                   & 1.15918E-01        & 0.9707
                      & 4.77803E-01        & 0.1060
                      & 5.99215E-01        & 0.0326
                      & 8.20493E-01        & -0.0378 \\
800                   & 5.35871E-02        & 1.1131
                      & 4.74317E-01        & 0.0106
                      & 5.50158E-01        & 0.1232
                      & 8.14650E-01        & 0.0103 \\
\hline
\space    &\multicolumn{4}{l}{\cellcolor{gray!35}{WENO-Z}}
          &\multicolumn{4}{l}{\cellcolor{gray!35}{LOP-GMWENO-Z}}\\
\cline{2-5}  \cline{6-9}
$N$                   & $L_{1}$ error      & $L_{1}$ order 
                      & $L_{\infty}$ error & $L_{\infty}$ order
              & $L_{1}$ error      & $L_{1}$ order 
                & $L_{\infty}$ error & $L_{\infty}$ order \\
\Xhline{0.65pt}
100                   & 5.47810E-01        & -
                      & 7.41658E-01        & -
                      & 5.56356E-01        & -
                      & 7.38792E-01        & - \\
200                   & 3.86995e-01        & 0.5014
                      & 6.85835e-01        & 0.1129
                      & 3.64352E-01        & 0.6107
                      & 7.09015E-01        & 0.0594\\
400                   & 2.02287e-01        & 0.9359
                      & 5.18993e-01        & 0.4021
                      & 1.74945E-01        & 1.0584
                      & 4.86425E-01        & 0.5436\\
800                   & 1.66703e-01        & 0.2791 
                      & 5.04564e-01        & 0.0407
                      & 6.10083E-02        & 1.5198
                      & 4.84660E-01        & 0.0052\\
\hline
\space    &\multicolumn{4}{l}{\cellcolor{gray!35}{WENO-Z$\eta(\tau_{81})$}}
          &\multicolumn{4}{l}{\cellcolor{gray!35}{LOP-GMWENO-Z$\eta(\tau_{81})$}}\\
\cline{2-5}  \cline{6-9}
$N$                   & $L_{1}$ error      & $L_{1}$ order 
                      & $L_{\infty}$ error & $L_{\infty}$ order
              & $L_{1}$ error      & $L_{1}$ order 
                & $L_{\infty}$ error & $L_{\infty}$ order \\
\Xhline{0.65pt}
100                   & 4.74455E-01        & -
                      & 7.86642E-01        & -
                      & 5.61608E-01        & -
                      & 7.46364E-01        & - \\
200                   & 2.42963e-01        & 0.9655
                      & 6.39818e-01        & 0.2980
                      & 3.95129E-01        & 0.5072
                      & 7.18477E-01        & 0.0549\\
400                   & 1.33752e-01        & 0.8612
                      & 6.01344e-01        & 0.0895
                      & 1.75622E-01        & 1.1699
                      & 4.85420E-01        & 0.5657\\
800                   & 5.89144e-02        & 1.1829
                      & 5.73819e-01        & 0.0676
                      & 6.07081E-02        & 1.5325
                      & 4.85758E-01        & -0.0010\\
\hline
\space    &\multicolumn{4}{l}{\cellcolor{gray!35}{WENO-A}}
          &\multicolumn{4}{l}{\cellcolor{gray!35}{LOP-GMWENO-A}}\\
\cline{2-5}  \cline{6-9}
$N$                   & $L_{1}$ error      & $L_{1}$ order 
                      & $L_{\infty}$ error & $L_{\infty}$ order
              & $L_{1}$ error      & $L_{1}$ order 
                & $L_{\infty}$ error & $L_{\infty}$ order \\
\Xhline{0.65pt}
100                   & 5.64755E-01        & -
                      & 7.27197E-01        & -
                      & 5.73967E-01        & -
                      & 7.37449E-01        & - \\
200                   & 5.31200e-01        & 0.0884
                      & 7.70910e-01        & -0.0842
                      & 4.06455E-01        & 0.4979
                      & 7.53060E-01        & -0.0302\\
400                   & 4.08352e-01        & 0.3794
                      & 6.93282e-01        & 0.1531
                      & 1.70772E-01        & 1.2510
                      & 5.21841E-01        & 0.5292\\
800                   & 2.95123e-01        & 0.4685
                      & 5.90637e-01        & 0.2312
                      & 6.18950E-02        & 1.4642
                      & 4.91304E-01        & 0.0870\\
\hline
\end{tabular*}
\end{myFontSize}
%\end{scriptsize}
%\end{footnotesize}
\end{table}

\begin{table}[!ht]
%\begin{footnotesize}
%\begin{scriptsize}
\begin{myFontSize}
\centering
\caption{Numerical errors and convergence orders of accuracy on 
solving BiCWP, $t = 2000$.}
\label{table:AccuracyTest:BiCWP:t2000}
\begin{tabular*}{\hsize}
{@{}@{\extracolsep{\fill}}cllllllll@{}}
\hline
\space    &\multicolumn{4}{l}{\cellcolor{gray!35}{WENO5-ILW}}  
          &\multicolumn{4}{l}{\cellcolor{gray!35}{WENO-JS}}  \\
\cline{2-5}  \cline{6-9}
$N$                   & $L_{1}$ error      & $L_{1}$ order 
                      & $L_{\infty}$ error & $L_{\infty}$ order
              & $L_{1}$ error      & $L_{1}$ order 
                & $L_{\infty}$ error & $L_{\infty}$ order \\
\Xhline{0.65pt}
100                   & 3.47945E-01        & -
                      & 5.14695E-01        & -
                      & 6.83328E-01        & -
                      & 5.82442E-01        & - \\
200                   & 1.96104E-01        & 0.8272
                      & 4.64745E-01        & 0.1473
                      & 5.89672E-01        & 0.2127
                      & 6.41175E-01        & -0.1386\\ 
400                   & 1.35386E-01        & 0.5345
                      & 4.74241E-01        & -0.0292
                      & 5.56639E-01        & 0.0832
                      & 5.94616E-01        & 0.1088\\
800                   & 7.96037E-02        & 0.7662
                      & 4.74182E-01        & 0.0002
                      & 4.72439E-01        & 0.2366
                      & 5.73614E-01        & 0.0519\\
\hline
\space    &\multicolumn{4}{l}{\cellcolor{gray!35}{WENO-Z}}
          &\multicolumn{4}{l}{\cellcolor{gray!35}{LOP-GMWENO-Z}}\\
\cline{2-5}  \cline{6-9}
$N$                   & $L_{1}$ error      & $L_{1}$ order 
                      & $L_{\infty}$ error & $L_{\infty}$ order
              & $L_{1}$ error      & $L_{1}$ order 
                & $L_{\infty}$ error & $L_{\infty}$ order \\
\Xhline{0.65pt}
100                   & 4.63184E-01        & -
                      & 5.26156E-01        & -
                      & 5.22837E-01        & -
                      & 4.91514E-01        & - \\
200                   & 3.13567E-01        & 0.5628
                      & 4.84876E-01        & 0.1179
                      & 2.69161E-01        & 0.9579
                      & 5.26700E-01        & -0.0997\\ 
400                   & 2.23255E-01        & 0.4901
                      & 5.10834E-01        & -0.0752
                      & 1.77055E-01        & 0.6043
                      & 5.00232E-01        & 0.0744\\
800                   & 1.74777E-01        & 0.3532
                      & 5.19528E-01        & -0.0243
                      & 9.20312E-02        & 0.9440
                      & 4.90299E-01        & 0.0289\\
\hline
\space    &\multicolumn{4}{l}{\cellcolor{gray!35}{WENO-Z$\eta(\tau_{81})$}}
          &\multicolumn{4}{l}{\cellcolor{gray!35}{LOP-GMWENO-Z$\eta(\tau_{81})$}}\\
\cline{2-5}  \cline{6-9}
$N$                   & $L_{1}$ error      & $L_{1}$ order 
                      & $L_{\infty}$ error & $L_{\infty}$ order
              & $L_{1}$ error      & $L_{1}$ order 
                & $L_{\infty}$ error & $L_{\infty}$ order \\
\Xhline{0.65pt}
100                   & 3.72310E-01        & -
                      & 5.87087E-01        & -
                      & 4.91059E-01        & -
                      & 5.27869E-01        & - \\
200                   & 2.39565E-01        & 0.6361
                      & 6.71171E-01        & -0.1931
                      & 2.72777E-01        & 0.8482
                      & 4.90307E-01        & 0.1065\\ 
400                   & 1.47989E-01        & 0.6949
                      & 5.78458E-01        & 0.2145
                      & 1.75297E-01        & 0.6379
                      & 5.01408E-01        & -0.0323\\
800                   & 8.57620E-02        & 0.7871
                      & 5.75825E-01        & 0.0066
                      & 9.05607E-02        & 0.9528
                      & 4.89516E-01        & 0.0346\\
\hline
\space    &\multicolumn{4}{l}{\cellcolor{gray!35}{WENO-A}}
          &\multicolumn{4}{l}{\cellcolor{gray!35}{LOP-GMWENO-A}}\\
\cline{2-5}  \cline{6-9}
$N$                   & $L_{1}$ error      & $L_{1}$ order 
                      & $L_{\infty}$ error & $L_{\infty}$ order
              & $L_{1}$ error      & $L_{1}$ order 
                & $L_{\infty}$ error & $L_{\infty}$ order \\
\Xhline{0.65pt}
100                   & 5.29370E-01        & -
                      & 5.34285E-01        & -
                      & 5.10558E-01        & -
                      & 5.07303E-01        & - \\
200                   & 4.14728E-01        & 0.3521
                      & 5.62824E-01        & -0.0751
                      & 2.80716E-01        & 0.8630
                      & 5.39137E-01        & -0.0878\\ 
400                   & 3.34332E-01        & 0.3109
                      & 5.47207E-01        & 0.0406
                      & 1.77541E-01        & 0.6610
                      & 5.06588E-01        & 0.0898\\
800                   & 2.67754E-01        & 0.3204
                      & 5.49441E-01        & -0.0059
                      & 8.95020E-02        & 0.9882
                      & 4.93665E-01        & 0.0373\\
\hline
\end{tabular*}
\end{myFontSize}
%\end{scriptsize}
%\end{footnotesize}
\end{table}

\subsubsection{With high-order critical points}
\label{subsubsec:LongCritical}
In order to show the advantage of the LOP-GMWENO-X schemes that they 
can obtain high resolutions for the problem with high-order critical 
points at long output times, we compute the 1D linear advection 
equation $u_{t} + u_{x} = 0, x \in (3.5, 6.5)$ with an initial 
condition given by $u(x, 0) = - \exp\Big(-\big(x-5.0\big)^{10}\Big)\cos^{10}\Big(\pi\big(x - 5.0\big)\Big)$. The CFL number is chosen to 
be $(\Delta x)^{2/3}$. For convenience, this test is dubbed HCP 
problem here.
 
It is trivial to know that the exact solution is $u(x,t) = -\exp\Big(-\big((x - t) - 5.0\big)^{10}\Big)\cos^{10}\Big(\pi\big((x - t) - 5.0\big)\Big)$. To compare the dissipations, we compute the increased 
errors, defined by $\chi_{1}=\frac{L_{1}^{\mathrm{Y}}(t)-L_{1}^{\mathrm{ILW}}(t)}{L_{1}^{\mathrm{ILW}}(t)}\times100\%$ and $\chi_{2}=\frac{L_{2}^{\mathrm{Y}}(t)-L_{2}^{\mathrm{ILW}}(t)}{L_{2}^{\mathrm{ILW}}(t)}\times100\%$, where $L_{j}^{\mathrm{ILW}}(t)$ and $L_{j}^{\mathrm{Y}}(t)$ stand for the $L_{j}$-norm ($j=1,2$) of numerical 
errors of the WENO5-ILW scheme and the scheme ``Y'' respectively. 

Table \ref{tab::long-run:N300} gives the $L_{1}$ and $L_{2}$ errors 
and the corresponding increased errors computed by various 
considered WENO schemes with $N=300$ and $t = 300, 600, 900, 1200$. 
For all output times, the WENO-JS scheme generates the 
largest numerical errors and it produces the highest increased 
errors. The WENO-Z and WENO-A schemes also have very large numerical 
errors that are slightly smaller than those of the WENO-JS scheme. 
Accordingly, the increased errors of these schemes are excessively 
large. In contrast, the numerical errors and the corresponding 
increased errors of the associated LOP-GMWENO-Z and LOP-GMWENO-A 
schemes are significantly reduced to a tolerable level. Actually, 
they can achieve much smaller numerical errors that get very close 
to those of the WENO-ILW scheme. In addition, the LOP-GMWENO-Z$\eta(\tau_{81})$ scheme is also able to maintain the numerical errors at 
an acceptable level to ensure that its increased errors are also at 
a tolerable level. It appears that the WENO-Z$\eta(\tau_{81})$ 
schemes gets solutions almost as accurate as, or even more accurate 
than, those of the WENO-ILW scheme. Of course, this is good, at 
least for this test. However, as discussed earlier, it suffers from 
lack of robustness as its dissipation is too small on solving 
problems with discontinuities, especially for long output time 
simulations. 

\begin{table}[!ht]
%\begin{footnotesize}
%\begin{scriptsize}
\begin{myFontSize}
\centering
\caption{$L_{1}$-, $L_{2}$-norm of numerical errors and their 
associated increased errors (in percentage) on solving the HCP 
problem with $N = 300$.}
\label{tab::long-run:N300}
\begin{tabular*}{\hsize}
{@{}@{\extracolsep{\fill}}rlrlrlrlr@{}}
\hline
\space    &\multicolumn{4}{l}{\cellcolor{gray!35}{WENO5-ILW}}  
          &\multicolumn{4}{l}{\cellcolor{gray!35}{WENO-JS}}  \\
\cline{2-5}  \cline{6-9}
Time, $t$             & $L_{1}$ error      & $\chi_{1}$
                      & $L_{2}$ error & $\chi_{2}$
                      & $L_{1}$ error      & $\chi_{1}$
                      & $L_{2}$ error & $\chi_{2}$ \\
\Xhline{0.65pt}
300                   & 6.30537E-03        & -
                      & 5.69915E-03        & -
                      & 8.88538E-02        & 1309\%
                      & 7.58151E-02        & 1230\%\\
600                   & 1.14068E-02        & -
                      & 9.89571E-03        & -
                      & 2.17193E-01        & 1804\%
                      & 1.78733E-01        & 1706\%\\
900                   & 1.58862E-02        & -
                      & 1.34968E-02        & -
                      & 2.84952E-01        & 1694\%
                      & 2.29984E-01        & 1604\%\\
1200                  & 1.98304E-02        & -
                      & 1.65947E-02        & -
                      & 3.32245E-01        & 1575\%
                      & 2.64896E-01        & 1496\%\\
\hline
\space    &\multicolumn{4}{l}{\cellcolor{gray!35}{WENO-Z}}
          &\multicolumn{4}{l}{\cellcolor{gray!35}{LOP-GMWENO-Z}}\\
\cline{2-5}  \cline{6-9}
Time, $t$             & $L_{1}$ error      & $\chi_{1}$         
                      & $L_{2}$ error & $\chi_{2}$     
              & $L_{1}$ error      & $\chi_{1}$         
                & $L_{2}$ error & $\chi_{2}$      \\
\Xhline{0.65pt}
300                   & 3.56061E-02        & 465\%
                      & 3.70325E-02        & 550\%
                      & 1.17228E-02        & 86\%
                      & 9.94127E-03        & 74\%\\
600                   & 8.21730E-02        & 620\%
                      & 8.66509E-02        & 776\%
                      & 2.06455E-02        & 81\%
                      & 1.86636E-02        & 89\%\\
900                   & 1.10948E-01        & 598\%
                      & 1.20020E-01        & 789\%
                      & 2.88134E-02        & 81\%
                      & 2.42690E-02        & 80\%\\
1200                  & 1.33598E-01        & 574\%
                      & 1.30955E-01        & 689\%
                      & 3.06112E-02        & 54\%
                      & 2.47195E-02        & 49\%\\
\hline
\space    &\multicolumn{4}{l}{\cellcolor{gray!35}{WENO-Z$\eta(\tau_{81})$}}
          &\multicolumn{4}{l}{\cellcolor{gray!35}{LOP-GMWENO-Z$\eta(\tau_{81})$}}\\
\cline{2-5}  \cline{6-9}
Time, $t$             & $L_{1}$ error      & $\chi_{1}$         
                      & $L_{2}$ error & $\chi_{2}$     
              & $L_{1}$ error      & $\chi_{1}$         
                & $L_{2}$ error & $\chi_{2}$      \\
\Xhline{0.65pt}
300                   & 6.41911E-03        & 2\%
                      & 6.41911E-03        & 13\%
                      & 1.10215E-02        & 75\%
                      & 9.92087E-03        & 74\%\\
600                   & 1.18436E-02        & 4\%
                      & 1.00019E-02        & 1\%
                      & 2.03525E-02        & 78\%
                      & 1.77994E-02        & 80\%\\
900                   & 1.65166E-02        & 4\%
                      & 1.34843E-02        & 0\%
                      & 2.74958E-02        & 73\%
                      & 2.29759E-02        & 70\%\\
1200                  & 2.01696E-02        & 2\%
                      & 1.63065E-02        & -2\%
                      & 2.92445E-02        & 47\%
                      & 2.45566E-02        & 48\%\\
\hline
\space    &\multicolumn{4}{l}{\cellcolor{gray!35}{WENO-A}}
          &\multicolumn{4}{l}{\cellcolor{gray!35}{LOP-GMWENO-A}}\\
\cline{2-5}  \cline{6-9}
Time, $t$             & $L_{1}$ error      & $\chi_{1}$         
                      & $L_{2}$ error & $\chi_{2}$     
              & $L_{1}$ error      & $\chi_{1}$         
                & $L_{2}$ error & $\chi_{2}$      \\
\Xhline{0.65pt}
300                   & 1.20889E-01        & 1817\%
                      & 1.13688E-01        & 1895\%
                      & 1.10467E-02        & 75\%
                      & 9.93024E-03        & 74\%\\
600                   & 1.84388E-01        & 1516\%
                      & 1.63149E-01        & 1549\%
                      & 2.03180E-02        & 78\%
                      & 1.81152E-02        & 83\%\\
900                   & 1.96326E-01        & 1136\%
                      & 1.70146E-01        & 1161\%
                      & 2.57805E-02        & 62\%
                      & 2.21596E-02        & 64\%\\
1200                  & 2.12758E-01        & 973\%
                      & 1.82152E-01        & 998\%
                      & 2.75976E-02        & 39\%
                      & 2.35188E-02        & 42\%\\
\hline
\end{tabular*}
\end{myFontSize}
%\end{scriptsize}
%\end{footnotesize}
\end{table}

\begin{figure}[!ht]
\centering
  \includegraphics[height=0.38\textwidth]
  {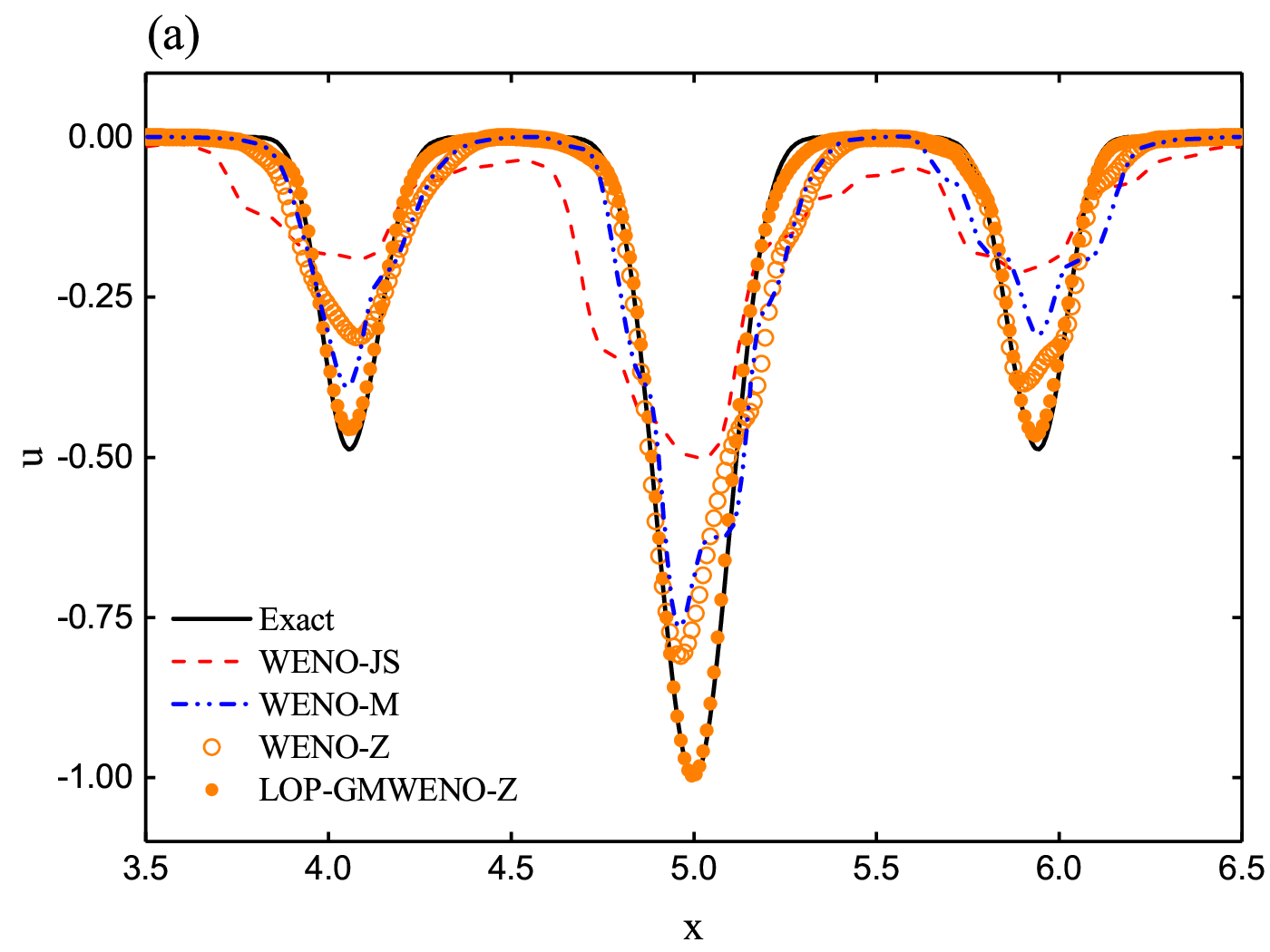}
  \includegraphics[height=0.38\textwidth]
  {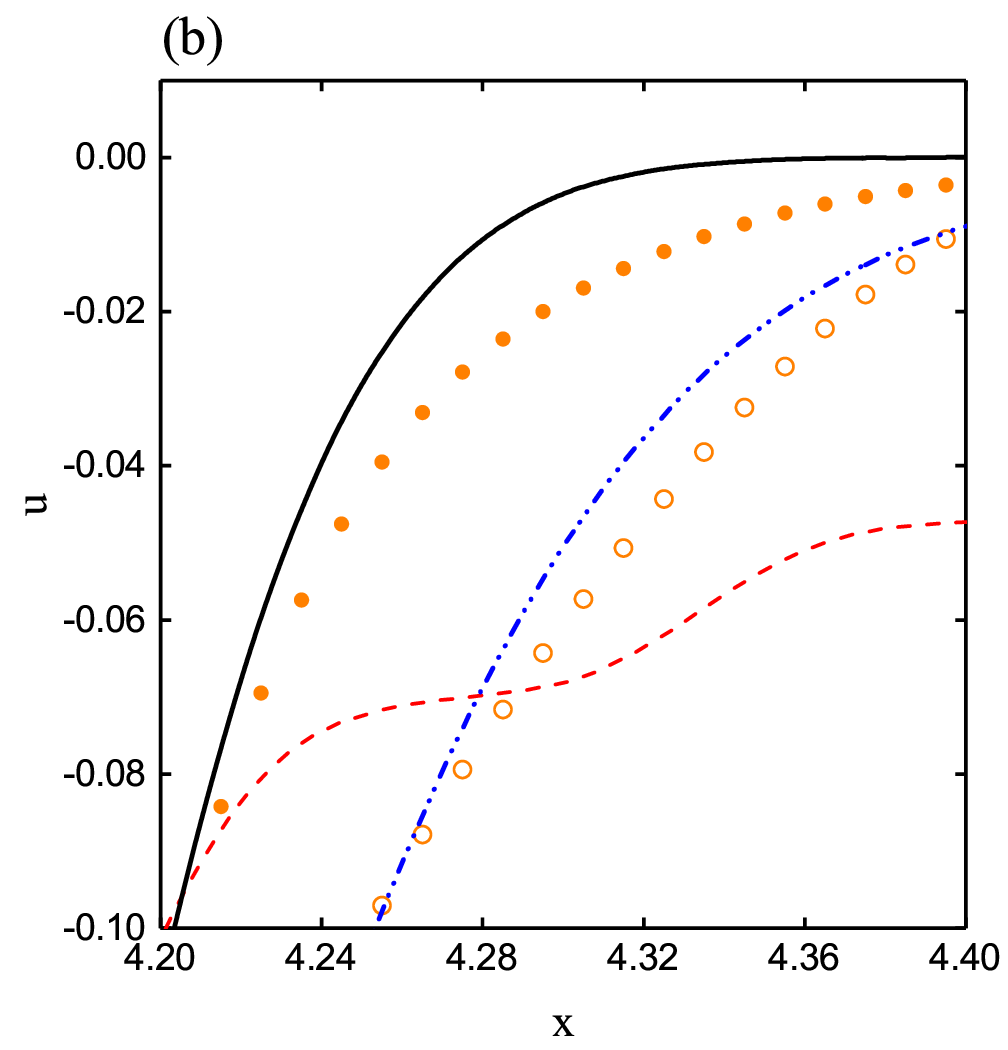}\\
  \vspace{0.5em}
  \includegraphics[height=0.38\textwidth]
  {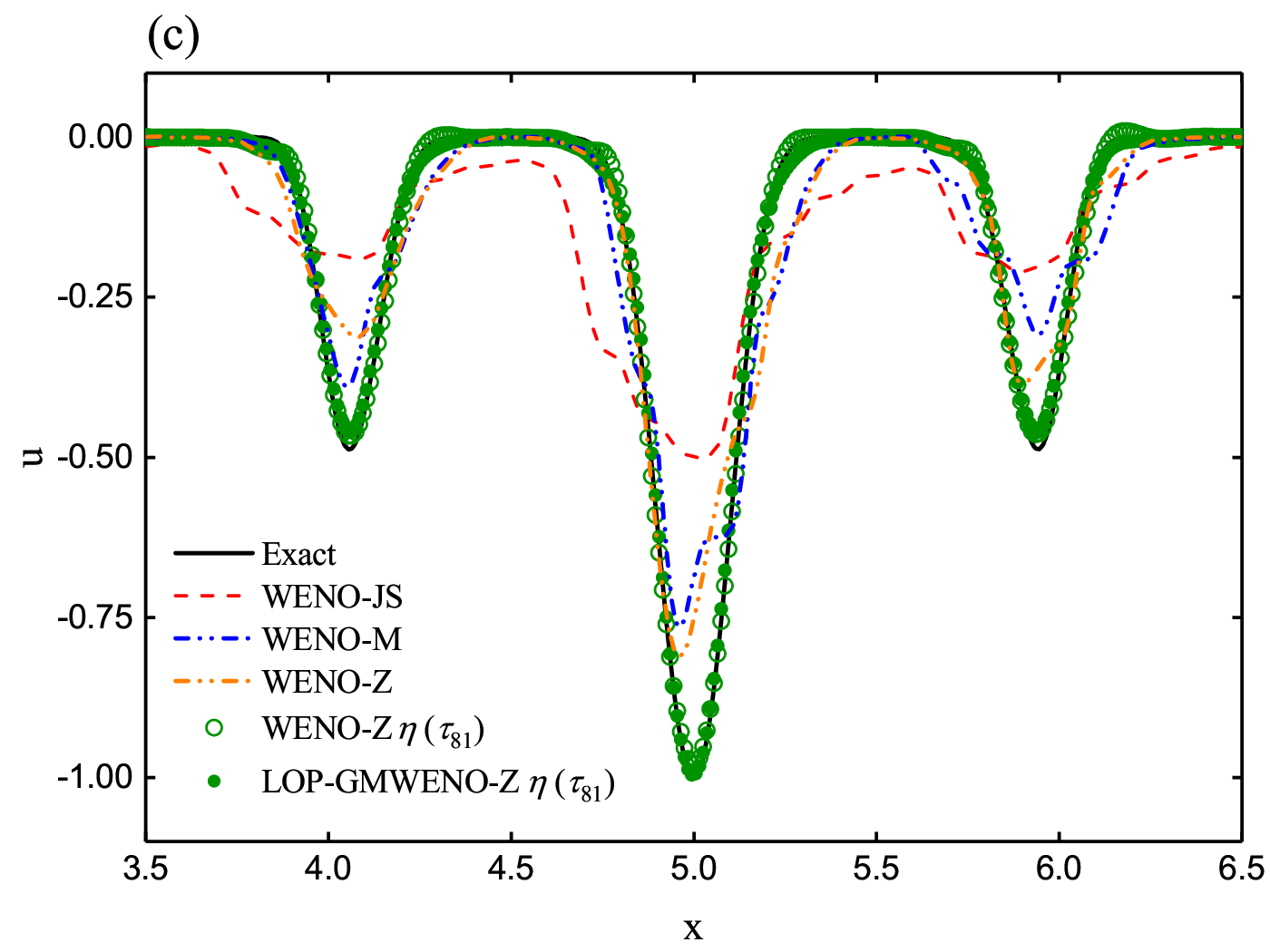}
  \includegraphics[height=0.38\textwidth]
  {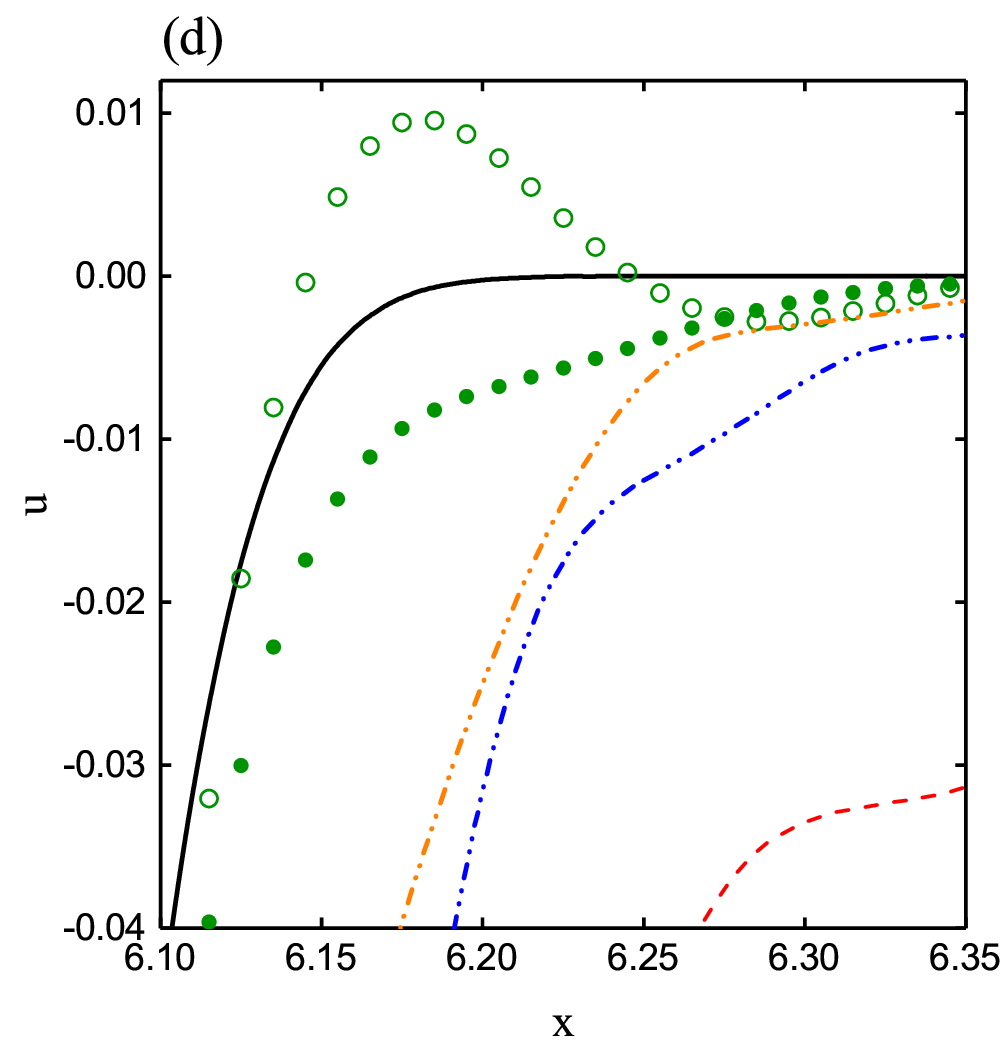} \\
  \vspace{0.5em}  
  \includegraphics[height=0.38\textwidth]
  {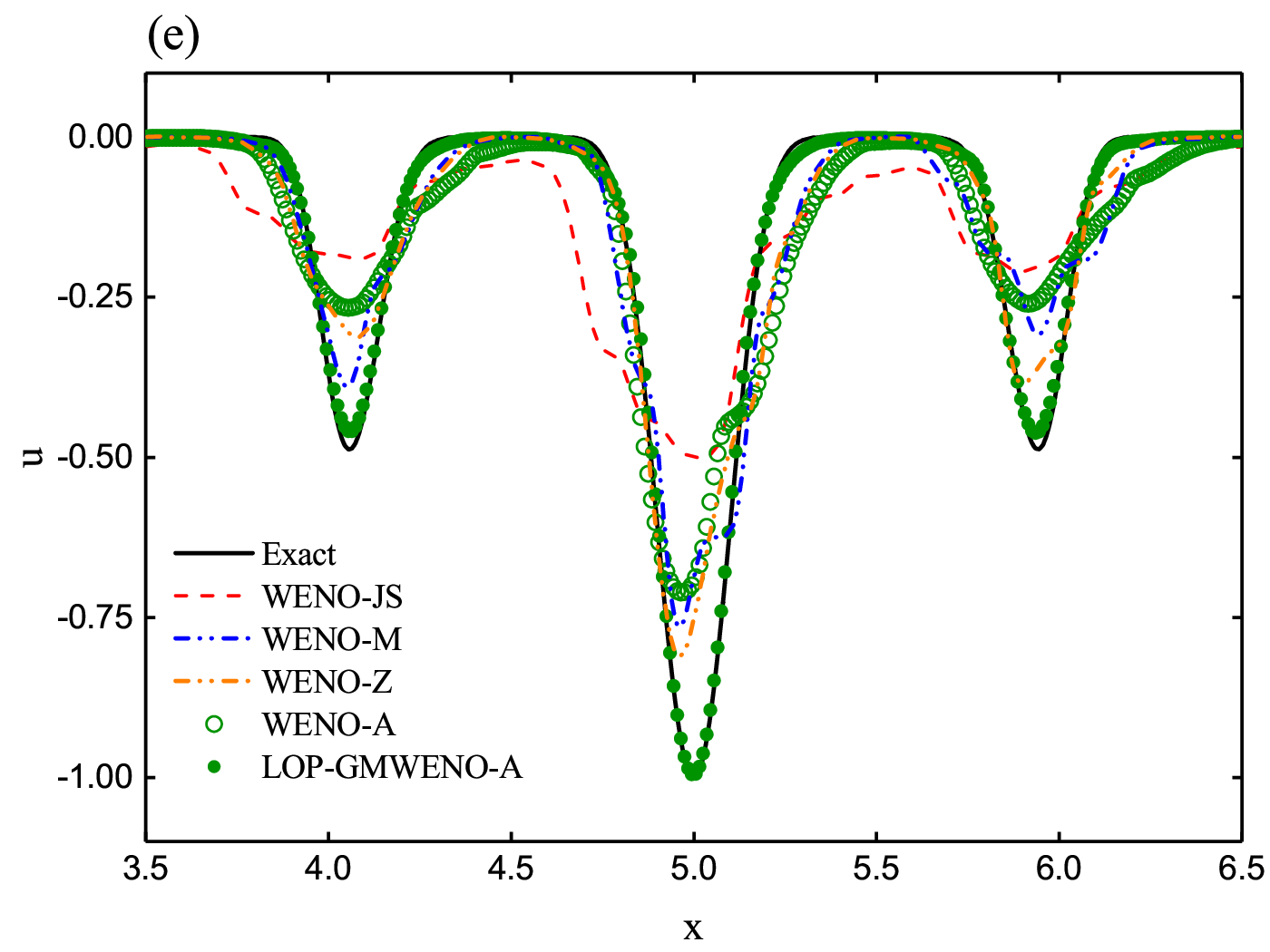}      
  \includegraphics[height=0.38\textwidth]
  {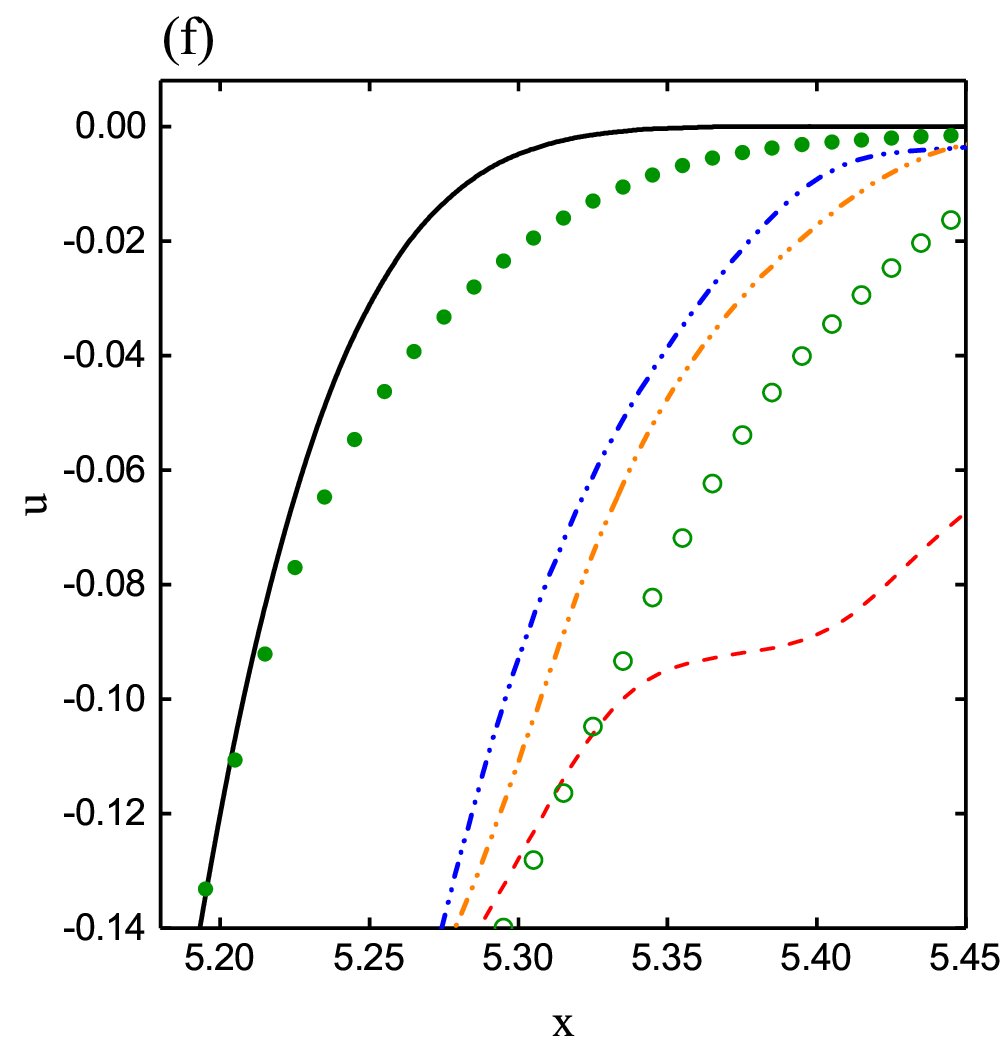}  
  \caption{Performance of various WENO schemes for the HCP problem at output time $t = 1200$ with a uniform mesh size of $N = 300$.}
\label{fig:ex:long-run:N300}
\end{figure}

In Fig. \ref{fig:ex:long-run:N300}, we plot the $x-u$ profiles of 
the LOP-GMWENO-X schemes and their associated WENO-X schemes at 
output time $t = 1200$. For the purpose of comparison, we also plot 
the $x-u$ profiles of the WENO-JS \cite{WENO-JS} and WENO-M 
\cite{WENO-M} schemes. From Fig. \ref{fig:ex:long-run:N300}, we can 
intuitively see that the WENO-JS scheme shows the lowest resolution, 
followed by the WENO-M scheme. As expected, the WENO-Z and WENO-A 
schemes show very low resolutions. However, the resolutions of the 
associated LOP-GMWENO-Z and LOP-GMWENO-A schemes have been improved 
significantly. Although the resolutions are slightly lower than 
those of its associated WENO-Z$\eta(\tau_{81})$ scheme, the 
LOP-GMWENO-Z$\eta(\tau_{81})$ scheme still shows far better 
resolutions than the WENO-JS, WENO-M, WENO-Z and WENO-A schemes.

%%% Local Variables:
%%% mode: latex
%%% TeX-master: "article"
%%% End:

\section{Numerical results}
\label{NumericalExperiments}
In this section, we compare the numerical performance of the 
considered WENO schemes by solving the system of hyperbolic 
conservation laws. The governing equations are the following Euler 
equations
\begin{equation}
\begin{array}{l}
\begin{aligned}
&\dfrac{\partial \rho}{\partial t} + \nabla \cdot (\rho \mathbf{u}) = 0,\\
&\dfrac{\partial (\rho \mathbf{u})}{\partial t} + \nabla \cdot (\rho \mathbf{u}\mathbf{u} + p\delta) = 0, \\
&\dfrac{\partial E}{\partial t} + \nabla \cdot \Big(\mathbf{u}(E + p)\Big) = 0,
\end{aligned}
\end{array}
\label{Numer:EulerEquations}
\end{equation}
where $\rho, \mathbf{u}, E, p$ are the density, the velocity vector, 
the total energy and the pressure. To close this set of equations, 
the equation of state (EOS) for an ideal polytropic gas $E = \frac{p}{\gamma - 1} + \frac{1}{2}\rho \mathbf{u}\mathbf{u}$ with $\gamma = 1.4$ as the ratio of specific heats is employed. 

The CFL number is set to be $0.5$ for both 1D and 2D problems below 
(unless otherwise noted), and the global Lax-Friedrichs flux 
splitting with the local characteristic decomposition \cite{WENO-JS} 
is employed. The WENO schemes are applied dimension-by-dimension to 
solve the two-dimensional Euler system. Zhang et al. \cite{FVMaccuracyProofs03} proposed two classes of finite volume WENO 
schemes in two dimension very carefully. For the same reason as 
discussed in \cite{MOP-WENO-ACMk}, the case of Class A is taken in 
this paper.

\subsection{One-dimensional Euler system}
\label{subsec:Euler1D}
In this subsection, we apply the LOP-GMWENO-X schemes with X = Z, Z$\eta(\tau_{81}), A$, as well as their associated WENO-X schemes, to 
solve several one-dimensional Euler problems. In order to compare 
the performance of these schemes with the recently-published 
low-dissipation shock-capturing ENO-family schemes, say, TENO 
schemes \cite{TENO-Fu-JCP2016,TENO-Fu-JCP2018,TENO-Fu-CPC2019-01,TENO-Fu-CPC2019-02,TENO-Fu-JSC2019,TENO-Fu-CiCP2019,TENO-Fu-CMAME2021,TENO-Haimovich-CF2017} in terms of the low-dissipation property, 
the classical TENO5 scheme with the threshold $C_{T}$ as $10^{-5}, 10^{-3}$ and $10^{-1}$ are considered. We also give the results of the 
WENO-JS scheme.

\subsubsection{Shock-tube problem}
To examine the shock-capturing capability of the LOP-GMWENO-X 
schemes, we compute two widely concerned shock-tube problems, that 
is, the Sod problem \cite{SodShock-tubeProblem} and the Lax problem \cite{LaxShock-tubeProblem}. The initial conditions are
\begin{equation*}
\begin{aligned}
&\text{Sod: } \qquad \big( \rho, u, p \big)(x, 0) =\left\{
\begin{array}{ll}
(1.0, 0.0, 1.0), & x \in [0.0, 0.5], \\
(0.125, 0.0, 0.1), & x \in [0.5, 1.0],
\end{array}\right. \\
&\text{Lax: } \qquad \big( \rho, u, p \big)(x, 0) =\left\{
\begin{array}{ll}
(0.445, 0.698, 3.528), & x \in [-5, 0], \\
(0.500, 0.000, 0.571), & x \in [0, 5].
\end{array}\right.
\end{aligned}
\label{ex:Shock-tube}
\end{equation*}
The transmissive boundary conditions are used, and the output times 
are 0.25 and 1.3, respectively. In order to compare the performance 
of the considered schemes on various coarse mesh resolutions, three 
different mesh sizes of $N = 300, 900, 1500$ are used.

\begin{figure}[!ht]
\centering
  \includegraphics[height=0.26\textwidth]
  {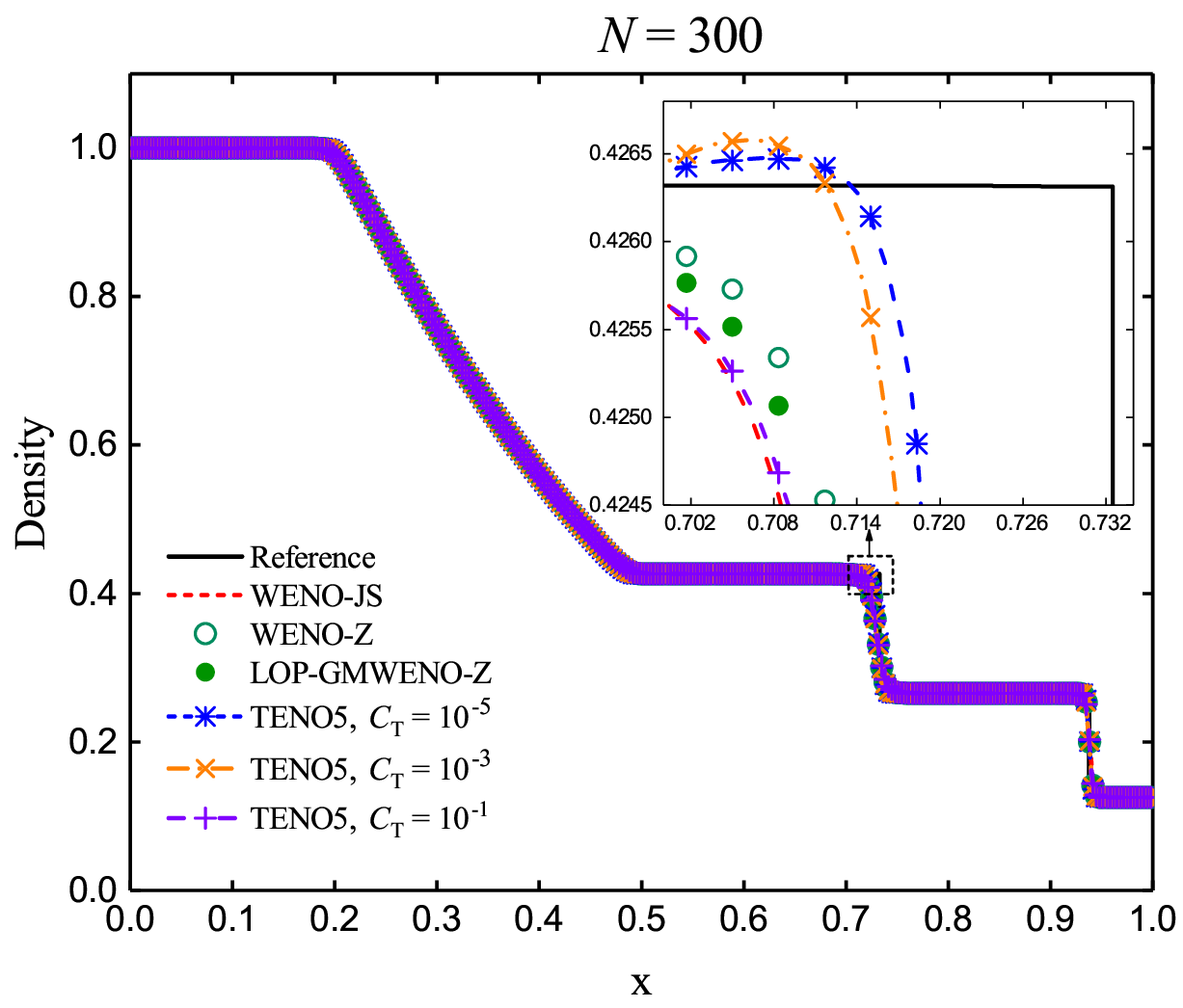}          \hspace{1.2ex}
  \includegraphics[height=0.26\textwidth]
  {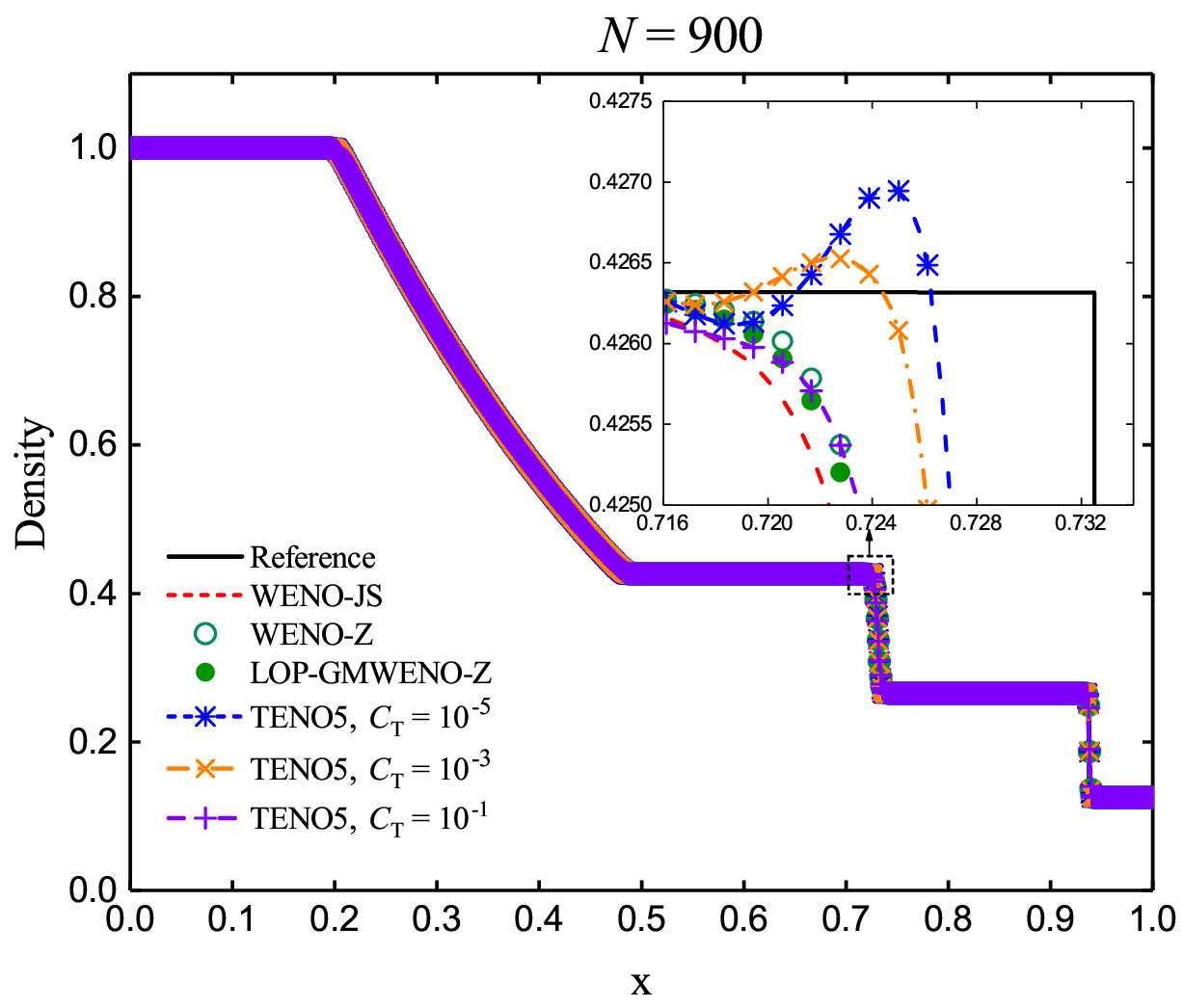}         \hspace{1.2ex}
  \includegraphics[height=0.26\textwidth]
  {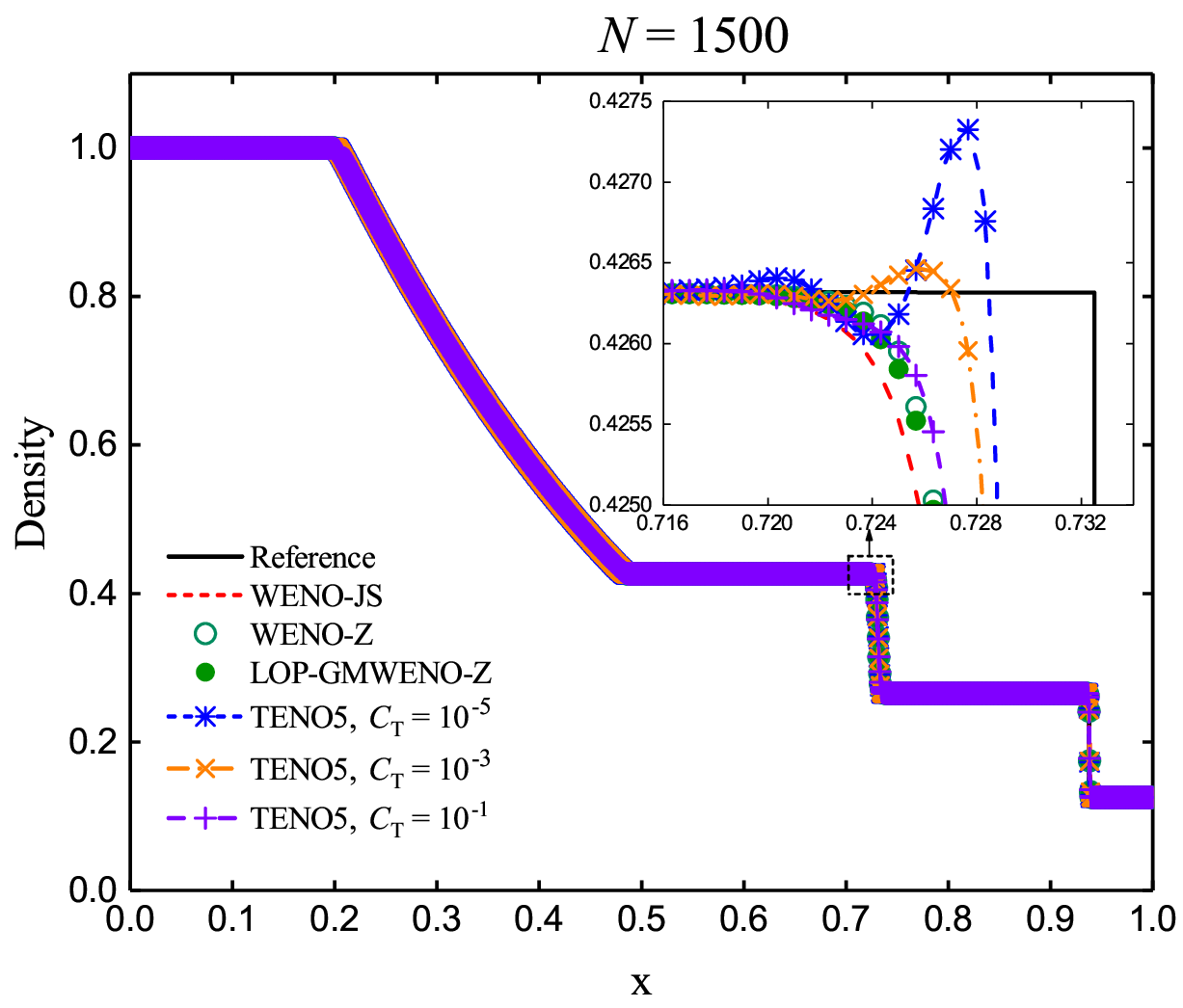}       \\ 
  \includegraphics[height=0.26\textwidth]
  {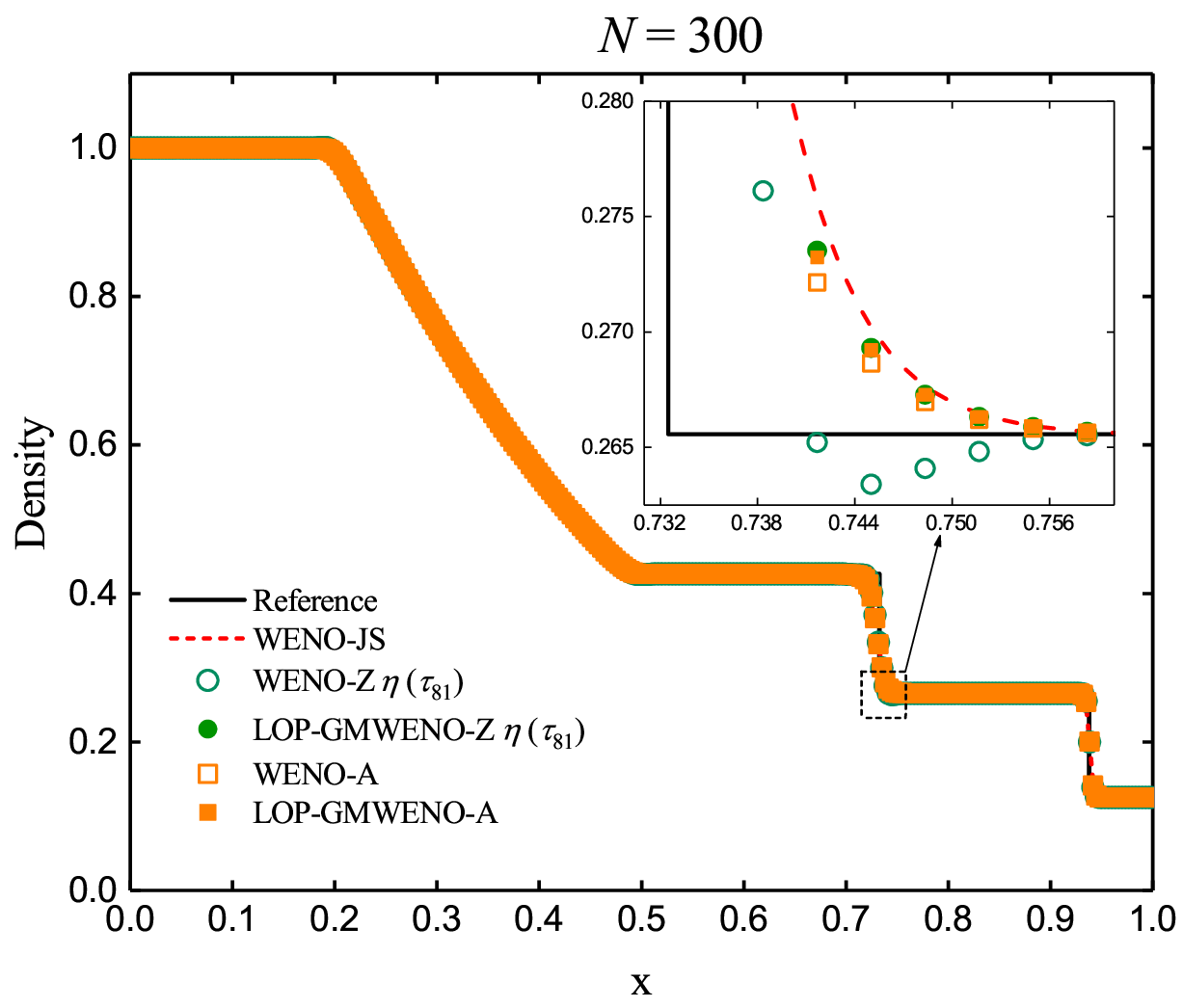}          \hspace{1.2ex}
  \includegraphics[height=0.26\textwidth]
  {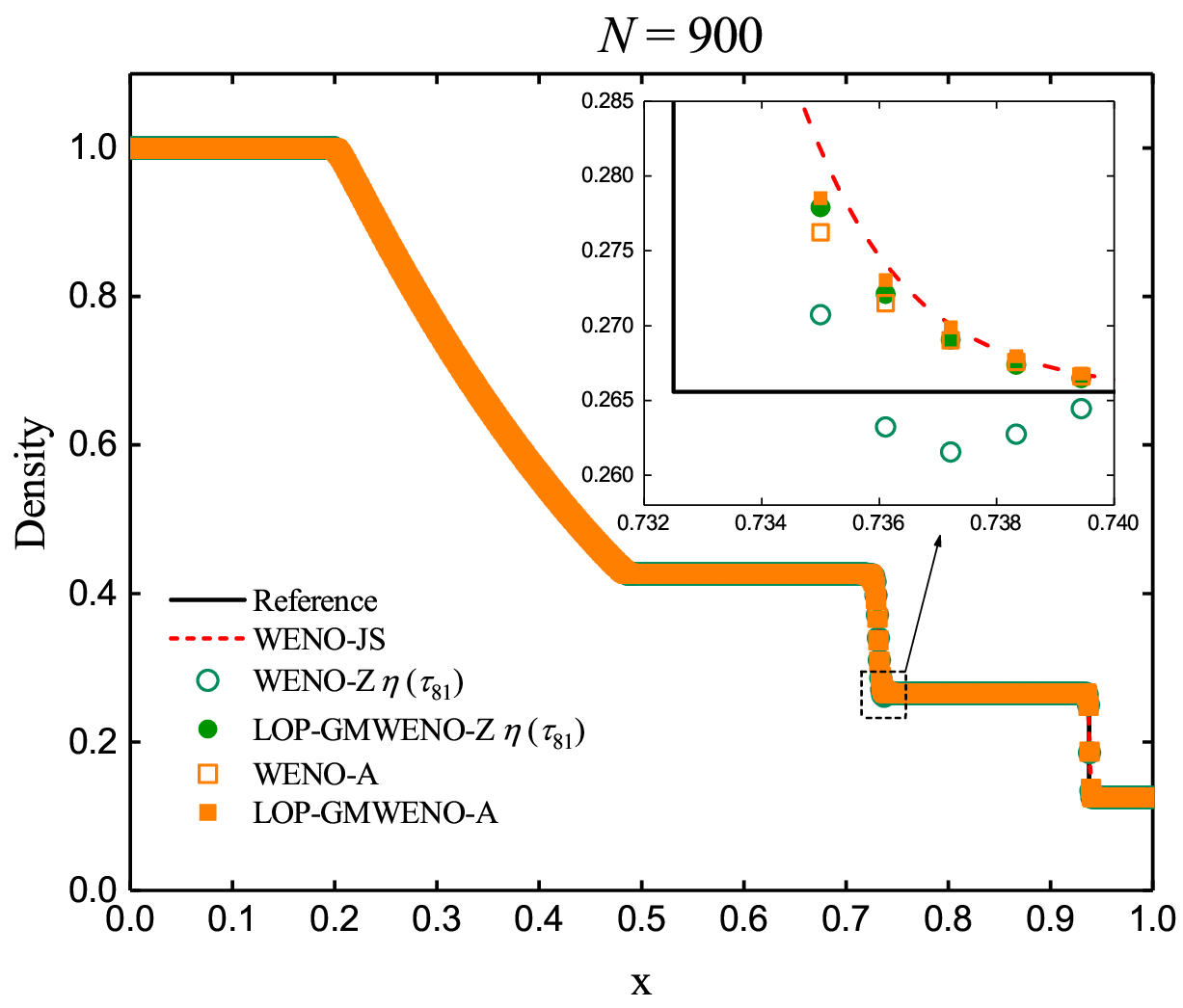}         \hspace{1.2ex}
  \includegraphics[height=0.26\textwidth]
  {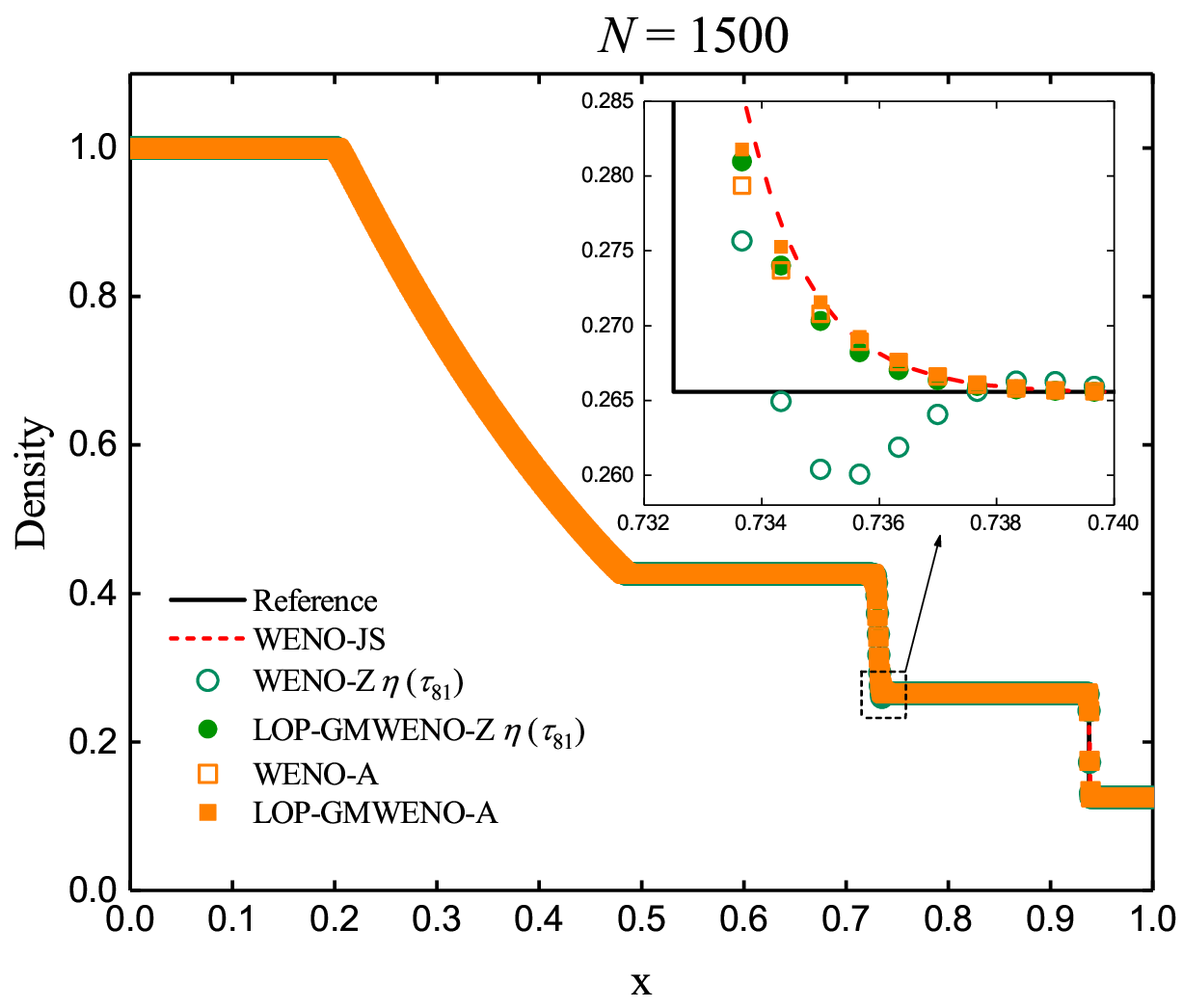}   
\caption{Density plots for the Sod shock-tube problem on various coarse mesh resolutions.}
\label{fig:ex:Sod}
\end{figure}

\begin{figure}[!ht]
\centering
  \includegraphics[height=0.26\textwidth]
  {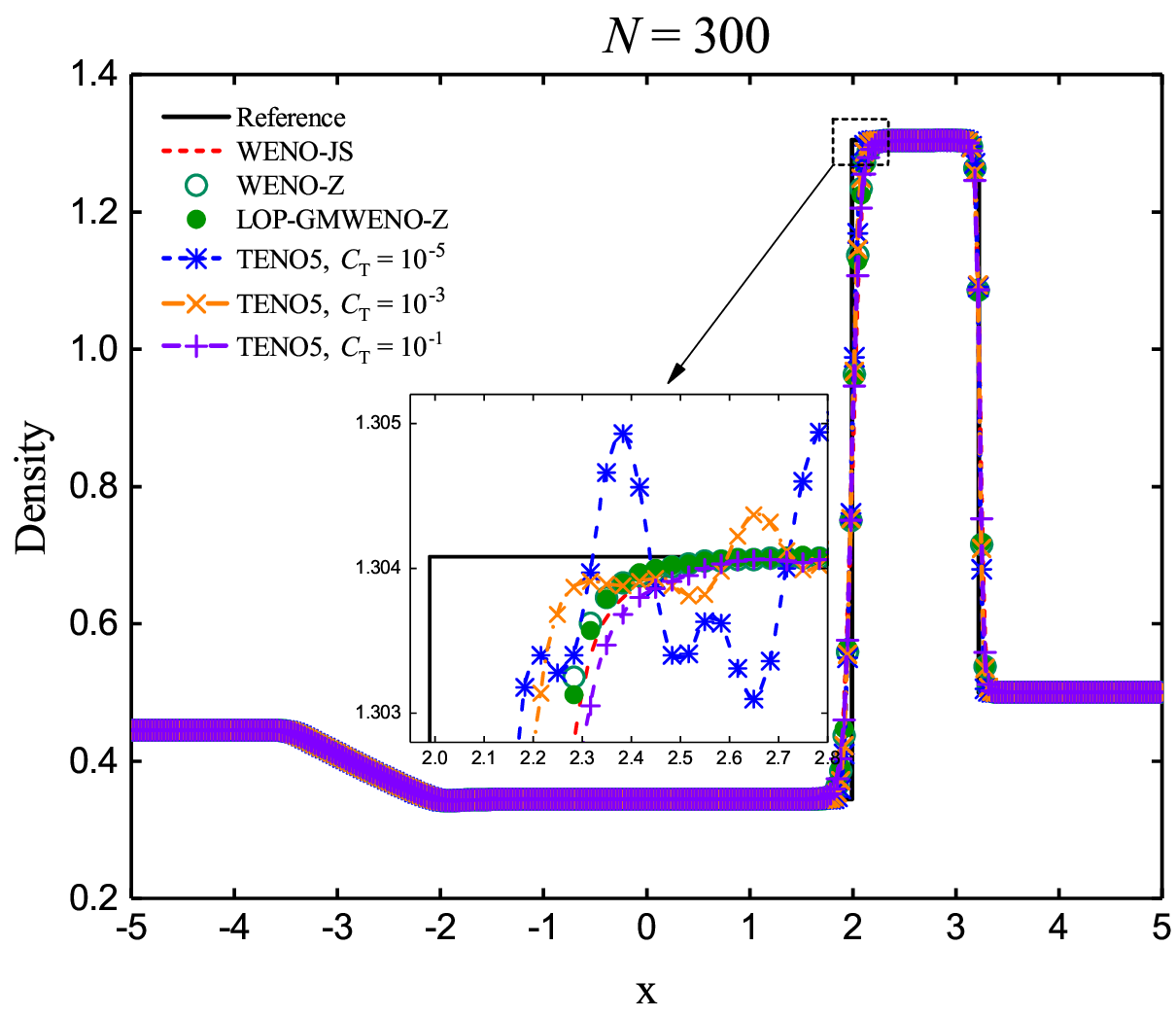}          \hspace{1.2ex}
  \includegraphics[height=0.26\textwidth]
  {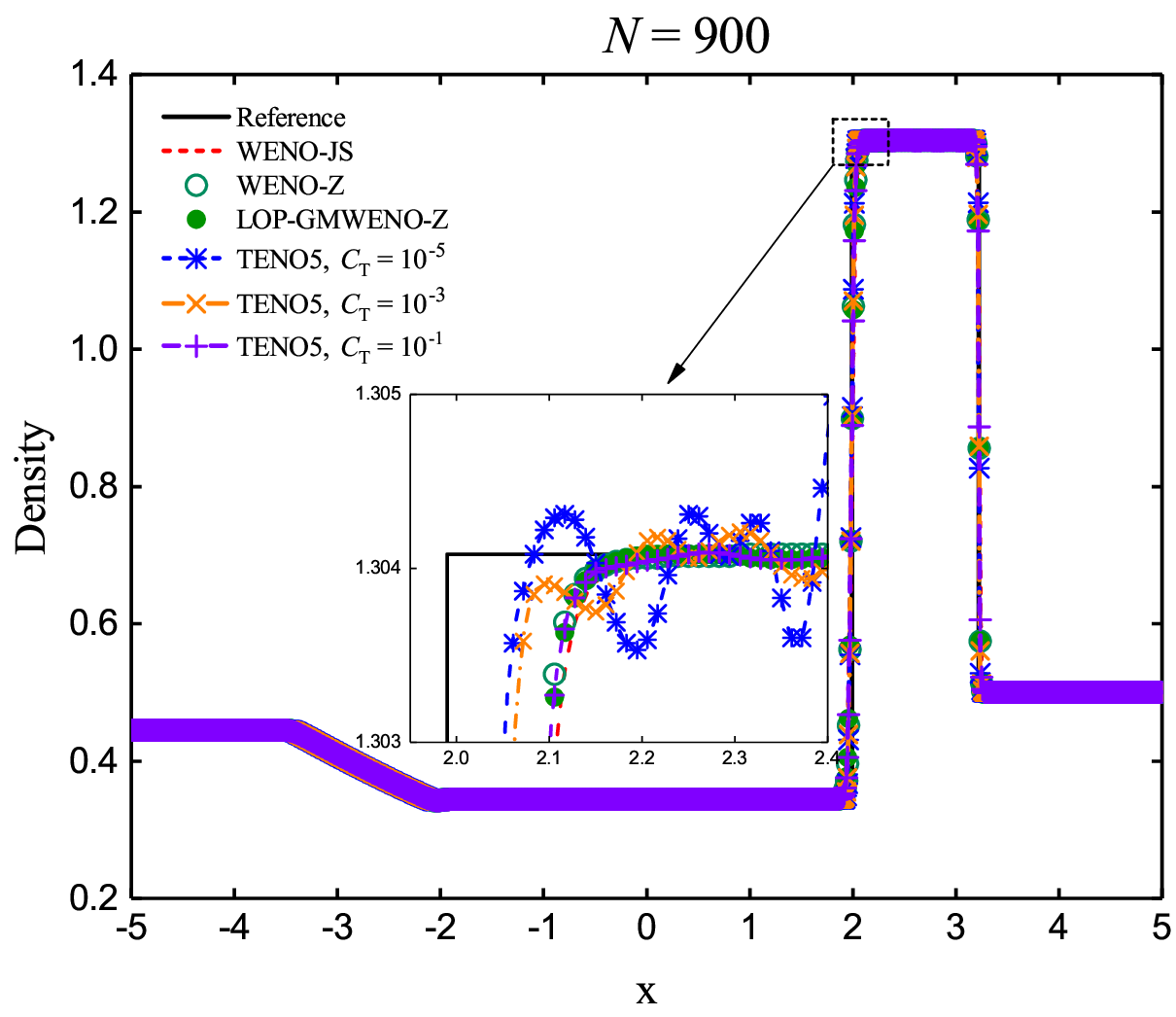}         \hspace{1.2ex}
  \includegraphics[height=0.26\textwidth]
  {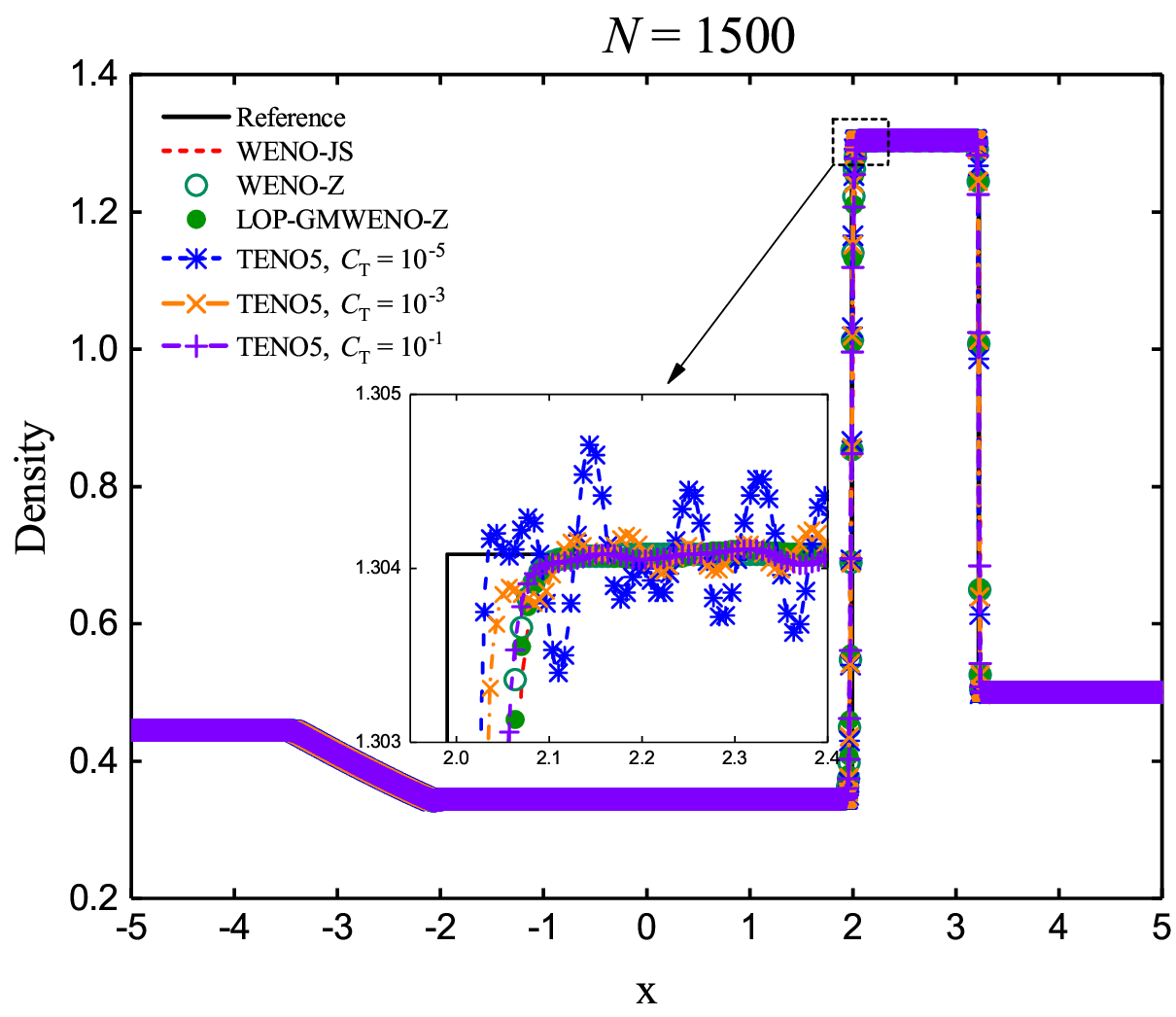}       \\ 
  \includegraphics[height=0.26\textwidth]
  {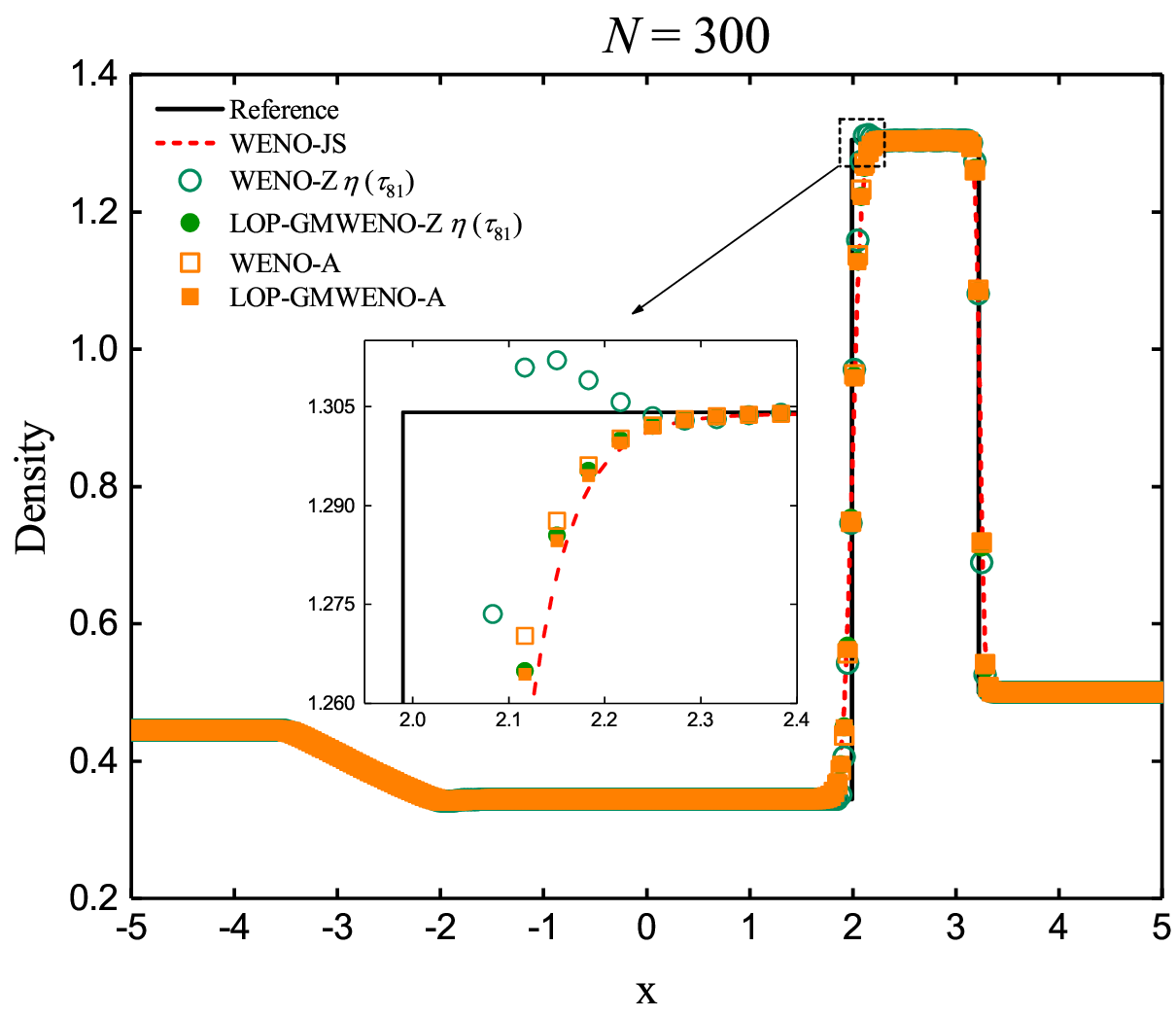}          \hspace{1.2ex}
  \includegraphics[height=0.26\textwidth]
  {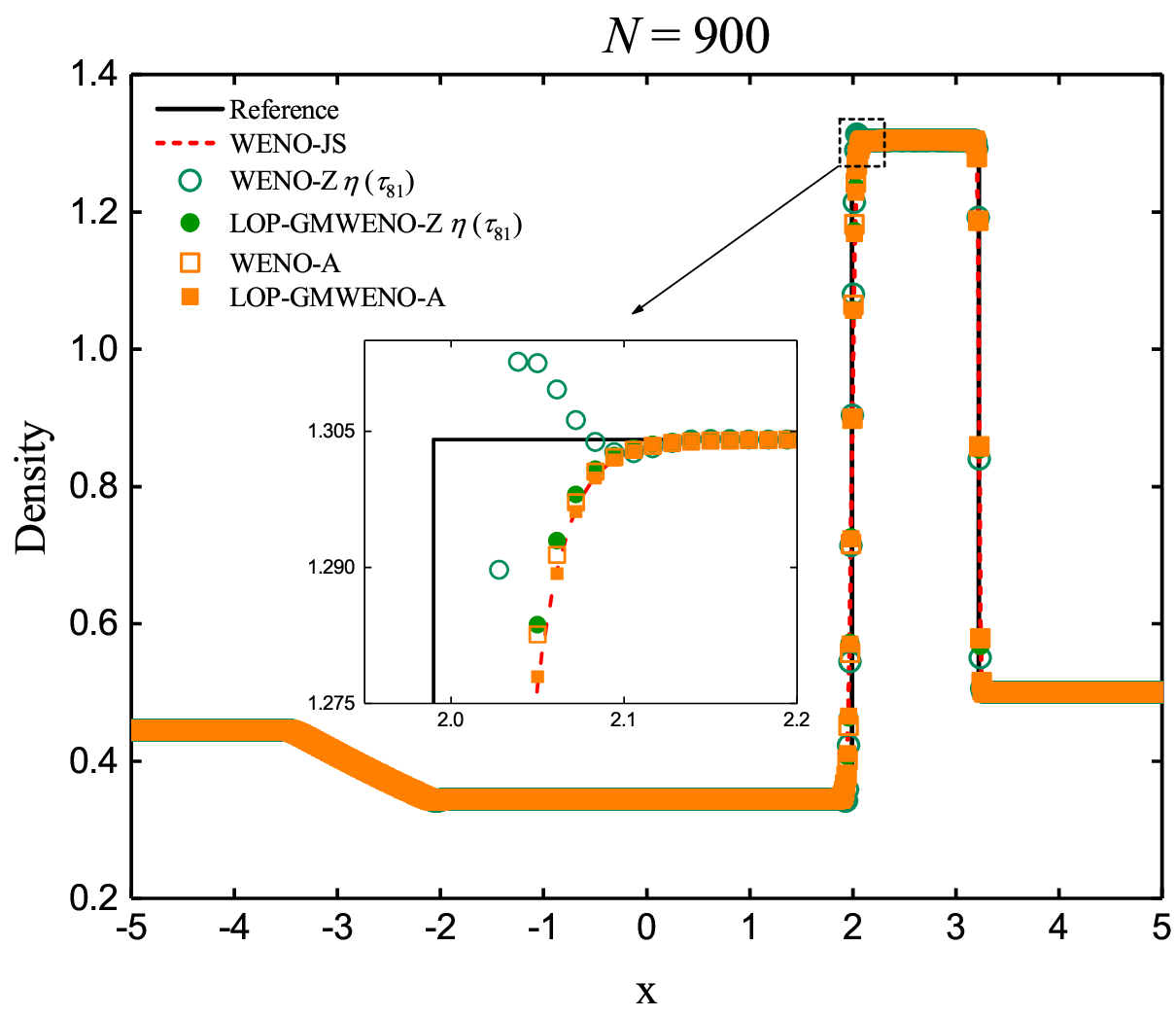}         \hspace{1.2ex}
  \includegraphics[height=0.26\textwidth]
  {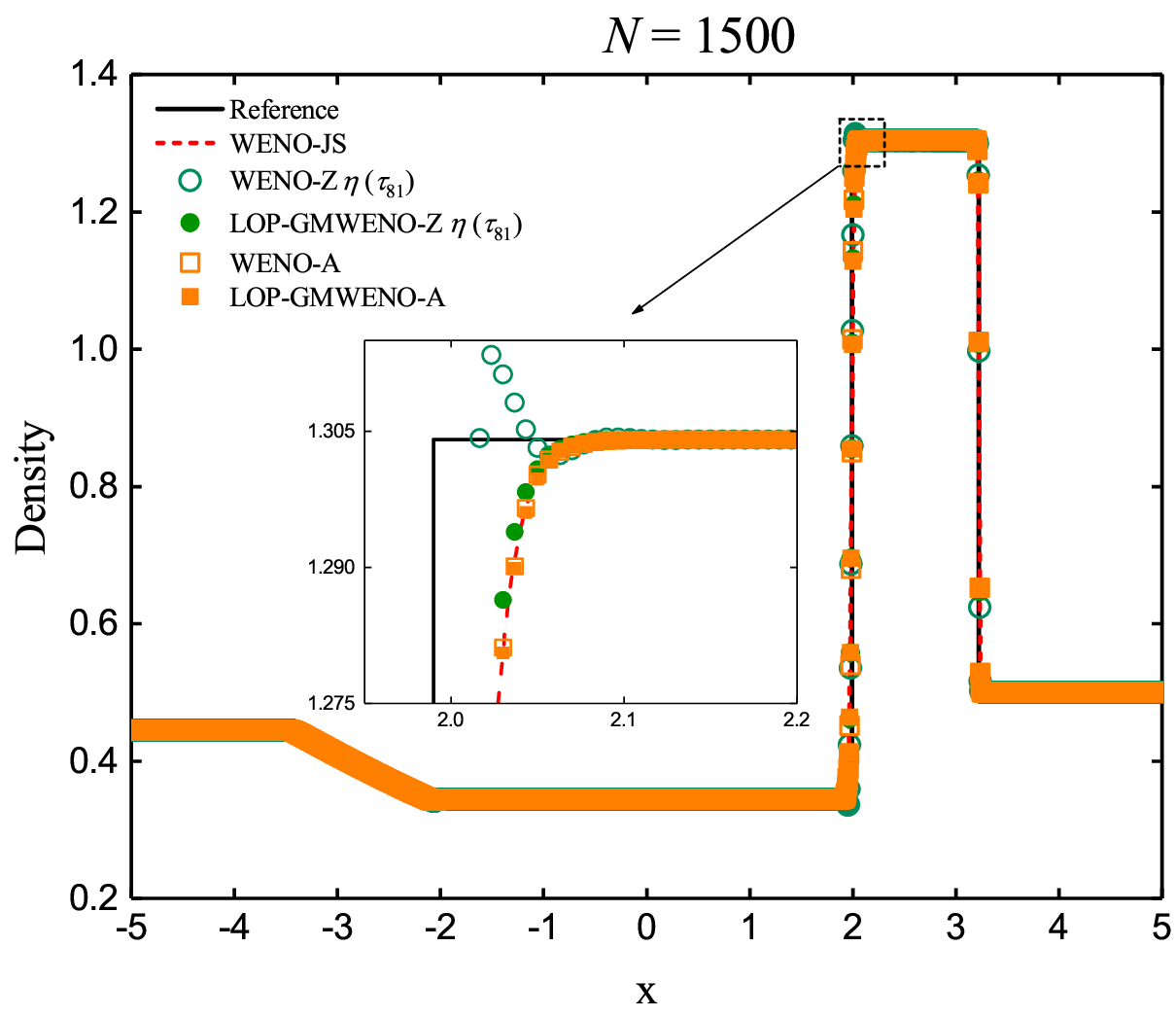}
\caption{Density plots for the Lax shock-tube problem on various coarse mesh resolutions.}
\label{fig:ex:Lax}
\end{figure}

Fig. \ref{fig:ex:Sod} and Fig. \ref{fig:ex:Lax} show the computed 
density profiles of the considered schemes. It is observed that the 
LOP-GMWENO-Z and LOP-GMWENO-A schemes give results with comparable 
or slightly lower resolutions than their associated WENO-Z and 
WENO-A schemes but their resolutions are better than the WENO-JS 
scheme. The resolution of the LOP-GMWENO-Z$\eta(\tau_{81})$ scheme 
appears to be lower than its associated WENO-Z$\eta(\tau_{81})$. 
However, if we take a closer look, we can see that the WENO-Z$\eta(\tau_{81})$ scheme generates spurious oscillations while the 
LOG-GMWENO-Z$\eta(\tau_{81})$ scheme can avoid the spurious 
oscillations successfully. Similarly, the TENO5 schemes with the 
threshold $C_{T}$ as $10^{-5}$ and $10^{-3}$ attain the much higher 
resolutions than other schemes while they produce very slightly 
spurious oscillations (a very close look is needed to observe this). 
One can remove these spurious oscillations by increasing the 
threshold $C_{T}$. Actually, when $C_{T} = 10^{-1}$ is taken, the 
spurious oscillations have been removed successfully at the price of 
the loss of resolutions.

\subsubsection{Shock-density wave interaction}
To demonstrate the excellent performance of the LOP-GMWENO-X schemes 
in the region with high-frequency but smooth waves, we simulate two 
typical shock-density wave interaction problems, that is, the 
Shu-Osher problem \cite{ENO-Shu1989} and the Titarev-Toro problem \cite{Titarev-Toro-1, Titarev-Toro-2,Titarev-Toro-3}. The initial 
conditions are
\begin{equation*}
\begin{aligned}
\text{Shu-Osher: } \qquad &\big( \rho, u, p \big)(x, 0) =\left\{
\begin{array}{ll}
(3.857143, 2.629369, 10.333333), & x \in [-5.0, -4.0], \\
(1.0 + 0.2\sin(5x), 0, 1), & x \in [-4.0, 5.0],
\end{array}\right. \\
\text{Titarev-Toro: } \qquad &\big( \rho, u, p \big)(x, 0) =\left\{
\begin{array}{ll}
(1.515695, 0.5233346, 1.80500), & x \in [-5.0, -4.5], \\
(1.0 + 0.1\sin(20\pi x), 0, 1), & x \in [-4.5, 5.0].
\end{array}\right.
\end{aligned}
\label{ex:shock-densityInteraction}
\end{equation*}
The transmissive boundary conditions are used at $x = \pm 5$.
The output times are set to be 1.8 and 5.0, and the uniform mesh 
sizes of $N = 300$ and $N = 1500$ are used for the Shu-Osher problem 
and the Titarev-Toro problem, respectively.

\begin{figure}[!ht]
\centering
  \includegraphics[height=0.35\textwidth]
  {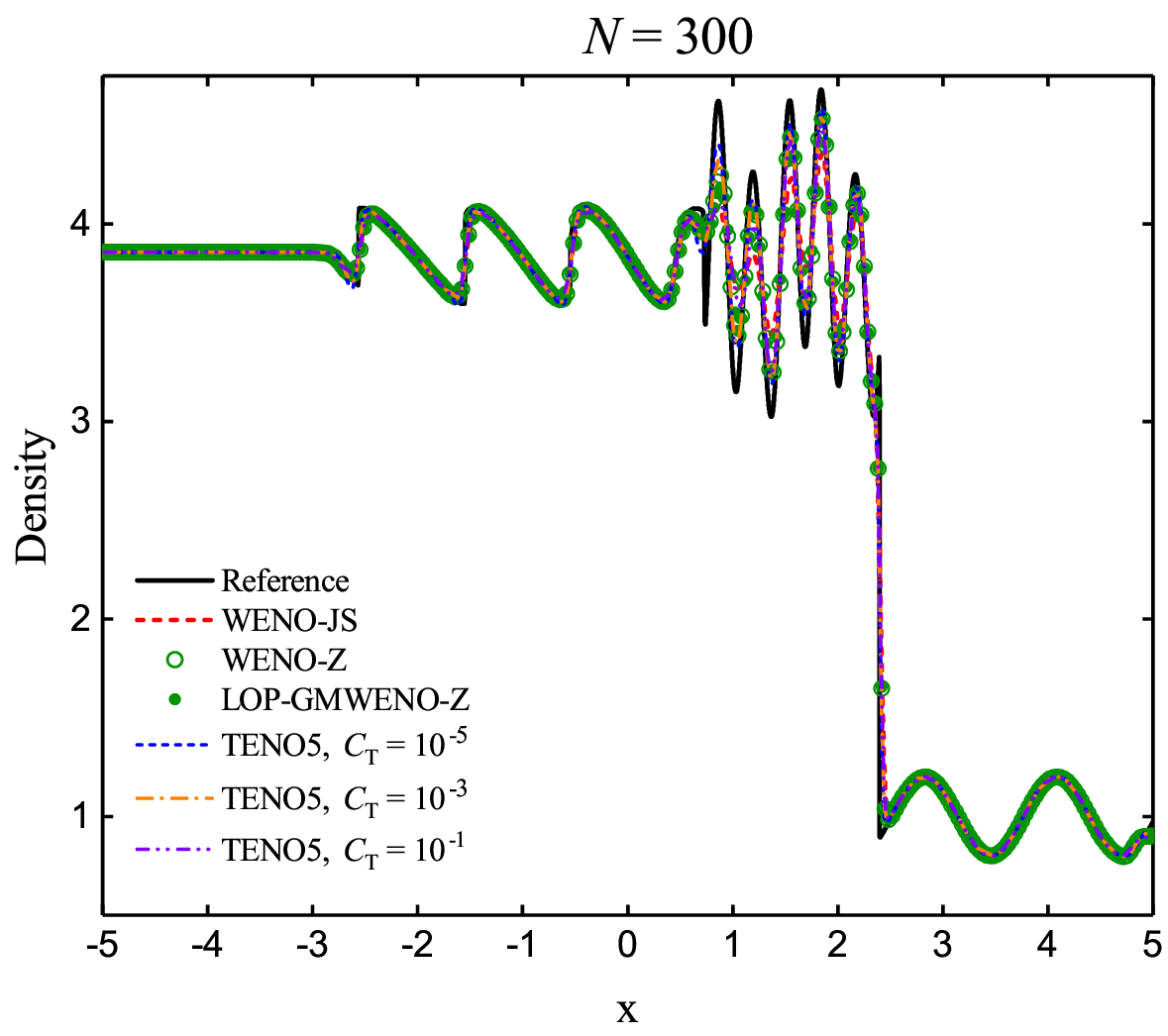}          \hspace{1.2ex}
  \includegraphics[height=0.35\textwidth]
  {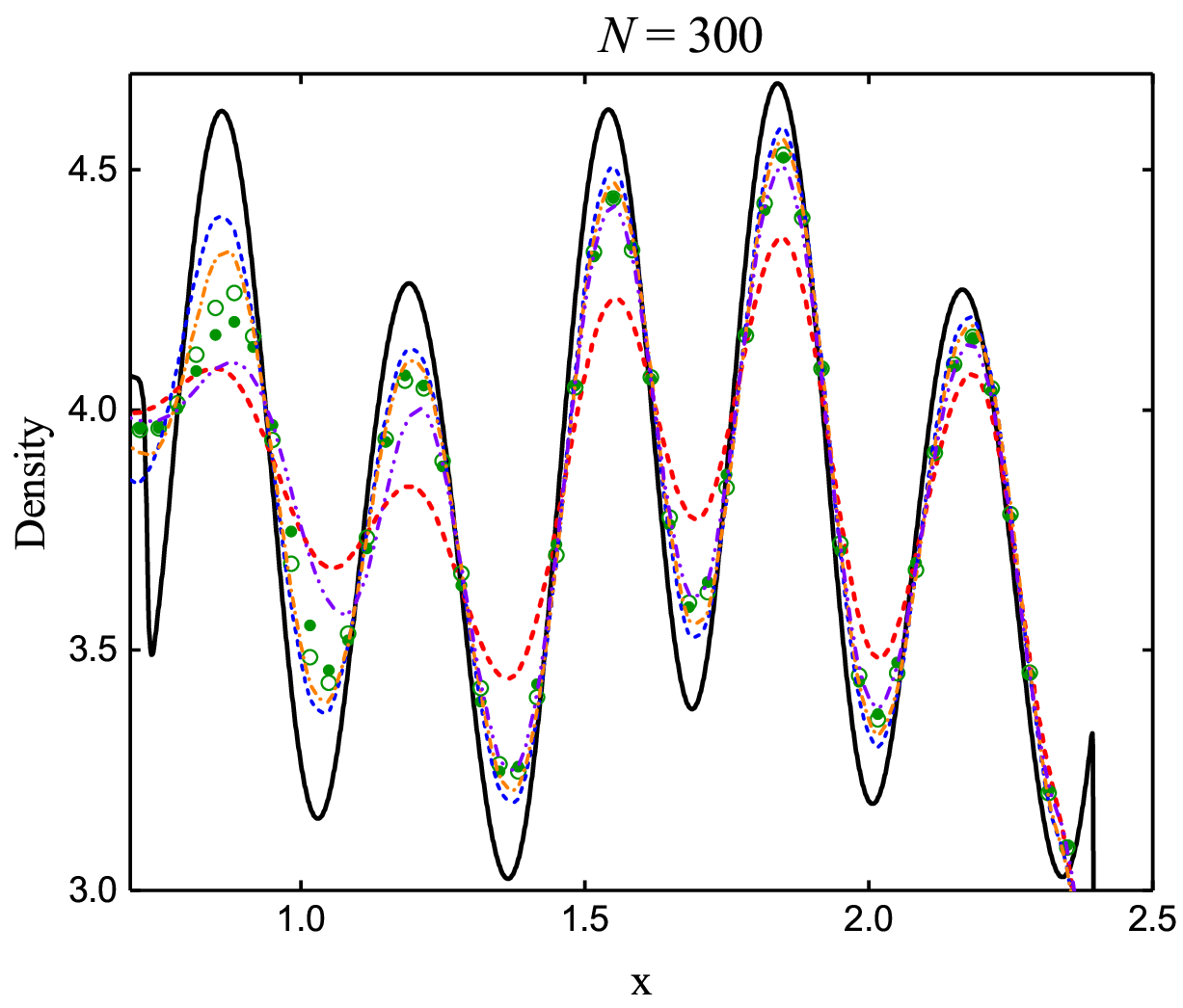}     \\ 
  \includegraphics[height=0.35\textwidth]
  {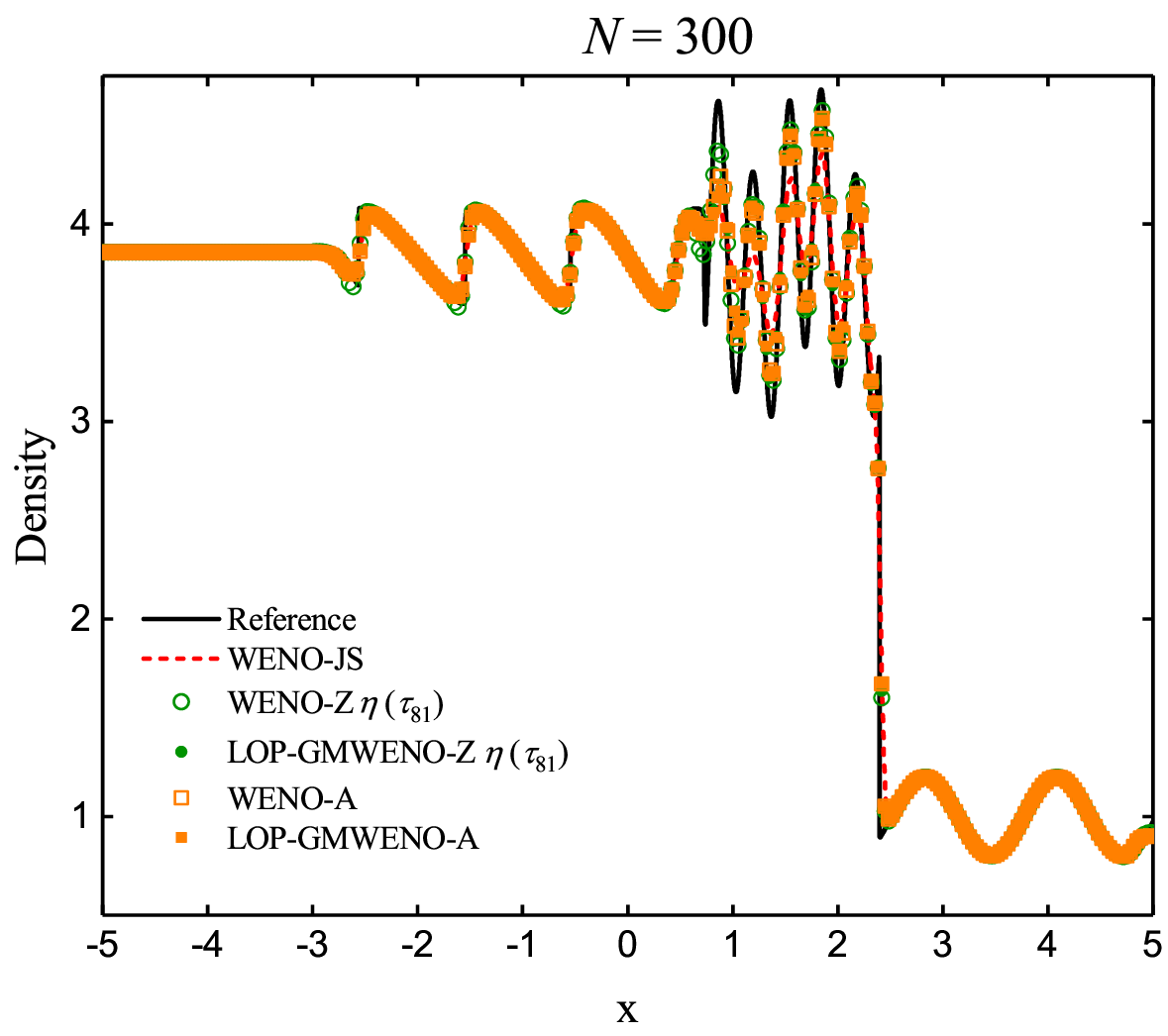}          \hspace{1.2ex}
  \includegraphics[height=0.35\textwidth]
  {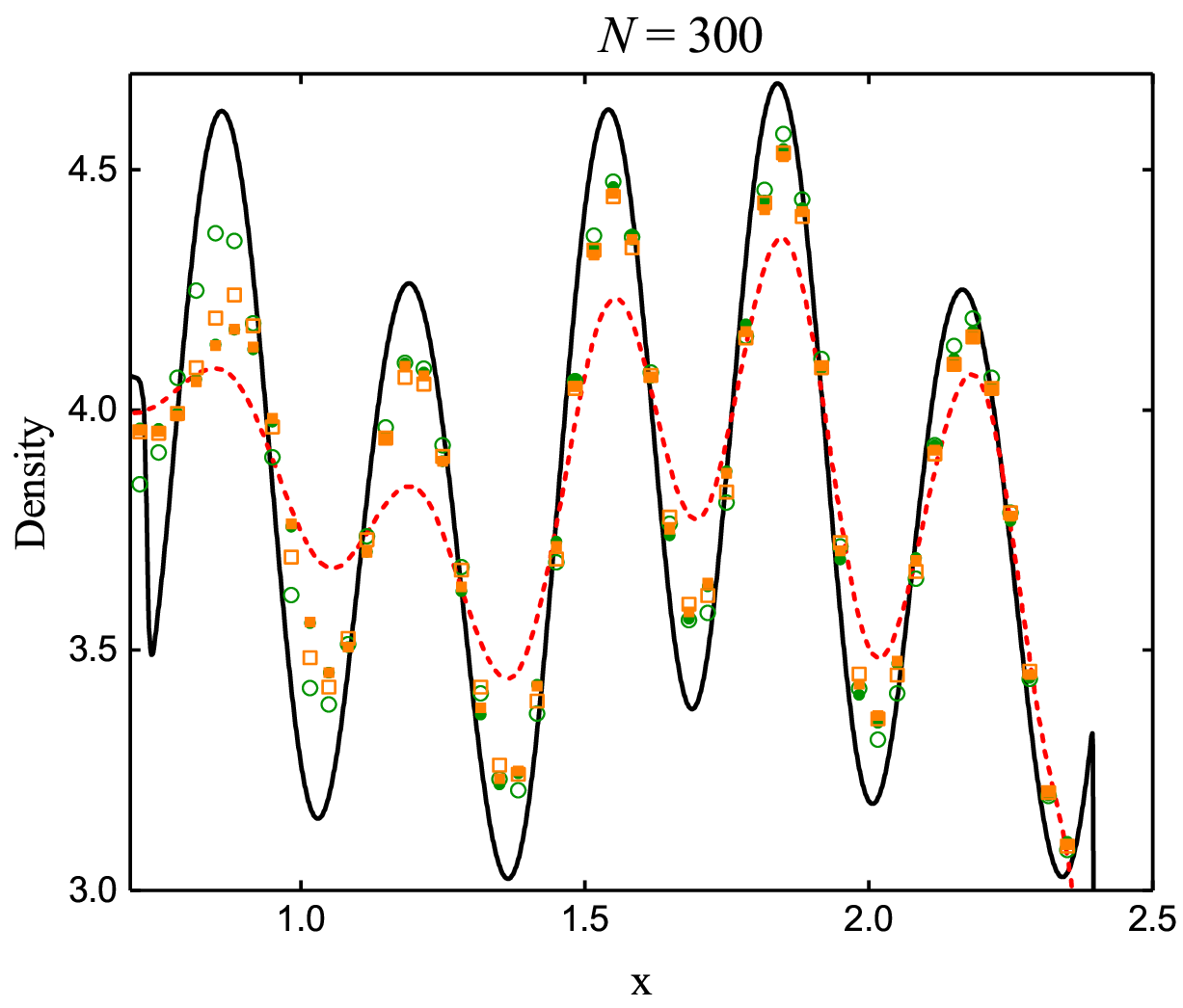}     \\  
\caption{Density plots for the Shu-Osher shock-density wave interaction problem.}
\label{fig:ex:Shu-Osher}
\end{figure}

\begin{figure}[!ht]
\centering
  \includegraphics[height=0.35\textwidth]
  {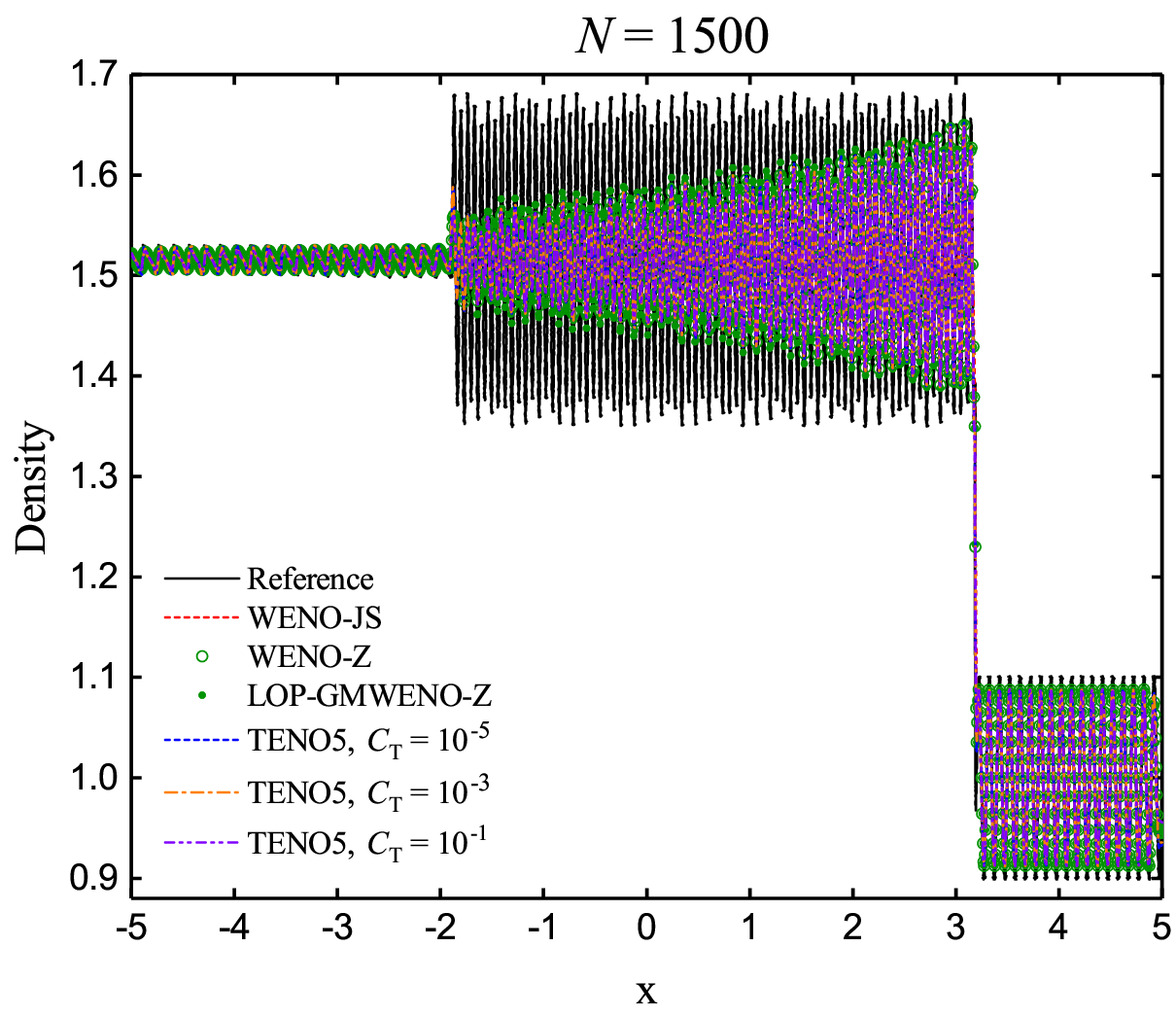}          \hspace{1.2ex}
  \includegraphics[height=0.35\textwidth]
  {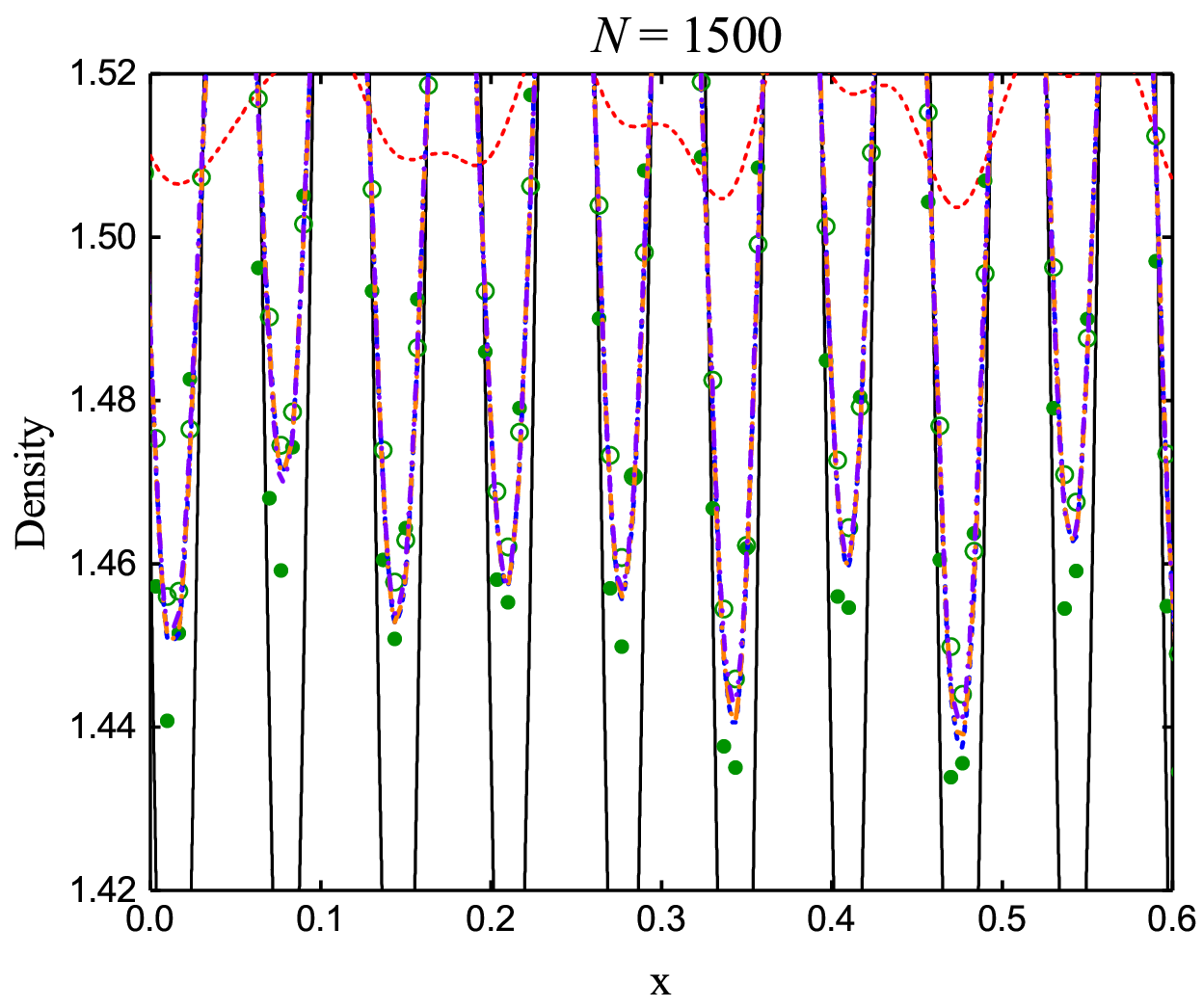}     \\ 
  \includegraphics[height=0.35\textwidth]
  {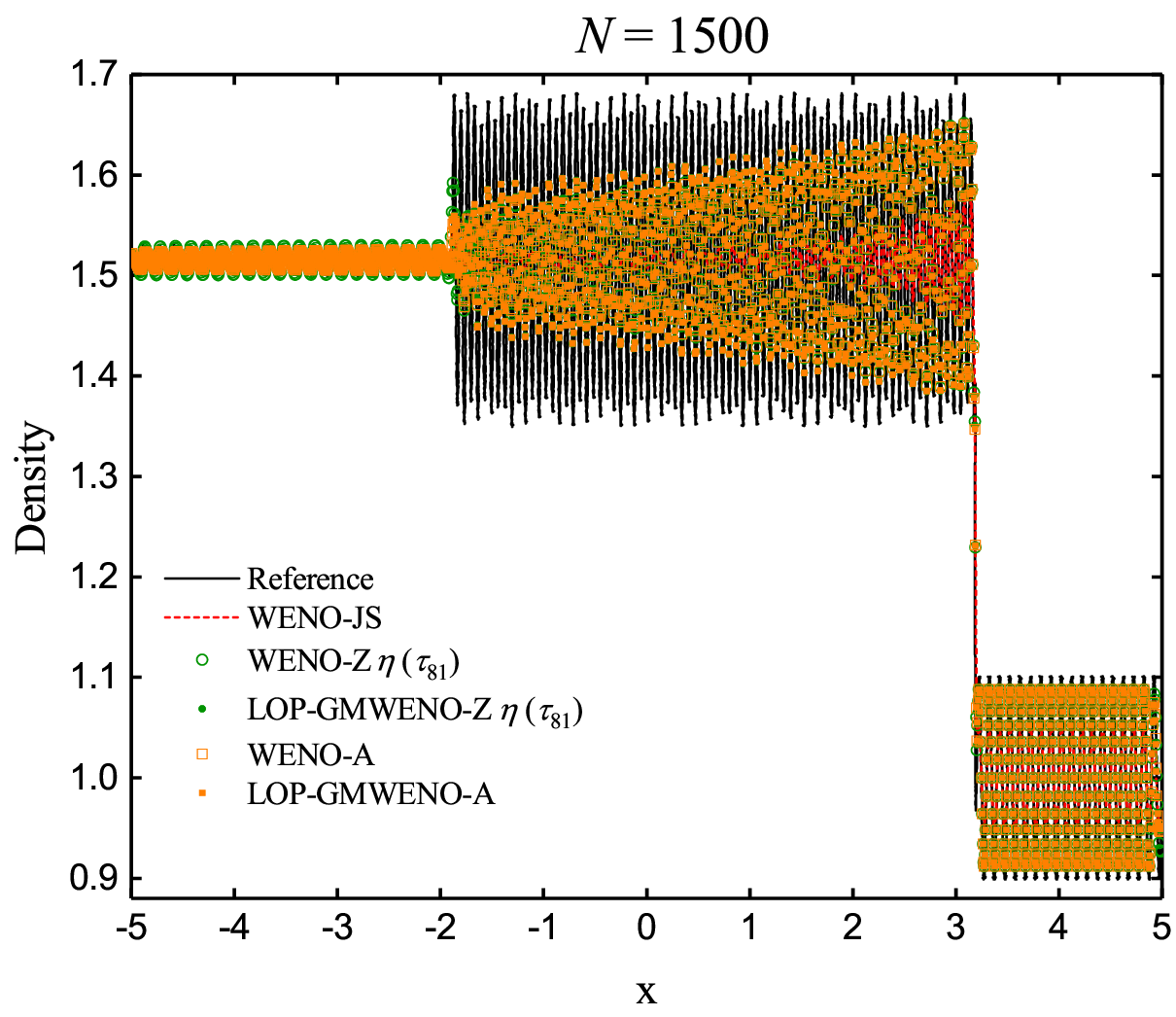}          \hspace{1.2ex}
  \includegraphics[height=0.35\textwidth]
  {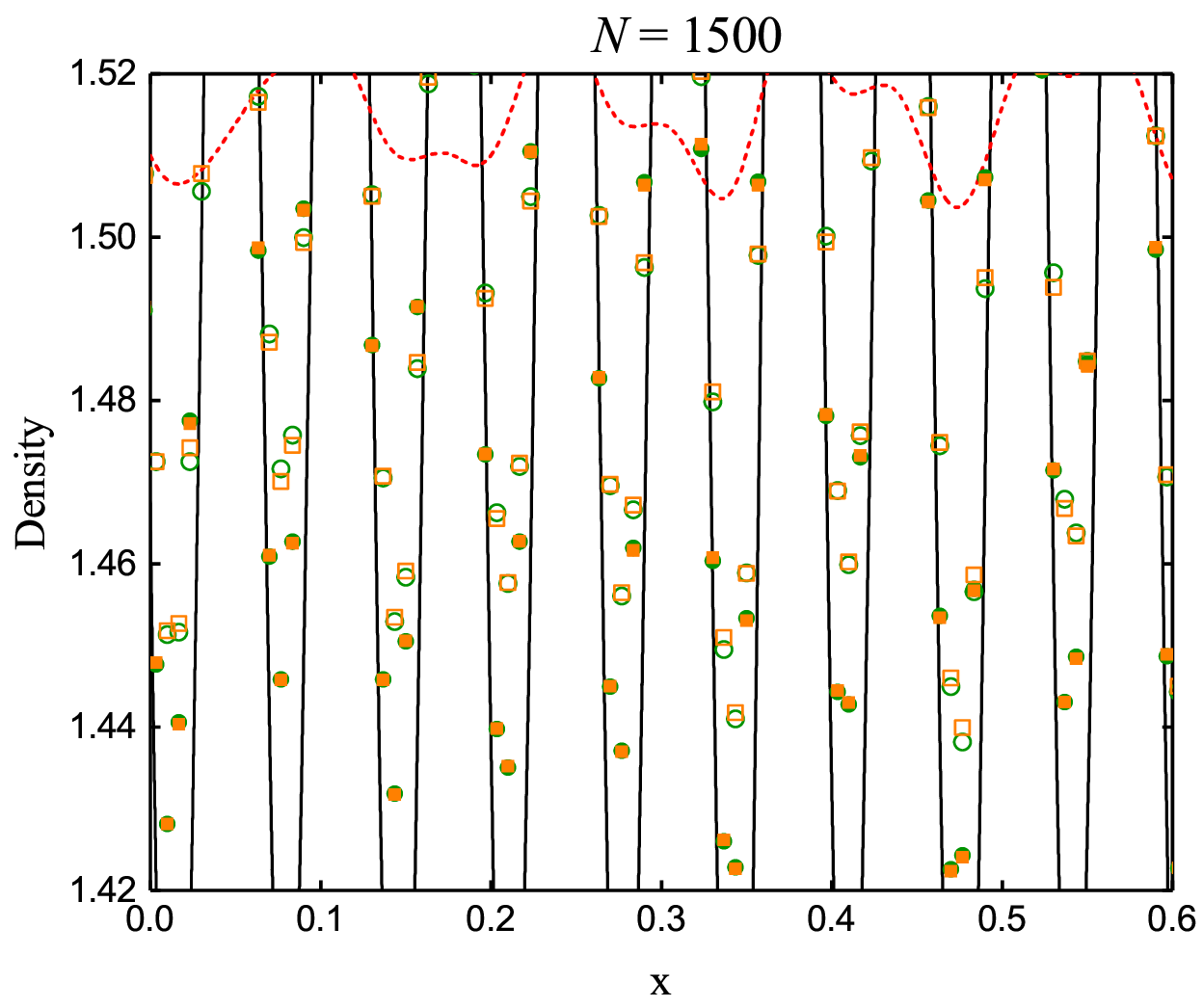}     \\  
\caption{Density plots for the Titarev-Toro shock-density wave interaction problem.}
\label{fig:ex:Titarev-Toro}
\end{figure}

The solutions of density of the LOP-GMWENO-Z, LOP-GMWENO-Z$\eta(\tau_{81})$, LOP-GMWENO-A schemes, as well as their associated WENO-Z, 
WENO-Z$\eta(\tau_{81})$, WENO-A schemes, are given in Fig. \ref{fig:ex:Shu-Osher} and Fig. \ref{fig:ex:Titarev-Toro}, where the 
reference solutions are computed by employing WENO-JS with $N = 10000$. Again, for comparison purpose, we show the solutions of the 
aforementioned TENO5 schemes and that of WENO-JS. As shown in Fig. \ref{fig:ex:Shu-Osher}, for the Shu-Osher problem, the LOP-GMWENO-Z, 
LOP-GMWENO-Z$\eta(\tau_{81})$ and LOP-GMWENO-A schemes provide 
comparable results with those of their associated WENO-Z, WENO-Z$\eta(\tau_{81})$ and WENO-A schemes which are far better than that of 
the WENO-JS scheme, and the TENO5 schemes with the threshold $C_{T}$ 
as $10^{-5}$ and $10^{-3}$ enjoy higher resolutions because of their 
low-dissipation property. As expected, when $C_{T} = 10^{-1}$ is 
taken, its resolution decreases to the same or slightly lower level 
of the LOP-GMWENO-Z scheme. For the Titarev-Toro problem, from Fig. \ref{fig:ex:Titarev-Toro}, we can find that the LOP-GMWENO-Z, 
LOP-GMWENO-Z$\eta(\tau_{81})$ and LOP-GMWENO-A schemes achieve much 
higher resolutions than their associated WENO-Z, WENO-Z$\eta(\tau_{81})$ and WENO-A schemes. Moreover, the considered TENO5 schemes can 
obtain better resolutions than the WENO-Z scheme, while their 
resolutions are evidently lower than those of all the considered 
LOP-GMWENO-X schemes.

\subsection{Two-dimensional Euler system}
\label{subsec:2DEuler}
In this subsection, we present numerical experiments with the 
two-dimensional Euler equations, such as the accuracy tests for 
smooth problems, 2D interacting blast waves problem, implosion and 
explosion problems, shock-vortex interaction (SVI) problem, regular 
shock reflection (RSR) problem, double Mach reflection (DMR) 
problem, forward facing step (FFS) problem and Rayleigh-Taylor 
instability (RTI) problem.

\subsubsection{Accuracy test for 2D smooth Euler problems}
We solve the 2D smooth Euler problem on the computational domain $[-1.0, 1.0]\times[-1.0, 1.0]$ to demonstrate the convergence property 
of the present scheme, and the density wave propagation problem \cite{AccuracyTest-Euler2D} with the following two initial conditions are 
calculated
\begin{equation}
\begin{aligned}
\text{Case 1: } \qquad &(\rho, u, v, p)(x, y, 0) = \Big( 1.0 + 0.2 \sin\big( \pi(x + y) \big), 0.7, 0.3, 1.0\Big), \\
\text{Case 2: } \qquad &(\rho, u, v, p)(x, y, 0) = \Bigg( 1.0 + 0.2 \sin\bigg( \pi(x + y) - \dfrac{\sin(\pi(x + y))}{\pi} \bigg), 0.7, 0.3, 1.0\Bigg).
\end{aligned}
\label{ex:Euler2D_smooth}
\end{equation}
The periodic boundary condition is used. We set $h = h_{x} = h_{y}$ 
in all calculations of this paper with $h_{x}, h_{y}$ as the uniform 
spatial step size in $x$- and $y$- direction. The CFL number is set 
to be $h^{2/3}$ so that the error for the overall scheme is a 
measure of the spatial convergence only. The calculation is advanced 
to $t = 2$.

To test the convergence orders, the following $L_{\infty}$-norm 
error of the density is computed
\begin{equation*}
\displaystyle
\begin{aligned}
L_{\infty} = \displaystyle\max_{\substack{1\leq i\leq N_{x}\\1\leq j
\leq N_{y}}} \Big\lvert \rho_{i,j}^{\mathrm{exact}} -(\rho_{h})_{i,j}
\Big\rvert,
\end{aligned}
\end{equation*}
where $N_{x}, N_{y}$ is the number of cells in $x-$ and $y-$ 
direction. $(\rho_{h})_{i,j}$ is the numerical solution of the 
density and $\rho_{i,j}^{\mathrm{exact}}$ is its exact solution.

\begin{figure}[!ht]
\centering
  \includegraphics[height=0.26\textwidth]
  {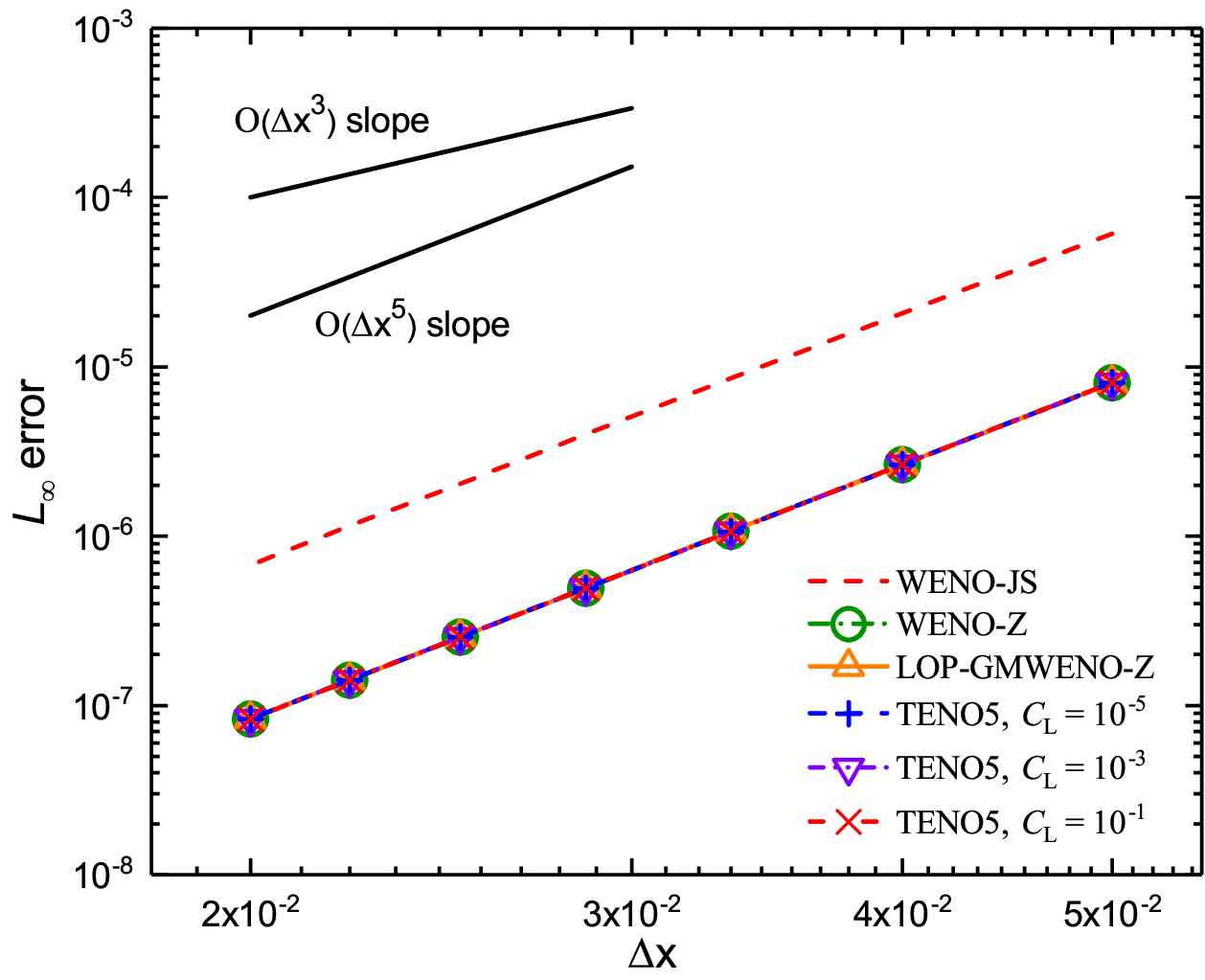} 
  \includegraphics[height=0.26\textwidth]
  {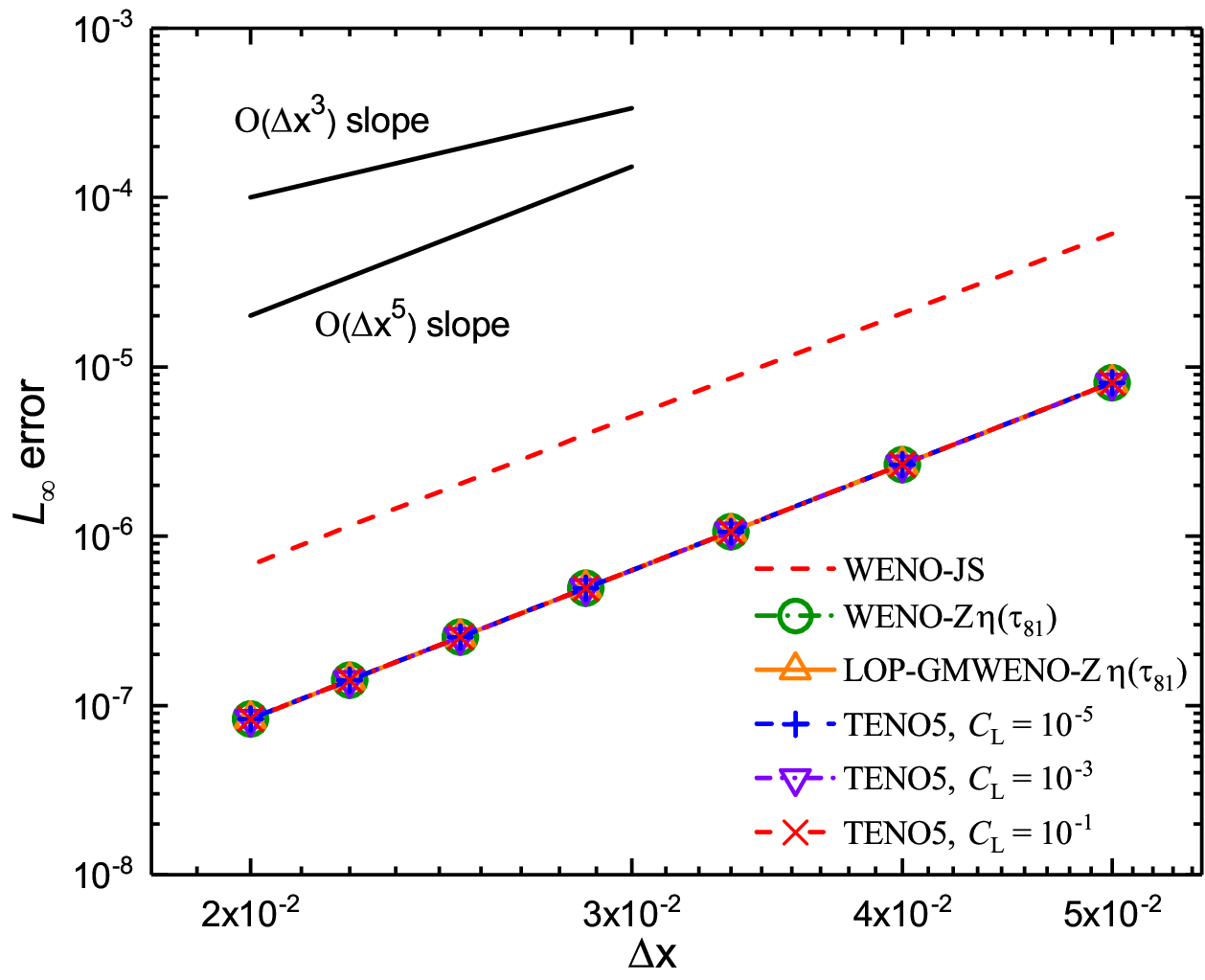}     
  \includegraphics[height=0.26\textwidth]
  {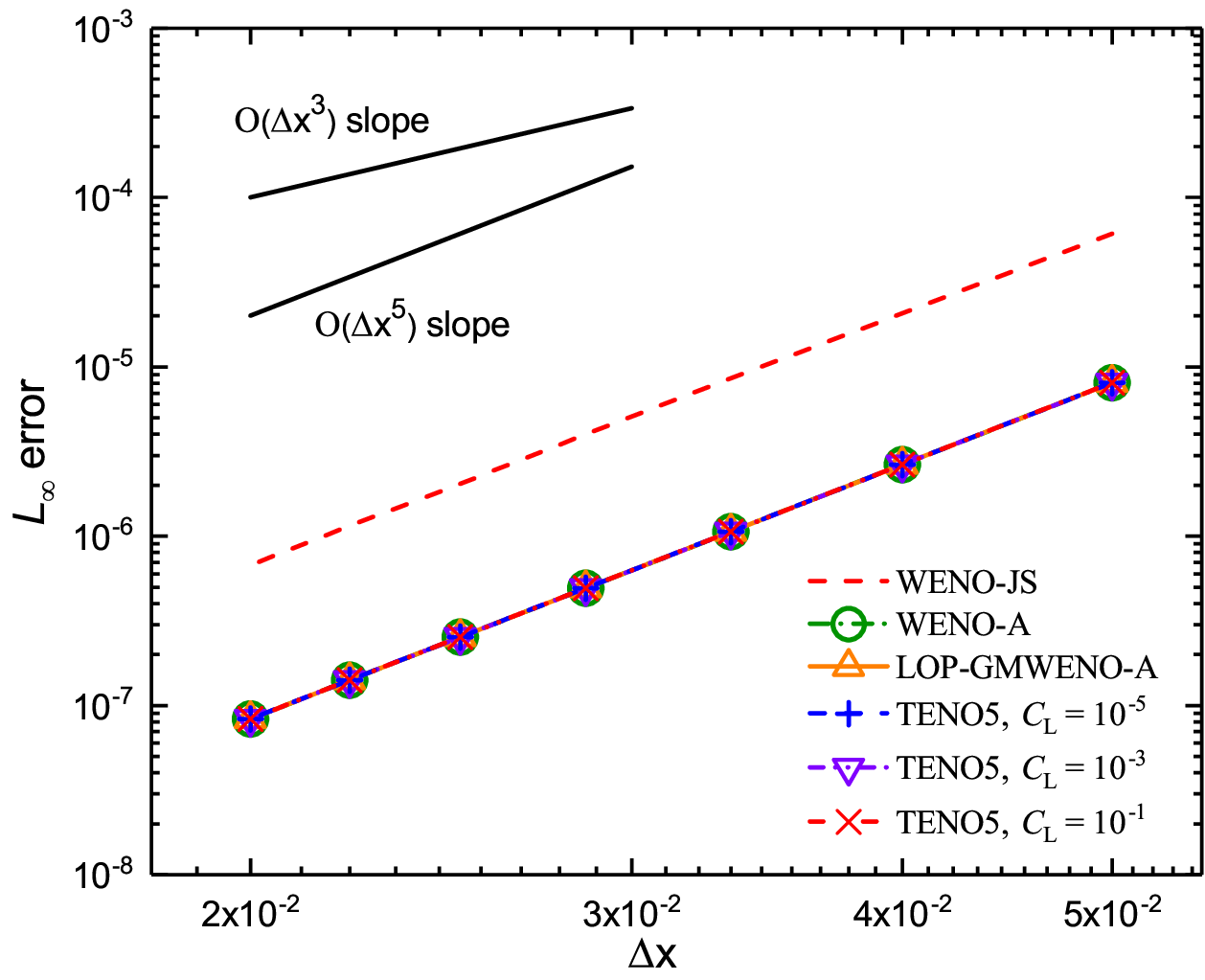} \\
  \includegraphics[height=0.26\textwidth]
  {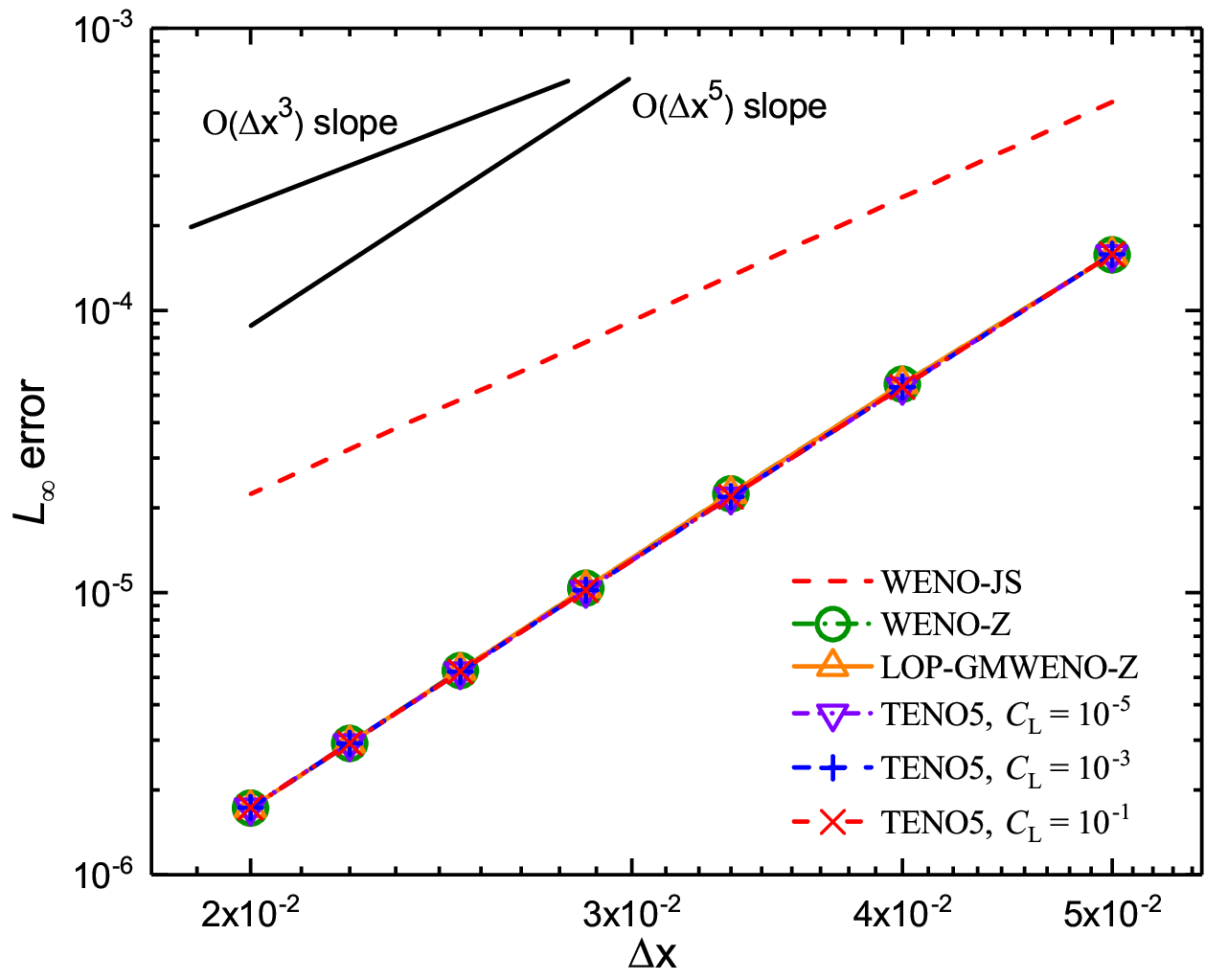} 
  \includegraphics[height=0.26\textwidth]
  {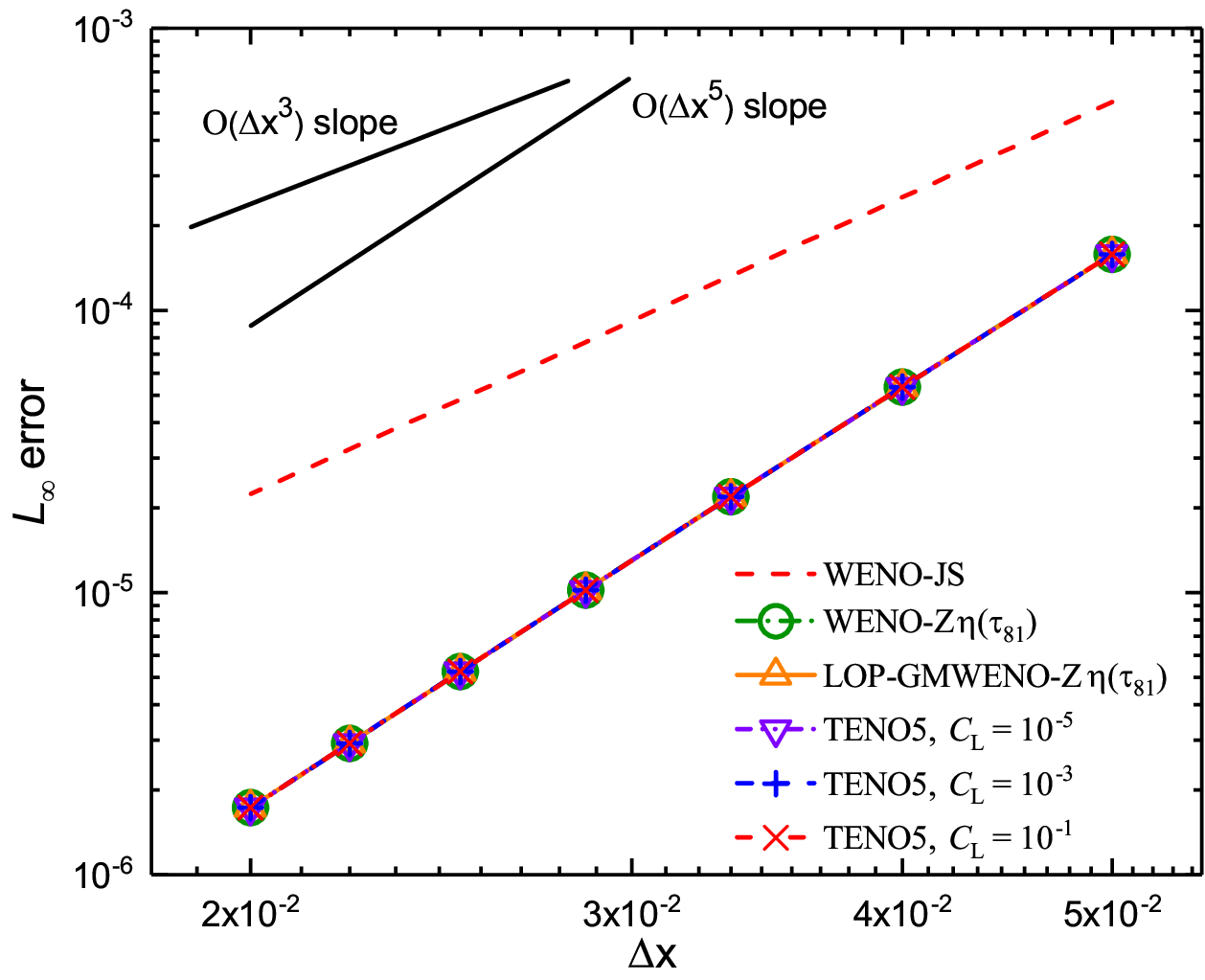}     
  \includegraphics[height=0.26\textwidth]
  {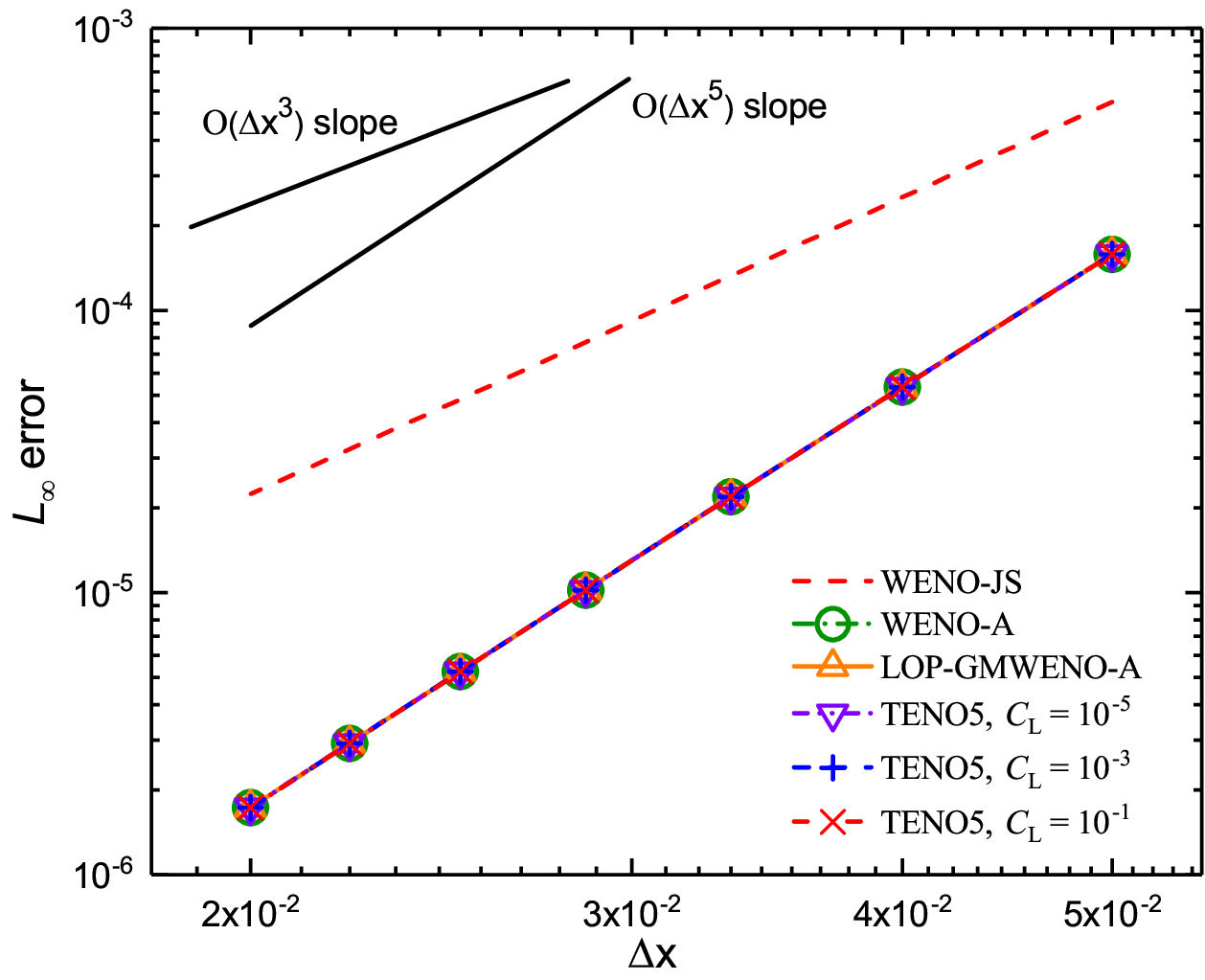}   
\caption{Convergence of the $L_{\infty}$-norm errors of various WENO schemes for 2D smooth Euler problems. Top: Case 1. Bottom: Case 2.}
\label{fig:ex:AccuracyTest:Convergence}
\end{figure}

For comparison purpose, we also apply the aforementioned TENO5 
schemes and the WENO-JS scheme to this accuracy test. In Fig. \ref{fig:ex:AccuracyTest:Convergence}, we provide the overall $L_{\infty}$
convergence behavior of all considered schemes. Obviously, as 
evidenced by the slope of the profiles, all the LOP-GMWENO-X schemes 
and their associated WENO-X schemes, as well as the TENO5 schemes 
with the threshold $C_{T}$ as $10^{-5}, 10^{-3}, 10^{-1}$, can 
recover the optimal convergence orders for both Case 1 and Case 2. 
Although the WENO-JS scheme can also obtain the optimal convergence 
order for Case 1, it can only achieve third-order convergence rate 
of accuracy for Case 2. Moreover, in terms of accuracy, the results 
of the LOP-GMWENO-X schemes are almost identical to those of the 
associated WENO-X schemes and the TENO5 schemes, and all these 
schemes have much smaller numerical errors than the WENO-JS scheme.

\begin{figure}[!ht]
\centering
  \includegraphics[height=0.27\textwidth]
  {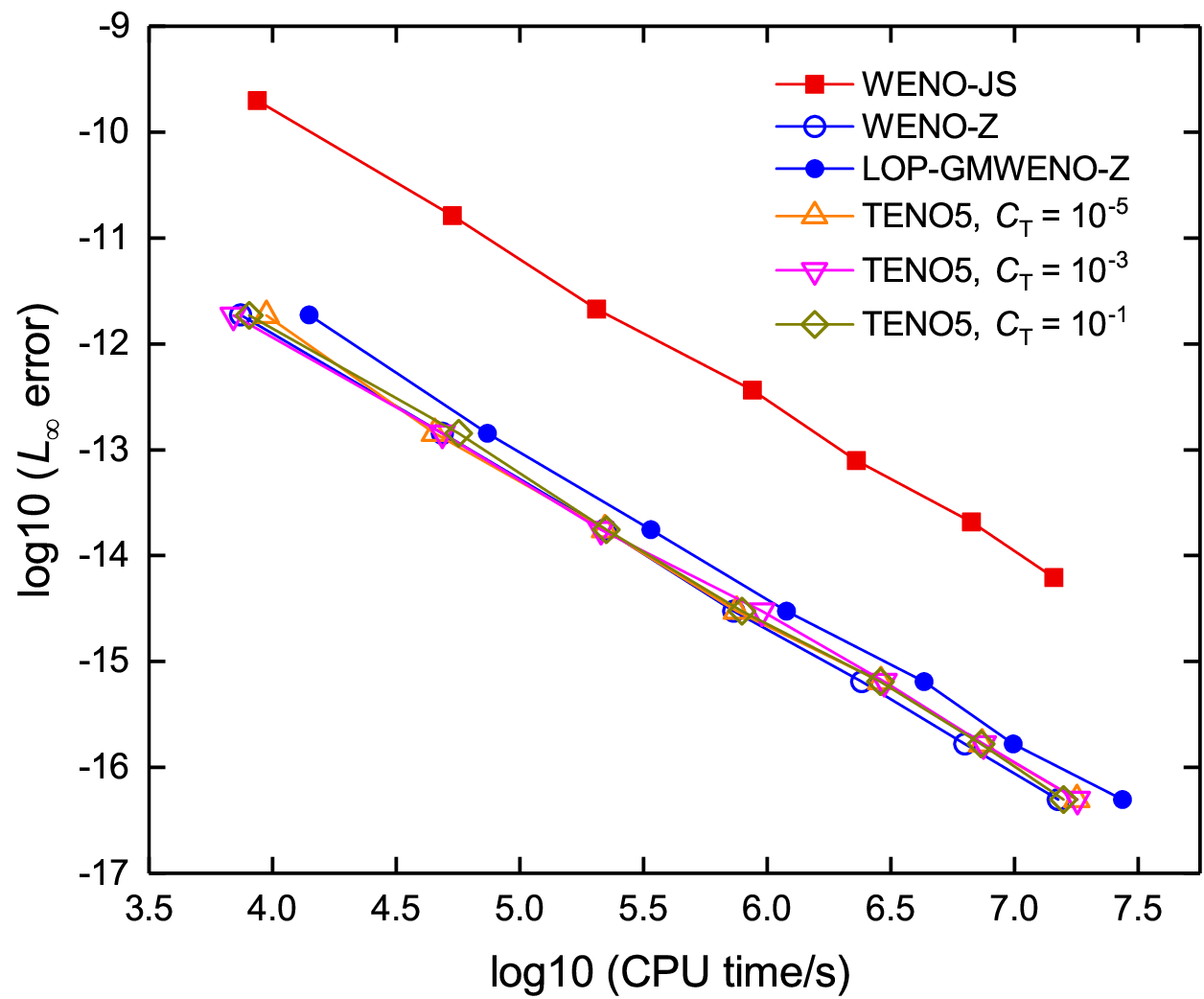} 
  \includegraphics[height=0.27\textwidth]
  {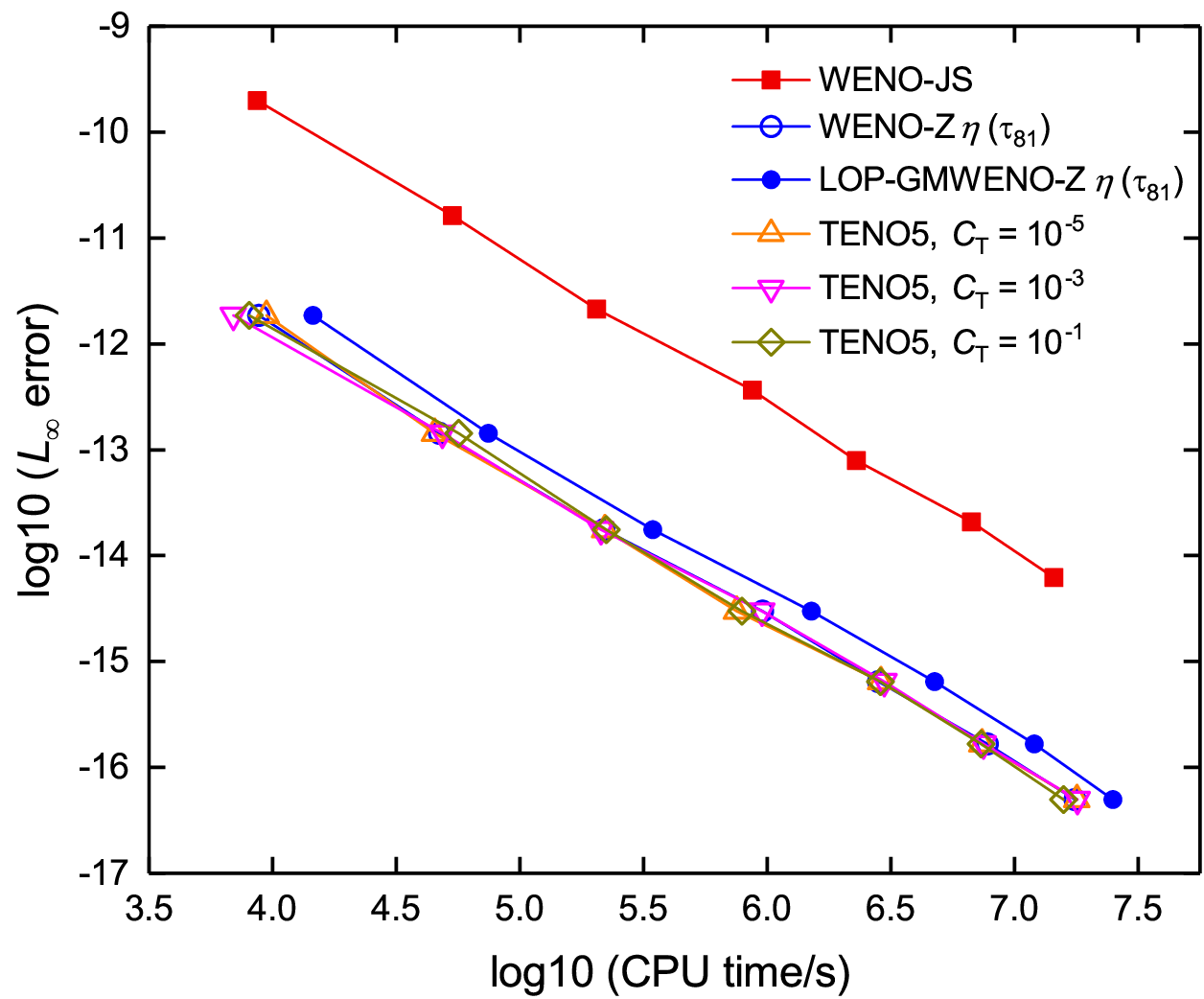}     
  \includegraphics[height=0.27\textwidth]
  {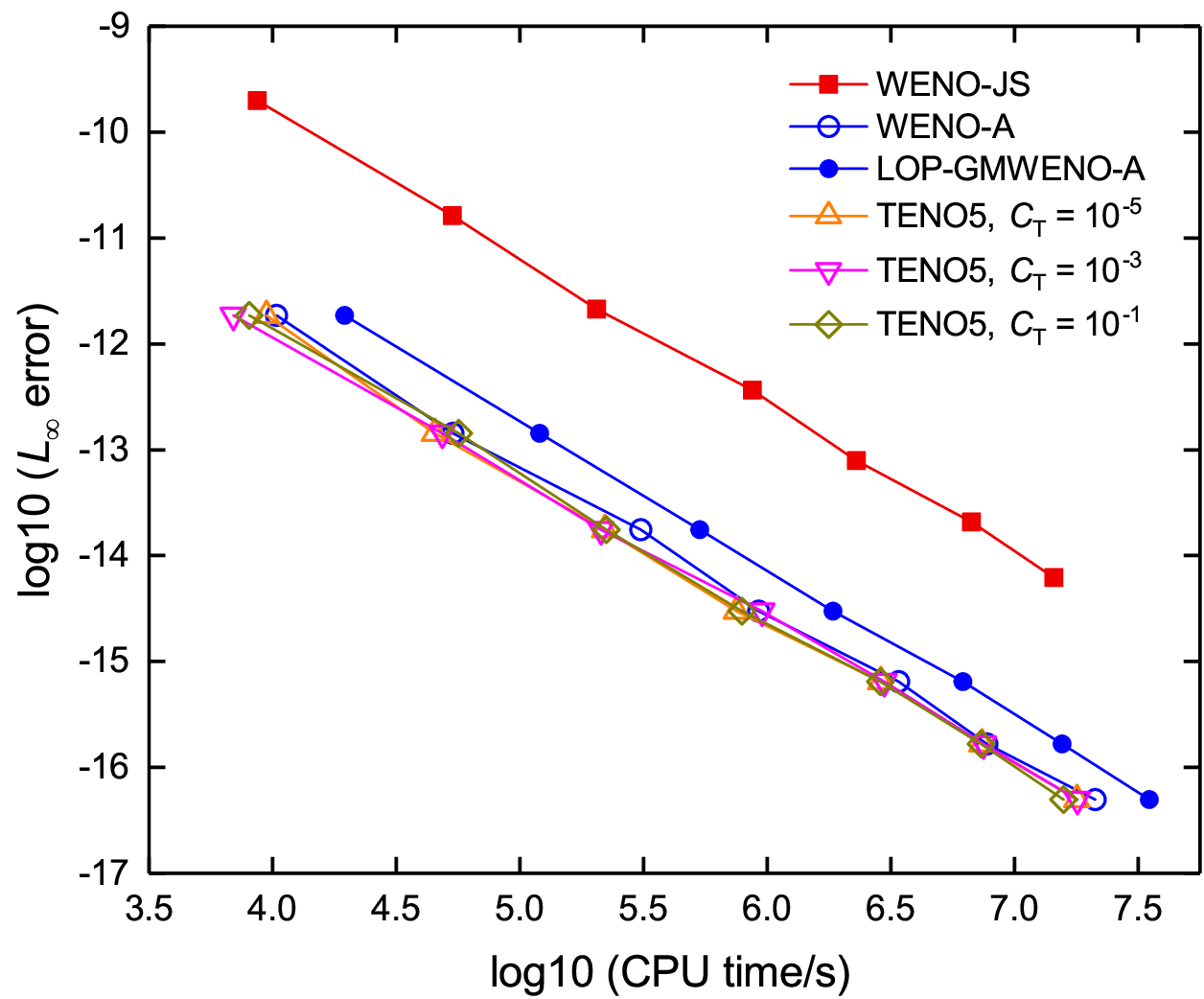} \\
  \includegraphics[height=0.27\textwidth]
  {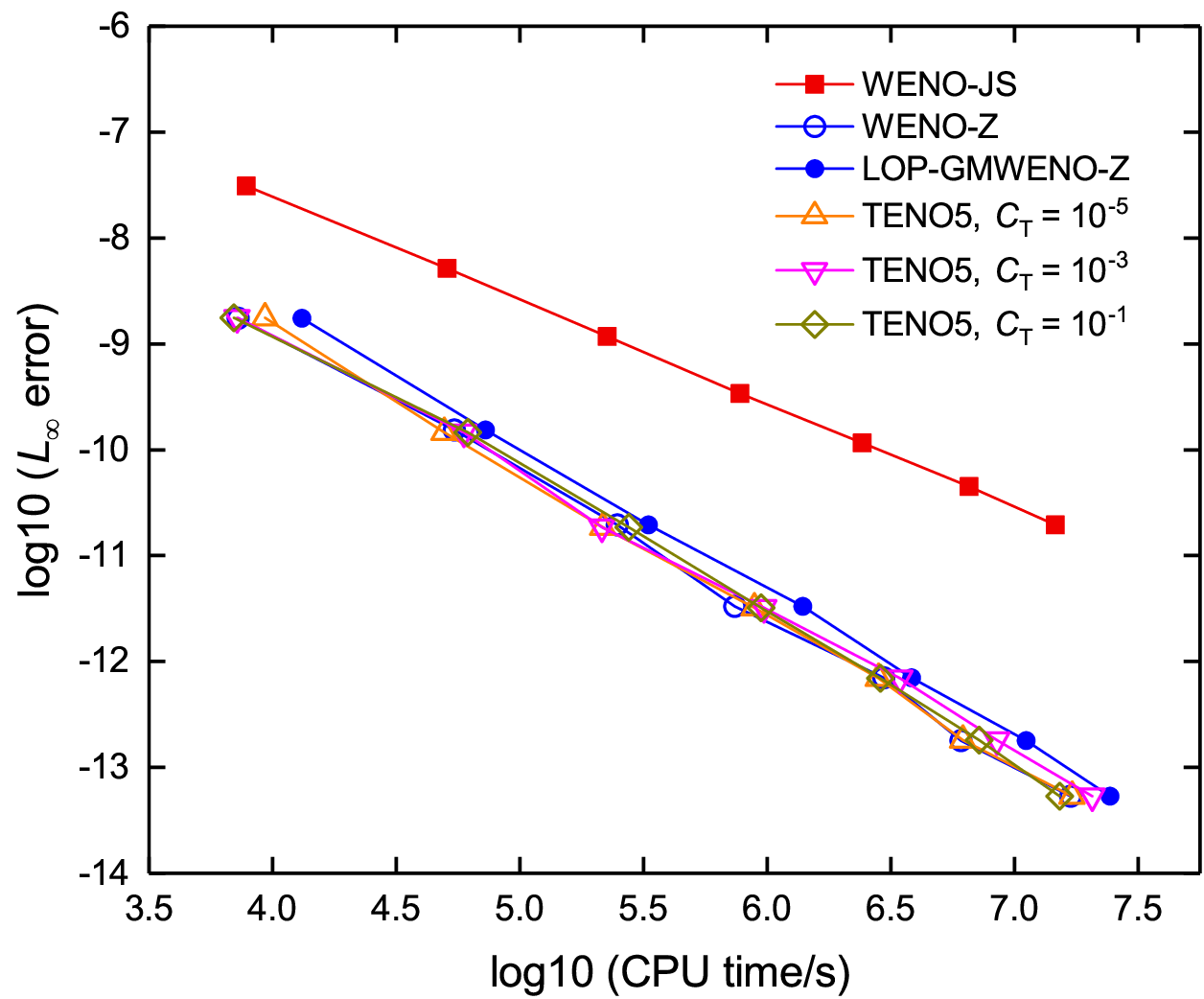} 
  \includegraphics[height=0.27\textwidth]
  {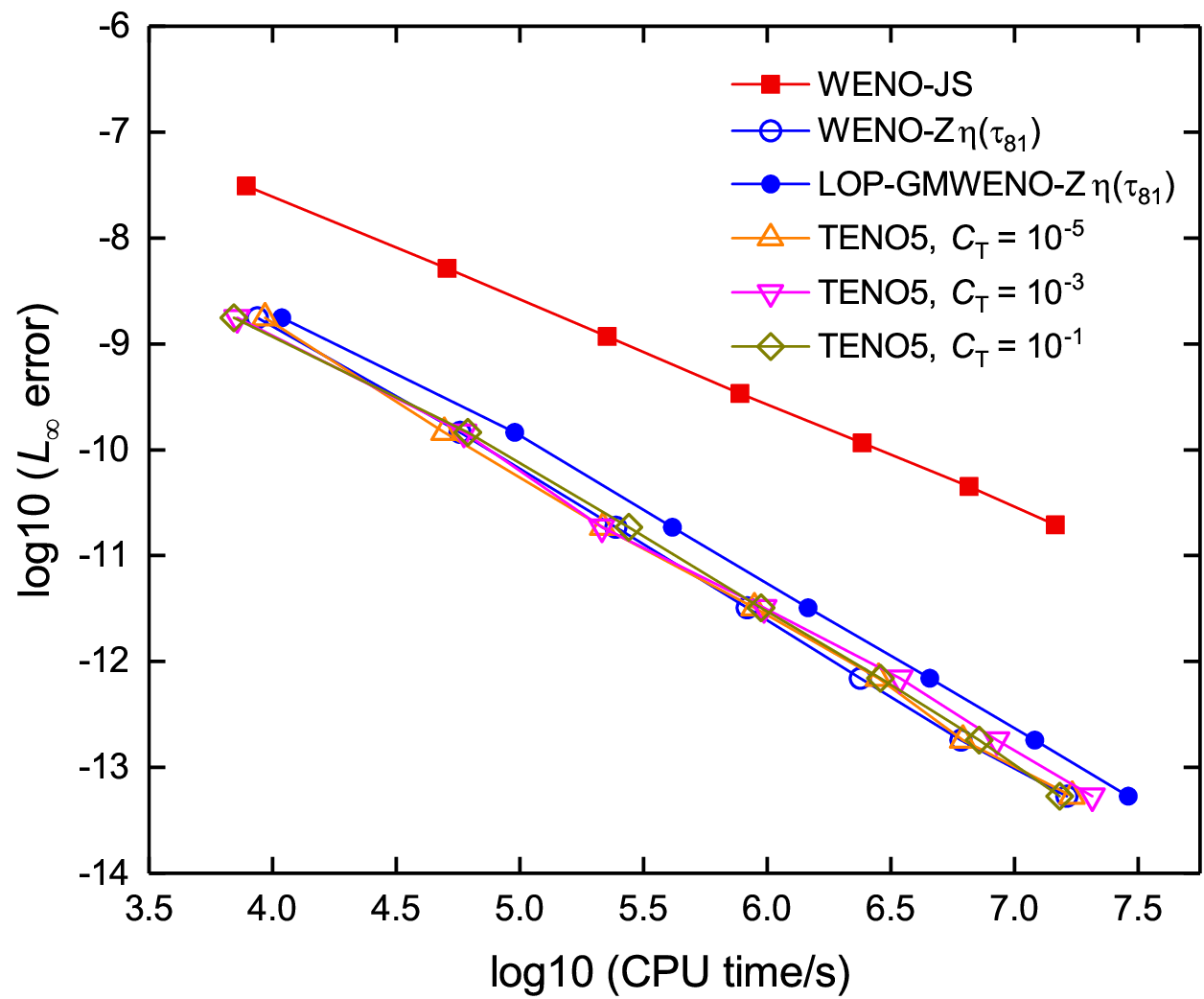}     
  \includegraphics[height=0.27\textwidth]
  {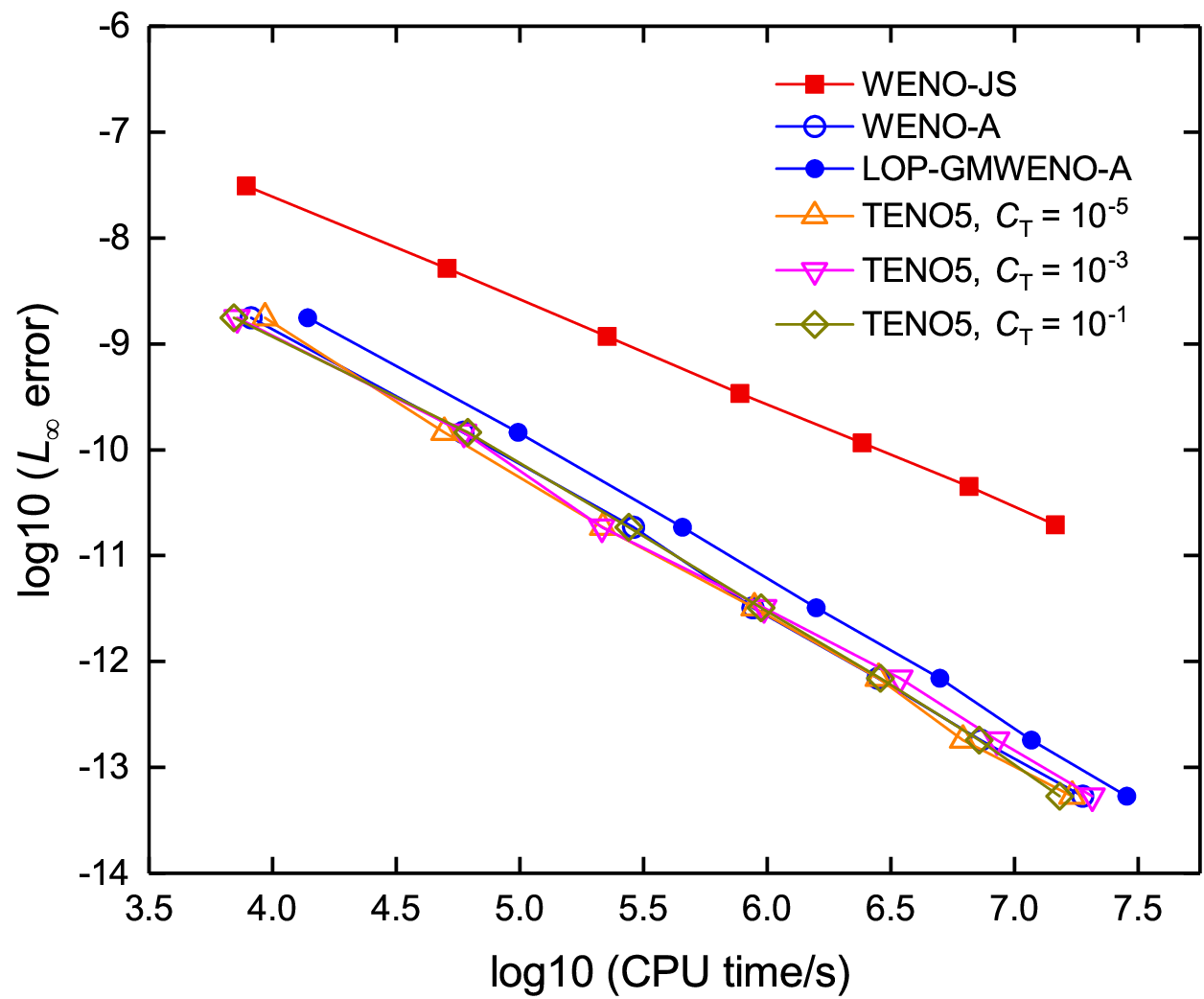}   
\caption{CPU-time versus $L_{\infty}$-norm errors of various WENO 
schemes for 2D smooth Euler problems. Top: Case 1. Bottom: Case 2.}
\label{fig:ex:AccuracyTest:efficiency}
\end{figure}

In order to examine the efficiency of Algorithm \ref{alg:posteriori}
, we plot the graphs for the CPU time versus the $L_{\infty}$-norm 
errors in Fig. \ref{fig:ex:AccuracyTest:efficiency}. Generally 
speaking, the LOP-GMWENO-X schemes have comparable or very slightly 
lower efficiency in comparison with their associated WENO-X schemes 
and the TENO5 schemes. Clearly, the efficiency of the LOP-GMWENO-X 
schemes are much higher than that of the WENO-JS scheme. 

\subsubsection{2D interacting blast waves}
We solve the 2D version of the standard blast-wave interaction problem first used by Woodward and Colella \cite{interactingBlastWaves-Woodward-Colella}. The computational domain $[0, 1]\times[0, 0.4]$ is initialized by
\begin{equation*}
\big( \rho, u, v, p \big)(x, y, 0) =\left\{
\begin{array}{ll}
(1, 0, 0, 1000),   & x \in [0, 0.1),  y \in [0, 0.4] \\
(1, 0, 0, 0.01),   & x \in [0.1, 0.9),  y \in [0, 0.4] \\
(1, 0, 0, 100),    & x \in [0.9, 1.0],  y \in [0, 0.4].
\end{array}\right.
\label{initial_2DEuler-blast}
\end{equation*}
The reflective boundary conditions are used at $x = 0, 1$ and the 
zero-gradient boundary conditions are used at $y = 0, 0.4$. The 
output time is 0.038, and the uniform mesh of $N_{x} \times N_{y} = 300 \times 120$ is used.

\begin{figure}[!ht]
\centering
  \includegraphics[height=0.26\textwidth]
  {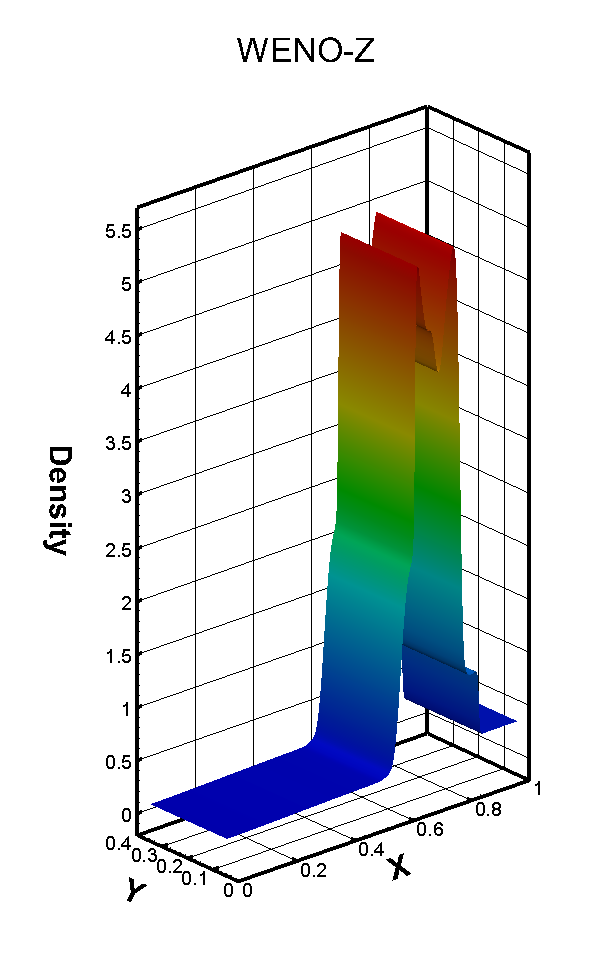}     
  \includegraphics[height=0.26\textwidth]
  {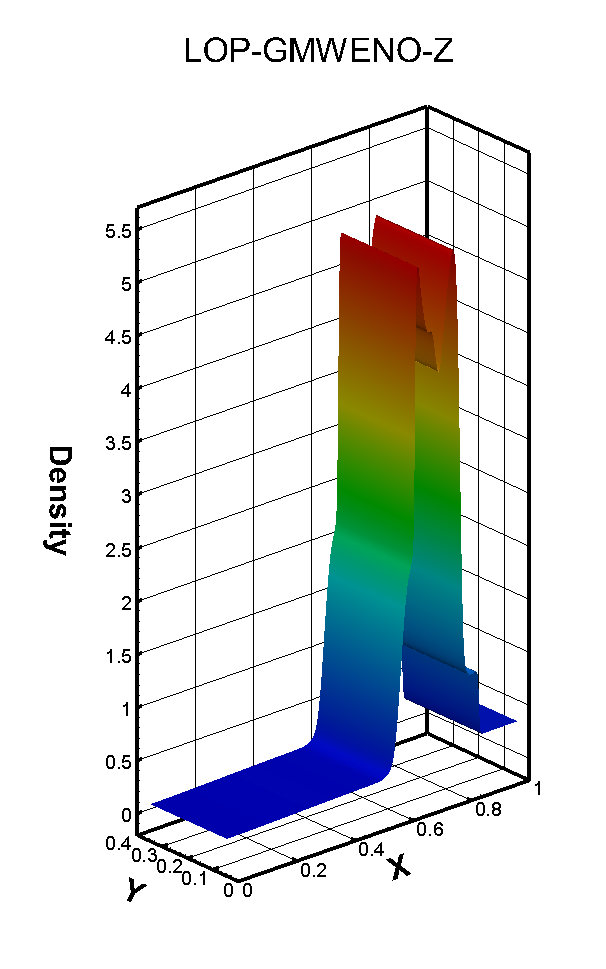} 
  \includegraphics[height=0.26\textwidth]
  {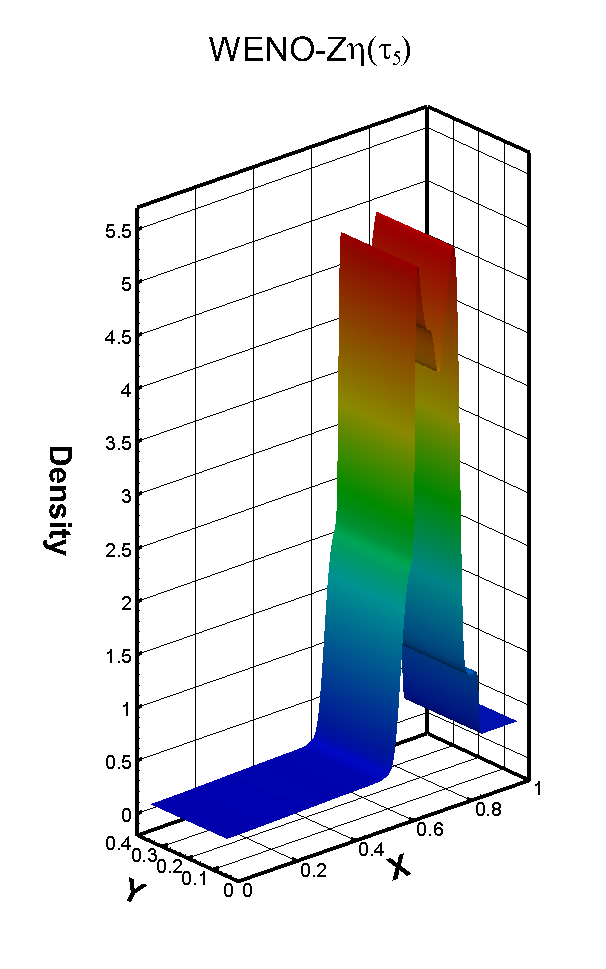}     
  \includegraphics[height=0.26\textwidth]
  {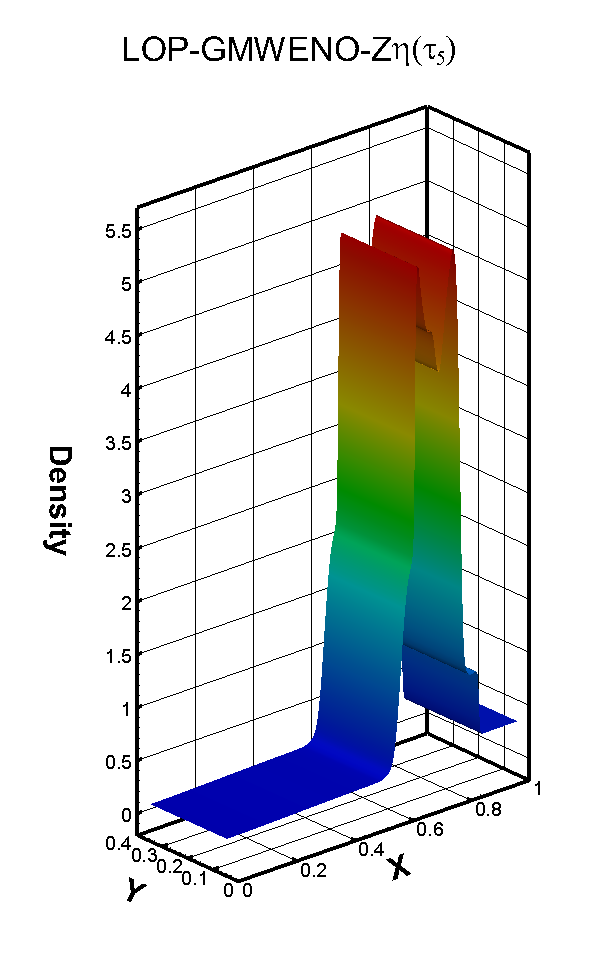} 
  \includegraphics[height=0.26\textwidth]
  {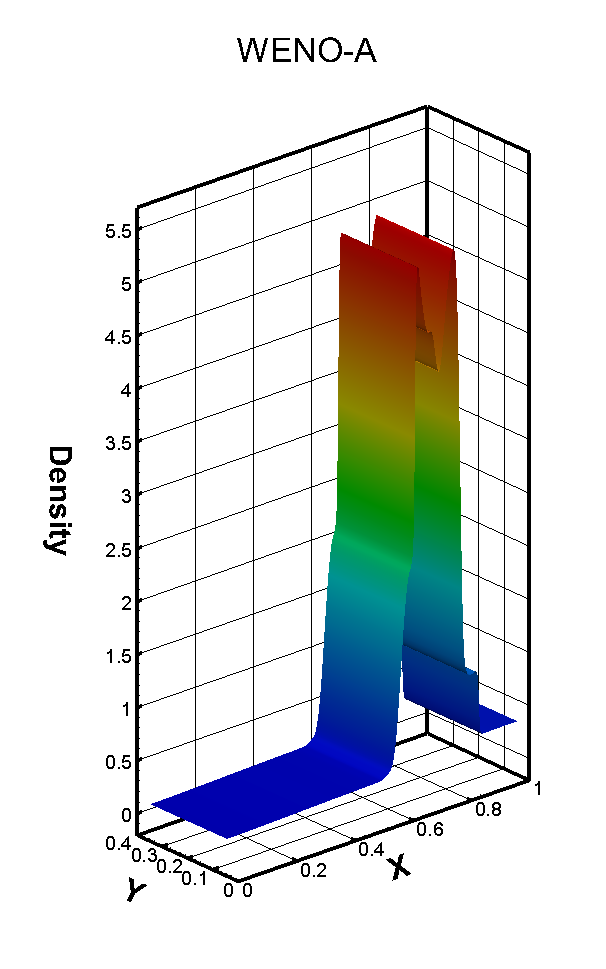}     
  \includegraphics[height=0.26\textwidth]
  {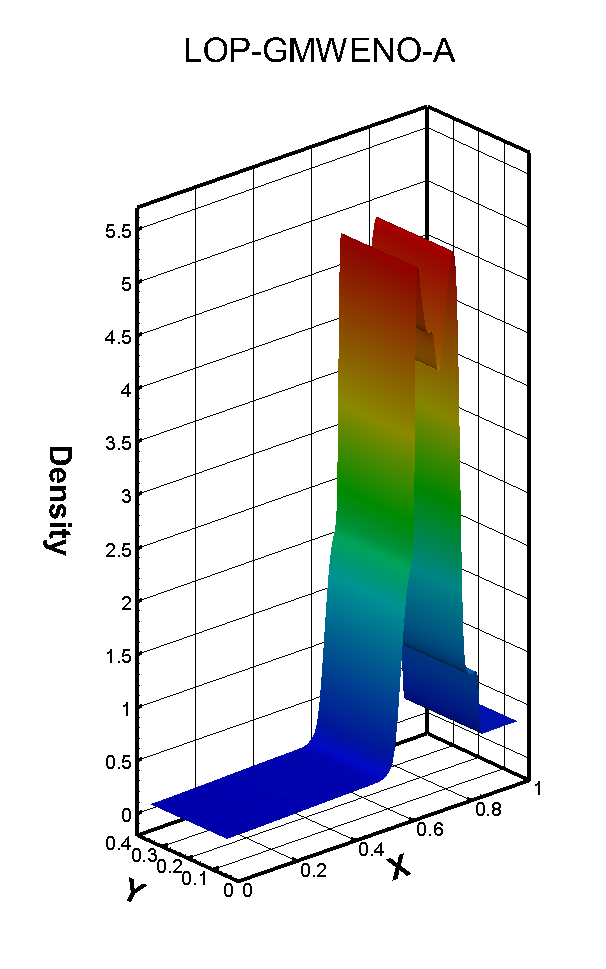} \\    
  \includegraphics[height=0.35\textwidth]
  {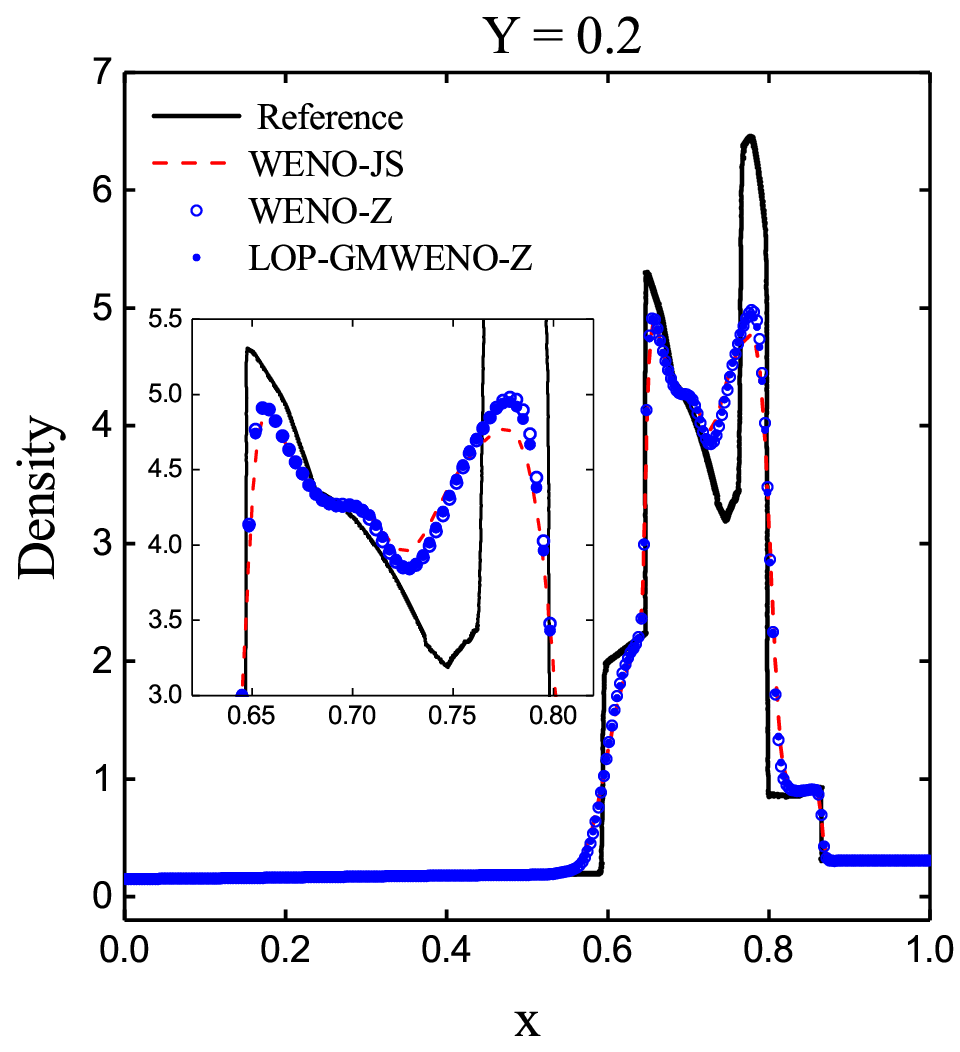} 
  \includegraphics[height=0.35\textwidth]
  {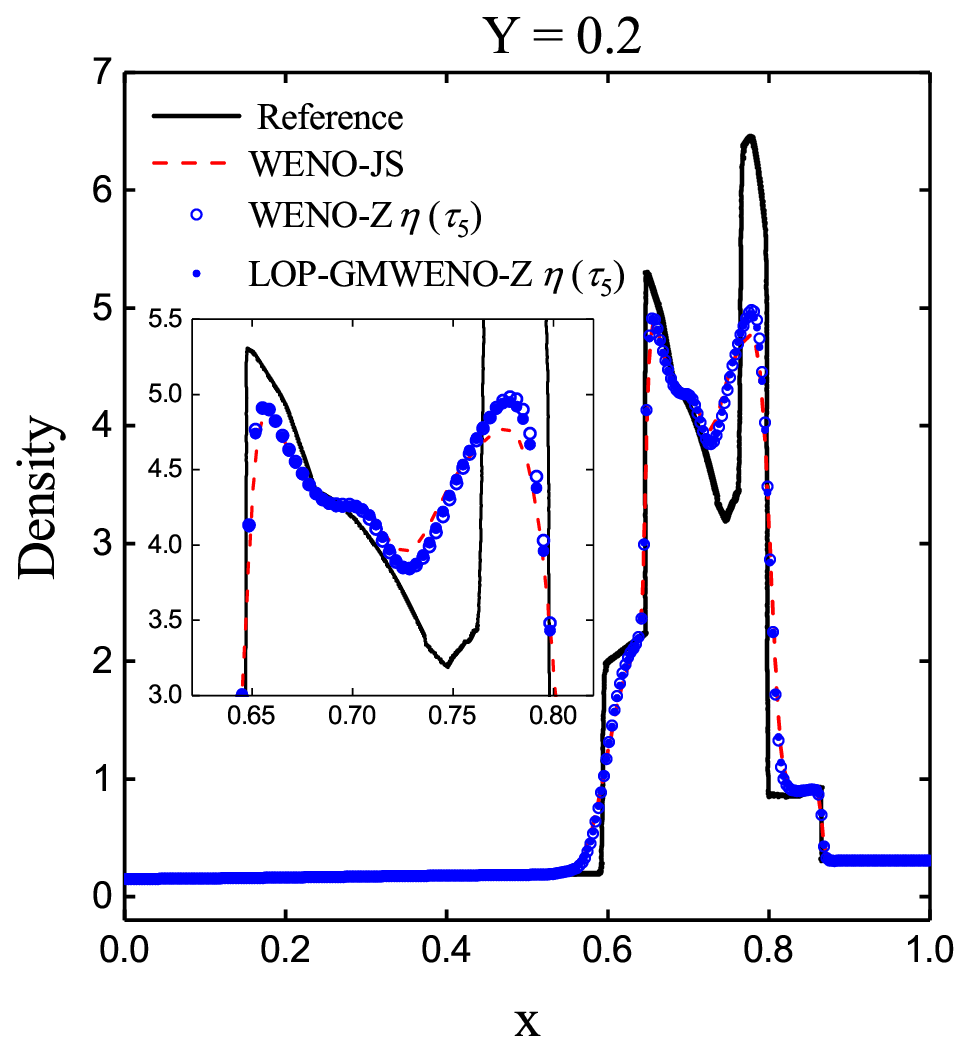}
  \includegraphics[height=0.35\textwidth]
  {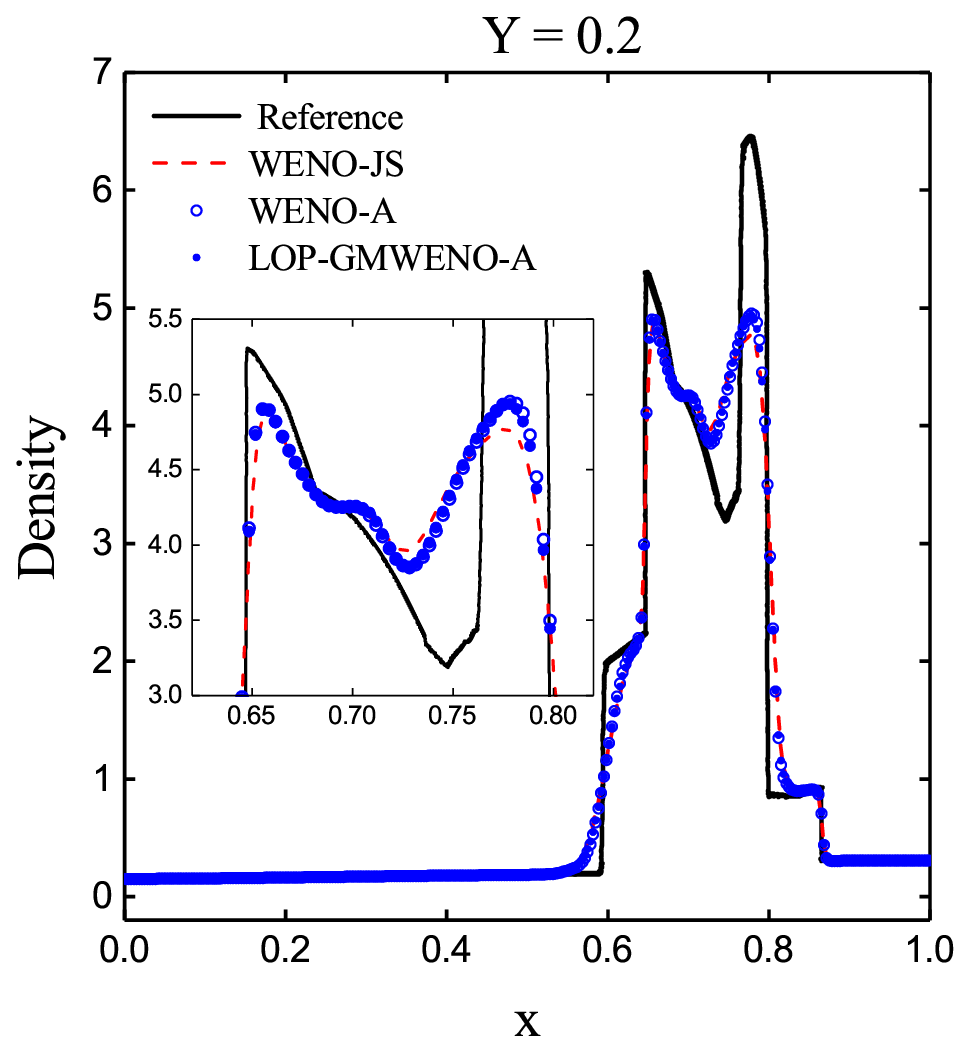}  \\    
\caption{Solutions of the interacting blast waves problem. Top: density contours. Bottom: density profiles at plane $y = 0.2$.}
\label{fig:ex:blast-wave}
\end{figure}

In Fig. \ref{fig:ex:blast-wave}, the density contours and the 
density profiles at plane $y = 0.2$ of the LOP-GMWENO-Z, LOP-GMWENO-Z$\eta(\tau_{5})$, LOP-GMWENO-A schemes and their associated WENO-Z, 
WENO-Z$\eta(\tau_{5})$, WENO-A schemes are presented. For 
comparison purpose, we also plot the density profiles of the WENO-JS 
scheme and the reference result is obtained by solving the 
corresponding 1D case using the WENO-JS scheme with $N = 10000$. It 
can be seen that the LOP-GMWENO-X schemes perform as well as their 
associated WENO-X schemes, and all of these schemes generate sharper 
results than the WENO-JS scheme. It should be noted that we present 
the result of the WENO-Z$\eta(\tau_{5})$ scheme but not that of the 
WENO-Z$\eta(\tau_{81})$ scheme in this test as it was indicated 
\cite{WENO-eta} that simulations under the use of WENO-Z$\eta(\tau_{81})$ scheme is found to be unstable.

\subsubsection{Implosion and explosion problems}
We simulate the implosion problem \cite{Implosion-SCIENCE-CHINA2010} 
and explosion problem \cite{WENO-ACM,Toro_RiemannSolvers}. In order 
to examine the performance for a long-time run, we reset the 
computational domain of the explosion problem used in \cite{WENO-ACM,Toro_RiemannSolvers} to be $[-10,10]\times[-10, 10]$ and the output 
time to be $t = 5.0$. For the implosion problem, we use the same 
computational domain of $[-20,20]\times[-20, 20]$ as in \cite{Implosion-SCIENCE-CHINA2010} where a large output time $t = 4.2$ was 
used. The initial conditions are given as
\begin{equation*}
\begin{aligned}
\text{Implosion: } \qquad &\big( \rho, u, v, p \big)(x, y, 0) = \left\{\begin{array}{ll}
%\begin{aligned}
(1, 0, 0, 1), & \mathrm{if} \sqrt{x^{2} + y^{2}} < 0.4,\\
(0.125, 0, 0, 0.1), & \mathrm{else},
%\end{aligned}
\end{array}\right. \\
\text{Explosion: } \qquad &\big( \rho, u, v, p \big)(x, y, 0) = \left\{\begin{array}{ll}
%\begin{aligned}
(2, 0, 0, 5), & \mathrm{if} \sqrt{x^{2} + y^{2}} < 4.0,\\
(1, 0, 0, 1), & \mathrm{else}.
%\end{aligned}
\end{array}\right.
\end{aligned}
\label{ex:Implosion-Explosion}
\end{equation*}
On all edges, the reflective boundary condition is used for the 
implosion problem, and the transmissive boundary condition is used 
for the explosion problem. $N_{x} \times N_{y} = 400 \times 400$ is 
used for these two problems.

\begin{figure}[!ht]
\centering
  \includegraphics[height=0.26\textwidth]
  {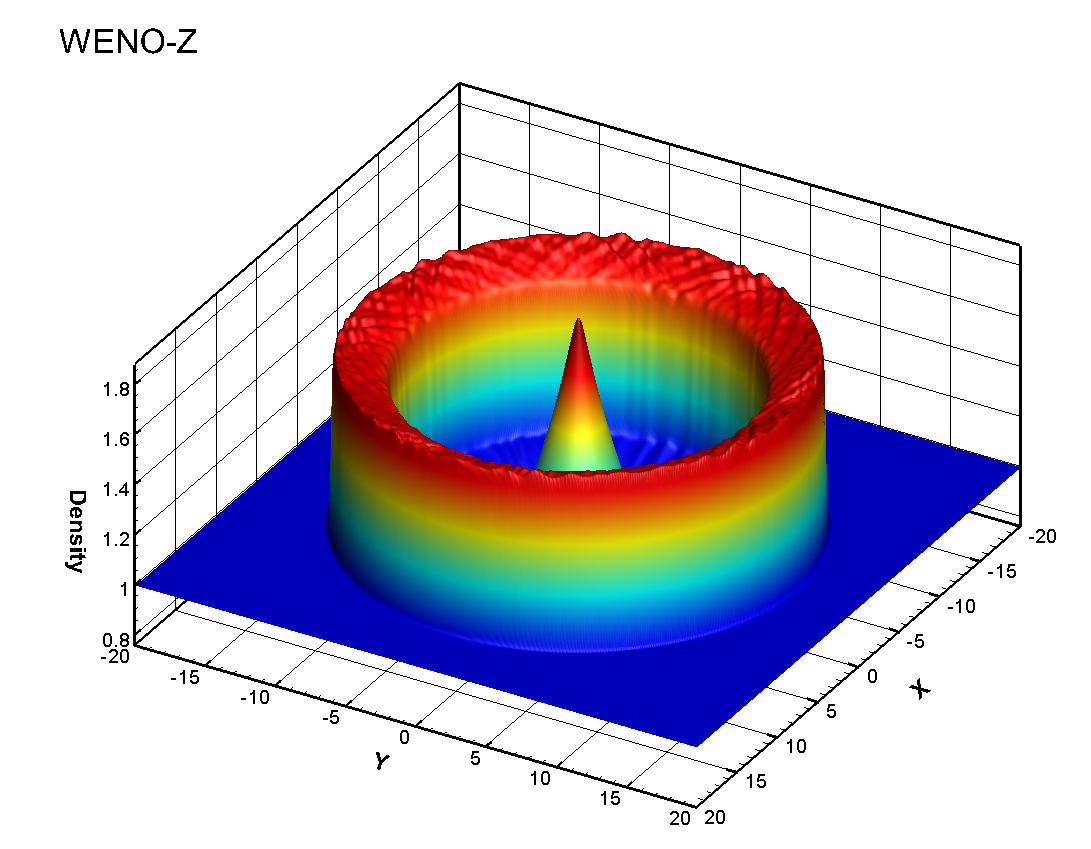}     
  \includegraphics[height=0.26\textwidth]
  {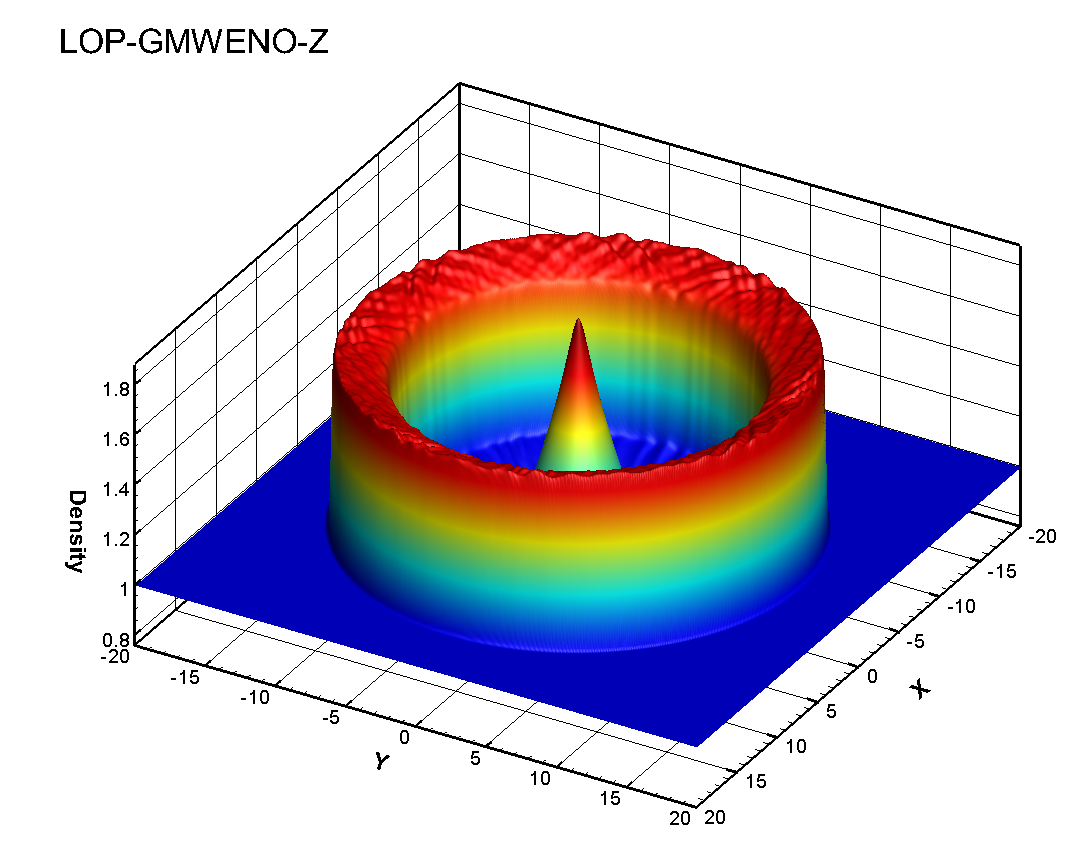} 
  \includegraphics[height=0.26\textwidth]
  {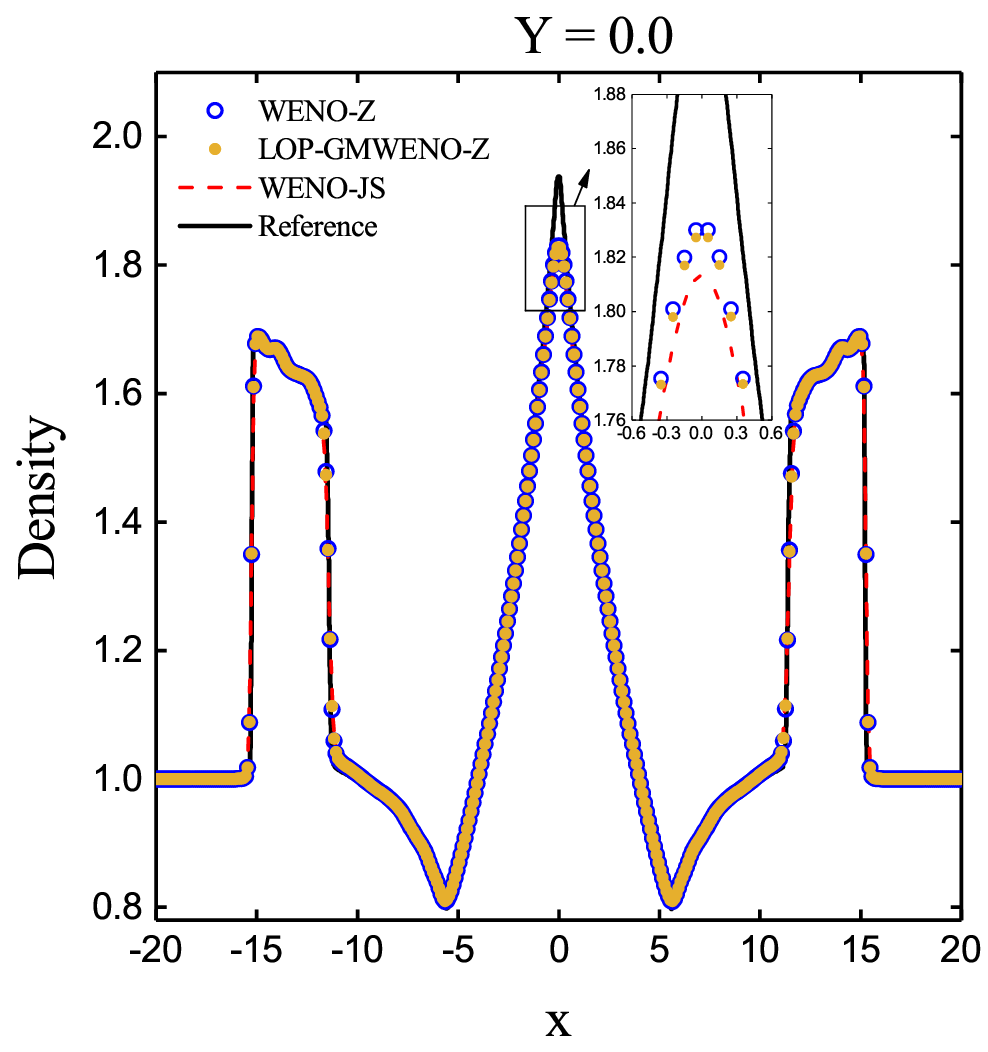}  \\   
  \includegraphics[height=0.26\textwidth]
  {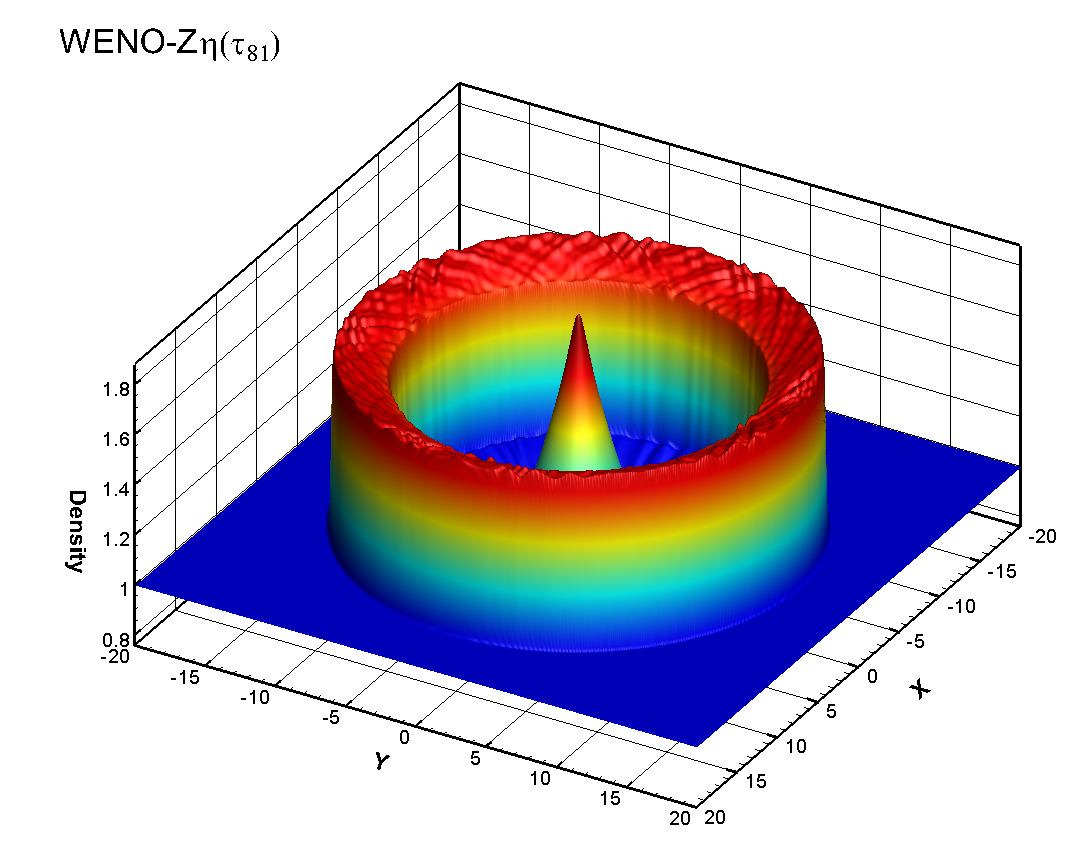}     
  \includegraphics[height=0.26\textwidth]
  {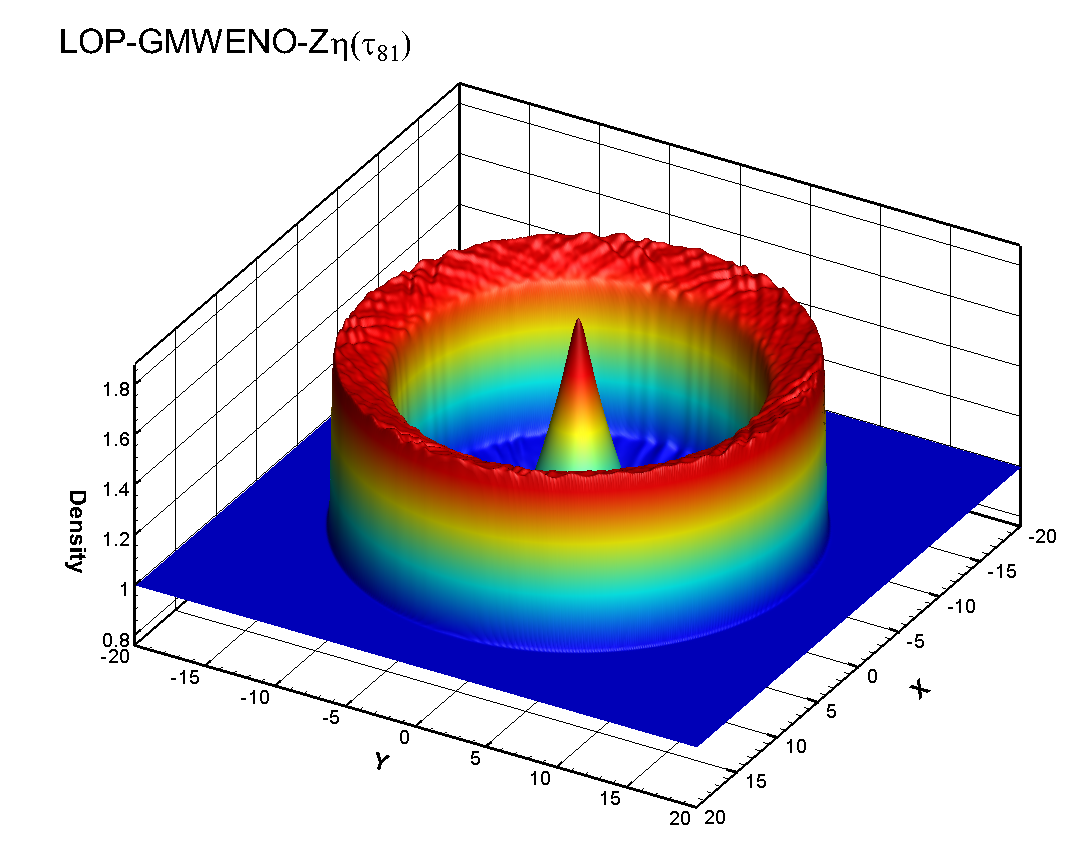} 
  \includegraphics[height=0.26\textwidth]
  {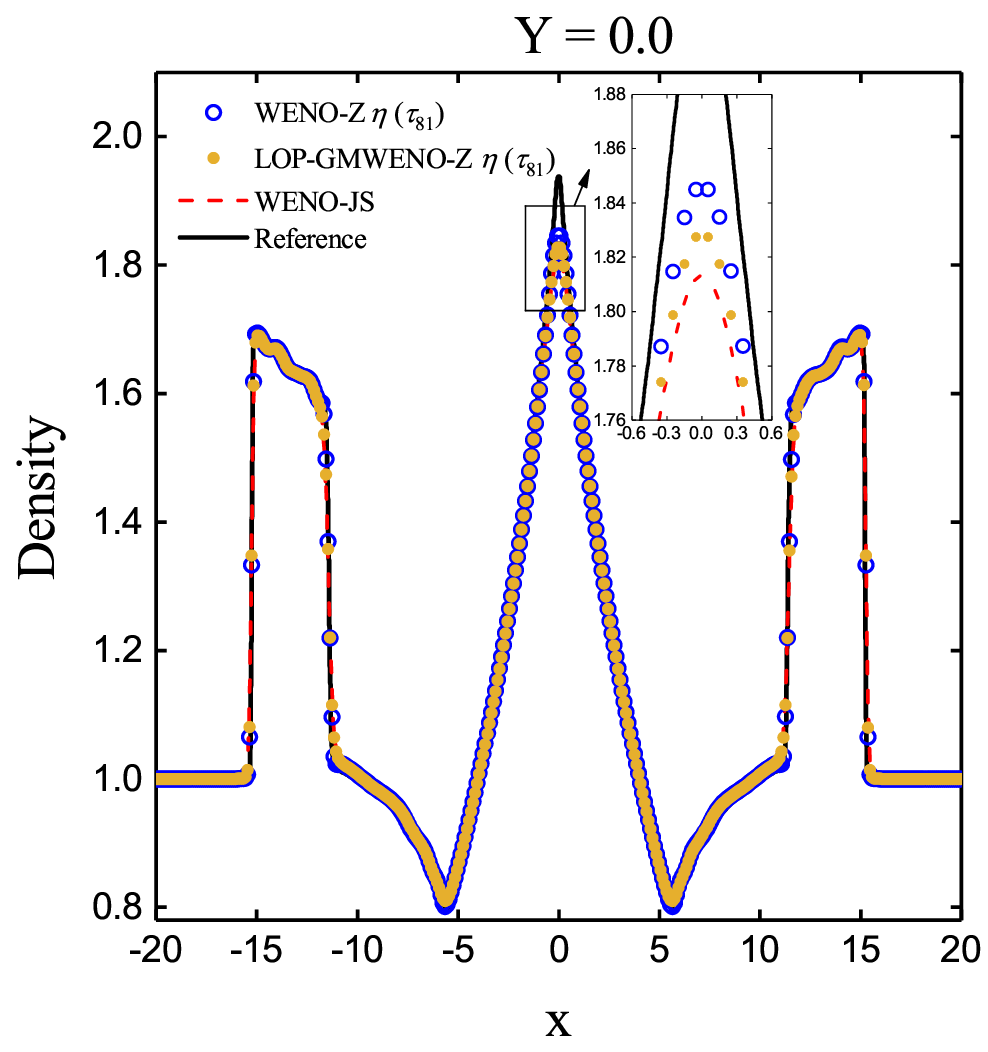}  \\   
  \includegraphics[height=0.26\textwidth]
  {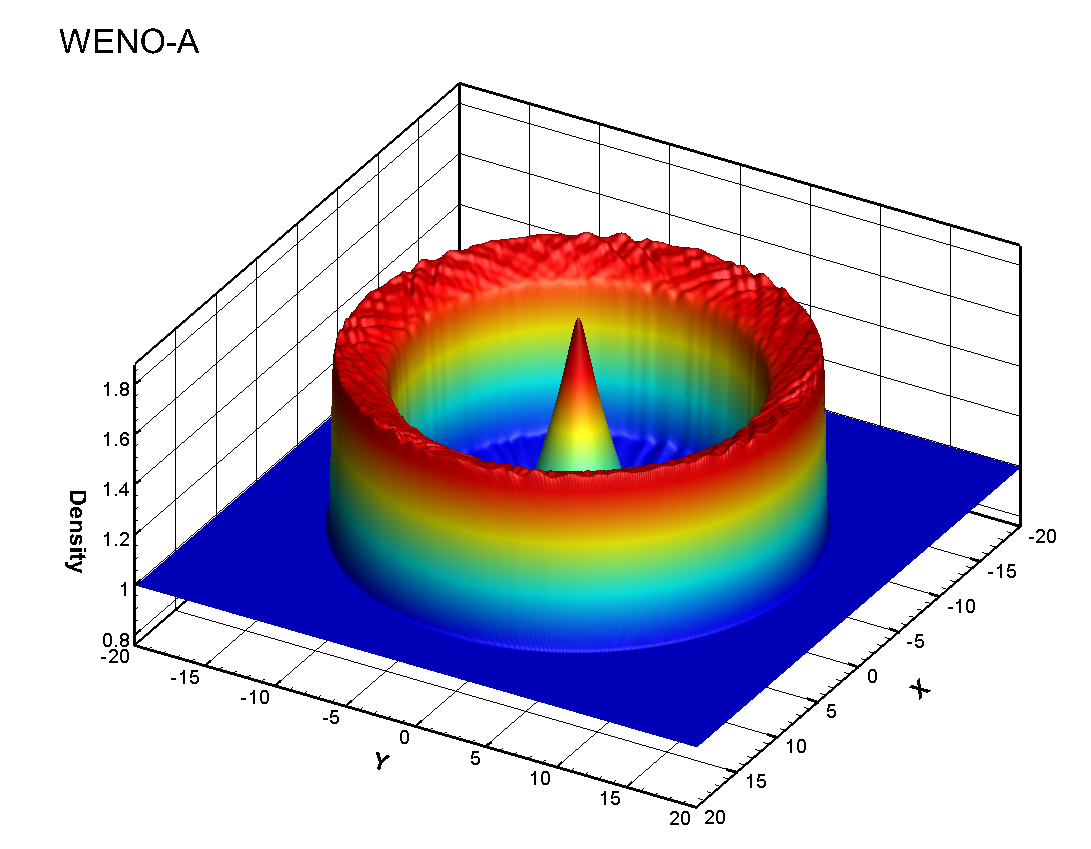}     
  \includegraphics[height=0.26\textwidth]
  {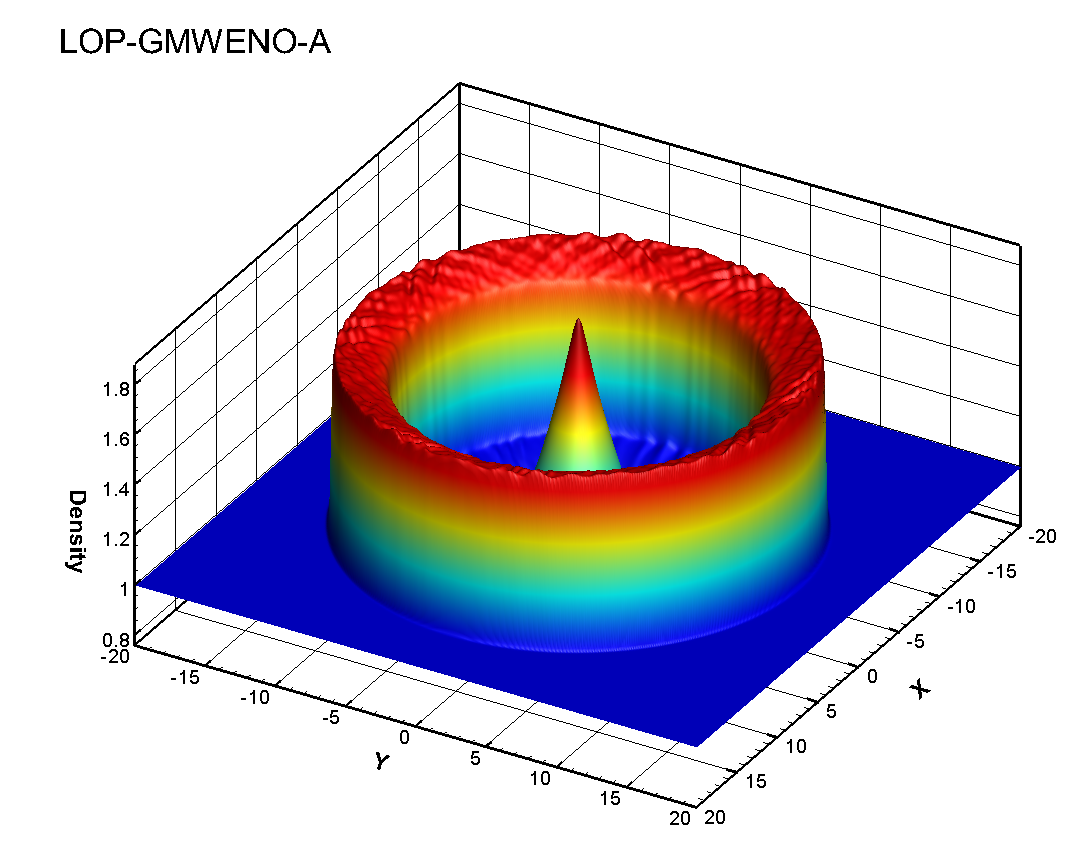} 
  \includegraphics[height=0.26\textwidth]
  {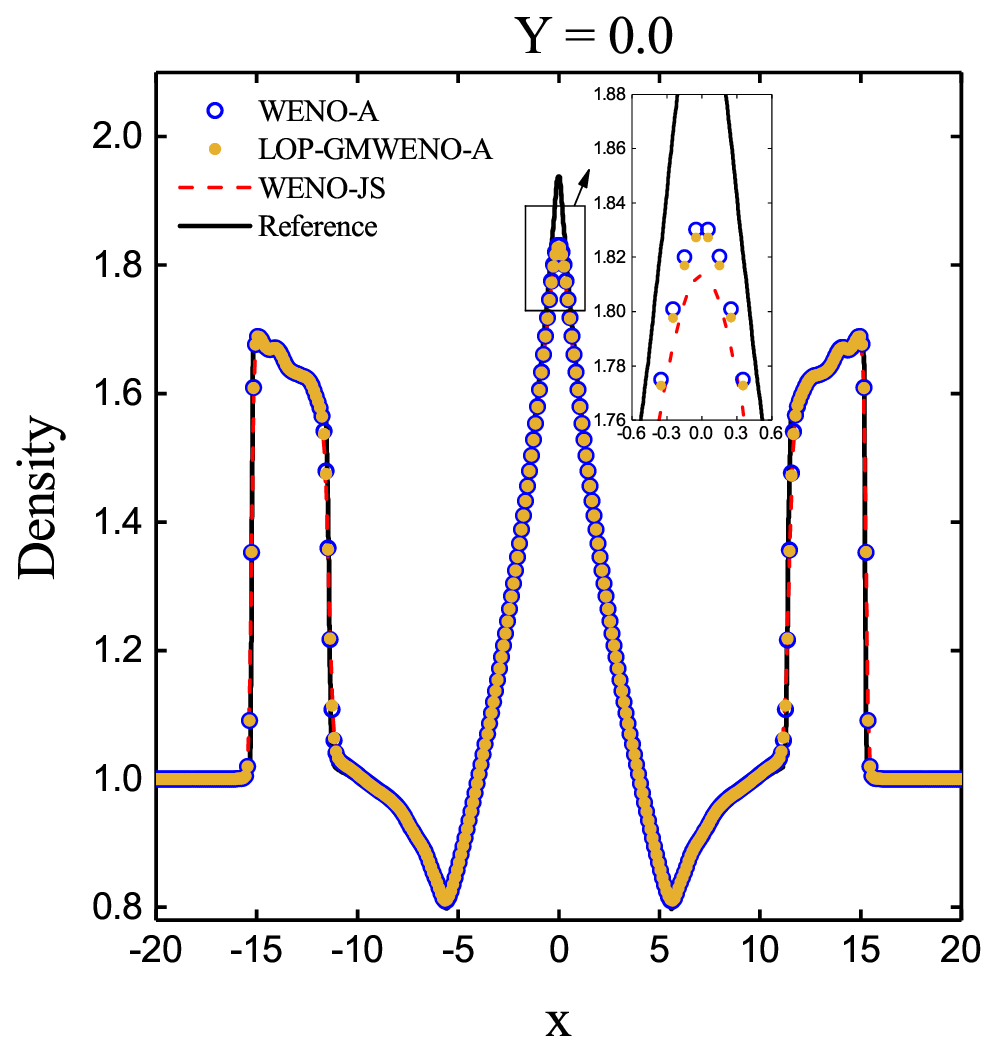}  \\         
\caption{Solutions of the 2D implosion problem. Left and middle: density contours. Right: density profiles at plane $y = 0.0$.}
\label{fig:ex:Implosion}
\end{figure}

\begin{figure}[!ht]
\centering
  \includegraphics[height=0.26\textwidth]
  {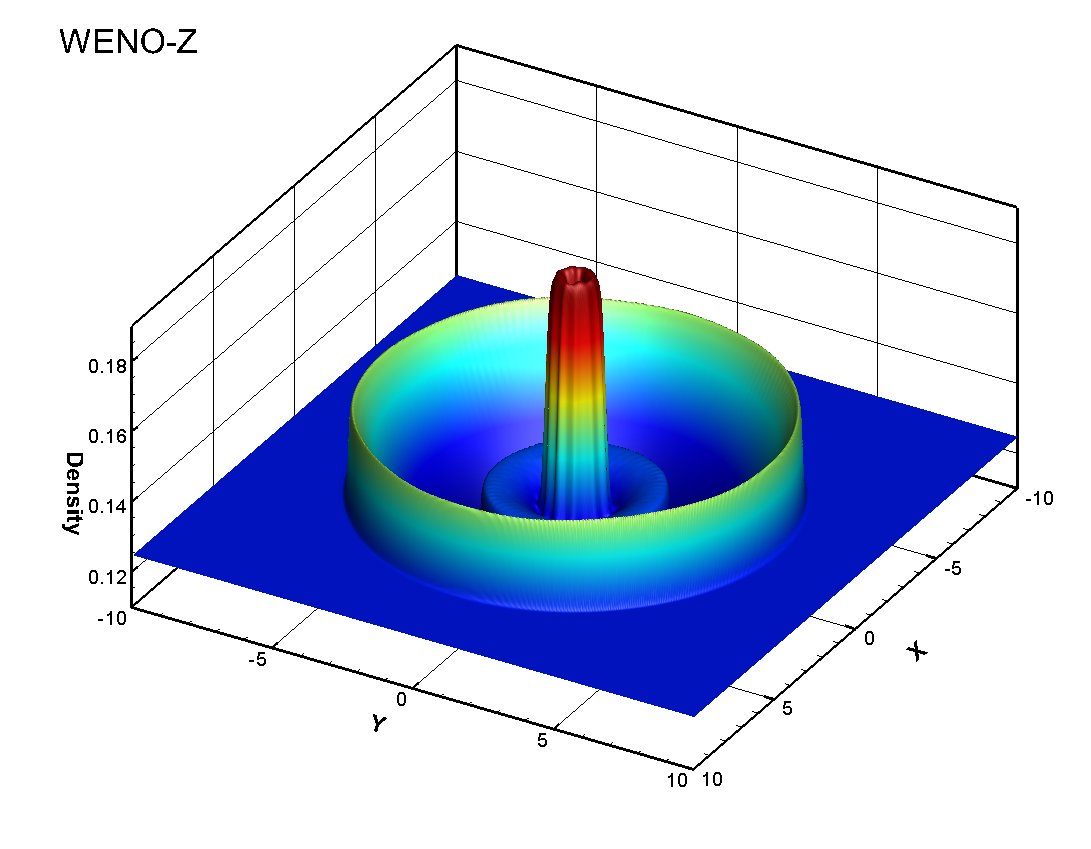}     
  \includegraphics[height=0.26\textwidth]
  {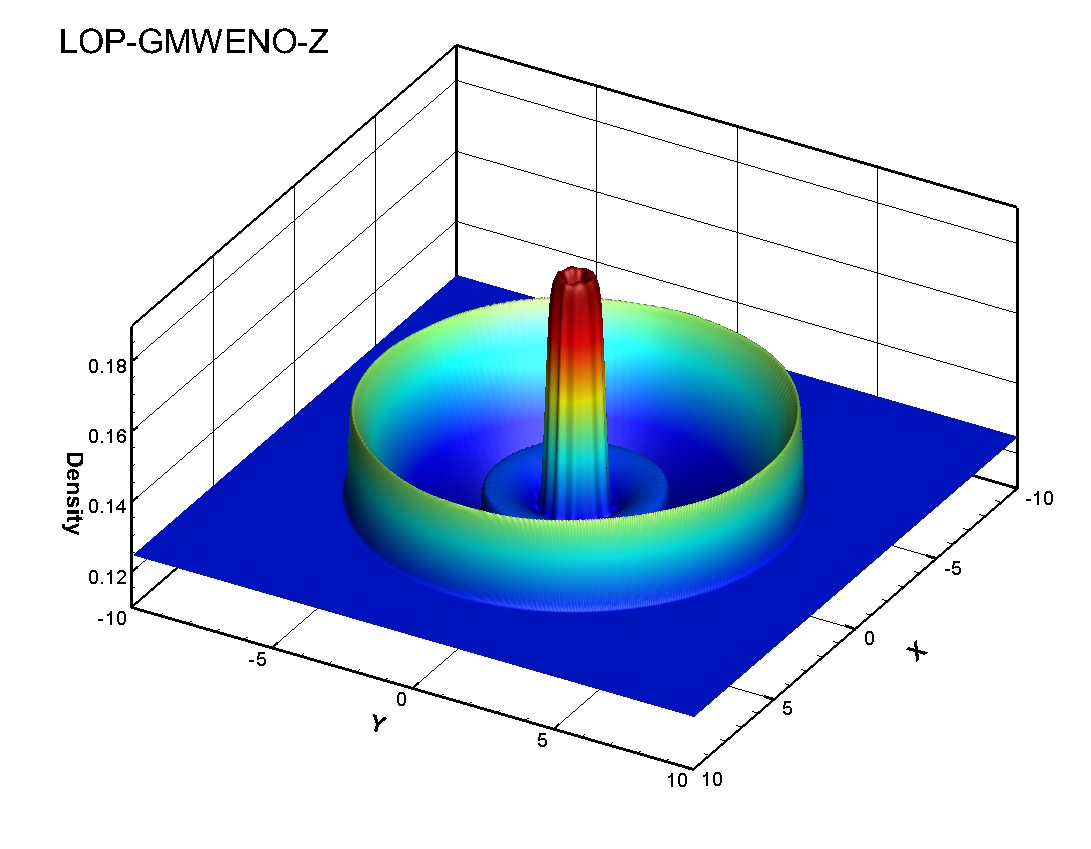} 
  \includegraphics[height=0.26\textwidth]
  {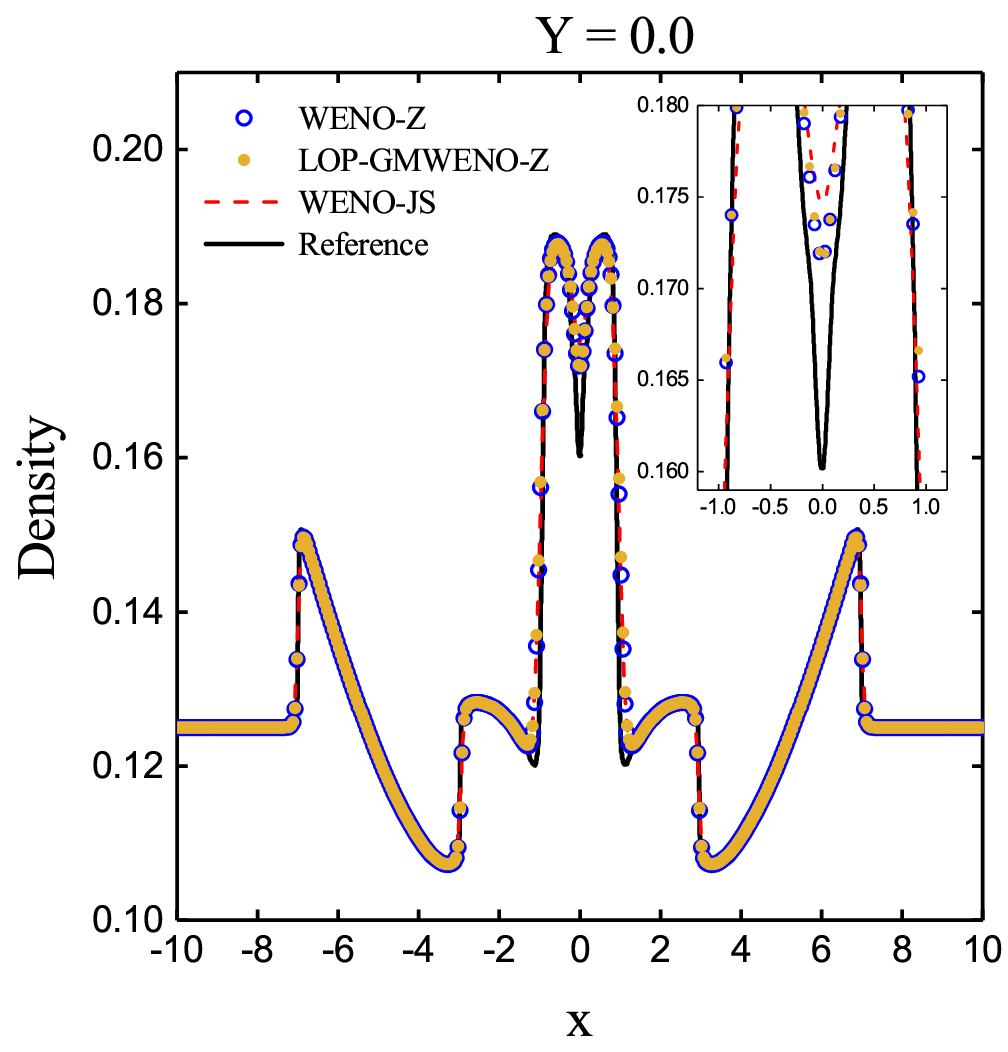}  \\   
  \includegraphics[height=0.26\textwidth]
  {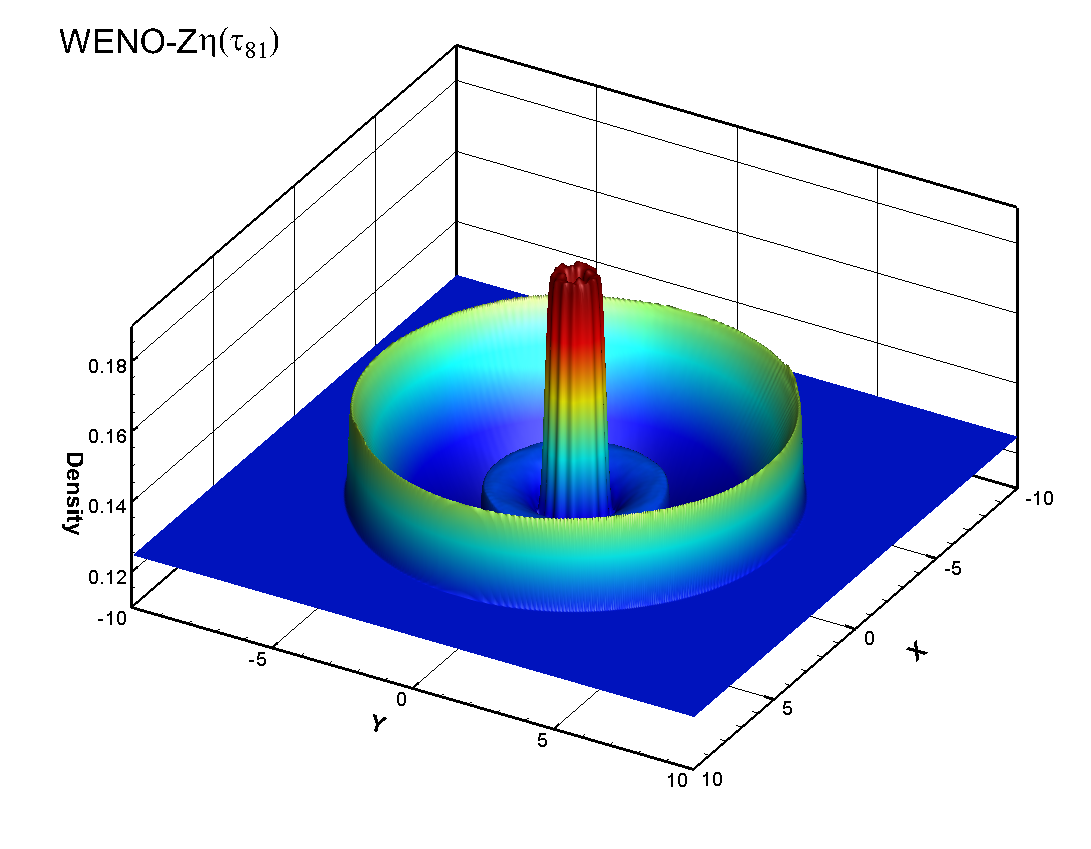}     
  \includegraphics[height=0.26\textwidth]
  {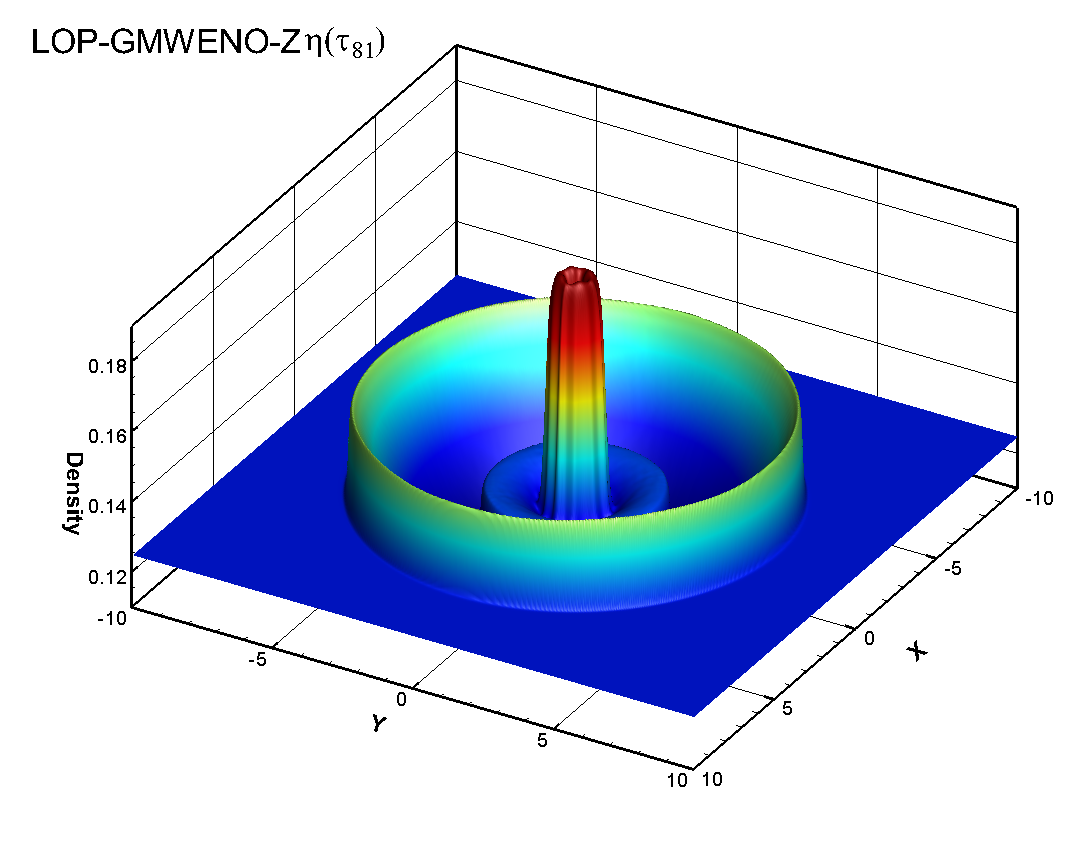} 
  \includegraphics[height=0.26\textwidth]
  {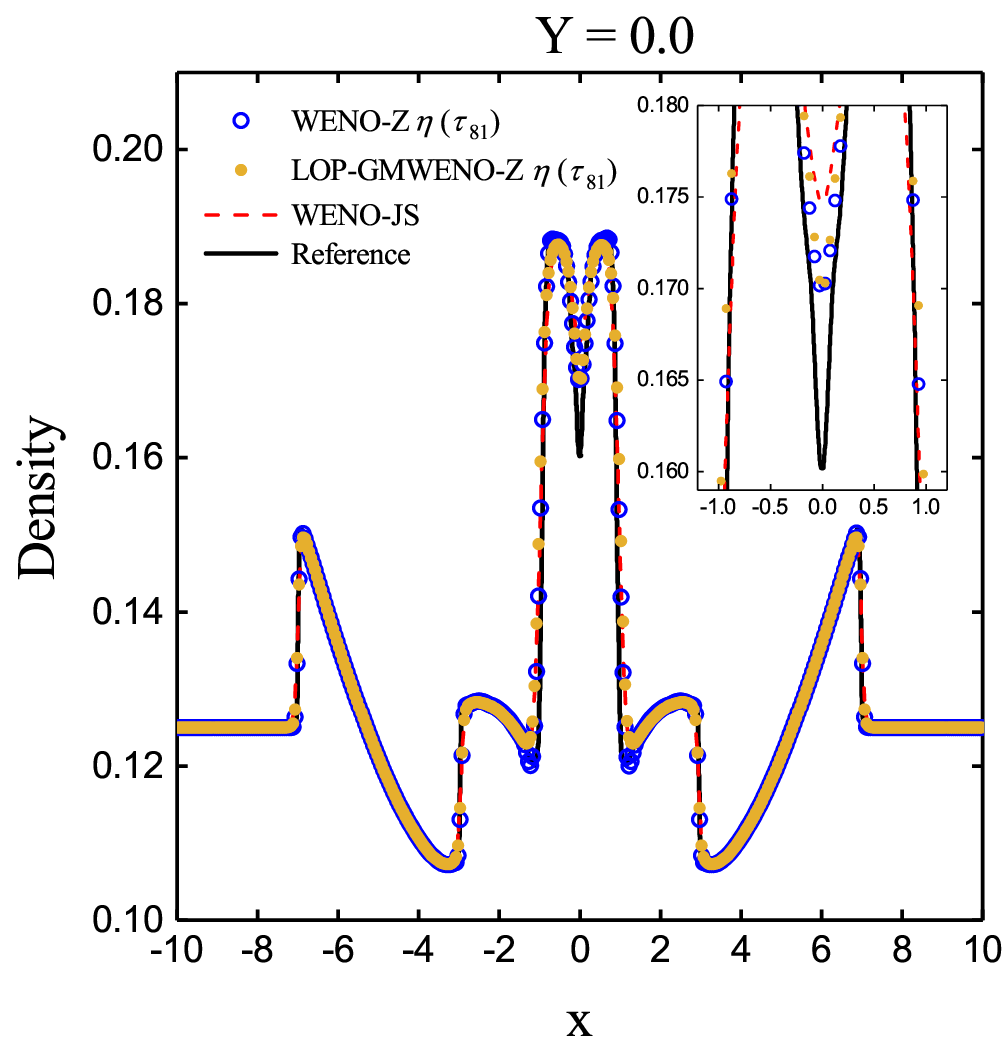}  \\   
  \includegraphics[height=0.26\textwidth]
  {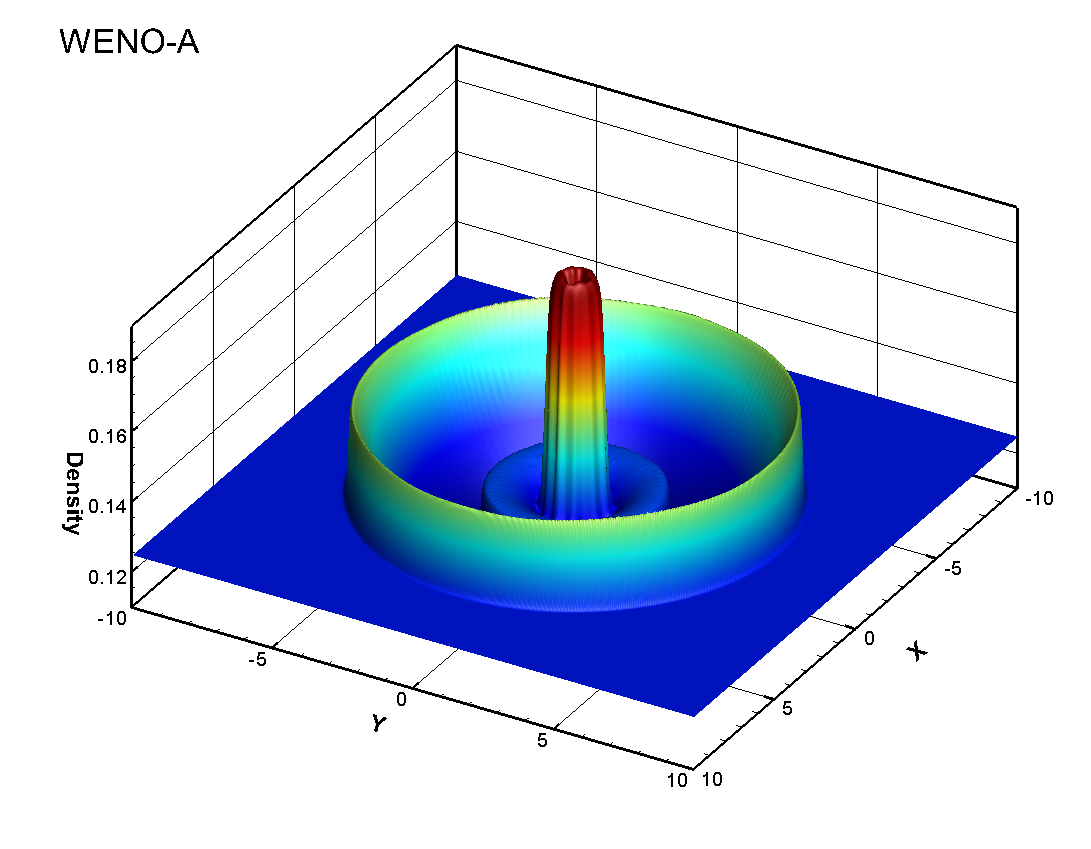}     
  \includegraphics[height=0.26\textwidth]
  {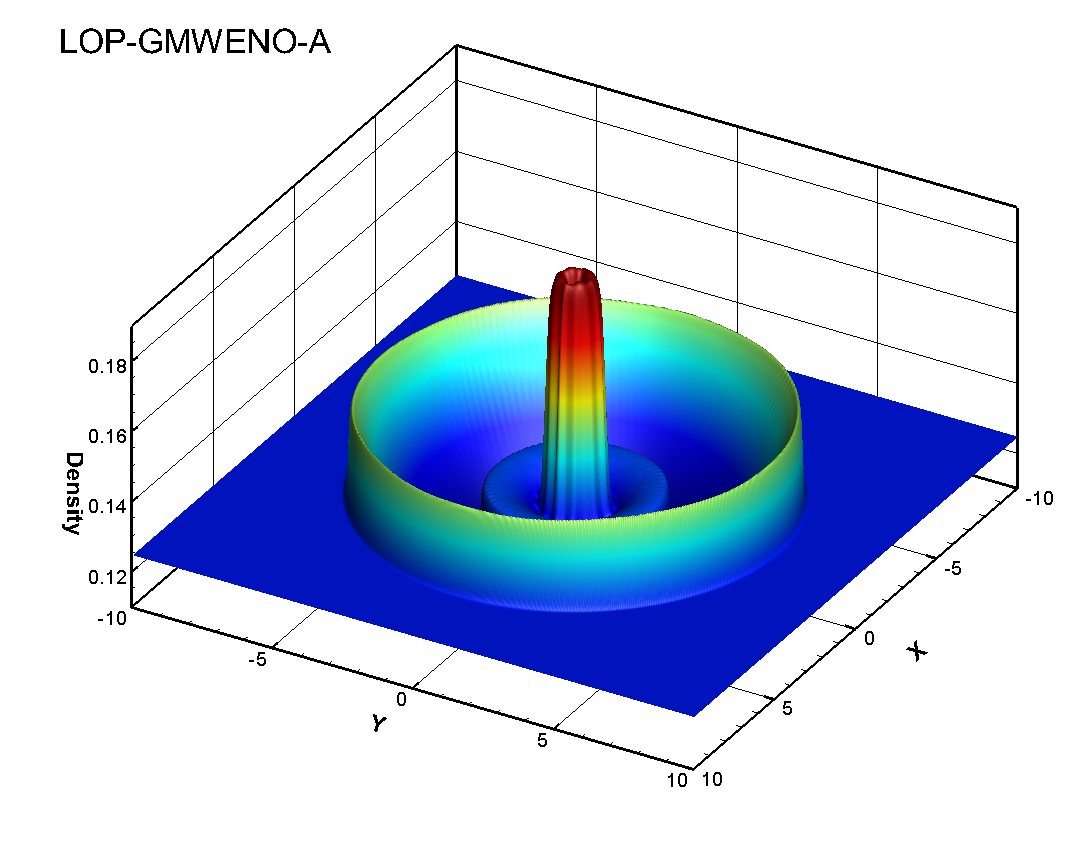} 
  \includegraphics[height=0.26\textwidth]
  {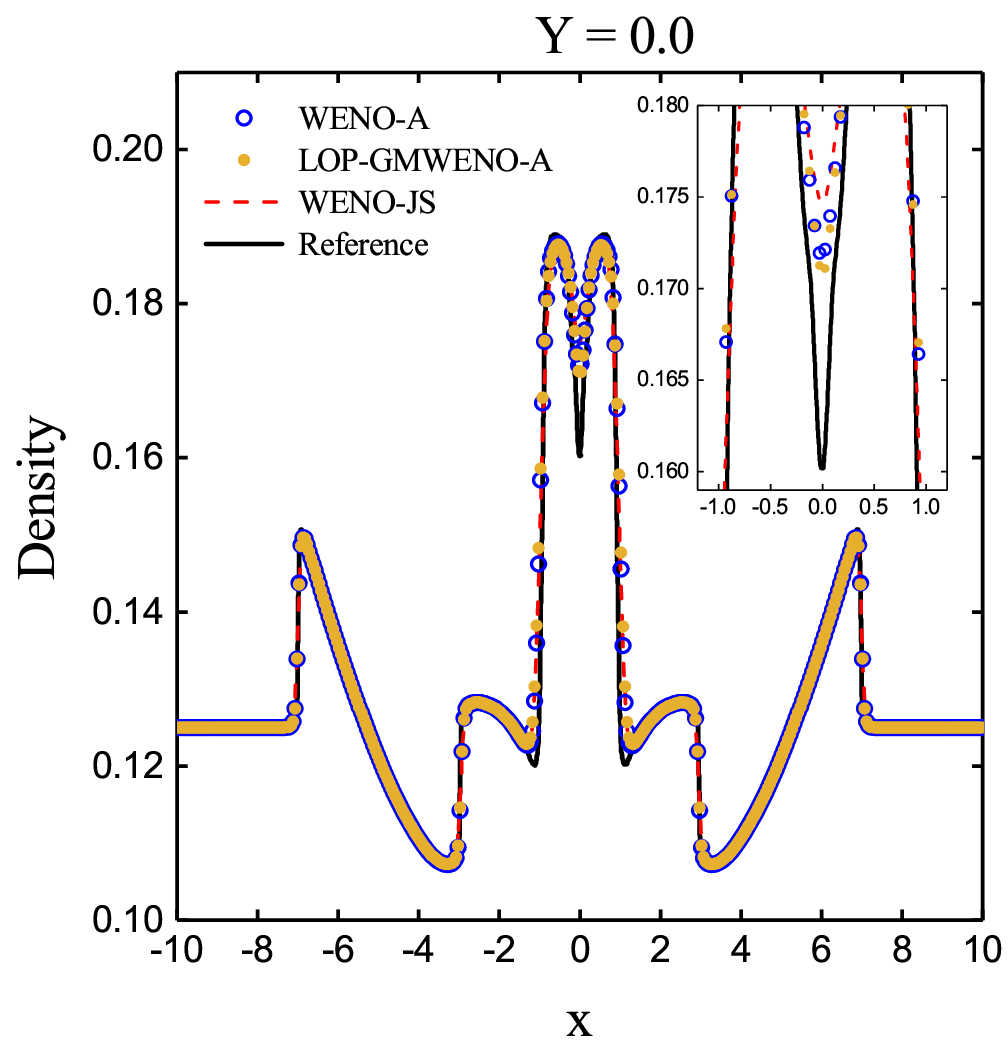}  \\         
\caption{Solutions of the 2D explosion problem. Left and middle: density contours. Right: density profiles at plane $y = 0.0$.}
\label{fig:ex:Explosion}
\end{figure}

In Fig. \ref{fig:ex:Implosion} and Fig. \ref{fig:ex:Explosion}, the 
density contours (the first two columns) and the desity profiles at 
plane $y = 0.0$ (the last column) of the LOP-GMWENO-Z, LOP-GMWENO-Z$\eta(\tau_{81})$ and LOP-GMWENO-A schemes, as well as their 
associated WENO-Z, WENO-Z$\eta(\tau_{81})$ and WENO-A schemes, are 
presented respectively. We also show the results of the WENO-JS 
scheme in the desity profiles (see the last column) where the 
reference solutions are computed by the WENO-JS scheme with $N_{x} \times N_{y} = 1000 \times 1000$. We can see that all considered 
schemes can successfully capture the main structures of the 
implosion and explosion problems. Moreover, for the implosion 
problem, as shown in the last column of Fig. \ref{fig:ex:Implosion}, 
the solutions of the LOP-GMWENO-Z and LOP-GMWENO-A schemes are very 
close to those of their associated WENO-Z and WENO-A schemes 
respectively, while the LOP-GMWENO-Z$\eta(\tau_{81})$ scheme has a 
lower resolution than the WENO-Z$\eta(\tau_{81})$ scheme. Similarly, 
for the explosion problem, as shown in the last column of Fig. \ref{fig:ex:Explosion}, the solutions of the LOP-GMWENO-Z and LOP-GMWENO-Z$\eta(\tau_{81})$ schemes are very close to those of their 
associated WENO-Z and WENO-Z$\eta(\tau_{81})$ schemes respectively, 
and the resolution of the LOP-GMWENO-A scheme appears to be higher 
than the WENO-A scheme. As expected, the WENO-JS scheme produces 
lower resolutions than all other considered schemes.

\subsubsection{Shock-vortex interaction}
The shock-vortex interaction problem is widely used to examine 
the performance of the high-resolution methods 
\cite{Shock-vortex_interaction-1,Shock-vortex_interaction-2,
Shock-vortex_interaction-3}. The computational domain is 
$[0,1]\times[0,1]$ and it is initialized by
\begin{equation*}
\big( \rho, u, v, p \big)(x, y, 0) = \left\{
\begin{aligned}
\begin{array}{ll}
\bf{U}_{\mathrm{L}}, & x < 0.5, \\
\bf{U}_{\mathrm{R}}, & x \geq 0.5, \\
\end{array}
\end{aligned}
\right.
\label{eq:initial_Euer2D:shock-vortex-interaction}
\end{equation*}
where ${\bf{U}}_{\mathrm{L}}=\big(\rho_{\mathrm{L}} + \delta \rho, u_{\mathrm{L}} + \delta u, v_{\mathrm{L}} + \delta v, p_{\mathrm{L}}+ 
\delta p\big)$ and ${\bf{U}}_{\mathrm{R}}=\big(\rho_{\mathrm{R}}, 
u_{\mathrm{R}}, v_{\mathrm{R}}, p_{\mathrm{R}}\big)$, with
\begin{equation*}
\begin{array}{l}
\begin{aligned}
&\rho_{\mathrm{L}} = 1, \qquad u_{\mathrm{L}} = \sqrt{\gamma}, 
\qquad v_{\mathrm{L}} = 0, \qquad p_{\mathrm{L}} = 1, \\
&\rho_{\mathrm{R}} = \rho_{\mathrm{L}}\bigg(\dfrac{\gamma - 1 + 
(\gamma + 1)p_{\mathrm{R}}}{\gamma + 1 + (\gamma - 1)p_{\mathrm{R}}} 
\bigg), u_{\mathrm{R}} = u_{\mathrm{L}}\bigg( \dfrac{1 - 
p_{\mathrm{R}}}{\sqrt{\gamma-1 + p_{\mathrm{R}}(\gamma + 1)}}\bigg), 
v_{\mathrm{R}} = 0, p_{\mathrm{R}} = 1.3,\\
&\delta \rho = \dfrac{\rho_{\mathrm{L}}^{2}}{(\gamma - 1)
p_{\mathrm{L}}}\delta T, 
\delta u = \epsilon \dfrac{y - y_{\mathrm{c}}}{r_\mathrm{c}}
\mathrm{e}^{\alpha(1-r^{2})}, 
\delta v = - \epsilon \dfrac{x - x_{\mathrm{c}}}{r_\mathrm{c}}
\mathrm{e}^{\alpha(1-r^{2})}, 
\delta p = \dfrac{\gamma \rho_{\mathrm{L}}^{2}}{(\gamma - 1)
\rho_{\mathrm{L}}}\delta T,
\end{aligned}
\end{array}
\end{equation*}
and $\epsilon = 0.3, r_{\mathrm{c}} = 0.05, \alpha = 0.204, 
x_{\mathrm{c}}=0.25, y_{\mathrm{c}}=0.5, r=\sqrt{((x-x_{\mathrm{c}}
)^{2} + (y - y_{\mathrm{c}})^{2})/r_{\mathrm{c}}^{2}}, \delta T=-(
\gamma-1)\epsilon^{2}\mathrm{e}^{2\alpha (1-r^{2})}/(4\alpha\gamma)$.
The transmissive boundary condition is used on all edges. The 
computational domain is discretized into $800 \times 800$ uniform 
cells, and the calculations are advanced in time up to $t = 0.6$.

\begin{figure}[!ht]
\centering
  \includegraphics[height=0.34\textwidth]
  {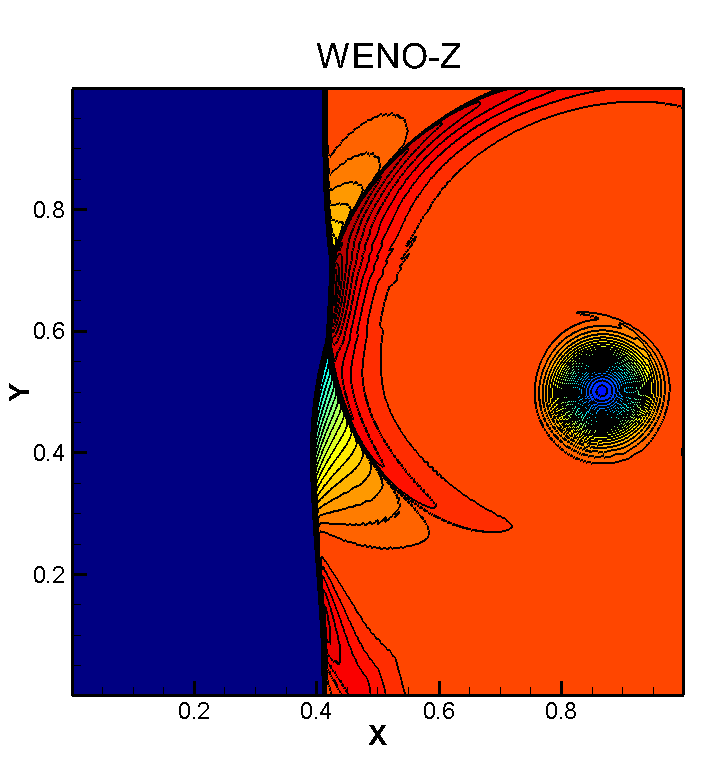}   
  \includegraphics[height=0.34\textwidth]
  {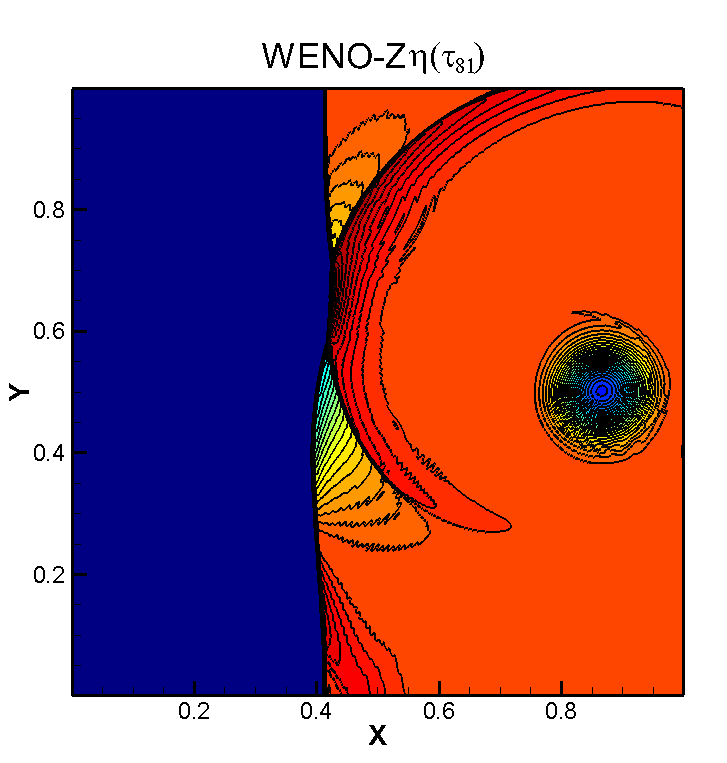}  
  \includegraphics[height=0.34\textwidth]
  {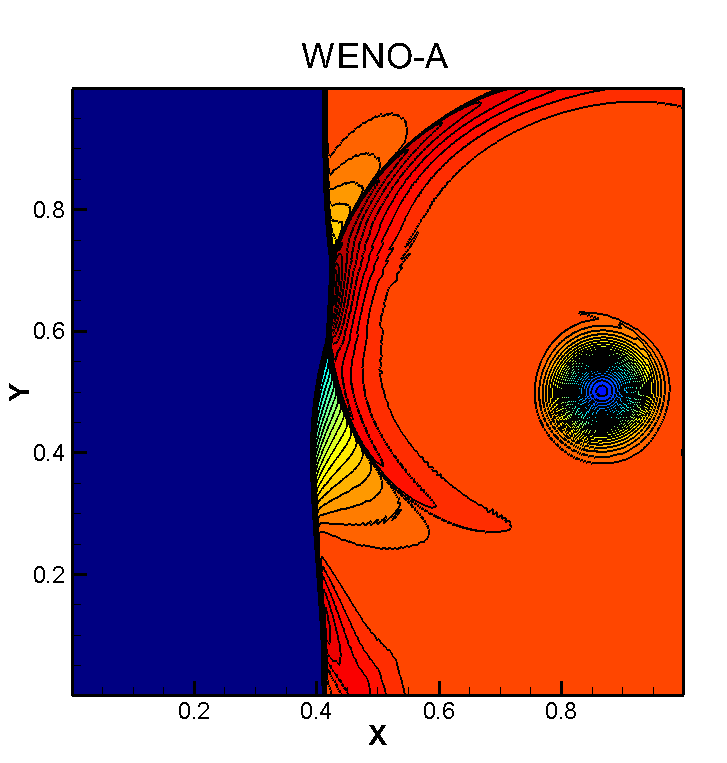}      \\
  \includegraphics[height=0.34\textwidth]
  {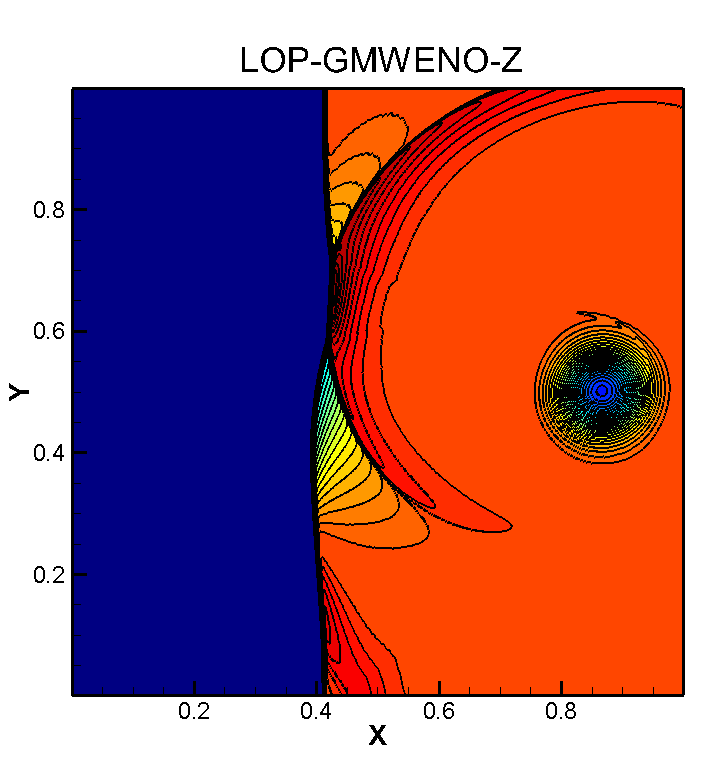}
  \includegraphics[height=0.34\textwidth]
  {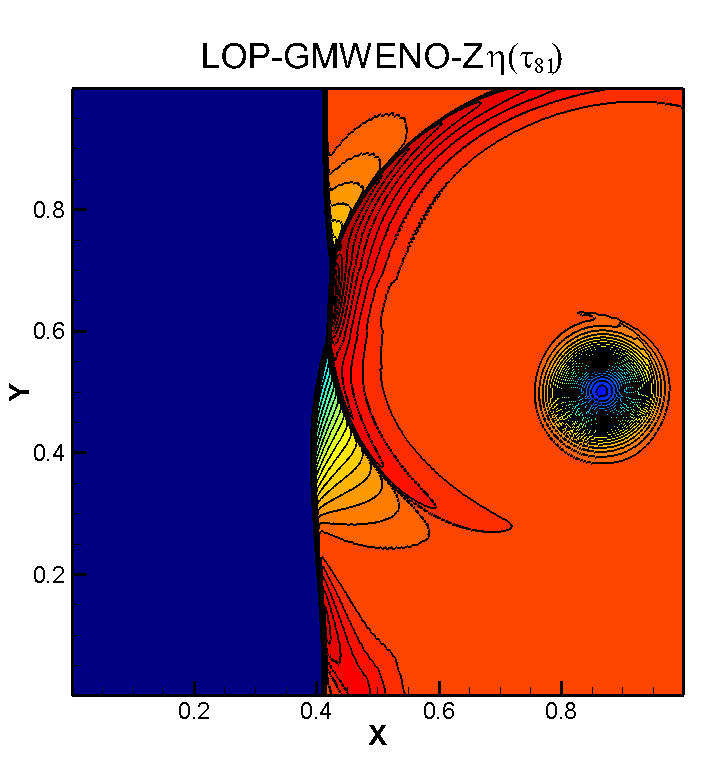}  
  \includegraphics[height=0.34\textwidth]
  {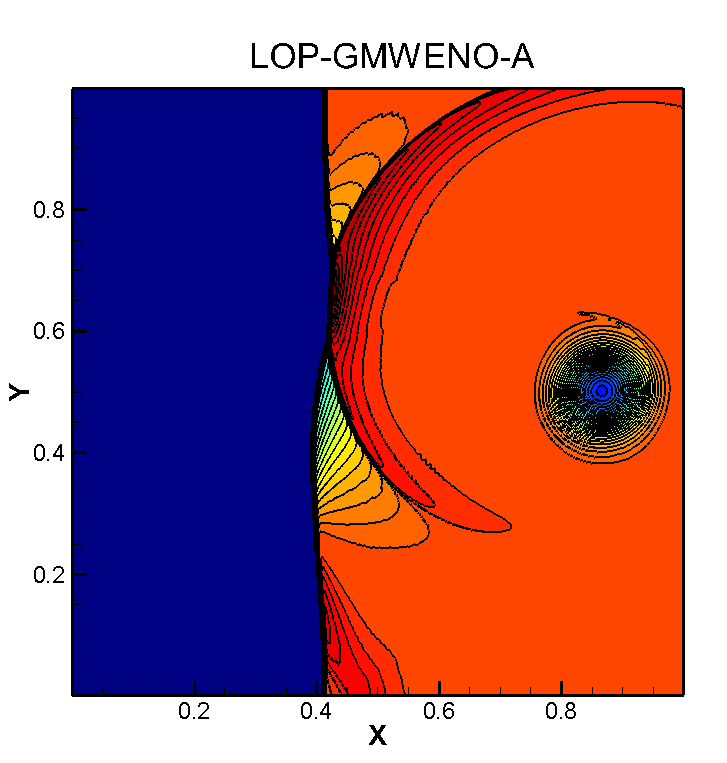} 
\caption{Density plots for the Shock-vortex interaction using $50$ 
contour lines with range from $1.0$ to $2.7$.}
\label{fig:ex:SVI:1}
\end{figure}

\begin{figure}[!ht]
\centering
  \includegraphics[height=0.30\textwidth]
  {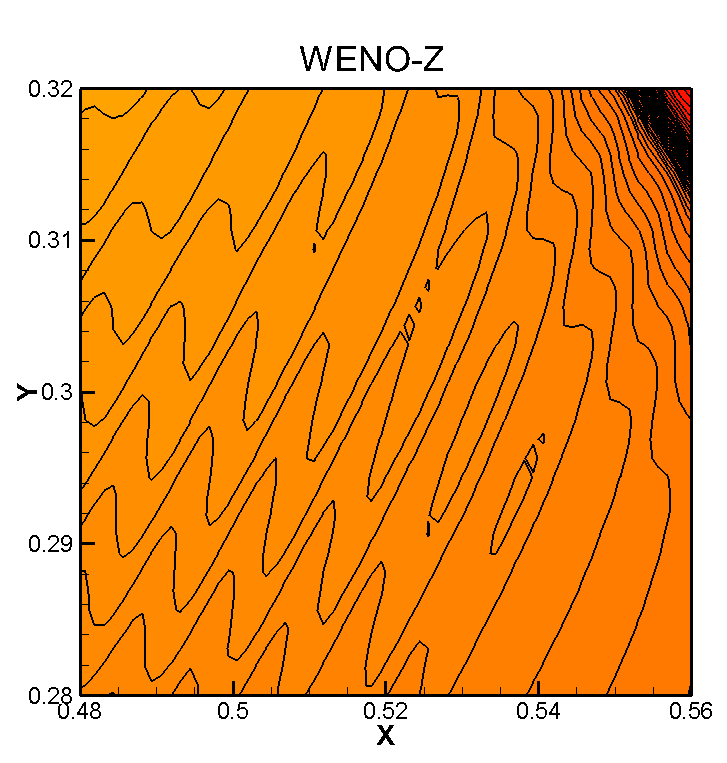}        \hspace{1.2ex}
  \includegraphics[height=0.30\textwidth]
  {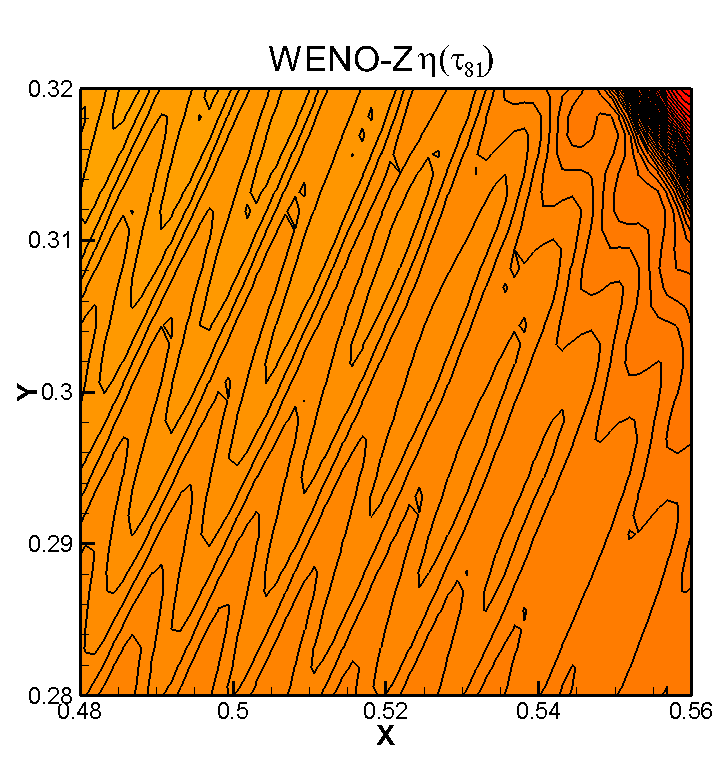}    \hspace{1.2ex}
  \includegraphics[height=0.30\textwidth]
  {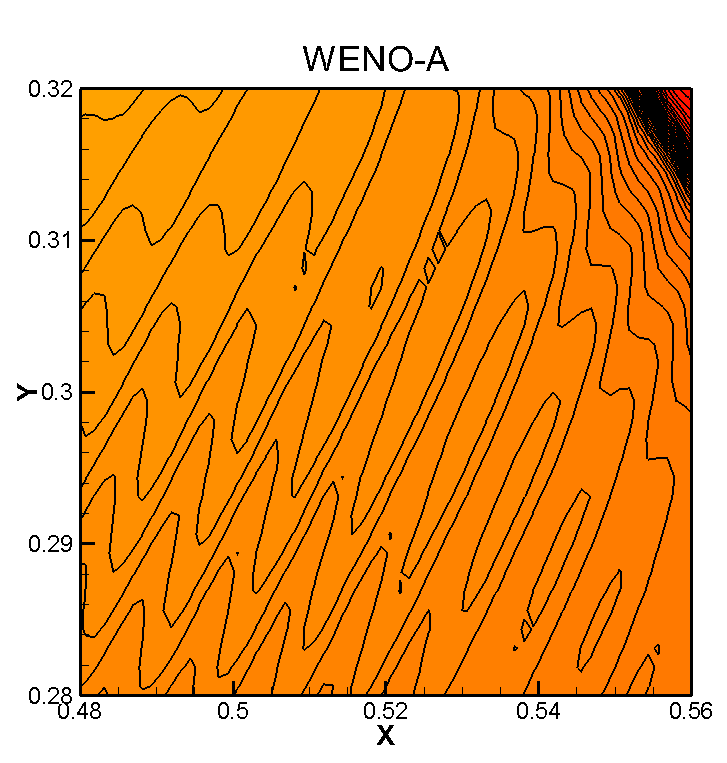}     \\
  \includegraphics[height=0.30\textwidth]
  {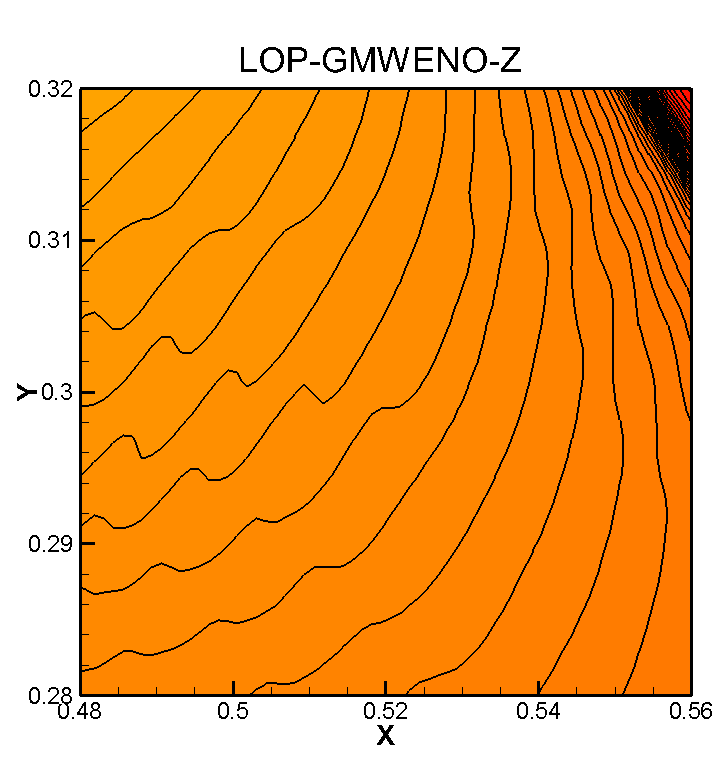}    \hspace{1.2ex}
  \includegraphics[height=0.30\textwidth]
  {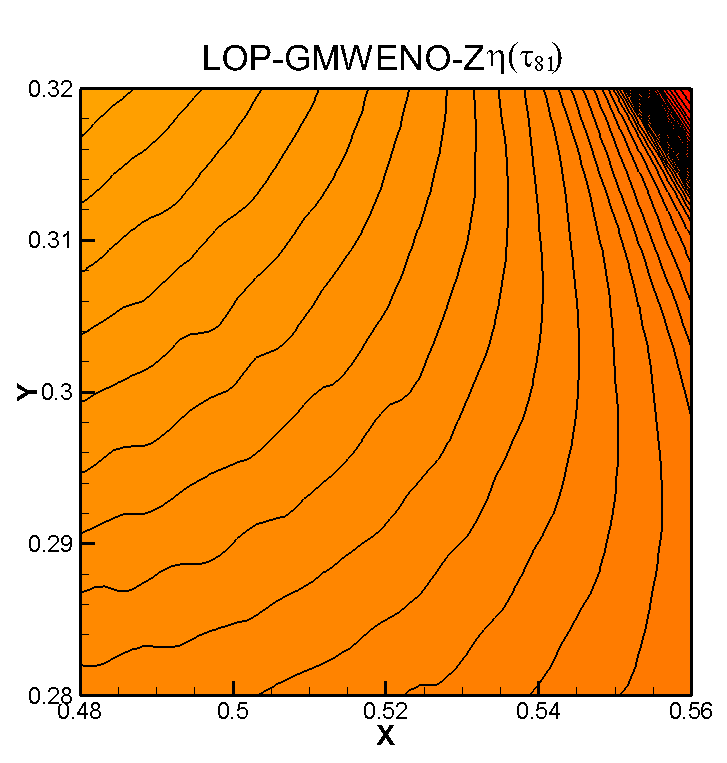} \hspace{1.2ex}
  \includegraphics[height=0.30\textwidth]
  {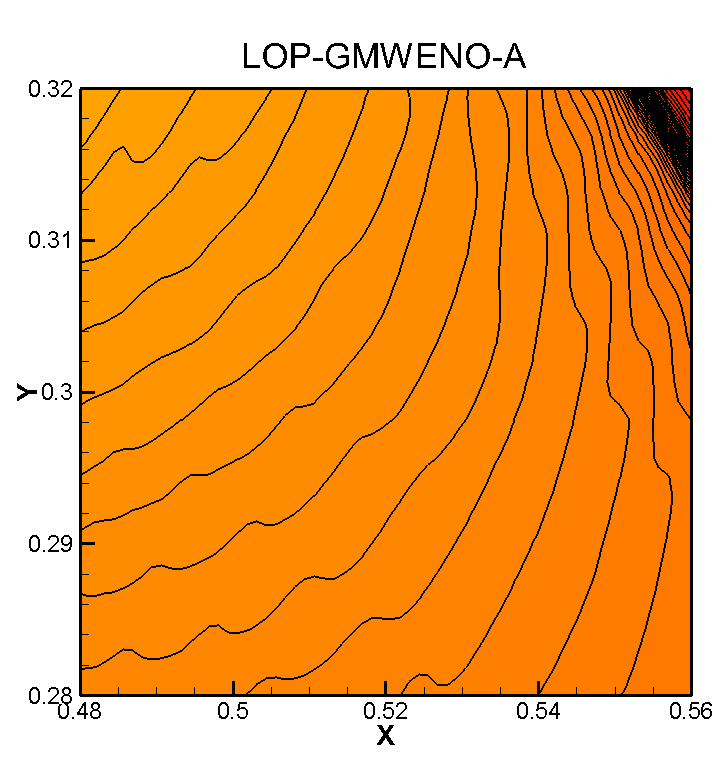}\\    
\caption{The zoomed-in density plots using $400$ contour lines with range from $1.0$ to $2.7$.}
\label{fig:ex:SVI:2}
\end{figure}

\begin{figure}[!ht]
\centering
  \includegraphics[height=0.26\textwidth]
  {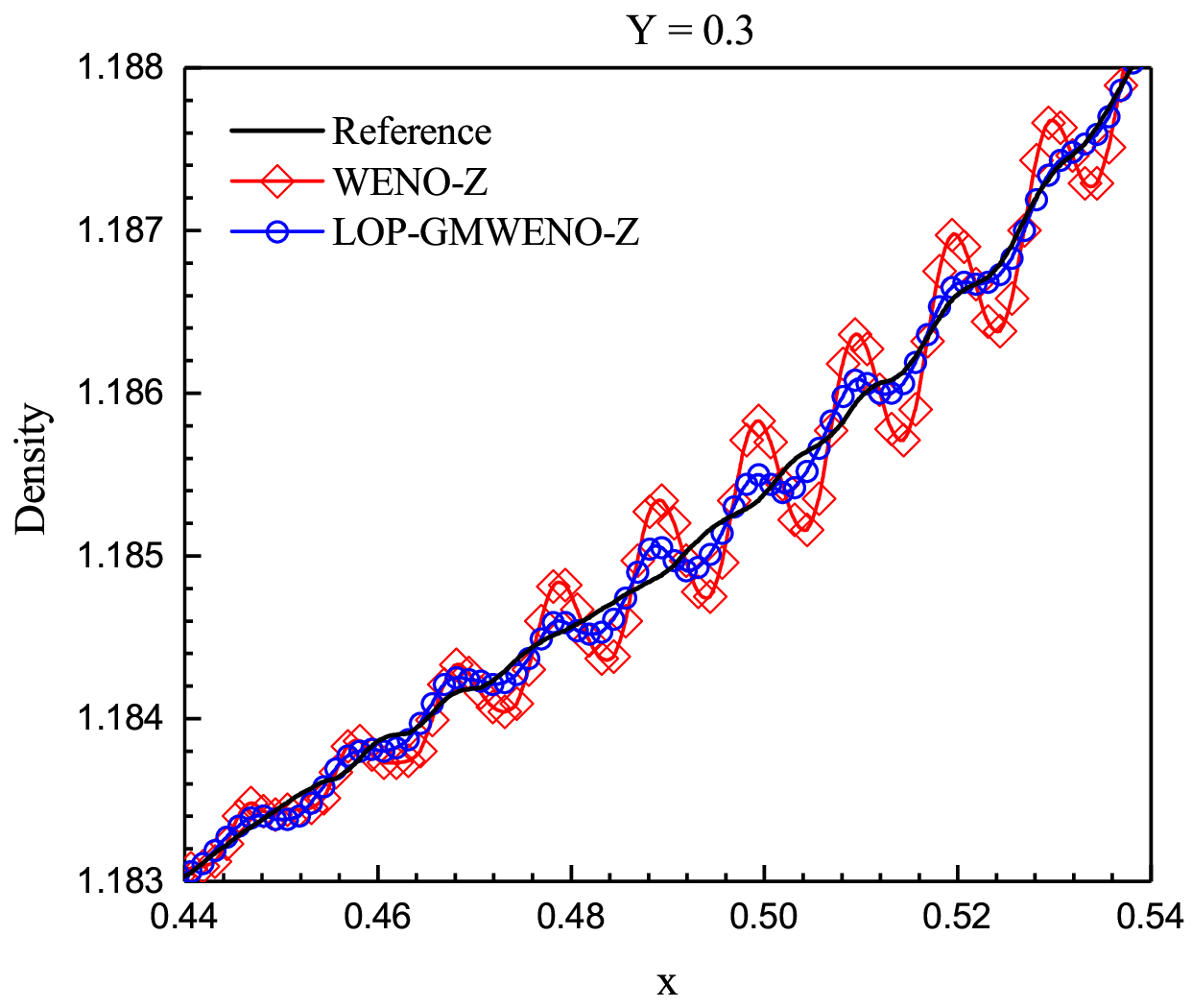}     
  \includegraphics[height=0.26\textwidth]
  {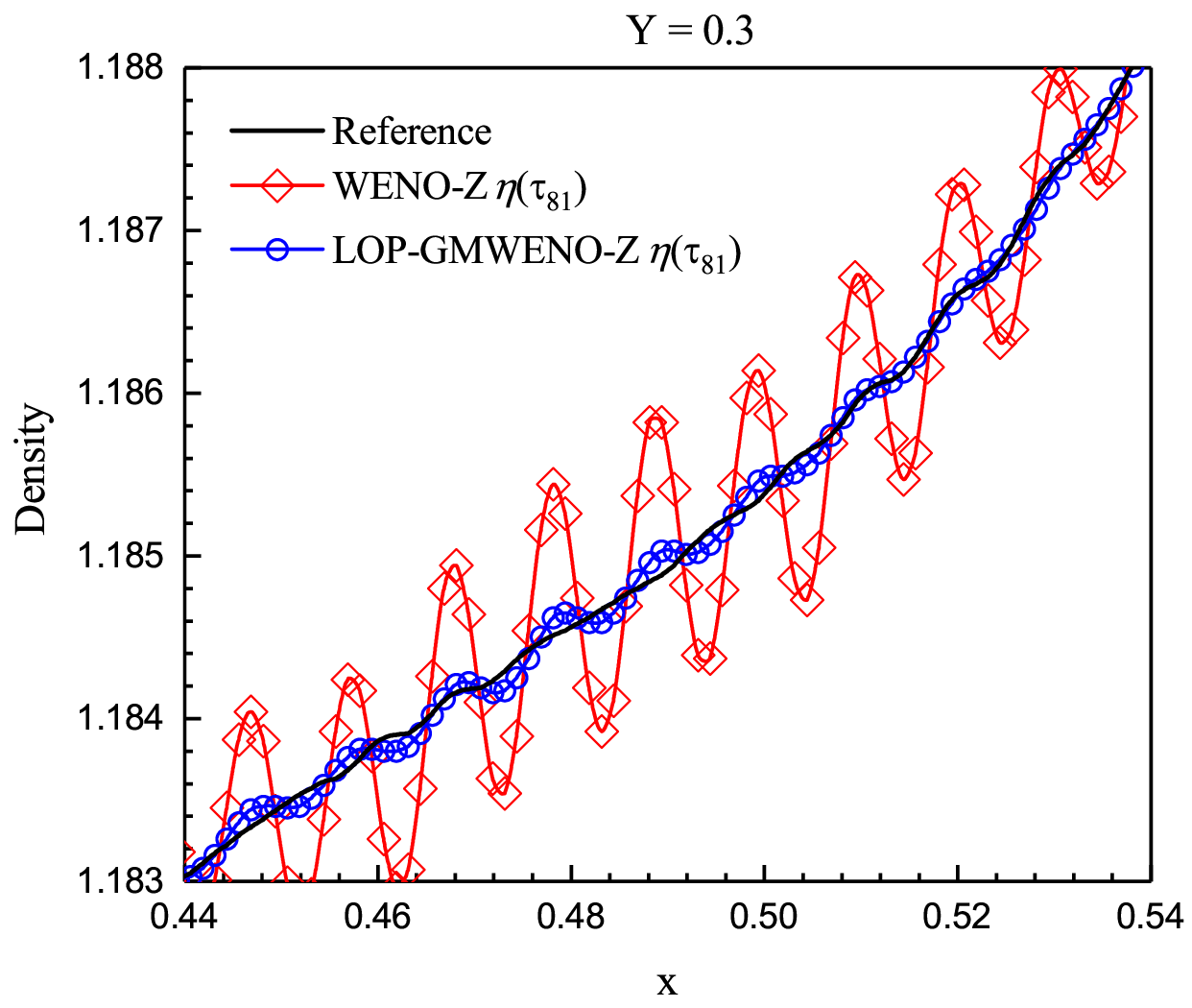} 
  \includegraphics[height=0.26\textwidth]
  {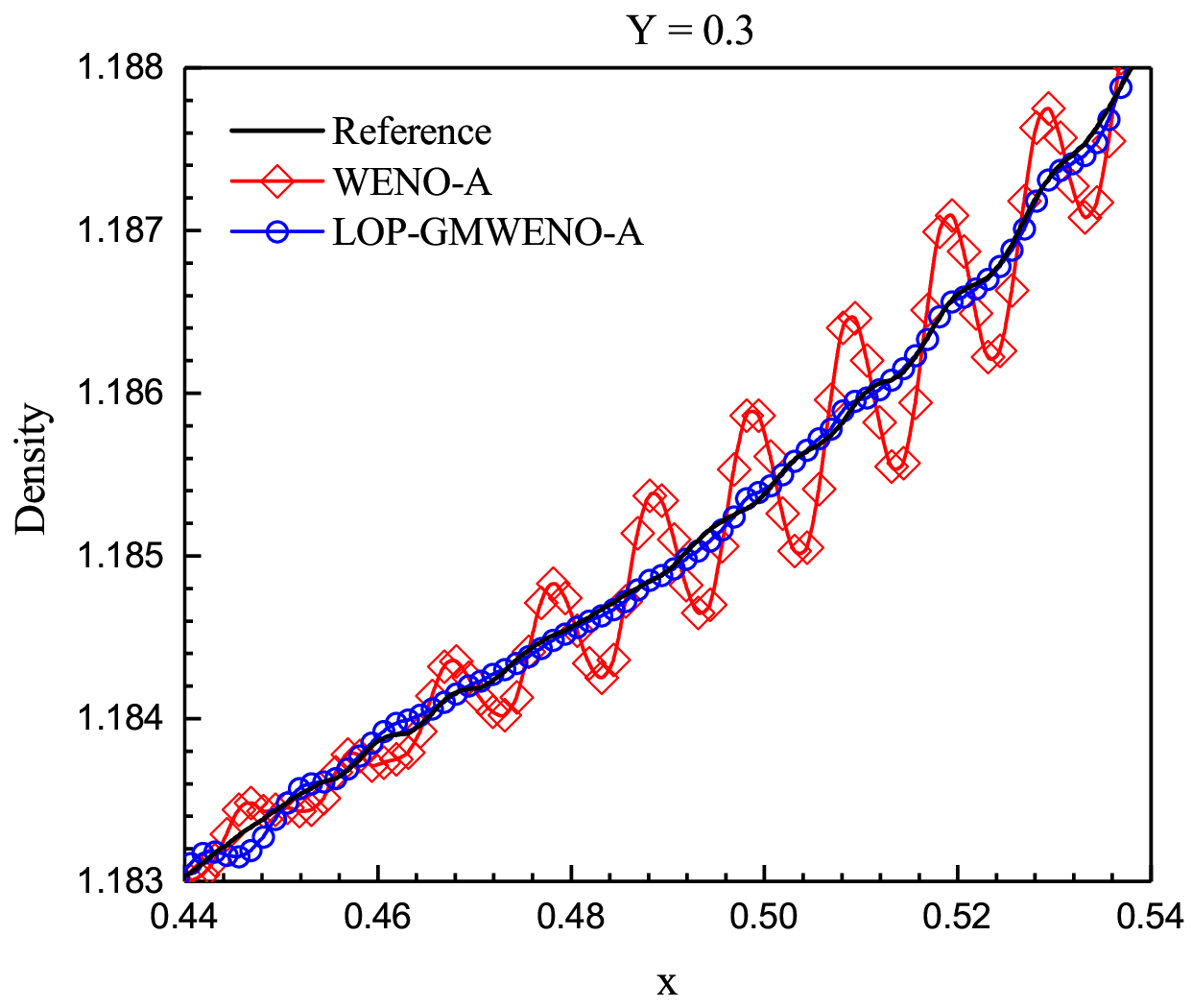}  
\caption{The cross-sectional slices of density plot along the plane 
$y = 0.30$ where $x \in [0.44, 0.54]$.}
\label{fig:ex:SVI:4}
\end{figure}

We consider the LOP-GMWENO-X schemes with X = Z, Z$\eta(\tau_{81})$, 
A and their associated WENO-X schemes (hereinafter the same). Fig. 
\ref{fig:ex:SVI:1} gives the solutions computed by all considered 
WENO schemes, and the main structure of the shock and vortex after 
the interaction are properly captured by all considered schemes. 
Fig. \ref{fig:ex:SVI:2} presents the zoomed-in view to show the 
numerical oscillations more clearly. It indicates that the WENO-X 
schemes generate severe post-shock oscillations while their 
associated LOP-GMWENO-X schemes can significantly decrease these 
oscillations. In order to demonstrate this more clearly, in Fig. 
\ref{fig:ex:SVI:4}, we give the cross-sectional slices of density 
plot along the plane $y = 0.30$ in $x \in (0.44, 0.54)$ where the 
reference solution is obtained using the WENO-JS scheme with a 
uniform mesh size of $1600 \times 1600$. Of course, the LOP-GMWENO-X 
schemes can not thoroughly remove the numerical oscillations. 
However, it is obvious that they can significantly reduce these 
oscillations compared to their associated WENO-X schemes. Indeed, 
this is a huge improvement for the WENO-Z$\eta(\tau_{81})$ and 
WENO-A schemes on reducing the post-shock oscillations, and the 
improvement for the WENO-Z scheme is also noticeable. Therefore, we 
conclude that this should be an additional advantage of the 
WENO-Z-type schemes with \textit{LOP} generalized mappings.

\subsubsection{Regular shock reflection}
As used in \cite{WENO-ZS}, we simulate the regular shock 
reflection. This is a typical 2D steady flow. The computational 
domain is a rectangle of 4 length units times 1 length unit, say, 
$[0, 4]\times[0, 1]$ here. It is initialized by
\begin{equation*}
\big( \rho, u, v, p \big)(x, y, 0) = (1.0, 2.9, 0.0, 1.0/\gamma),
\label{eq:IC_Euer2D:RSR}
\end{equation*}
where $\gamma = 1.4$. Along the bottom and the right edges, the 
reflection and transmissive boundary conditions are used 
respectively. The following Dirichlet boundary conditions are used 
along the top and left edges
\begin{equation*}
\big( \rho, u, v, p \big)(x, y, t) = \left\{
\begin{array}{ll}
(1.69997, 2.61934, -0.50632, 1.52819), & \text{where} \quad y = 1,\\
(1.0, 2.9, 0.0, 1.0/\gamma), & \text{where} \quad x = 0.
\end{array}
\right.
\label{eq:BC_Euer2D:RSR}
\end{equation*}
A uniform mesh size of $800 \times 200$ is used and the final time 
is $t = 2.5$.

\begin{figure}[!ht]
\centering
  \includegraphics[height=0.17\textwidth]
  {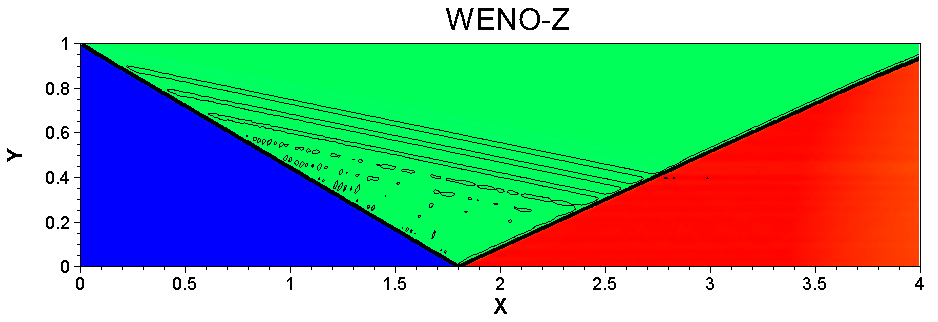}  
  \includegraphics[height=0.17\textwidth]
  {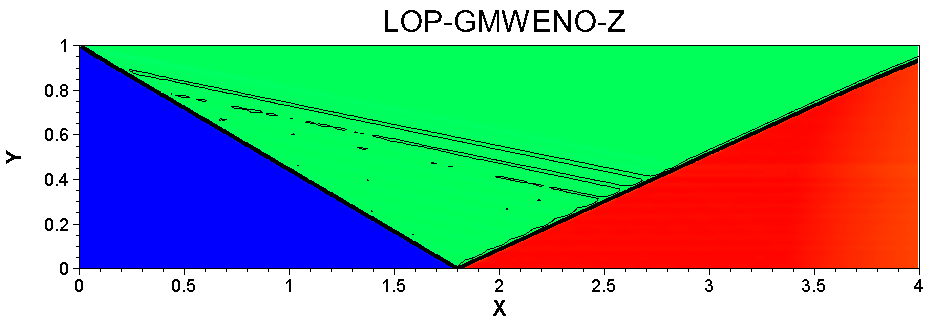}\\
  \vspace{4.0ex}
  \includegraphics[height=0.17\textwidth]
  {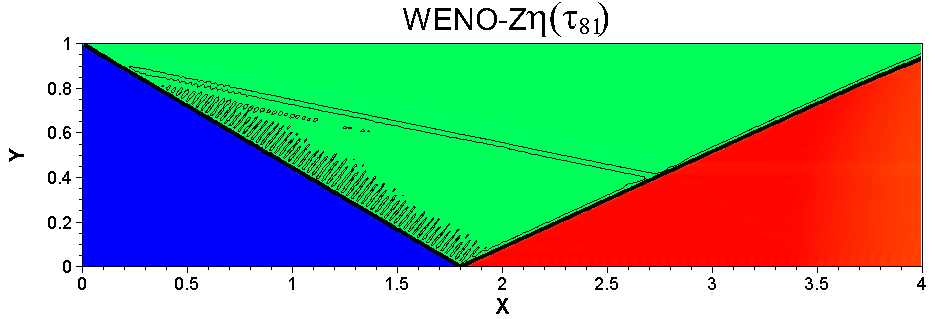}    
  \includegraphics[height=0.170\textwidth]
  {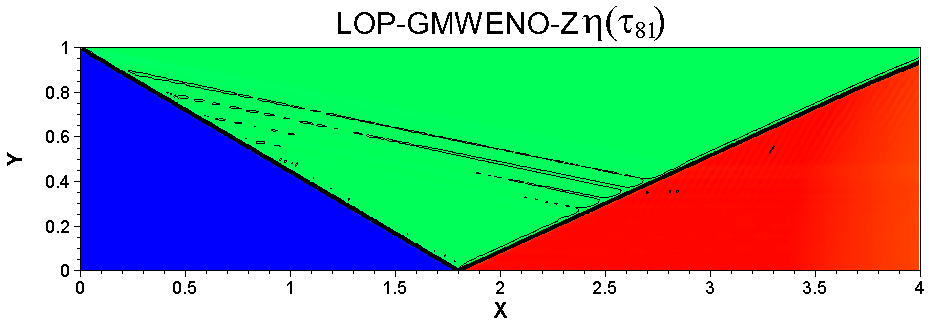}\\
  \vspace{4.0ex}
  \includegraphics[height=0.17\textwidth]
  {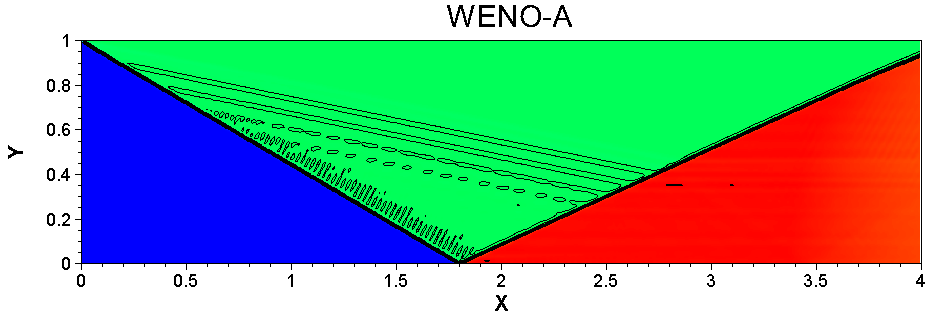}       
  \includegraphics[height=0.17\textwidth]
  {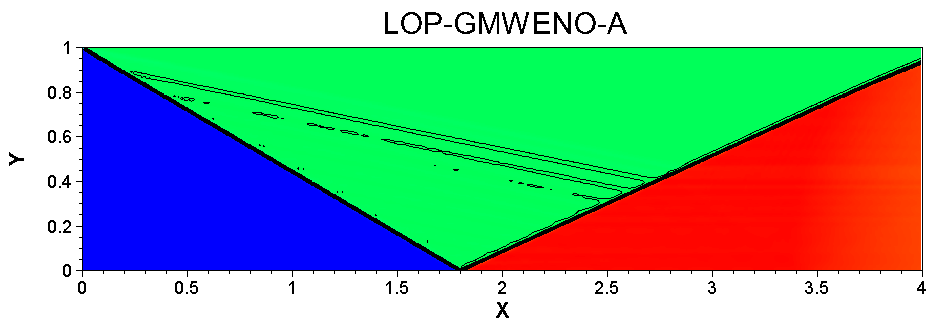}
\caption{Density plots for the regular shock reflection using $30$ 
contour lines with range from $1.05$ to $2.7$.}
\label{fig:ex:RSR:1}
\end{figure}

\begin{figure}[!ht]
\flushleft
  \includegraphics[height=0.26\textwidth]
  {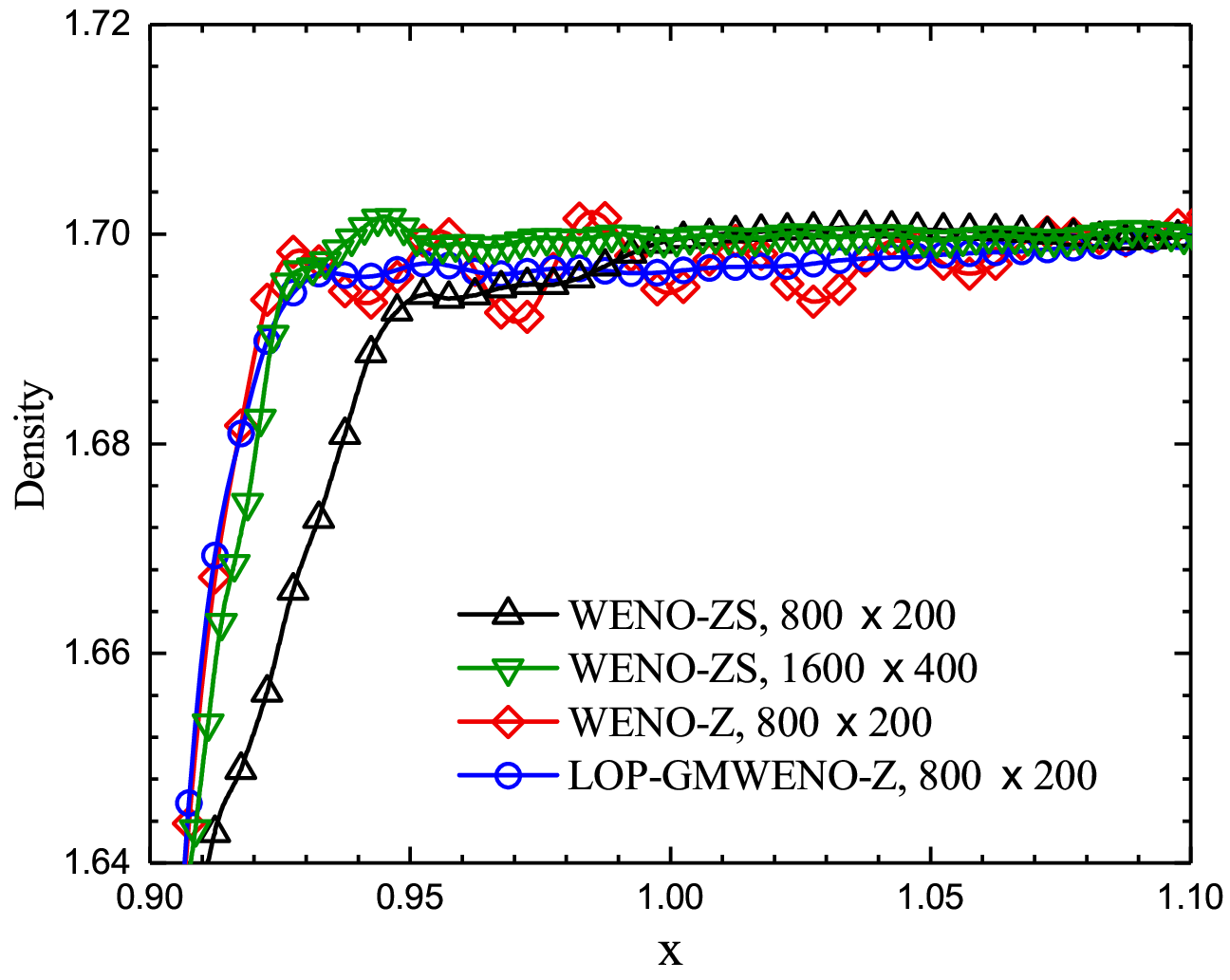}     
  \includegraphics[height=0.26\textwidth]
  {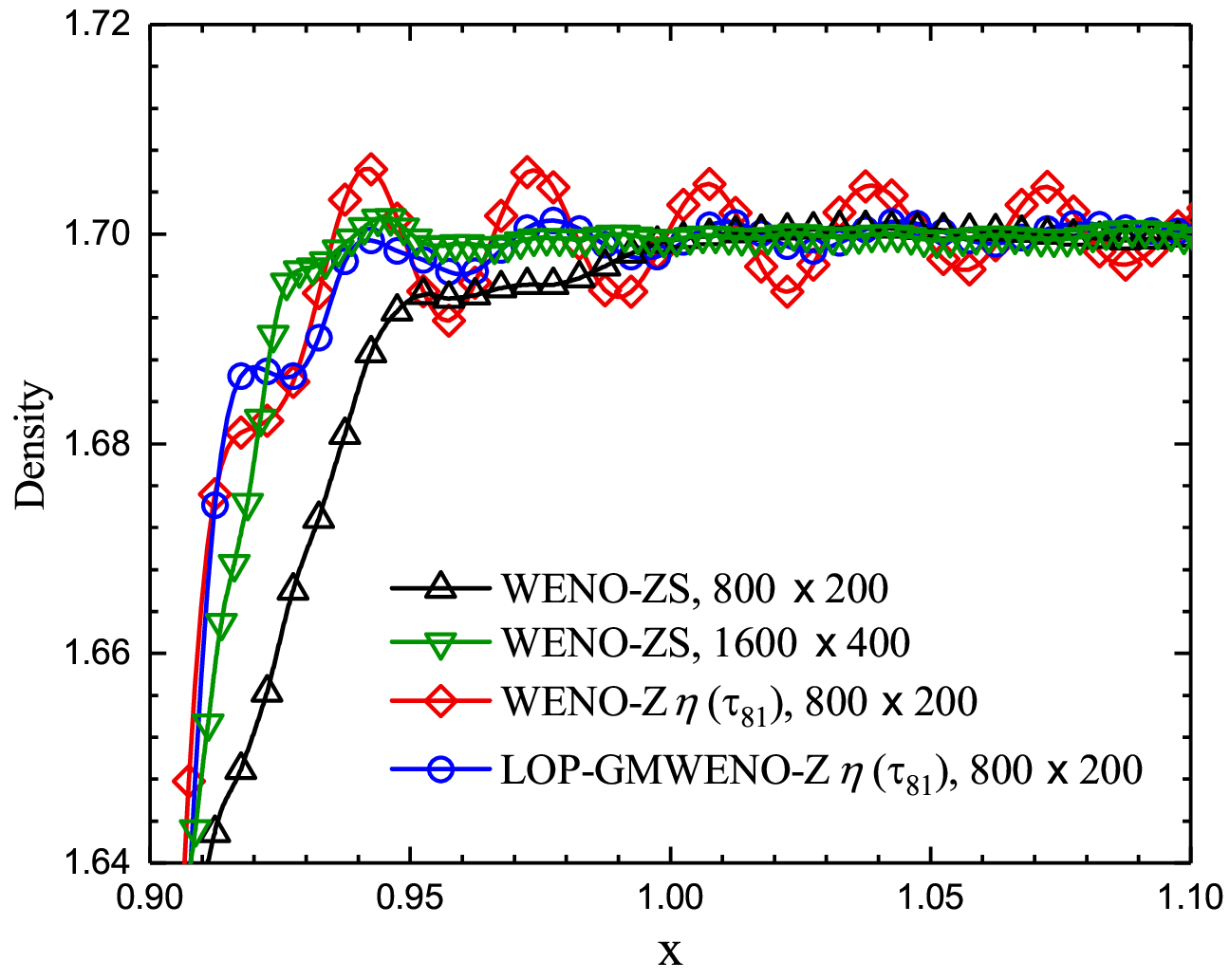} 
  \includegraphics[height=0.26\textwidth]
  {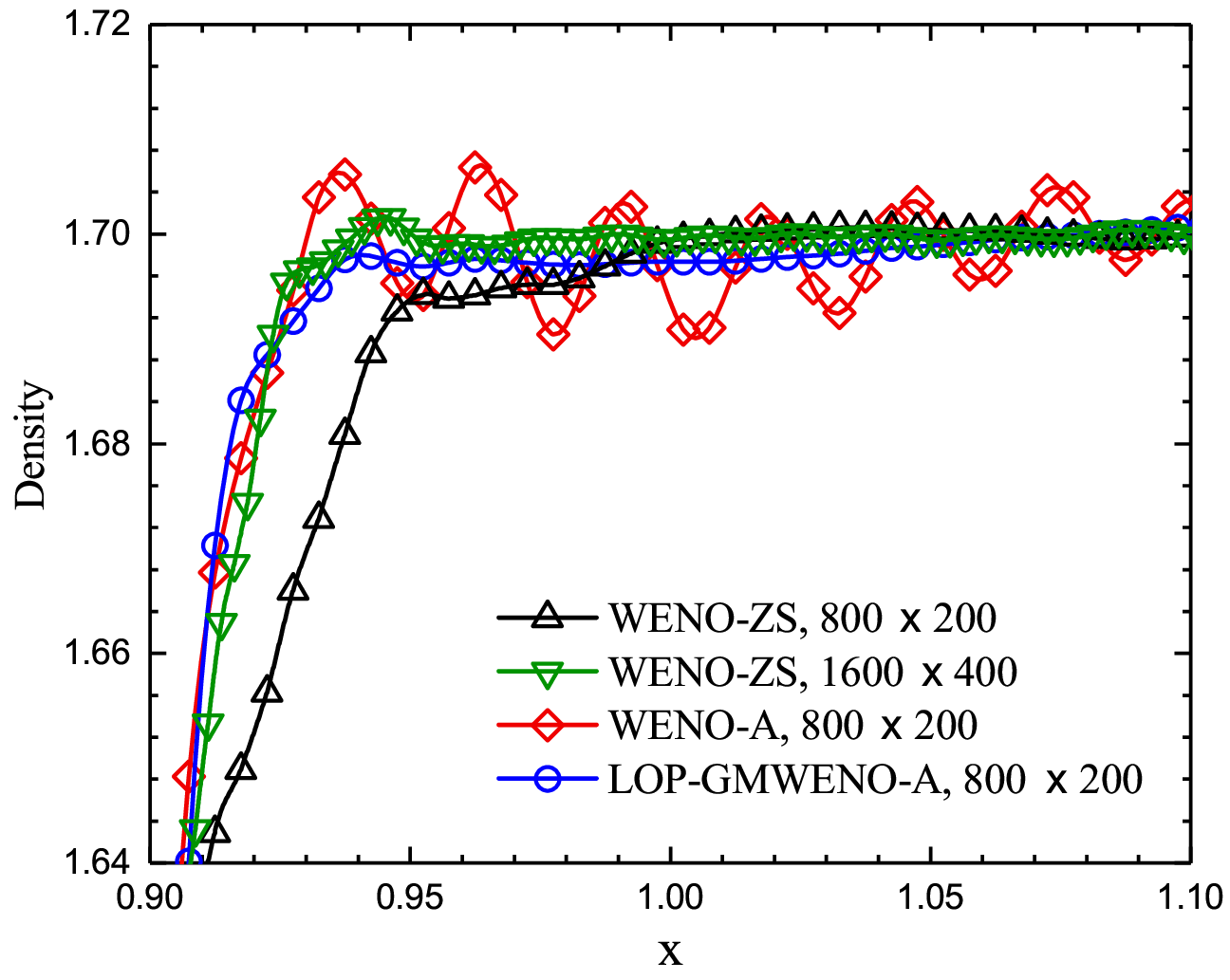}  
\caption{The cross-sectional slices of density plot along the plane 
$y = 0.50$ where $x \in [0.9, 1.1]$.}
\label{fig:ex:RSR:2}
\end{figure}

In Fig. \ref{fig:ex:RSR:1}, we show the density contours computed by 
the LOP-GMWENO-X schemes (the right column) and their associated 
WENO-X schemes (the left column). In general, all the considered 
schemes can capture the main structure of the shock transitions for 
this problem. Unfortunately, all these schemes produce the 
post-shock numerical oscillations. However, in comparison with the 
WENO-X schemes without LOP, their associated LOP-GMWENO-X schemes 
can significantly reduce these oscillations. In order to more 
clearly manifest this, we show the cross-sectional slices of density 
plot along the plane $y = 0.5$ in $x \in (0.9, 1.1)$ in Fig. 
\ref{fig:ex:RSR:2}. For comparison purpose, we also plot the results 
of the WENO-ZS \cite{WENO-ZS} scheme using two different mesh sizes 
of $N_{x} \times N_{y} = 800 \times 200, 1600 \times 400$. It is 
well known that WENO-ZS is able to either remove or significantly 
reduce the post-shock oscillations. It can easily be found that the 
post-shock oscillations produced by the WENO-X schemes are much 
severer than those of their associated LOP-GMWENO-X schemes. In 
other words, the post-shock oscillations of the LOP-GMWENO-X schemes 
are considerably reduced compared to those of their associated 
WENO-X schemes. As mentioned before, this should be an advantage of 
the WENO-Z-type schemes with \textit{LOP} mappings. Furthermore, 
although the post-shock oscillations generated by the LOP-GMWENO-X 
schemes are slightly severer than that of the WENO-ZS scheme, the 
LOP-GMWENO-X schemes achieve higher resolutions than the WENO-ZS 
scheme under the same mesh resolution. Indeed, the resolutions of 
the LOP-GMWENO-X schemes with $N_{x} \times N_{y} = 800 \times 200$ 
is comparable to that of WENO-ZS with $N_{x} \times N_{y} = 1600 
\times 400$. In summary, this might be another additional advantage 
of the LOP-GMWENO-X schemes in some ways.

\subsubsection{Double Mach reflection of a strong shock}
This problem was initially proposed by Woodward and Colella \cite{interactingBlastWaves-Woodward-Colella}. The computational domain is set to be $[0, 4]\times[0, 1]$ and it is initialized by
\begin{equation*}
\big( \rho, u, v, p \big)(x, y, 0) = \left\{
\begin{aligned}
\begin{array}{ll}
(8.0, 8.25\cos \theta, -8.25\sin \theta, 116.5),         & x < x_{0} + \dfrac{y}{\sqrt{3}}, \\
(1.4, 0.0, 0.0, 1.0),    & x \geq x_{0} + \dfrac{y}{\sqrt{3}}, \\
\end{array}
\end{aligned}
\right.
\label{eq:initial_Euer2D:DMR}
\end{equation*}
where $\theta = \frac{\pi}{6}$ and $x_{0}=\frac{1}{6}$. The inflow 
boundary condition with the post-shock values as stated above and 
the outflow boundary condition are used at $x = 0$ and $x = 4$ 
respectively. At $y = 0$, the post-shock values are imposed at $[0, x_{0})$, while the reflective boundary condition is applied to $[x_{0}, 4]$. At $y = 1$, the fluid variables are defined as to the exact 
solution of the Mach 10 moving oblique shock. The uniform meshes of $N_{x} \times N_{y} = 2000 \times 500$ and the output time $t = 0.2$ 
are used.

\begin{figure}[!ht]
\flushleft
  \includegraphics[height=0.2\textwidth]
  {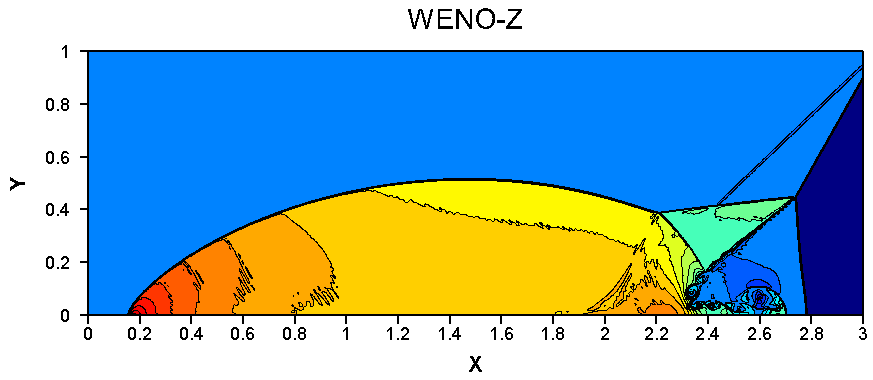} \hspace{1.8em}
  \includegraphics[height=0.2\textwidth]
  {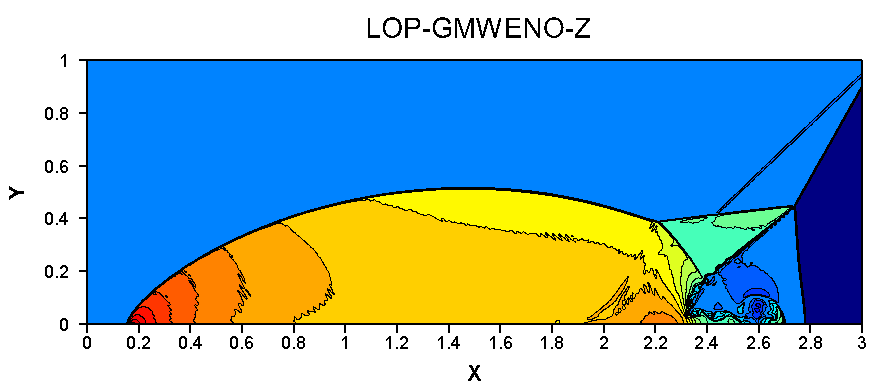}   \\
  \includegraphics[height=0.2\textwidth]
  {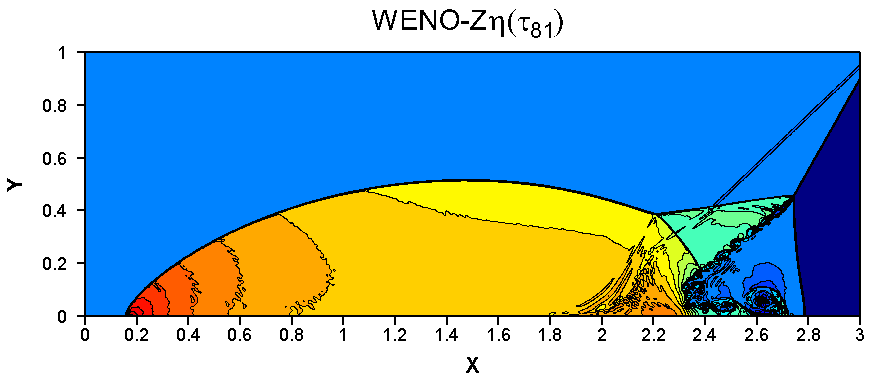} \hspace{1.8em} 
  \includegraphics[height=0.2\textwidth]
  {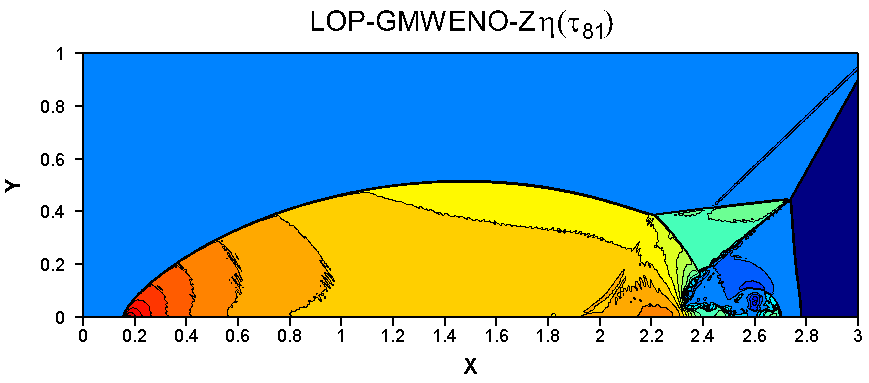}   \\         
  \includegraphics[height=0.2\textwidth]
  {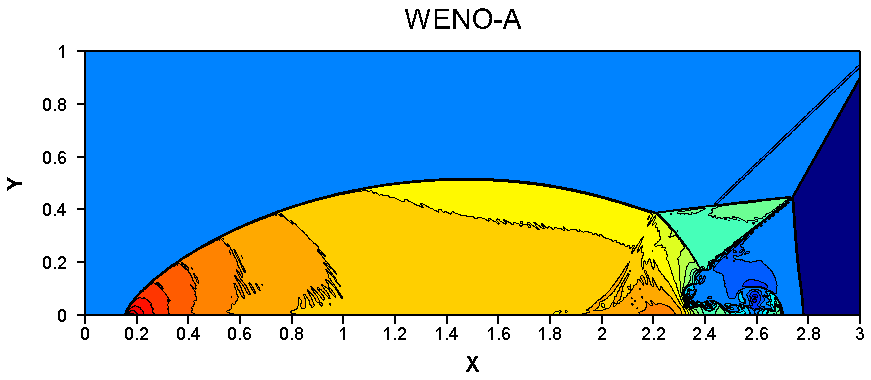} \hspace{1.8em} 
  \includegraphics[height=0.2\textwidth]
  {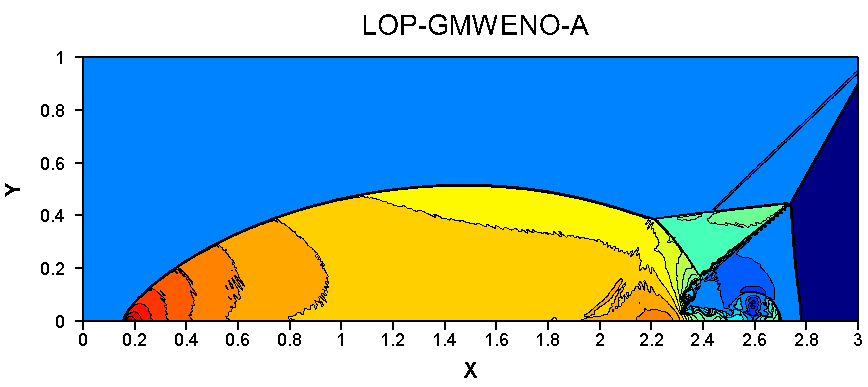}\\
  \includegraphics[height=0.305\textwidth]
  {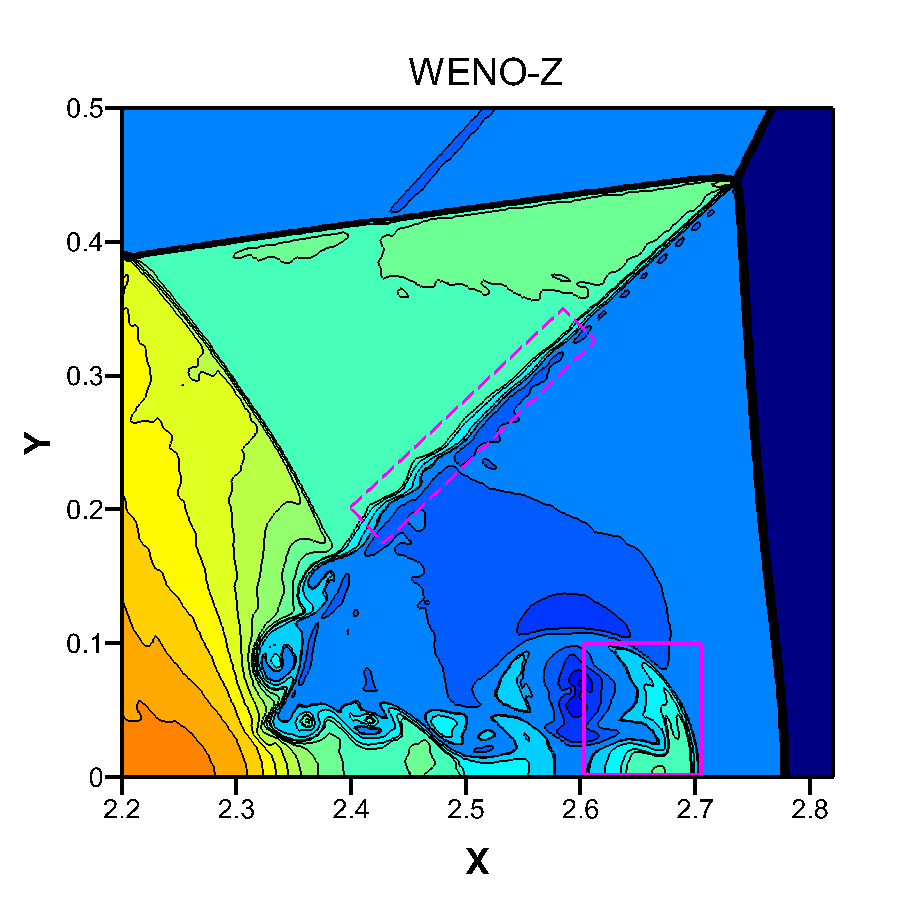} \hspace{0.9em}  
  \includegraphics[height=0.305\textwidth]
  {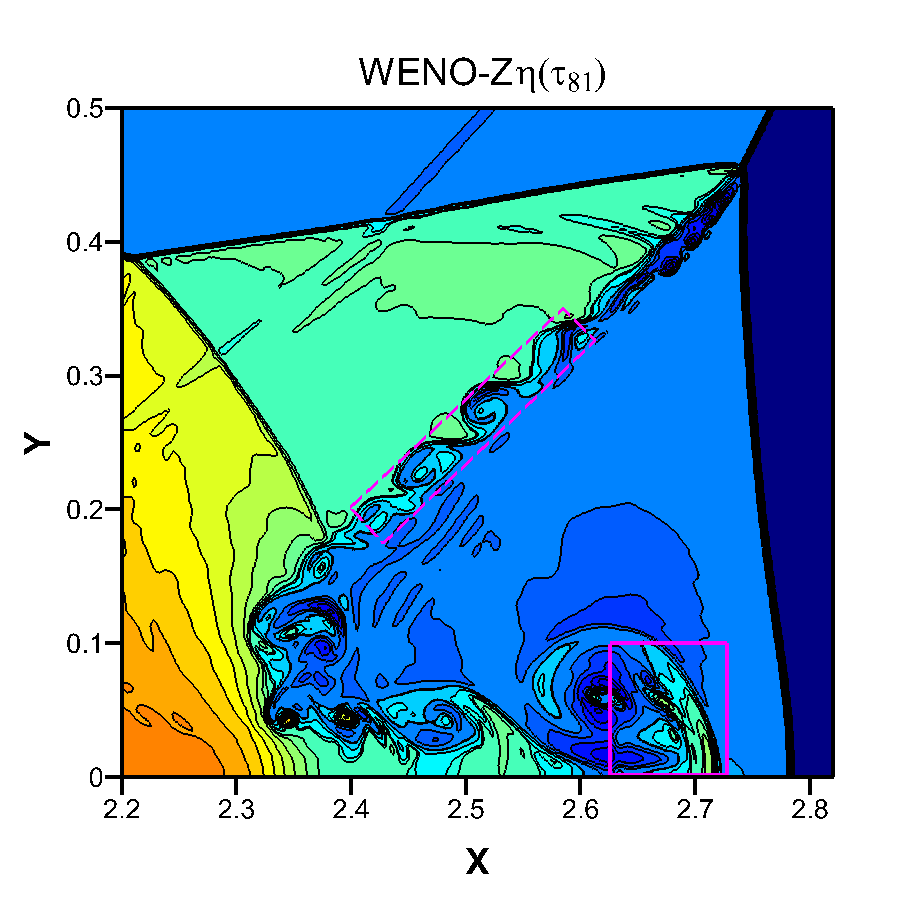} \hspace{0.9em}  
  \includegraphics[height=0.305\textwidth]
  {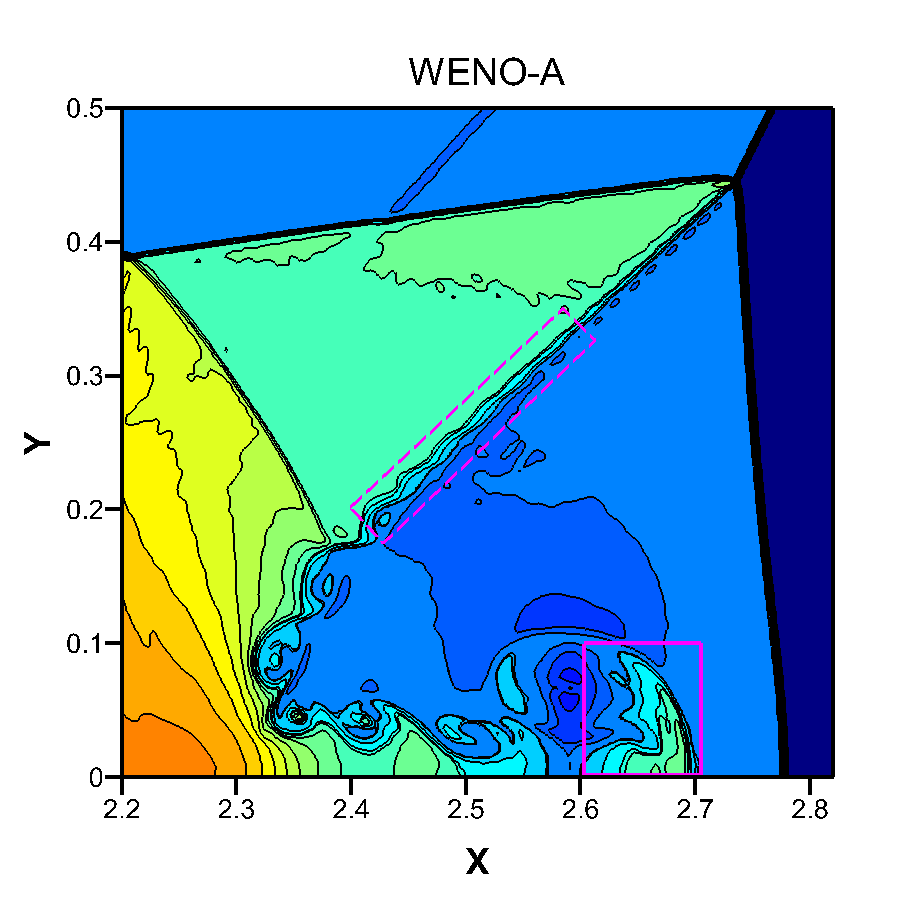} \\
  \includegraphics[height=0.305\textwidth]
  {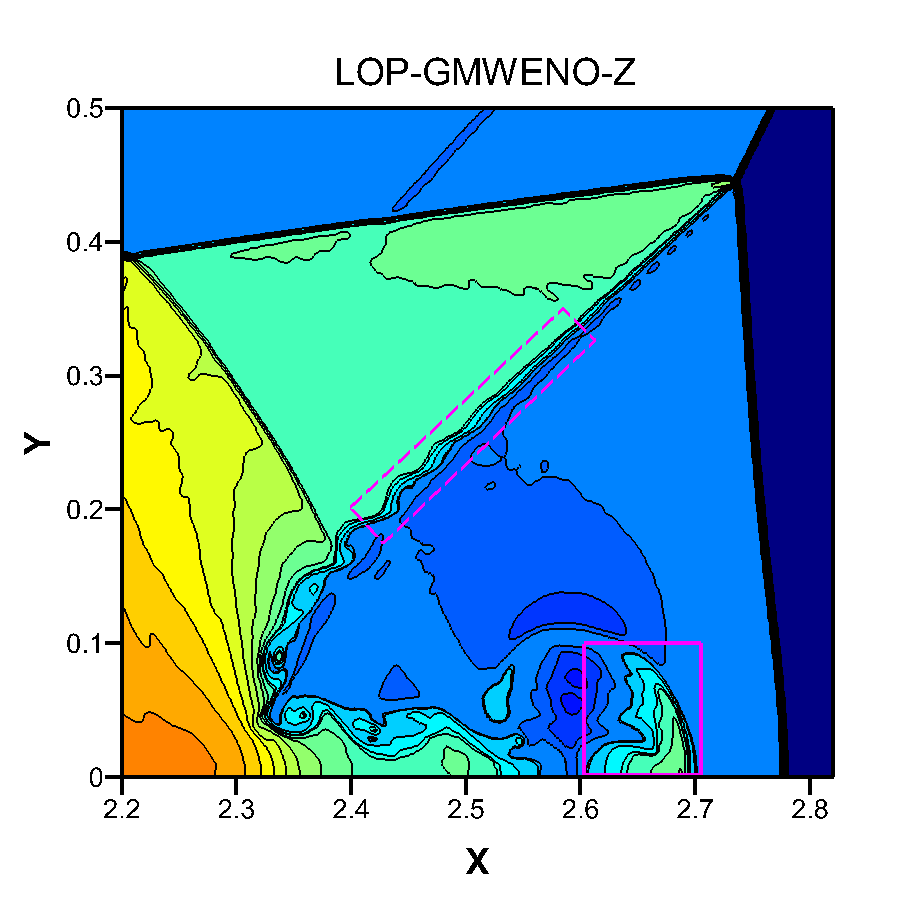} \hspace{0.9em}  
  \includegraphics[height=0.305\textwidth]
  {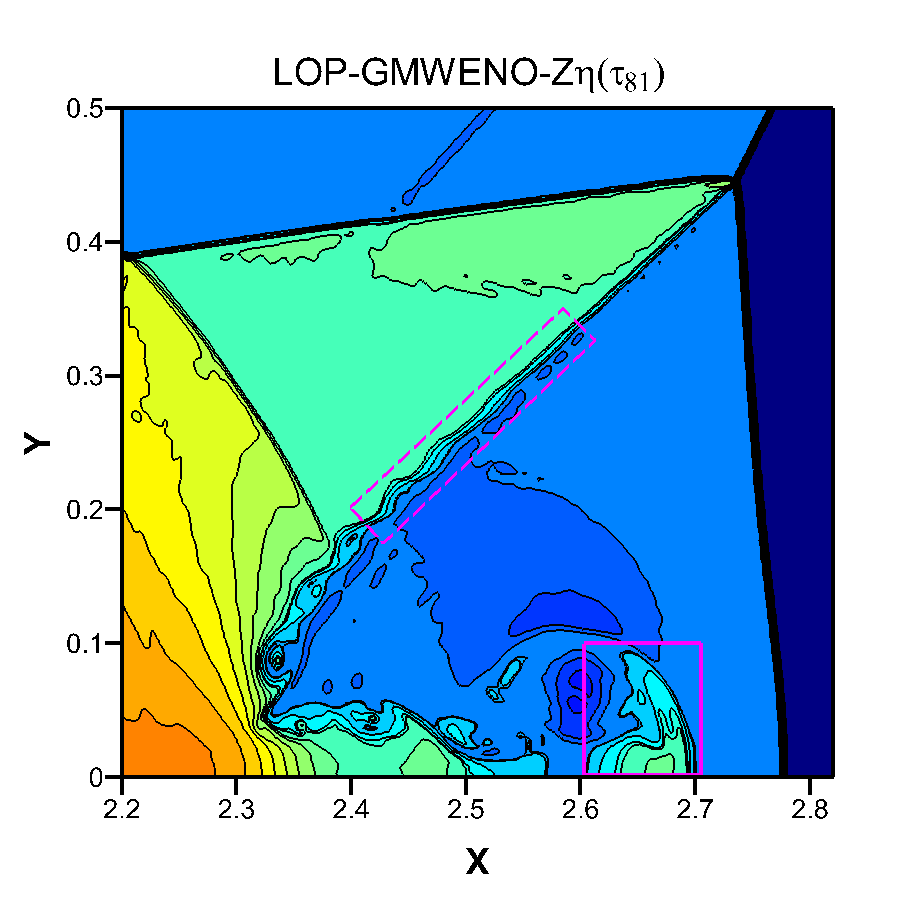} \hspace{0.9em}  
  \includegraphics[height=0.305\textwidth]
  {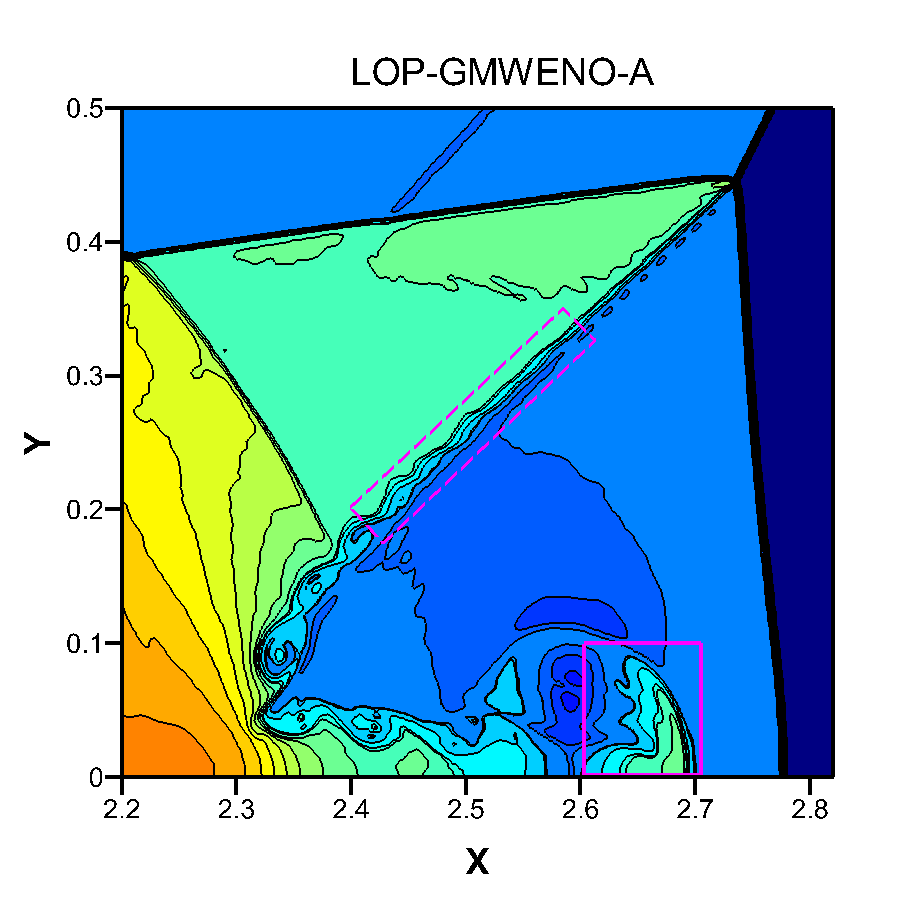}  
\caption{Full and zoomed-in density contours for the DMR problem at output time $t = 0.2$ with a uniform mesh size of $2000 \times 500$.}
\label{fig:ex:DMR}
\end{figure}

Fig. \ref{fig:ex:DMR} illustrates the results from the LOP-GMWENO-Z, 
LOP-GMWENO-Z$\eta(\tau_{81})$, LOP-GMWENO-A schemes and their 
associated WENO-Z, WENO-Z$\eta(\tau_{81})$, WENO-A schemes. All 
these schemes can properly give the main structure of the flow field 
and successfully capture the companion structure and the small 
vortices generated along the slip lines. In addition, they have 
resolved much richer vortical structures than the WENO-JS schemes 
whose solution at the same grid space can be found in \cite{WENO-ACM}. Closer inspection of Fig. \ref{fig:ex:DMR} shows that the 
LOP-GMWENO-Z and LOP-GMWENO-A schemes can slightly decrease the 
numerical oscillations of their associated WENO-Z and WENO-A schemes.

\subsubsection{Forward facing step problem}
In recent studies, this benchmark problem originally presented by 
Woodward and Colella \cite{interactingBlastWaves-Woodward-Colella} 
have been widely used to test the performances of various high order 
schemes \cite{ADER-WENO-1,ADER-WENO-2, WENO-eta,WENO-MAIMi,WENO-ACM}. Its setup is as follows: a step with a height of $0.2$ length 
units located $0.6$ length units from the left-hand end of a wind 
tunnel, which is $1$ length unit wide and $3$ length units long. The 
computational domain of this problem is $\Omega = [0,0.6]\times[0,1]\cup[0.6,3]\times[0.2,1]$ and it is initialized by
\begin{equation*}
\big( \rho, u, v, p \big)(x, y, 0) = (1.4, 3.0, 0.0, 1.0), \quad
(x,y) \in \Omega.
\label{eq:initial_Euer2D:FFS}
\end{equation*}
Reflective boundary conditions are applied along the walls of the 
wind tunnel and the step. At the left and right boundaries, inflow 
and outflow conditions are applied respectively. The mesh resolution 
is $900\times300$.

\begin{figure}[!ht]
\centering
  \includegraphics[height=0.23\textwidth]
  {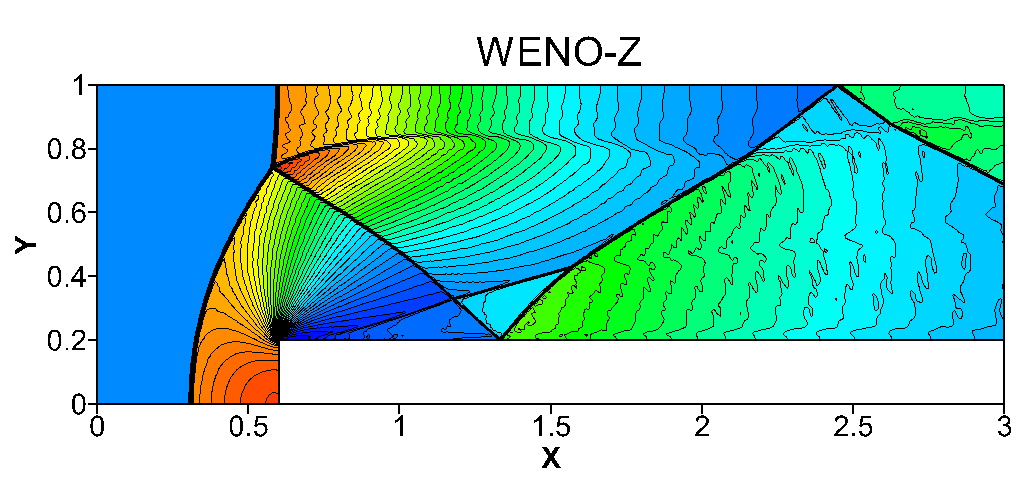} 
  \includegraphics[height=0.23\textwidth]
  {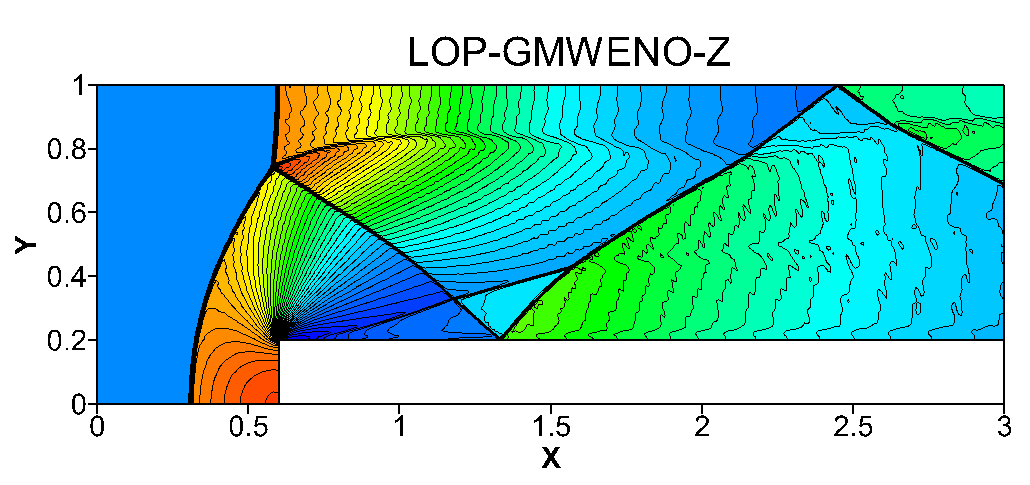}   \\  
  \includegraphics[height=0.23\textwidth]
  {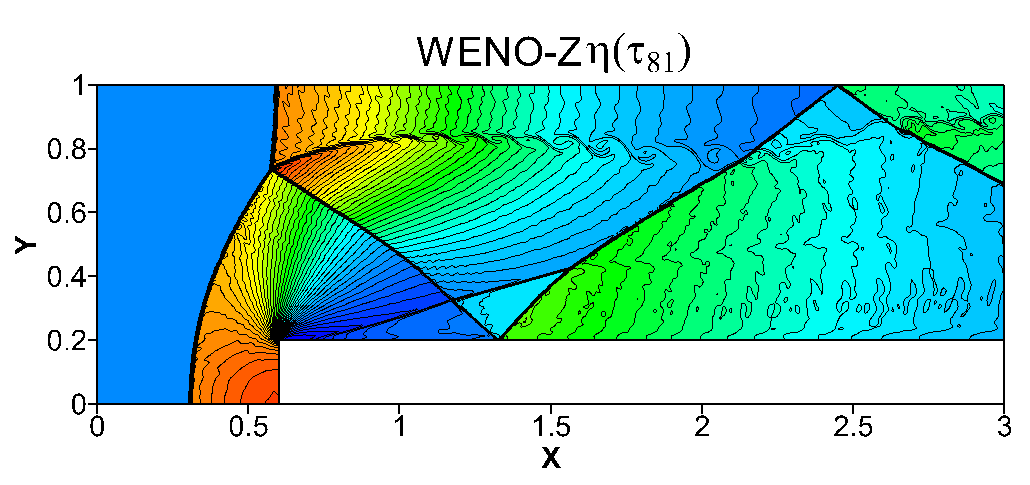} 
  \includegraphics[height=0.23\textwidth]
  {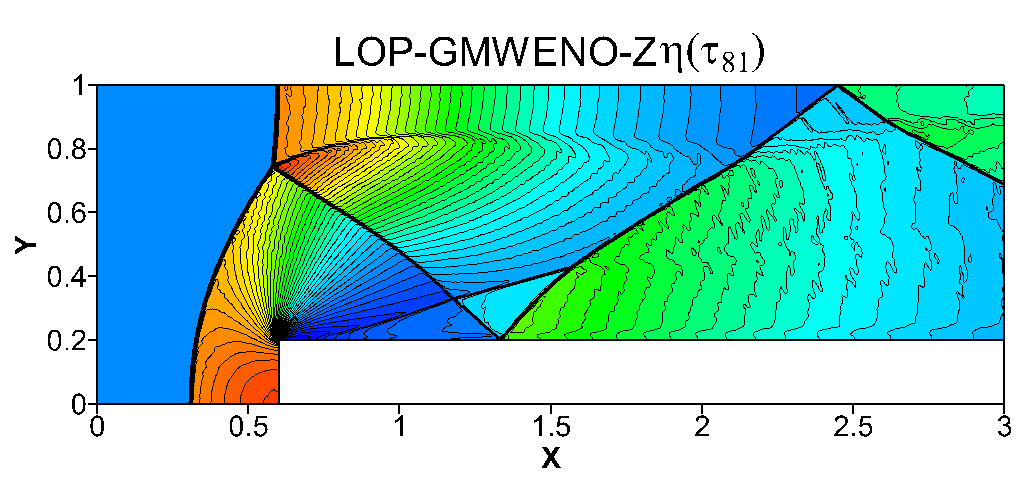} \\    
  \includegraphics[height=0.23\textwidth]
  {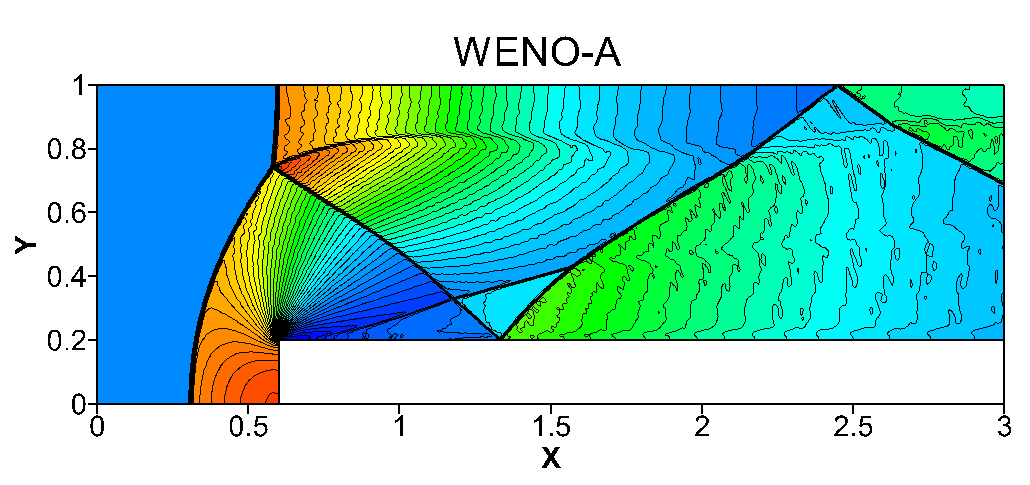} 
  \includegraphics[height=0.23\textwidth]
  {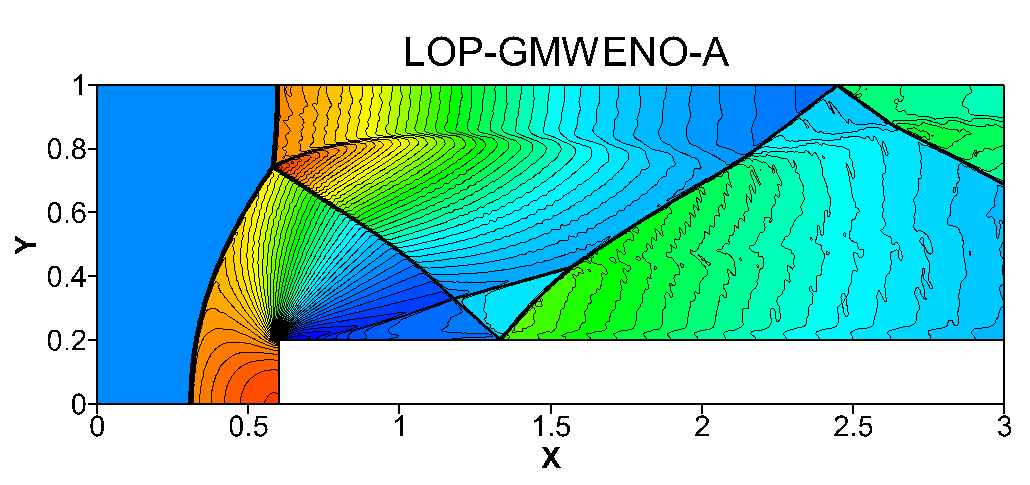}  
\caption{Density contours for the forward facing step problem at output time $t = 4.0$ with a uniform mesh size of $900 \times 300$.}
\label{fig:ex:FFS}
\end{figure}

Fig. \ref{fig:ex:FFS} shows the density contours of the 
LOP-GMWENO-Z, LOP-GMWENO-Z$\eta(\tau{81})$, LOP-GMWENO-A schemes and 
their associated WENO-Z, WENO-Z$\eta(\tau{81})$, WENO-A schemes at 
the final time $t = 4.0$. We can see that all these schemes can 
simulate the complicated structure of the flow field successfully.

\subsubsection{Rayleigh-Taylor instability}
The inviscid Rayleigh-Taylor instability problem used in \cite{RTI-01_JingShi_2003,RTI-02_XuZhengfu_2005} is computed here. The 
computational domain is set to be $[0,1/4]\times[0,1]$ and it is 
initialized by
\begin{equation*}
\big( \rho, u, v, p \big)(x, y, 0) = \left\{
\begin{aligned}
\begin{array}{ll}
(2,0,-0.025c\cdot \cos(8\pi x),2y + 1),         & y \leq 0.5, \\
(1,0,-0.025c\cdot \cos(8\pi x),y + 1.5),        & y > 0.5,
\end{array}
\end{aligned}
\right.
\label{eq:initial_Euer2D:RTI}
\end{equation*}
where the speed of sound is $c=\sqrt{\gamma p / \rho}$ with the 
ratio of specific heats of $\gamma = 5/3$. At $x = 0$ and $x =0.25$, 
the reflective boundary conditions are used. At $y = 0$ and $y = 1$, 
the following Dirichlet boundary conditions are imposed
\begin{equation*}
\begin{aligned}
\begin{array}{l}
\big( \rho, u, v, p \big)(x, y, t) = \left\{
\begin{array}{ll}
(2, 0, 0, 1),   & y = 0, \\
(1, 0, 0, 2.5), & y = 1.
\end{array}
\right.
\end{array}
\end{aligned}
\label{eq:boundary:Euler2D:RTI}
\end{equation*}

\begin{figure}[!ht]
\centering
  \includegraphics[height=0.45\textwidth]
  {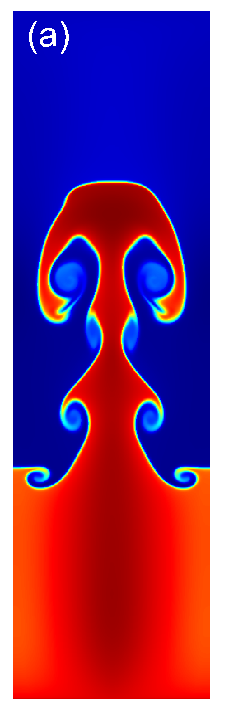} 
  \includegraphics[height=0.45\textwidth]
  {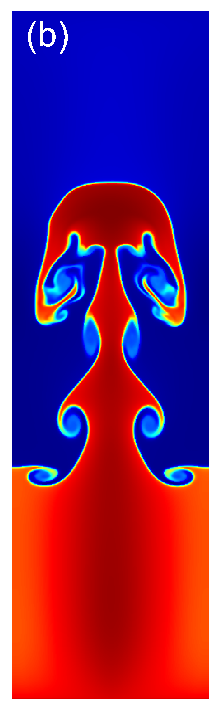}     
  \includegraphics[height=0.45\textwidth]
  {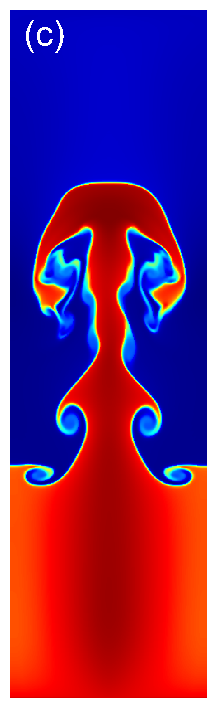} 
  \includegraphics[height=0.45\textwidth]
  {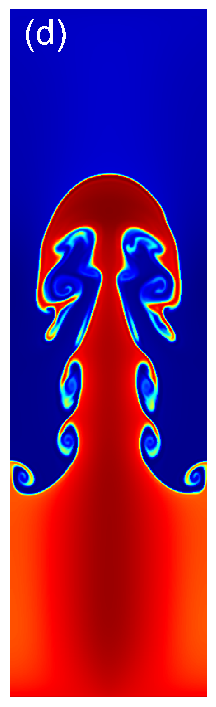}     
  \includegraphics[height=0.45\textwidth]
  {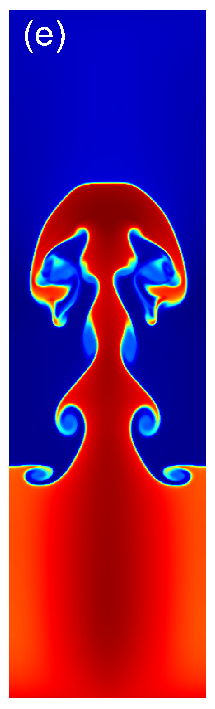} 
  \includegraphics[height=0.45\textwidth]
  {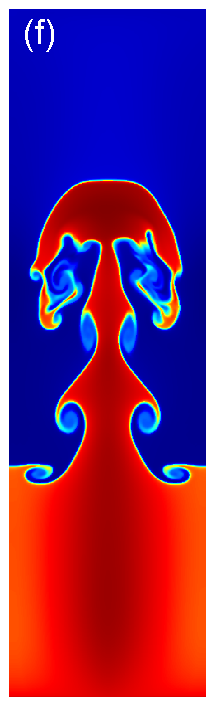}     
  \includegraphics[height=0.45\textwidth]
  {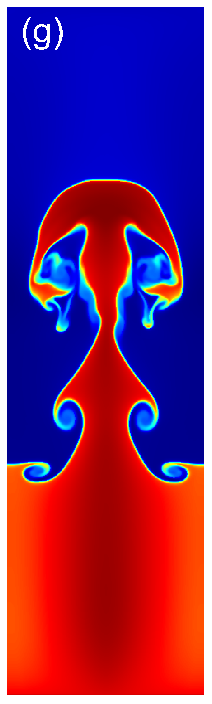}            
\caption{Solutions of the Rayleigh–Taylor instability problem: (a) WENO-JS, (b) WENO-Z, (c) LOP-GMWENO-Z, (d) WENO-Z$\eta(\tau_{81})$, (e) LOP-GMWENO-Z$\eta(\tau_{81})$, (f) WENO-A, (g) LOP-GMWENO-A.}
\label{fig:ex:RTI}
\end{figure}

Fig. \ref{fig:ex:RTI} shows the numerical results of all comsidered 
WENO schemes at a spatial resolution $N_{x} \times N_{y} = 240 \times
960$ with the output time $t = 1.98$. Clearly, the LOP-GMWENO-Z, 
LOP-GMWENO-Z$\eta(\tau_{81})$, LOP-GMWENO-A schemes and their 
associated WENO-Z, WENO-Z$\eta(\tau_{81})$, WENO-A schemes produce 
finer structures than the WENO-JS scheme. Furthermore, the 
LOP-GMWENO-Z, LOP-GMWENO-Z$\eta(\tau_{81})$, LOP-GMWENO-A schemes 
include more noticeable flow asymmetry than their associated WENO-Z, 
WENO-Z$\eta(\tau_{81})$, WENO-A schemes, and it has been indicated \cite{TENO-Fu-CPC2019-01} that the low-dissipation property tends to 
break the flow symmetry, and in other words \cite{TENO-Fu-JCP2016}, 
a more dissipative scheme is more likely to prevent this asymmetry.

%%% Local Variables:
%%% mode: latex
%%% TeX-master: "article"
%%% End:

\section{Conclusions}
\label{secConclusions} 

In this paper, we investigate extending the order-preserving (OP) 
criterion to the WENO-Z-type schemes. The locally order-preserving 
(LOP) mapping is introduced resulting in the improved WENO-Z-type 
schemes, dubbed LOP-GMWENO-X. The major advantages of the new 
LOP-GMWENO-X schemes are in fourfolds. Firstly and most 
importantly, they can successfully remove spurious oscillations and 
meanwhile obtain high resolutions for long simulations of hyperbolic 
problems with discontinuities while their associated WENO-Z-type 
schemes can not. Secondly, on solving hyperbolic problems with 
high-order critical points for long output times, they can achieve 
considerable high resolution. Thirdly, in the region with 
high-frequency but smoothwaves, they can get evidently higher 
resolution than their associated WENO-Z-type schemes. And lastly, 
their post-shock oscillations of the 2D Euler problems with strong 
shock waves are much less and smaller than those of their 
counterpart WENO-Z-type schemes. Extensive numerical experiments 
have been performed to demonstrate these advantages. The numerical 
tests also show that the new schemes can achieve designed 
convergence orders in smooth regions even in the presence of 
critical points.

%%% Local Variables:
%%% mode: latex
%%% TeX-master: "article"
%%% End:

%\input{article_appendix}

\bibliographystyle{model1b-shortjournal-num-names}
\bibliography{refs}

\end{document}